\documentstyle[epsf]{book}

\newcommand{\ee}{\end{equation}}
\newcommand{\be}{\begin{equation}}
\def\adots{\cdot^{\displaystyle{~\cdot}^{\displaystyle{~\cdot}} }}
\def\bin#1#2{\pmatrix{ #1 \cr #2 \cr}}  

\begin{document}

\font\aglio=cmcsc10
\font\cipolla=cmbx10 at 15pt

\null

\vskip 3cm

\centerline{\bf \huge Inverse Problem for   }
\vskip 0.4cm
\centerline{\bf \huge Semisimple Frobenius Manifolds }
\vskip 0.4cm
\centerline{\bf  \huge Monodromy Data and the Painlev\'e VI Equation}
\vskip 2cm
\centerline{\aglio Davide Guzzetti, SISSA-ISAS, October 2000}
 \centerline{SISSA Preprint: 101/2000/FM}

\vskip 3cm

 \centerline{\bf Abstract}
\vskip 0.5 cm
This work is a part the  Ph.D. thesis of Davide Guzzetti,   with 
the supervision of professor B. Dubrovin.  

\vskip 0.2 cm 
We study the inverse problem for   
semisimple Frobenius manifolds of dimension three. We 
 explicitly compute a parametric form of  the solutions of the WDVV
equations of associativity in terms of
solutions of a special Painlev\'e VI equation and we 
show that the solutions are
labelled by a set of monodromy data. The procedure is a relevant application of
the theory of isomonodromic deformations.  We use the parametric form to 
construct  polynomial and algebraic solutions of the WDVV
equations. We also apply the parametric form to construct  
  the generating function of Gromov-Witten invariants
corresponding to the Frobenius manifold called quantum cohomology of the
two dimensional projective space. 

\vskip 0.2 cm 

As a necessary step, we give a contribution to the analysis of the Painlev\'e
 VI equation. We find a class of solutions that   
covers almost all the values of the  monodromy data associated to the
equation, except one point in the space of the data.   
We describe the asymptotic behavior close to
the critical points in terms of two parameters and we 
find the relation  among the parameters at the
different critical points (connection problem). In this way we  
unify and extend   pre-existing results. In particular, we study the elliptic
 representation of Painlev\'e VI and the critical behavior it implies. 
\vskip 0.2 cm

\frontmatter

\chapter{Introduction}
\vskip 3 cm 

 This work is devoted to the construction of solutions of the
 WDVV equations of associativity in 2-D topological field theory as an
 application of the theory of Frobenius manifolds, isomonodromic deformations
 and Painlev\'e equations. It is a part of my  Ph.D. thesis in mathematical
 physics submitted at SISSA-ISAS in October 2000. 

\vskip 0.3 cm 
The WDVV equations of associativity were introduced by Witten \cite{Witten},
Dijkgraaf, Verlinde E., Verlinde H. \cite{DVV}.  
They are  
differential equations satisfied by the {\it primary free energy} $F(t)$ in 
 two-dimensional topological field theory. $F(t)$  is a function 
of the coupling constants $t:=(t^1,t^2,...,t^n)$ $t^i\in {\bf C}$. Let
$\partial_{\alpha}:={\partial \over \partial t_{\alpha}}$.   Given a 
non-degenerate  symmetric matrix $\eta_{\alpha \beta}$, $\alpha,\beta=1,...,n$,
 and numbers  $q_1=0,q_2,...,q_n,d$,   the WDVV equations are
\be
\partial_{\alpha}\partial_{\beta} \partial_{\lambda} F~\eta^{\lambda \mu} 
~\partial_{\mu}\partial_{\gamma}\partial_{\delta} F ~=~ \hbox{ the same with 
$\alpha$, $\delta$ exchanged},
\label{WDVVintroduzione}
\ee
$$
\partial_1 \partial_{\alpha}\partial_{\beta} F = \eta_{\alpha\beta},~~~~~
F(\lambda^{1-q_1} t^1,...,\lambda^{1-q_n}t^n)=\lambda^{3-d} F(t^1,....t^n),~~~
 \lambda \in {\bf C}\backslash \{0\}$$

 The theory of  Frobenius manifolds was introduced  
 by B. Dubrovin \cite{Dub4} to formulate in
 geometrical terms the  
   WDVV equations. It has  links to many branches of
mathematics  like singularity theory and 
 reflection groups \cite{Saito} \cite{SYS}  \cite{Dub6}  \cite{Dub1},
algebraic and 
enumerative geometry \cite{KM} \cite{Manin}, isomonodromic deformations
theory, boundary value problems 
 and Painlev\'e equations \cite{Dub2}.  

If we define $ c_{\alpha \beta \gamma}(t):= \partial_{\alpha}\partial_{\beta} 
\partial_{\gamma} F(t)$, $c_{\alpha\beta}^{\gamma}(t) := \eta^{\gamma \mu} 
c_{\alpha \beta \mu}(t)$   (sum omitted),
and we consider a vector space $A$=span($e_1,...,e_n$), then we obtain a 
family of algebras $A_t$  with the  
multiplication 
$ 
  e_{\alpha} \cdot e_{\beta} := c_{\alpha \beta}^{\gamma}(t) e_{\gamma}. 
$ 
Equation (\ref{WDVVintroduzione}) is equivalent to associativity.

 A {\it Frobenius manifold} 
is a smooth/analytic manifold $M$ 
over ${\bf C}$ whose tangent space $T_t M$ at any $t \in M$ 
is an {\it associative, commutative algebra}. Moreover, there 
exists a non-degenerate bilinear form $<~~, ~~>$ 
 defining a {\it flat metric}.  
The variables $t^1,..,t^n$ are the 
flat coordinates for a point  $t\in M$ and $\eta_{\alpha
 \beta}=<\partial_{\alpha},\partial_{\beta}>$, $\alpha,\beta
=1,2,...,n$. The structure constants in $T_t M$ with respect to the basis
$\partial_1,...,\partial_n$ are   $
c_{\alpha \beta \gamma}(t)= \partial_{\alpha}\partial_{\beta}  
\partial_{\gamma} F(t)$.  

  The manifold is characterized by a family of 
 flat connections $\tilde{\nabla}(z)$, 
 parametrized by a complex number $z$, such that for $z=0$
 the connection is associated to $<~~,~~>$.  To find a flat coordinate $\tilde{t}(t,z)$ we impose $\tilde{\nabla}(z) d\tilde{t}
=0$, which becomes the linear system
$$
    \partial_{\alpha} \xi = z C_{\alpha}(t) \xi,~~~~~~
     \partial_z \xi = \left[{\cal U}(t) +{\hat{\mu} \over z} \right] \xi,
$$
where $
\hat{\mu} := \hbox{diag}(\mu_1,...,\mu_n)$, $\mu_{\alpha}:=q_{\alpha} -{d 
\over 2}$, 
$ 
   \xi$ is a column vector of components $\xi^{\alpha}=\eta^{\alpha\mu}
\partial \tilde{t} /\partial t^{\mu}$, $\alpha=1,...,n$  
and $ C_{\alpha}(t):=
\bigl(
                                  c_{\alpha \gamma}^{\beta}(t)\bigr)$, $
{\cal U}:= \bigl( (1-q_{\mu})t^{\mu} c_{\mu \gamma}^{\beta}(t)  \bigr)
$ are $n\times n$ matrices (sum over repeated indices is omitted). 
 
We restrict to  {\it semisimple} Frobenius manifolds, 
namely analytic Frobenius manifolds  such that the 
matrix ${\cal U}$ can be diagonalized with distinct eigenvalues on an 
open dense subset ${\cal M} \subset M$. 
Then, there exists an invertible matrix $\Psi(t)$ such that 
$
   \Psi {\cal U} \Psi^{-1} = \hbox{diag} (u_1,...,u_n)=:U$, $u_i\neq u_j 
$ for  $i\neq j$ 
 on ${\cal M}$. The  second equation of the above system   becomes:
\be
  {dY\over dz}=\left[U+{V(u)\over z} \right]Y,~~~~u:=(u_1,...,u_n),~~~u_i\in
  {\bf C},
~~~~~V=\Psi^{-1} \hat{\mu} \Psi. 
\label{sistemaitroduzione}
\ee

As it is proved in \cite{Dub1} \cite{Dub2}, 
 $u_1$,...,$u_n$ are local coordinates on ${\cal M}$. 
 Locally we obtain a  change of coordinates, $t^{\alpha} = 
t^{\alpha}(u)$, then $\Psi=\Psi(u)$, 
$V=V(u)$. A local  chart of ${\cal M}$ 
is reconstructed by  
parametric formulae: 
\be
t^{\alpha}=t^{\alpha}(u),~~~~~F=F(u)
\label{parametricintroduzione}
\ee
where $t^{\alpha}(u)$, $F(u)$ are certain meromorphic functions  of
$(u_1,...,u_n)$, $u_i\neq u_j$, which  can be obtained from 
 the coefficients of the system (\ref{sistemaitroduzione}). 

\vskip 0.2 cm 

  The dependence of the system  on $u$ is {\it
  isomonodromic} \cite{JMU}. 
This means that the monodromy data of the system, to
  be  introduced below, do not change for a small
deformation of $u$.   Therefore, the coefficients of the system 
 in every
local chart of ${\cal M}$   are 
 naturally labeled by the monodromy data. To calculate the functions
  (\ref{parametricintroduzione})  in every local chart one has to reconstruct
  the  system  ({\ref{sistemaitroduzione})  
from its {\it monodromy data}. This is  the 
{\it inverse problem}.

\vskip 0.2 cm 

The inverse problem can be formulated as a 
 {\it Riemann-Hilbert boundary value problem}. 
 It can be proved  \cite{Miwa} \cite{Malgrange} \cite{Dub2} 
that if the boundary value problem has solution 
at $u=u^0$ (such that
$u_i^0\neq u_j^0$) for given monodromy data, then 
the solution  is unique and 
it defines $V(u)$, $\Psi(u)$ and 
\ref{parametricintroduzione}   as 
 analytic functions in a neighborhood of $u^0$. 

 Moreover, $V(u)$, $\Psi(u)$  and (\ref{parametricintroduzione}) 
can be continued analytically as 
meromorphic functions on the universal covering of ${\bf C}^n \backslash 
\hbox{diagonals}$, where ``diagonals'' stands for the union of all the 
sets $\{ u \in {\bf C}^n ~|~ u_i=u_j, ~~i\neq j\}$. 
Since $(u_1,...,u_n)$ are local 
coordinates on $M$ they are defined up to permutation. Thus, the 
analytic continuation of the local structure of  ${\cal M}$ 
is described by the 
{\it braid group}, namely the 
fundamental group of 
$({\bf C}^n \backslash 
\hbox{diagonals})/ {\cal S}_n $ (${\cal S}_n$ is the symmetric group of $n$ 
elements). Since every local chart of the atlas covering the manifold is 
labelled by monodromy data, then there exists an action of the 
braid group itself on the monodromy data corresponding to the change of
coordinate chart. This action is described 
in \cite{Dub2} and chapter \ref{cap1}.

\vskip 0.2 cm 

In order to understand the global structure of the manifold $M$ we have to
study  the solution of the inverse problem 
and  (\ref{parametricintroduzione}) 
when two  or more distinct coordinates $u_i$, $u_j$, etc,  merge. $\Psi(u)$, 
$V(u)$ and 
 (\ref{parametricintroduzione}) are multi-valued meromorphic functions of
$u=(u_1,...,u_n)$ and the branching occurs when $u$ goes around a 
 loop around
the set of diagonals $\bigcup_{ij}\{ u \in {\bf C}^n ~|~ u_i=u_j, ~~i\neq
j\}$.    $\Psi(u)$, $V(u)$ and 
 (\ref{parametricintroduzione}) have singular behaviour if $u_i\to u_j$
$(i\neq j$). We call such
behavior {\it critical behaviour}. Although it is impossible to
solve  the boundary value problem exactly,  
except for special cases occurring for $2\times 2$ systems, 
we may hopefully compute the
asymptotic/critical 
behaviour of the  solution,  
using the isomonodromic deformation method 
\cite{JMU} \cite{IN}. 
We will face the problem in the first non-trivial case, namely for three
dimensional Frobenius manifolds. 

Instead of analyzing the boundary value
problem directly, we exploit the isomonodromic dependence  of  the system  
(\ref{sistemaitroduzione}) on $u$, which implies that 
  the solution of the inverse problem must
satisfy  the nonlinear equations 
\be
{\partial V \over \partial u_k} = \left[V_k,V \right]
\label{isomonodromyintroduzione1}
\ee
 \be
{\partial \Psi, \over \partial u_k}= V_k \Psi,~~~~~k=1,...,n
\label{isomonodromyintroduzione2}
\ee
where   $V_k$ is a $n\times n$ matrix whose entries are  
$(V_k)_{ij}={\delta_{ki}-\delta_{kj}\over u_i-u_j}~V_{ij}$.   The WDVV
equations are equivalent to (\ref{isomonodromyintroduzione1})
(\ref{isomonodromyintroduzione2}).  
 For 3-dimensional Frobenius manifolds, (\ref{isomonodromyintroduzione1})
(\ref{isomonodromyintroduzione2})  are
reduced to a special case of the Painlev\'e 6 equation \cite{Dub2}:  
$$
{d^2y \over dx^2}={1\over 2}\left[ 
{1\over y}+{1\over y-1}+{1\over y-x}
\right]
           \left({dy\over dx}\right)^2
-\left[
{1\over x}+{1\over x-1}+{1\over y-x}
\right]{dy \over dx}$$
\be
+{1\over 2}
{y(y-1)(y-x)\over x^2 (x-1)^2}
\left[
(2\mu-1)^2+{x(x-1)\over (y-x)^2}
\right]
,~~~~~\mu \in {\bf C},~~~~x={u_3-u_1\over u_2-u_1}
\label{painleveintroduzione}
\ee
The parameter $\mu$ appears in the matrix $\hat{\mu}=$diag$(\mu,0,-\mu)$ of
(\ref{sistemaitroduzione}). 
 In chapter \ref{cap1} and \ref{inverse reconstruction} we show that the 
 entries of $V(u)$ and $\Psi(u)$  are rational
functions of $x$,  $y(x)$, ${dy\over dx}$.
  The critical behaviour of $V(u)$, $\Psi(u)$ and
(\ref{parametricintroduzione}) is a consequence of the critical behaviour of
the transcendent $y(x)$ close to the {\it critical points}
$x=0,1,\infty$.

Let  $F_0(t):= 
{1\over 2}\left[ (t^1)^2 t^3+t^1(t^2)^2\right]$. 
 We prove in chapter \ref{inverse reconstruction} that for generic
$\mu$ the parametric representation (\ref{parametricintroduzione}) becomes  
\be
 t^2(u)= \tau_2(x,\mu) ~(u_2-u_1)^{1+\mu},~~~~
  t^3(u)=\tau_3(x,\mu)~(u_2-u_1)^{1+2\mu}
\label{tintroduzione11}
\ee
\be
F(u)=F_0(t)+ {\cal F}(x,\mu)~(u_2-u_1)^{3+2\mu}
\label{Fintroduzione11},~~~~x={u_3-u_1\over u_2-u_1}
\ee
where $ \tau_2(x,\mu), \tau_3(x,\mu), {\cal F}(x,\mu)$ are 
 certain 
rational functions of $x$, $y(x)$, ${d y \over dx}$ and $\mu$, which we
 computed explicitly. More details
will be given below.

\vskip 0.2 cm 

The two integration constants in $y(x)$ -- and thus in the corresponding
 solution 
 of (\ref{isomonodromyintroduzione1})
(\ref{isomonodromyintroduzione2}) -- and the parameter  $\mu$ are contained
 in the three entries $(x_0,x_1,x_{\infty})$ of the {\it Stokes's matrix}
$$
     S=\pmatrix{ 1&x_{\infty}&x_0\cr
              0&1&x_1\cr
              0&0&1 \cr
             }, ~~~~~\hbox{ such that }
              x_0^2+x_1^2+x_{\infty}^2-x_0x_1x_{\infty}=4 \sin^2(\pi \mu) 
$$ 
 of  the system (\ref{sistemaitroduzione}). 
 The Stokes matrix is part of the {\it monodromy data} of the 
system (\ref{sistemaitroduzione}). They will
 be discussed in chapter \ref{cap1}. Here we briefly introduce the Stokes
 matrix. 
At $z=\infty$ there is a formal solution  of (\ref{sistemaitroduzione})
$
   {Y}_F = \left[ I + {F_1 \over z}+{F_2 \over z^2}+
... \right]~e^{z~U}$ 
where $F_j$'s are $n \times n$ matrices. It is a well known   result  
that fundamental matrix  solutions exist which have ${Y}_F$ as 
asymptotic expansion for $z \to \infty$ \cite{BJL1} . 
 Let $l$ be a 
  generic oriented line  passing through the origin.  Let  $l_{+}$ be the
 positive half-line and  $l_{-}$ the negative one. 
 Let  ${\cal S}_L$ and ${\cal S}_R$ be two sectors in the complex plane 
to the
 left and to the right of $l$ respectively. 
 There exist unique fundamental matrix solutions  $Y_L$ and $Y_R$ having the 
asymptotic expansion $Y_F$ for $x\to \infty$ 
in ${\cal S}_L$ and ${\cal S}_R$ 
respectively \cite{BJL1}.  They are related by a connection matrix $S$,
called {\it 
Stokes matrix}, such that $
    {Y}_L(z)={Y}_R(z)S
$  
for $z\in l_{+}$.

\vskip 0.2 cm

 A further step is to invert 
 (\ref{parametricintroduzione}) in order to obtain a closed form  
 $F=F(t^1,...,t^n)$.  
The final purpose of the inversion is to understand the analytic properties of
the solution $F(t)$ of the WDVV equations.

\vskip 0.3 cm 
The entire procedure described above is
 an application of the isomonodromic deformation theory to the WDVV equations
 and is the object of the thesis. 

 The first step is to properly choose the monodromy data in order to arrive at
 physically or geometrically interesting Frobenius manifolds.
 We consider the quantum cohomology of
 projective spaces as an important example of semisimple 
 Frobenius manifold. For this example we compute the monodromy data.

\vskip 0.2 cm
 The
{\it quantum cohomology} 
of ${\bf CP}^d$, denoted by $QH^{*}({\bf CP}^d)$, is a $(d+1)$-dimensional
semisimple Frobenius manifold (\cite{KM} \cite{M} and chapter \ref{cap2}
below).  
 It has relations to enumerative geometry. 
A well known example is the quantum cohomology of ${\bf CP}^2$. It corresponds
to the solution of
the WDVV equations for $n=3$  which generates 
 the numbers $N_k$  of rational curves ${\bf
CP}^1\longrightarrow{\bf CP^2}$ of degree $k$ passing through $3k-1$ generic
points. Namely 
\be
F(t^1,t^2,t^3)= 
{1\over 2}\left[ (t^1)^2 t^3+t^1(t^2)^2\right]+ \sum_{k=1}^{\infty} {N_k \over
(3k-1)!} (t^3)^{3k-1} e^{k t^2}, ~~~\hbox{ for } (t^3)^3 e^{t^2} \to 0
\label{kontsevichintroduzione}
\ee
 The global analytic properties of this function are unknown, except for
 $(t^3)^3 e^{t^2} \to 0$,   and the 
inverse reconstruction of the corresponding Frobenius manifold starting from
 its monodromy data may shed some
light on these properties.

\vskip 0.5 cm 

\noindent
{\bf $\bullet$ Results }

\vskip 0.3 cm
{\bf i)} The monodromy data of $QH^{*}({\bf CP}^d)$.  The result 
is already pubblished in \cite{guz}, therefore we omit the entire
part of the thesis devoted to it  and we refer to the paper. We just recall
the main theorem of \cite{guz}

\vskip 0.2 cm
\noindent
{\bf Theorem [chapter \ref{stokes}]:} { \it There exists a chart of
$QH^{*}({\bf CP}^d)$ where the Stokes's matrix $S=(s_{ij})$,
$i,j=1,2,...,d+1$,   has the canonical form:
$$ 
   s_{ii}=1,~~~~~~s_{ij}=(-1)^{j-i}\bin{d+1}{j-i},~~~~s_{ji}=0,~~~~i<j
$$
For any other local chart of $QH^{*}(\bf{CP}^{d})$ the Stokes matrix is
obtained from the canonical form by the action of the braid group.}

\vskip 0.2 cm 
We recall that such computation is the first step in the process of inverse
reconstruction of a semisimple Frobenius manifold starting from its monodromy
data. 

We also recall the there is a further motivation for this
computation. Namely it proves, for projective spaces, a long-lasting
conjecture aboute links between topological field theory and the theory of
derived categories of coherent sheaves. We refer to \cite{guz} for details,
(and to \cite{CV} \cite{Zas}, \cite{Dub3} for the origin of the conjecture). 

\vskip 0.2 cm 
In \cite{guz} 
 we also studied  the structure of the monodromy group of $QH^{*}({\bf CP}^d)$ 
and we proved that it is related to the hyperbolic triangular groups (the case
 $d=2$ was already studied in \cite{Dub2}). The
 concept of {\it monodromy group of a Frobenius manifold} is explained in
 chapter \ref{cap1}.

\vskip 0.4 cm 

{\bf ii)}
 Once the monodromy data, in particular the Stokes' matrix, are known, we
have to solve the inverse problem for the system (\ref{sistemaitroduzione}) 
 in order to   
  obtain  $V(u)$, $\Psi(u)$ and the parametric representation
 (\ref{parametricintroduzione}).  We face the problem for the first
 non-trivial case $n=3$, which is
reduced to the  Painlev\'e 6 equation. 
Therefore, we devote  chapter \ref{PaInLeVe} to the investigation of the
behaviour of  $y(x)$ close to the {\it critical points} 
 $x=0,1,\infty$ and to the
connection problem between the parameter governing that behaviour at different
critical points. Moreover we give the explicit dependence of the parameters
on the entries $(x_0,x_1,x_{\infty})$ of the Stokes' matrix.  

\vskip 0.2 cm

The six classical Painlev\'e  equations were discovered by Painlev\'e \cite{pain} 
 and 
Gambier \cite{gamb}, 
who classified all the second order ordinary differential equations
of the type 
$$
         {d^2 y \over dx^2}= {\cal R}\left(x,y,{dy\over dx}\right)
$$
where ${\cal R}$ is rational in ${dy\over dx}$, meromorphic in $x$ and
 $y$. The Painlev\'e  equations 
satisfy the {\it Painlev\'e property} of  absence of movable critical
 singularities. The 
 general solution of the $6^{\hbox{th}}$ Painlev\'e equation  
 can be analytically continued to a meromorphic function on the universal
   covering of ${\bf CP}^1\backslash \{ 0,1 ,\infty \}$. 
 For generic values of the integration constants and of the parameters in
   the equation, the solution can not be expressed via elementary or classical
   transcendental functions. For this reason, the solution is called a {\it
   Painlev\'e transcendent}.

\vskip 0.2 cm 
The connection problem for a class of solutions to  
the Painlev\'e 6 equation was solved by
Jimbo \cite{Jimbo}  
for the general Painlev\'e  equation  with generic values of its 
coefficients $\alpha$, $\beta$, $\gamma$ $\delta$ (in standard notation of
\cite{IN}),  using the
isomonodromic deformation theory developed in \cite{JMU} \cite{JM1}.  
 Later, Dubrovin-Mazzocco \cite{DM} applied Jimbo's procedure to
 (\ref{painleveintroduzione}), with the restriction  $2\mu \not \in {\bf Z}$.  
The connection problem was solved by the authors above for the class of
transcendents 
having the following local behaviour at the critical points $x=0,1,\infty$: 
\be
y(x)= a^{(0)} x^{1-\sigma^{(0)}}(1+O(|x|^{\delta})),~~~~x\to 0,
\label{loc1introduzione}
 \ee
\be
y(x)= 1-a^{(1)}(1-x)^{1-\sigma^{(1)}} (1+O(|1-x|^{\delta})),~~~~x\to 1,
\label{loc2introduzione}
 \ee
\be
y(x)= a^{(\infty)}
 x^{-\sigma^{(\infty)}}(1+O(|x|^{-\delta})),~~~~x\to \infty,
\label{loc3introduzione}
\ee
where $\delta$ is a small positive number, $a^{(i)}$ and
 $\sigma^{(i)}$ are complex numbers such that $a^{(i)}\neq 0$ and 
 $$0\leq \Re \sigma^{(i)}<1.$$ 
This behaviour is true if $x$ converges
 to the critical points inside a sector with vertex
 on the corresponding critical point, along a radial direction 
in the $x$-plane.  
The {\it connection problem}, i.e.  the problem of 
finding  the relation among the three pairs $(\sigma^{(i)},a^{(i)})$,
$i=0,1,\infty$, was solved  thanks to  the link 
 between the Painleve' equation and a Fuchsian system of
 differential equations
$$
   {dY\over dz}=\left[ {A_0(x)\over z}+{A_x(x) \over z-x}+{A_1(x)\over
z-1}\right] Y$$
where the  $2\times 2$ matrices  $A_i(x)$ ($i=0,x,1$ are labels) 
have isomonodromic dependence on $x$ and satisfy Schlesinger equations.  
 The local behaviours (\ref{loc1introduzione}), (\ref{loc2introduzione}),
(\ref{loc3introduzione})  were
proved  using  
a result on the asymptotic
 behaviour of a class of  solutions of Schlesinger equations 
 proved by Sato, Miwa, Jimbo in \cite{SMJ}.  The
 connection problem was solved because the parameters
 $\sigma^{(i)}$, $a^{(i)}$ were expressed as functions of 
 the monodromy data of the fuchsian system
 associated to (\ref{painleveintroduzione}). For  studies on the 
 asymptotic behaviour of the
coefficients of  Fuchsian systems and Schlesinger
equations  see also \cite{Boli}.

The  monodromy data of the Fuchsian system turn out to be expressed in terms
of  the  triple  $(x_0,x_1,x_{\infty})$ 
of entries of the Stokes matrix \cite{DM}.  There
 exists a one-to-one correspondence between  triples and  branches of 
 the Painlev\'e transcendents.\footnote{ 
There are only some exceptions to the one-to-one correspondence above, which
are already treated in \cite{M}. In order to rule them out we require 
that at most one of the entries $x_i$ of the triple may be zero and that
$(x_0,x_1,x_{\infty})\not \in \{(2,2,2)$ $(-2,-2,2)$, $(2,-2,-2)$, $(-2,2,-2)
\}$. See \cite{M}.}  
In other words, any branch $y(x)$ is parametrized by
 a triple, namely  
 $y(x)=y(x;x_0,x_1,x_{\infty})$.   As it is proved in \cite{DM}, the
 transcendents (\ref{loc1introduzione}),  (\ref{loc2introduzione}),
 (\ref{loc3introduzione}) are parametrized by a triple according to the
 formulae 
$$ 
   x_i^2= 4 \sin^2\left({\pi\over 2}
   \sigma^{(i)}\right),~~~i=0,1,\infty,~~~~~~~~0\leq \Re \sigma^{(i)}<1.
$$
A more complicated expression gives $a^{(i)}=a^{(i)}(x_0,x_1,x_{\infty})$.

\vskip 0.2 cm

Due to the  restriction $0\leq \Re \sigma^{(i)}<1$, the formulae of
Dubrovin-Mazzocco do not work  
if at least one $x_i$ ($i=0,1,\infty$) is real and $|x_i|\geq
 2$. This is the case of $QH^{*}({\bf
CP}^2)$, because $(x_0,x_1,x_{\infty})=(3,3,3)$. 
 To overcome this limitation, in  chapter \ref{PaInLeVe} 
 we find the critical behaviour and we solve the
 connection problem for  all the triples satisfying
$$
    x_i\neq \pm 2 ~~\Longrightarrow ~~ \sigma^{(i)} \neq 1,~~~~i=0,1,\infty
$$
 The method used is an extension of Jimbo and Dubrovin-Mazzocco's method and
 relies on the isomonodromy deformation theory.  We prove that 

\vskip 0.2 cm 
\noindent
{\bf Theorem 1 [chapter \ref{PaInLeVe}]: } {\it Let $\mu \neq 0$.  For any $\sigma^{(0)} \not \in
(-\infty,0)\cup [1,+\infty)$, for
any $a^{(0)} \in {\bf 
C}$, $a\neq 0$, for any $\theta_1,\theta_2  \in {\bf R}$  and for any
$0<\tilde{\sigma }<1$, 
there exists a sufficiently small positive $\epsilon^{(0)}$ 
such that  the  equation (\ref{painleveintroduzione}) has a solution   
$y(x;\sigma^{(0)},a^{(0)})$ with the behaviour  
$$
    y(x;\sigma^{(0)},a^{(0)})=a^{(0)} x^{1-\sigma^{(0)}} \left(
1+O(|x|^{\delta})
\right)~~~~,0<\delta<1,
$$
as $x\to 0$ in the domain $
D(\sigma^{(0)})=D(\epsilon^{(0)};\sigma^{(0)};\theta_1,\theta_2,\tilde{\sigma})
  $ 
defined by  
$$|x|<\epsilon^{(0)},~~~~~ 
        \Re \sigma^{(0)} \log|x|+\theta_2 \Im \sigma^{(0)} 
\leq \Im \sigma^{(0)} \arg(x)
        \leq (\Re \sigma^{(0)} - \tilde{\sigma}) \log|x| + \theta_1 \Im
        \sigma^{(0)} 
$$
For $\sigma^{(0)}=0$ the domain is simply $|x|<\epsilon^{(0)}$. 
}
\vskip 0.2 cm
We note that $\epsilon^{(0)}$ depends on the choice of $\theta_1$ and $a$. 
The critical behaviour in theorem 1 coincides 
 with (\ref{loc1introduzione}) for $0\leq \Re \sigma^{(0)}<1$, 
but for $\Re\sigma^{(0)}<0$ and
$\Re \sigma^{(0)}\geq 1$ it holds true if $x\to 0$ along a
spiral, according to the shape of  $
D(\epsilon^{(0)};\sigma^{(0)};\theta_1,\theta_2,\tilde{\sigma})$. 
Instead, the behaviour when $x\to 0$ along a radial path may be  more
complicated and no indication is given by the theorem.

By symmetries of (\ref{painleveintroduzione}), 
we also 
prove the existence of solutions with  local behaviour at 
 $x=1$ 
$$y(x,\sigma^{(1)},a^{(1)})= 1-a^{(1)} (1-x)^{1-\sigma^{(1)}}
(1+O(|1-x|^{\delta }))~~~~x\to 1 
$$
$$
   a^{(1)}\neq 0,~~~~\sigma^{(1)}\not \in (-\infty,0)\cup [1,+\infty)
$$
and $$
    y(x; \sigma^{(\infty)},
 a^{(\infty)}) = a^{(\infty)} {x}^{ \sigma^{(\infty)} }\left(
 1+O({1\over |x|^{\delta}}) \right)~~~~~\tilde{x} \to \infty
$$
$$
  a^{(\infty)}\neq 0,~~~~\sigma^{(\infty)}\not \in (-\infty,0)\cup [1,+\infty)
$$
in suitable domains $D(\sigma^{(1)})$, $D(\sigma^{(\infty)})$ 
which will be described in chapter (\ref{PaInLeVe}).

\vskip 0.2 cm 
We also prove that 
the analytic continuation of a branch 
$y(x;x_0,x_1,x_{\infty})$ to the domains of theorem 1 
 is governed by parameters $\sigma^{(i)}$, $a^{(i)}$ 
given by the following 

\vskip 0.2 cm
\noindent
{\bf Theorem 2 [chapter \ref{PaInLeVe}]:} {\it 
 For any set of monodromy data ($x_0,~x_1,~x_{\infty}$)
 such that $x_0^2+x_1^2+x_{\infty}^2-x_0x_1x_{\infty}=4 \sin^2(\pi
 \mu)$ and $x_i\neq \pm2$ 
there exist a unique solution 
$y(x;\sigma^{(i)},a^{(i)})$ in $D(\sigma^{(i)})$
  with parameters   $\sigma^{(i)}$  and $a^{(i)}$
 obtained as follows:  
 $$
   x_i^2= 4 \sin^2 \left( {\pi \over 2} \sigma^{(i)}  \right) ,~~~~~~~~~
\sigma^{(i)}\in {\bf C}\backslash \{(-\infty,0)\cup [1,+\infty)
\},~~~~~~i=0,1,\infty 
$$
$$                  %
  a^{(0)}={iG(\sigma^{(0)},\mu)^2\over 2 \sin(\pi \sigma^{(0)})} 
      \Bigl[
 2(1+e^{-i\pi\sigma^{(0)}})-f(x_0,x_1,x_{\infty})(x_{\infty}^2+e^{-i\pi
\sigma^{ (0)}} x_1^2)  
\Bigr]  ~ f(x_0,x_1,x_{\infty})
$$
where  
$$
   f(x_0,x_1,x_{\infty}):={4-x_0^2\over 2-x_0^2-
2\cos(2\pi\mu)},
~~~~G(\sigma^{(0)},\mu)= {1\over 2}{4^{\sigma^{(0)}}\Gamma({\sigma^{(0)}+1\over 2})^2\over
                  \Gamma(1-\mu+ {\sigma^{(0)}\over
                  2})\Gamma(\mu+{\sigma^{(0)}\over 2})}$$
The parameters $a^{(1)}$, $a^{(\infty)}$ 
 are obtained like $a^{(0)}$,  provided that we do the
substitutions $(x_0,x_1,x_{\infty})\mapsto ( x_1,x_0,x_0
   x_1-x_{\infty})$ and  $(x_0,x_1,x_{\infty})\mapsto (x_{\infty},
-x_1,x_0-x_1x_{\infty})$ respectively (and  $\sigma^{(0)} \mapsto
\sigma^{(1)}$ and  $\sigma^{(0)} \mapsto
\sigma^{(\infty)}$  respectively) in the above formulae. 
}  
\vskip 0.2 cm

We remark that the theorem is actually a bit more complicated, we need to
 distinguish some sub-cases and to be careful about the definition of the
 branch cuts for $y(x;x_0,x_1,x_{\infty})$; 
 we refer to chapter
 \ref{PaInLeVe} for details.

 The {\it connection problem} 
 for the transcendents $y(x;\sigma^{(i)}, 
a^{(i)})$ is now solved, because we are able to compute  $
(\sigma^{(i)}, 
a^{(i)})$ for $i=0,1,\infty$ in terms of a fixed triple $(x_0,x_1,x_{\infty})$.

We also discuss the problem of the analytic continuation of the branch
$y(x;x_0,x_1,x_{\infty})$, namely we discuss 
how $(x_0,x_1,x_{\infty})$ change when $x$ describes loops around
$x=0,1,\infty$.

\vskip 0.2 cm

The above theorem implies that we can always restrict to the case $0\leq \Re
\sigma^{(i)}\leq 1$, $\sigma^{(i)}\neq 1$, so the critical behaviours
$y(x;\sigma^{(i)}, a^{(i)})$ coincide with (\ref{loc1introduzione}),
(\ref{loc2introduzione}), (\ref{loc3introduzione}), except for the case $\Re
\sigma^{(i)} =1$, where the critical behaviour holds true only if $x$
converges to a critical point along a spiral. 

We use  the elliptic representation of the transcendents in order to
investigate this last case and its critical behaviour along radial paths. We
can restrict here to $x=0$ because the symmetries of
(\ref{painleveintroduzione}) yield the behaviour close to the other critical
points.  The elliptic representation was introduced by R.Fuchs in \cite{fuchs}:
$$ 
  y(x)=\wp\left({u(x)\over 2};\omega_1(x),\omega_2(x)\right)
$$
Here  $u(x)$ solves a non-linear second order differential equation and
  $\omega_1(x)$, $\omega_2(x)$ are two elliptic integrals, expanded for
  $|x|<1$ in terms of hypergeometric functions: 
$$
\omega_1(x)= {\pi \over 2} 
\sum_{n=0}^{\infty}{ \left[\left({1 \over 2}\right)_n\right]^2
  \over (n!)^2 } x^n
$$
$$
\omega_2(x)= -{i\over 2}\left\{ 
\sum_{n=0}^{\infty}{ \left[\left({1 \over 2}\right)_n\right]^2
  \over (n!)^2 } x^n ~\ln(x) +   
\sum_{n=0}^{\infty}{ \left[\left({1 \over 2}\right)_n\right]^2
  \over (n!)^2 } 2\left[ \psi(n+{1\over 2}) - \psi(n+1)\right]
x^n \right\}
$$
where $\psi(z):= {d\over dz} \ln \Gamma(z)$. 

 We study  the critical behaviour implied by this representation. 
  We show that the representation also provides the critical
 behaviour along radial paths for $\Re \sigma^{(i)}=1$. More precisely we
 prove the following

\vskip 0.2 cm
\noindent
{\bf Theorem 3 [chapter \ref{PaInLeVe}]:} 
 {\it For any complex $\nu_1$, $\nu_2$ such that 
$$ 
  \nu_2\not \in (-\infty,0]\cup [2,+\infty)
$$
there exists a sufficiently small $r$ such that 
$$
  y(x)= \wp(\nu_1 \omega_1(x)+\nu_2 \omega_2(x)
  +v(x);\omega_1(x),\omega_2(x))
$$
in the domain $ {\cal D}(r;\nu_1,\nu_2)$ defined as 
$$
|x|<r,~~~~   \Re \nu_2 \ln|x|+ C_1-\ln r < \Im \nu_2
   \arg x < (\Re \nu_2 -2)\ln|x| +C_2 + \ln r,
$$
$$
C_1:= -\bigl[4 \ln 2 \Re \nu_2 + \pi \Im \nu_1\bigr],~~~~C_2:=C_1+8\ln 2.
$$
The function $v(x)$ is holomorphic  in $ {\cal D}(r;\nu_1,\nu_2)$ and has
convergent expansion 
$$
  v(x)= \sum_{n\geq 1} a_n x^n +\sum_{n\geq 0,~m\geq 1} b_{nm} x^n
  \left({e^{-i\pi \nu_1}\over 16^{2-\nu_2}} x^{2-\nu_2}\right)^m +\sum_{n\geq
  0,~m\geq 1}c_{nm} x^n \left( 
{e^{i\pi \nu_1} \over 16^{\nu_2}} x^{\nu_2}\right)^m     
$$
where $a_n$, $b_{nm}$, $c_{nm}$ are certain rational functions of  $\nu_2$. 
Moreover, there exists a constant $M(\nu_2)$ depending on $\nu_2$ such that 
 $v(x)\leq M(\nu_2) \left(|x|+\left|{e^{-i\pi \nu_1}\over 16^{2-\nu_2}} x^{2-\nu_2} \right|+\left| 
{e^{i\pi \nu_1} \over 16^{\nu_2}} x^{\nu_2}\right| \right)$ in  ${\cal
D}(r;\nu_1,\nu_2)$ .
}

\vskip 0.2 cm

We note that for $\mu={1\over 2}$ the function $v(x)$ vanishes, and we
obtain {\it Piccard solutions} \cite{Picard} whose critical behaviour is
studied in \cite{M}. 

The  transcendent of 
theorem 3 coincides with $y(x;\sigma^{(0)},a^{(0)})$ of theorem
1 on the domain  $D(\epsilon^{(0)},\sigma^{(0)})
\cap  {\cal D}(r;\nu_1,\nu_2)$ with critical behaviour 
specified by  $\sigma^{(0)}=1-\nu_2$ and $a^{(0)}= -{1\over 4} \left[
{e^{i\pi\nu_1} \over 16^{\nu_2-1}}\right]$.  The identification of $a^{(0)}$
and $\sigma^{(0)}$  makes it possible 
 to connect $\nu_1$ and $\nu_2$ to the monodromy data 
$(x_0,x_1,x_{\infty})$ according to theorem
2.

On the other hand,
 we will prove that 
the behaviour  implied by the elliptic representation is oscillatory along
 paths contained in $ {\cal D}(r;\nu_1,\nu_2)$ which
 are parallel to the boundaries of the domain  in the $(\ln|x|,\arg(x))$
 plane, namely $ \Im \nu_2
   \arg x = (\Re \nu_2 -2)\ln|x| +C_2 + \ln r$ and $ \Im \nu_2
   \arg x= \Re \nu_2 \ln|x|+ C_1-\ln r$. 
This follows from the Fourier expansion of the Weierstrass elliptic function
 which will be discussed in section \ref{beyond}. 

In particular the case $\Re \nu_2=0$ ($\nu_2\neq 0$) coincides
 with  $\Re \sigma^{(0)}=1$ ($\sigma^{(0)}\neq 1$)  and the 
paths 
parallel to the  boundary 
$ \Im \nu_2  \arg x=  C_1-\ln r$ are radial paths. 
As a consequence of theorem 3, the critical behaviour along a
radial path (equivalently, inside a sector) is 
\be
 y(x)= O(x)+{ 1+O(x)\over \sin^2\left({\nu\over 2} \ln x - \nu \ln 16 +{\pi
\nu_1\over 2} + \sum_{m=1}^{\infty} c_{om}(\nu) \left[\left({e^{i\pi \nu_1}
\over 16^{i\nu}}\right)x^{i\nu}\right]^m\right)},  ~~~~x\to 0.
\label{fuchsshimomuraintroduzione}
\ee
The number $\nu$ is real, $ \nu \neq 0$ and  $\sigma^{(0)}=1-i\nu$.  
The series $ \sum_{m=1}^{\infty} c_{om}(\nu)
  \left[\left({e^{i\pi \nu_1} 
\over e^{i\nu}}\right)x^{i\nu}\right]^m$ converges and defines a
holomorphic and bounded function in the domain ${\cal D}(r;\nu_1, i\nu)$ 
$$ 
|x|<r,~~~~  C_1 - \ln r < \nu \arg x < -2 \ln |x| +C_2 + \ln r 
$$ 
 Note that not all the values of $\arg x$
are allowed, namely 
$ C_1-\ln r<\nu \arg(x)$.   
 Our belief is that if we extend the range of $\arg x$, then
 $y(x)$ may have  (movable) poles. We are
 not able to prove it in general, but we will produce 
an example in section \ref{beyond}.
 
\vskip 0.2 cm 

 We finally remark that the critical behaviour of Painlev\'e transcendents can
 also be investigated using  a representation due to S.Shimomura \cite{Sh}
 \cite{IKSY} 
$$
    y(x) = {1\over \cosh^2\left({u_{s}(x)\over 2}\right)} 
$$ 
where $u_{s}(x)$ 
solves a non linear differential equation of the second order. We will discuss
it in chapter \ref{PaInLeVe}. However, the connection problem in this
representation was not solved.

 In the thesis  we give an extended and unified 
picture of both elliptic and Shimomura's representations  and
Dubrovin-Mazzocco's  works, 
and we solve the
 connection problem for elliptic and Shimomura's representations.

\vskip 0.4 cm

{\bf iii)}
 In Chapter \ref{inverse reconstruction} we apply the results on the
critical behaviour of Painlev\'e transcendents to obtain the behaviour of the
parametric representation (\ref{parametricintroduzione}) for $n=3$. 

Let  $F_0(t):= 
{1\over 2}\left[ (t^1)^2 t^3+t^1(t^2)^2\right]$. 
 We prove that for generic
$\mu$ the parametric representation is 
\be
 t^2(u)= \tau_2(x,\mu) ~(u_2-u_1)^{1+\mu},~~~~
  t^3(u)=\tau_3(x,\mu)~(u_2-u_1)^{1+2\mu}
\label{tintroduzione}
\ee
\be
F(u)=F_0(t)+ {\cal F}(x,\mu)~(u_2-u_1)^{3+2\mu}
\label{Fintroduzione},~~~~x={u_3-u_1\over u_2-u_1}
\ee
where $ \tau_2(x,\mu), \tau_3(x,\mu), {\cal F}(x,\mu)$ will be computed 
explicitly as
rational functions of $x$, $y(x)$, ${d y \over dx}$ and $\mu$.
 The ratio $
{t^2\over
(t^3)^{{1+\mu\over 1+2\mu}}}
$ is independent of $(u_2-u_1)$. This is actually the crucial point, because
now  
the closed form $F=F(t)$ must be:
$$
F(t)=F_0(t)+ (t^3)^{{3+2\mu\over 1+2\mu}}\varphi\left({t^2\over
(t^3)^{{1+\mu\over 1+2\mu}}}\right)
$$
where the function $\varphi$ has to be determined by the inversion of
(\ref{tintroduzione}) (\ref{Fintroduzione}).

In the case of $QH^{*}({\bf CP}^2)$ we prove that we just need to take the
limit of the above $t^3$ and $F$ for $\mu\to -1$, but such a limit does not
exist  for $t^2$ and the correct form is 
\be
 t^2(u)=3 \ln(u_2-u_1) +3\int^x
 d\zeta{1\over \zeta +f(\zeta) }
\label{t2introduzione}
\ee
where $f(x)$ is again computed explicitly as a  
rational functions of $x$, $y(x)$, ${d y \over dx}$.  This time  $
   e^{t^2} (t^3)^3 
$ is independent of $(u_2-u_1)$ and so  
$$
 F(t)=F_0(t)+ {1\over t^3} ~\varphi\left( e^{t^2} (t^3)^3 \right)
$$

\vskip 0.2 cm
To our knowledge, this is the first time the {\it explicit}  parameterization
(\ref{tintroduzione}) (\ref{Fintroduzione}) (\ref{t2introduzione}) is given;  
although its proof is  mainly a computational problem (the theoretical problem
being already solved by the reduction to the Painlev\'e 6 eq. \cite{Dub1}),
it is very hard. Moreover, the knowledge of this explicit 
form  is necessary to proceed
to the inversion of the parametric formulae close to the diagonals.

 When the
transcendent behaves like (\ref{loc1introduzione}), or
(\ref{loc2introduzione}), or  (\ref{loc3introduzione}) with rational
exponents,  then $t$ and $F$ in 
(\ref{tintroduzione}) (\ref{Fintroduzione})                 
 are expanded in Puiseux series in  
$x$, $1-x$ of ${1\over x}$ . The expansion can be 
inverted, in order to obtain $F=F(t)$ in closed form 
as an expansion in $t$.  We apply the
procedure staring from the algebraic solutions \cite{DM} 
of (\ref{painleveintroduzione})
and we obtain the polynomial solutions of the WDVV equations. 

We also apply the procedure for $QH^{*}({\bf CP}^2)$. This time, $\Re
\sigma^{(i)}=1$, the transcendent has oscillatory behaviour and 
therefore the reduction of (\ref{t2introduzione}) (\ref{tintroduzione})
(\ref{Fintroduzione}) to closed form is
hard. Hence,  we 
 expand the transcendent in Taylor series close to a regular point
$x_{\hbox{reg}}$, we plug it into  (\ref{t2introduzione}) (\ref{tintroduzione})
(\ref{Fintroduzione})  and we obtain
$t$ and $F$ as a 
Taylor series in $(x-x_{\hbox{reg}})$. We invert the series and
we get a closed form $F=F(t)$. We prove that it is   precisely the 
solution 
(\ref{kontsevichintroduzione}). Thus, our procedure is an alternative way to
compute the numbers $N_k$ as an application of the isomonodromic deformations
theory.  

 We have just started to investigate the possibilities offered by the formulae
  (\ref{tintroduzione})
(\ref{Fintroduzione})  (\ref{t2introduzione}). 
We believe they will be a good tool to
 understand some analyticity properties of $F(t)$  in future
 investigations. Particularly, we hope to better 
understand the connection between
 the monodromy data of the quantum cohomology and the number of rational
 curves.   
This problem will be the
 object of further investigations.

\vskip 0.3 cm 

 The entire procedure at points i), ii), iii) is a significant application of
 the theory of isomonodromic deformations to a problem of mathematical
 physics:  solve the WDVV equations or, at least,  investigate the
 analytic properties of $F(t)$.

\vskip 0.4 cm

This work is organized as follows.  Chapter \ref{cap1} is a review
on Frobenius manifolds where 
we discuss in detail  the parameterization of the
manifold through monodromy data and the reduction to a Painlev\'e 6
equation. This chapter is mainly a synthesis of \cite{Dub1} \cite{Dub2}.

In chapter \ref{cap2} we introduce the quantum
cohomology of projective spaces and its connections to enumerative
geometry. We also propose  a numerical computation which  we did in order to
investigate the nature of the singular point which determines the radius of
convergence of the solution (\ref{kontsevichintroduzione}). 

Chapter \ref{rh2} is a didactic exposition of the reconstruction of a
2-dimensional Frobenius manifold starting from the isomonodromic deformation of
the associated linear system. The 2-dim case is exactly solvable, but it is a
good model for the general procedure in any dimension.

Chapters \ref{stokes}, \ref{PaInLeVe}, \ref{inverse reconstruction} 
contain the main results i), ii), iii) of the
thesis. 

\vskip 0.3 cm
 We would like to spend a word to explain the nature of this work.  
Although many 
 objects studied here come from enumerative and 
algebraic geometry, we never make use of  tools from those
 fields. We only require some knowledge of differential geometry and elementary
 topology. 
All the work
 is mainly analytical. We  use concepts and tools and we face problems 
 of   
 complex analysis,    
  asymptotic expansions, theory of linear systems of differential
 equations, isomonodromic deformations theory, Riemann-Hilbert  boundary
value problems,  Painlev\'e equations.   

\vskip 1 cm
{\it Acknowledgements.} I am grateful to B.Dubrovin for suggesting me the
 problem of this work and for many  discussions, guidance  and advice. 
  I would like to thank many mathematicians who talked to me during these
 years. I'll mention  some of them in different points of this work. 
  In particular,  I would like to thank 
 M. Bertola, S.Bianchini,  A. Its, 
  M.Mazzocco for many discussions.

\tableofcontents

\mainmatter


\chapter{ Introduction to Frobenius Manifolds}\label{cap1}

 This chapter is a review on the theory of Frobenius manifolds. The
  connection between Frobenius manifolds and the theory of isomonodromic
 deformations will be studied, in view of the inverse reconstruction of a
 Frobenius structure starting from a set of monodromy data of a system of
  linear differential equations.

\section{ The WDVV Equations of Associativity}

  The theory of Frobenius manifolds was introduced by B. Dubrovin in 
\cite{Dub4} for the Witten-Dijkgraaf-Verlinde-Verlinde equations of 
associativity (WDVV) \cite{Witten} \cite{DVV}. The WDVV equations are  
differential equations satisfied by the {\it primary free energy} $F(t)$ of  
 two-dimensional topological field theory. $F(t)$  is a function 
of the coupling constants $t:=(t^1,t^2,...,t^n)$ $t^i\in {\bf C}$. Given a 
non-degenerate  symmetric matrix $\eta^{\alpha \beta}$, $\alpha,\beta=1,...,n$,
  $F(t)$ satisfies: 
\be 
\partial_{\alpha}\partial_{\beta} \partial_{\lambda} F~\eta^{\lambda \mu} 
~\partial_{\mu}\partial_{\gamma}\partial_{\delta} F ~=~ \hbox{ the same with 
$\alpha$, $\delta$ exchanged},
\label{WDVV1}
\ee
where $\partial_{\alpha}:={\partial \over \partial t^{\alpha}}$. Sum
 over repeated indices is omitted. The relation 
between $\eta^{\alpha \beta}$ and $F(t)$ is given by 
\be
\partial_1 \partial_{\alpha}\partial_{\beta} F = \eta_{\alpha\beta},
\label{WDVV2}
\ee
where the matrix $
 (\eta_{\alpha \beta})$ is the inverse of the matrix $(\eta^{\alpha \beta})$. 
Finally, $F(t)$ must satisfy the {\it quasi-homogeneity condition}: given 
numbers $q_1,q_2,...,q_n,d$, $r_1,...,r_n$ ($r_{\alpha}=0$ if $q_{\alpha} 
 \neq 1$) 
  we require: 
\be
  E (F(t))= (3-d) F(t) + \hbox{ (at most) quadratic terms}
\label{WDVV3b},
\ee
where $E(F(t))$ means a differential operator $E$ applied to $F(t)$ 
and defined as
follows:  
$$
E= \sum_{\alpha=1}^n E^{\alpha} \partial_{\alpha},~~~~ E^{\alpha}=
       (1-q_{\alpha})t^{\alpha} +r_{\alpha},~~~ \alpha=1,...,n,
$$  
It is  called {\it Euler vector field}.
The equation (\ref{WDVV3b}) is
       the differential form of 
$F(\lambda^{1-q_1} t^1,...,\lambda^{1-q_n}t^n)=\lambda^{3-d} F(t^1,....t^n)
$, where $\lambda\neq 0$ is any complex number. 
If $q_{\alpha}=1$ we must read $r_{\alpha}
       \ln\lambda +t^{\alpha}$ instead of $\lambda^{1-q_{\alpha}}t^{\alpha}$
       in the argument of $F$.

 The equations (\ref{WDVV1}), 
(\ref{WDVV2}), (\ref{WDVV3b}) are the {\it WDVV equations}. 

\vskip 0.2 cm
 
If we define $ c_{\alpha \beta \gamma}(t):= \partial_{\alpha}\partial_{\beta} 
\partial_{\gamma} F(t)$, $c_{\alpha\beta}^{\gamma}(t) := \eta^{\gamma \mu} 
c_{\alpha \beta \mu}(t)$   (sum omitted),
and we consider a vector space $A$=span($e_1,...,e_n$), then we obtain a 
family of algebras $A_t$  with the  
multiplication 
$ 
  e_{\alpha} \cdot e_{\beta} := c_{\alpha \beta}^{\gamma}(t) e_{\gamma}. 
$ 
The algebra is {\it commutative} by definition of $c_{\alpha\beta}^{\gamma}$. 
Equation (\ref{WDVV1}) is equivalent to {\it associativity}.
The vector $e_1$ is the {\it unit element}  because  (\ref{WDVV2}) implies
$c_{1\beta}^{\gamma} =\delta^{\gamma}_{\beta}$. The bilinear form $<.,.>$
defined by
$$ 
   <e_{\alpha},e_{\beta}>= \eta_{\alpha\beta}
$$
is symmetric, non degenerate and {\it invariant}, namely
$<e_{\alpha}\cdot e_{\beta},e_{\gamma}>= <e_{\alpha},
e_{\beta}\cdot e_{\gamma}>  $.

\vskip 0.2 cm

The equation (\ref{WDVV2}) is integrated and yields:
$$
  F(t)= {1\over 6} \eta_{11} (t^1)^3 +{1\over 2} [\sum_{\alpha\neq 1}
  \eta_{1\alpha} t^{\alpha}] (t^1)^2+{1\over 2} 
[\sum_{\alpha\neq 1,\beta\neq 1} \eta_{\alpha\beta} t^{\alpha}t^{\beta} ]~t^1+
f(t^2,...,t^n)
$$
The function $f(t^2,...,t^n)$ is so far arbitrary and it is 
 defined {\it up to quadratic and linear terms} which do not
affect the third derivatives of $F(t)$.

By a linear change of coordinates ${t^{\prime}}^{\alpha} = A^{\alpha}_{\beta}
t^{\beta}$  which does not modify (\ref{WDVV1}), (\ref{WDVV2})
and 
(\ref{WDVV3b}) we can reduce $\eta:=(\eta_{\alpha\beta})$ to the form:
$$
     \eta=\pmatrix{ \eta_{11} & 0 & 0&       & 0 & 0 \cr
                    0         &   &  &       &   & 1 \cr
                    0         &   &  &       & 1 &   \cr
                    \vdots    &   &  &\adots &   &   \cr
                    0         &   &1 &       &   &   \cr
                    0         & 1 &  &       &   &   \cr},
     ~~~~\hbox{ if } 
                     \eta_{11}\neq 0          
$$
\vskip 0.3 cm
$$
     \eta=\pmatrix{   &   &         & & 1 \cr
                      &   &         &1& \cr 
                      &   & \adots  & & \cr
                      & 1 &         & & \cr
                    1 &   &         & &  \cr
                   }  ~~~~\hbox{ if } 
                     \eta_{11}= 0
$$ 
\noindent
All the other entries are zero. 
If $n\geq 2$, we  are going to deal with Frobemius manifolds such that
$$
   q_1=0,~~~~, q_{\alpha}+q_{n-\alpha+1}=d,~~~
     \eta=\pmatrix{   &   &         & & 1 \cr
                      &   &         &1& \cr 
                      &   & \adots  & & \cr
                      & 1 &         & & \cr
                    1 &   &         & &  \cr
                   }
$$
The condition $q_1=0$ is not very restrictive, because if $q_1\neq 1$ we can
 always reduce to $q_1=0$ by rescaling all the $q_{\alpha}$'s and $d$  in the
 quasi-homogeneity condition.  The free energy in this case is 
$$ 
  F(t)= {1\over 2} (t^1)^2 t^n +{1\over 2} t^1 \sum_{\alpha=2}^{n-1}
  t^{\alpha} t^{n-\alpha+1} ~+f(t^2,...,t^n)
$$

\subsection{ Examples} 

\noindent
$\bullet$ $n=1$
$$
  F(t)= {1\over 6} \eta~ (t)^3 +A~(t)^2+B~t+C
$$

\vskip 0.3 cm
\noindent
$\bullet$ $n=2$, $\eta_{11}=0$
$$
  F(t^1,t^2)= {1\over 2} (t^1)^2t^2+f(t^2)
$$
(\ref{WDVV1}) is automatically satisfied.   (\ref{WDVV3b}) is:
\vskip 0.2 cm
\noindent 
- If $d\neq 1$
$$
   (1-d) t^2 {df\over dt^2}= (3-d) f ~~\Rightarrow~~ f(t^2)= C~(t^2)^{3-d\over
   1-d} 
$$
where $C$ is a constant. 
\vskip 0.2 cm
\noindent 
- If $d=1$
$$
   r {df\over dt^2}=2 f~~\Rightarrow~~ f(t^2)=C ~
                             \left[e^{t^2}\right]^{2\over r}
$$
\vskip 0.2 cm
\noindent 
- If $d=1$, $r=0$, then $f=0$.
\vskip 0.2 cm
\noindent 
- For $d=-1$, there is also the solution 
                  $$ f(t^2)= C (t^2)^2 \ln(t^2)$$
\vskip 0.2 cm
\noindent 
- For $d=3$,there is also the solution 
                  $$ f(t^2)= C \ln(t^2)$$

\vskip 0.3 cm
\noindent
$\bullet$ $n=3$, $\eta_{11}=0$
$$
  F(t^1,t^2,t^3)= F_0(t^1,t^2,t^3)+ f(t^2,t^3),~~\hbox{ where } 
  F_0(t^1,t^2,t^3)= {1\over 2} \left[(t^1)^2t^3+t^1(t^2)^2\right]
$$
(\ref{WDVV1}) becomes
\be
   f_{222}f_{233}+f_{333}=(f_{223})^2
\label{WDVV12}
\ee
where the subscripts mean derivatives w.r.t. $t^2$ and $t^3$. The charges are:
$$
 q_1=0,~~~q_2={d\over 2},~~~q_3=d
$$
and (\ref{WDVV3b})
is translated into the following quasi-homogeneity conditions:
\vskip 0.2 cm
\noindent
- If $d\neq 1,2$
$$ 
   f(\lambda^{1-{d\over 2}}t^2,\lambda^{1-d} t^3) =\lambda^{3-d} f(t^2,t^3)
$$
thus
$$
   f(t^2,t^3)= (t^2)^{3-d\over 1-d/2} \varphi(t^2(t^3)^q),~~~q={1-d/2\over
   d-1}
$$
where $\varphi$ is an function of $t^2(t^3)^q$ such that (\ref{WDVV12}) is
satisfied . 

\vskip 0.2 cm
\noindent
- If $d = 2$ 
  $$
      f(t^2+r_2\ln\lambda,{1\over \lambda} t^3) = \lambda f(t^2,t^3)
$$
which implies
$$ 
  f(t^2,t^3)= {1\over t^3} \varphi((t^3)^{r_2} e^{t^2})
$$
 For $r=3$ this is the case of the {\it Quantum Cohomology} of the projective
 space $CP^2$ (see  chapter \ref{cap2}).

\vskip 0.2 cm
\noindent
- If $d = 1$ 
$$ 
  f(\lambda^{1\over 2} t^2,t^3+r_3 \ln \lambda) = \lambda^2 f(t^2,t^3)
$$
thus
 $$
  f(t^2,t^3)= (t^2)^4 \varphi((t^2)^{-2r_3} e^{t^3})
$$

In all the above cases, 
the equation  (\ref{WDVV12}) becomes an O.D.E. for $\varphi$.  
In particular, there are four polynomial solutions. Let $a\in {\bf
C}\backslash \{0\}$: 
\be
  F(t)= F_0(t)+a~ (t^2)^2(t^3)^2 +{4\over 15} a^2~ (t^3)^{5},~~~~d={1\over 2},
\label{A3}
\ee
\be
   F(t)= F_0(t)+a~(t^2)^3t^3+6a^2~(t^2)^2(t^3)^3+{216\over 35} a^4~
   (t^3)^7,~~~~d={2\over 3}, 
\label{B3}
\ee
\be
   F(t)= F_0(t)+a~(t^2)^3(t^3)^2+{9\over 5}a^2~(t^2)^2(t^3)^5 +{18\over 55}
   a^4~ 
   (t^3)^{11},~~~~d={4\over 5}, 
\label{H3}
\ee
$$
  F(t) = F_0(t) + a ~(t^2)^4,~~~~d=1.
$$
 Note that we can fix $a$ (it is an integration constant). This means that
each of the above solutions is considered as {\it one} solution and not as a
one-parameter family. 

For any  positive integer $m$ 
  there are analytic solutions at $t=0$ consisting in   
a one-parameter family if $d=2{m+1\over m+2}$, one solution 
if $d=2{m+2\over m+4}$, one solution if $d=2{m+3\over m+6}$. Note that $d\neq
  1,2$. 

 For $d=1$ there are three analytic solutions at $t=0$ (plus the polynomial
 solution above). We give them for a fixed value of the integration constant
 $a$: 
$$
  F(t)= F_0(t) -{1\over 24} (t^2)^4 +t^2 e^{t^3}, ~~~r_3={3\over 2},
$$
$$
 F(t) =F_0(t)  + {1\over 2} (t^2)^4 +(t^2)^2 e^{t^3} -{1\over 48} e^{3 t^3},
 ~~~ r_3= 1, 
$$
$$
  F(t)=F_0(t) -{1\over 72} (t^2)^4 +{2\over 3} (t^2)^3 e^{t^3} +{2\over 3}
  (t^2)^2 e^{2t^3} +{9 \over 16} e^{4t^3},~~~ r_3={1\over 2}
$$
For $r_3=0$ we also have
$$
  F(t)=F_0(t) -{(t^2)^4\over 16} \gamma(t^3)
$$
where $\gamma(\zeta)$ satisfies
$$ 
    \gamma^{\prime\prime\prime} = 6 \gamma \gamma^{\prime\prime} - 9
    (\gamma^{\prime})^2,~~~(\hbox{here}
    \gamma^{\prime}:=d\gamma(\zeta)/d\zeta). 
$$
 This is the {\it Chazy equation} and
 it has a solution analytic at $\zeta=i\infty$:
$$ 
\gamma(\zeta) =\sum_{n\geq 0}~ a_n e^{2\pi i n \zeta}
$$ 
The Chazy equation determines the $a_n$'s.

For $d=2$, $r_2=3$  there is a solution analytic at $t^2=-\infty$, $t^3=0$,
discovered 
by Kontsevich \cite{KM}:
$$
  F(t)=F_0(t) +{1\over t^3} \sum_{k\geq 0} A_k~ (t^3)^{3k} e^{k t^2}
$$
The $A_k$ are uniquely determined (once the integration constant 
$A_1$ is chosen) and the series converges around $(t^3)^{3} e^{ t^2}=0$. 

\vskip 0.2 cm
We will return later to some of the above solutions.

\vskip 0.3 cm


\section{ Frobenius Manifolds}\label{definizione FM}

Let's consider a smooth/analytic manifold $M$ of dimension $n$ 
over ${\bf C}$, whose tangent space $T_t M$ at any $t \in M$ 
is an {\it associative, commutative algebra} with {\it unit element} $e$, 
equipped
with  a non-degenerate bilinear form $<., .>$. Let denote by $\cdot$ the
product of two vectors;  the bilinear form  is {\it invariant}  w.r.t. the
product, namely $<u\cdot v ,w>=<u,v\cdot w>$ for any $u,v,w\in T_tM$. $T_tM$ is
called a {\it Frobenius algebra} (the name comes from a similar structure
studied by Frobenius in group theory).  

\vskip 0.2 cm 
\noindent
 We further suppose that 

\vskip 0.2 cm
1)  $<.,.>$ is {\it flat}. Therefore, there exist flat 
coordinates $t^1$,...,$t^n$ such that 
        $$
               <\partial_{\alpha},  \partial_{\beta}>=:\eta_{\alpha \beta}~~\hbox{ constant}
$$
where $\partial_{\alpha}={\partial \over \partial t^{\alpha}}$ is a basis. 
We denote by $\nabla$ the Levi-Civita connection. In particular
$\nabla_{\alpha}=\partial_{\alpha}$. 

\vskip 0.2 cm
2) $\nabla e=0$. So  we can choose 
 $e=\partial_1$. 

\vskip 0.2 cm
3)  the tensors  
$c(u,v,w):=<u\cdot v , w>$ and $\nabla_y c(u,v,w)$,  $u,v,w,y\in T_t M$,  are 
symmetric. 

We define
$ c_{\alpha \beta \gamma}(t):=<\partial_{\alpha} \cdot \partial_{\beta},
\partial_{\gamma}>$; the 
symmetry becomes the complete symmetry of $\partial_{\delta}
c_{\alpha \beta \gamma}(t)$ in the indices. This implies the existence 
of a function $F(t)$ such that
$\partial_{\alpha}\partial_{\beta}\partial_{\gamma} F(t)  
= c_{\alpha \beta \gamma}(t)$. $F$  satisfies the equation  (\ref{WDVV1})
because of the associativity of the algebra $T_t M$. 

The 
equation  
(\ref{WDVV2}) follows from  the axiom $\nabla e=0$  and the choice 
 $e=\partial_1$.

\vskip 0.2 cm
4) There exist an {\it Euler vector field} $E$ such that
   $$
      \hbox{ i)}~~~~ \nabla\nabla E=0
$$
   $$  
      \hbox{ ii)} ~~~~ \hbox{Lie}_E ~c = c
$$
  $$  \hbox{ iii})~~~~  
 \hbox{Lie}_E~ e = -e
$$
  $$
     \hbox{iv)}~~~~  
 \hbox{Lie}_E~ <.,.>=(2-d)<.,.>
$$
i) implies that $E^{\gamma}= \sum_{\beta} a^{\gamma}_{\beta} t^{\beta}
+r_{\beta}$ in flat coordinates.
 Assuming {\it diagonalizability} of $\nabla E$ we reduce to 
   $$ 
        E=\sum_{\alpha=1}^n ~\bigl((1-q_{\alpha}) t^{\alpha}+r_{\alpha}\bigr)
        \partial_{\alpha}
$$
iii) implies $q_1=0$; iv) implies
$(q_{\alpha}+q_{\beta}-d)\eta_{\alpha\beta}=0$ and ii) implies
$E(c_{\alpha\beta}^{\gamma})= (q_{\alpha}+q_{\beta}+q_{\gamma}-d)
c_{\alpha\beta}^{\gamma}$ (no sums), and then  $E(F)=(3-d)F +$ quadratic
terms. The condition $(q_{\alpha}+q_{\beta}-d)\eta_{\alpha\beta}=0$ can be put
in the form  
$q_{\alpha}+q_{n-\alpha+1} =d$ if $\eta_{11}=0$.

\vskip 0.2 cm
\noindent{\bf Definition:} The manifold $M$ equipped with such a structure is
called a {\it Frobenius manifold}. 
\vskip 0.2 cm 
 
  In this way, the  WDVV equations are reformulated as geometrical conditions
  on $M$. Frobenius manifolds  arise as geometric 
structures in many branches of mathematics, like enumerative geometry and 
quantun cohomology \cite{KM} \cite{Manin}. Moreover  there is a FM structure 
on the 
space of orbits of a Coxeter group (the solutions (\ref{A3}), (\ref{B3}),
  (\ref{H3}) correspond to the groups $A_3$, $B_3$, $H_3$ respectively) 
and on the universal unfolding of simple 
singularities \cite{Saito} \cite{SYS}  \cite{Dub6}  \cite{Dub1} (see also 
\cite{Bert}). We will return later to the examples.


\section{ Deformed flat connection}

 The connection $\nabla$ can be deformed by a complex parameter $z$. We
introduce  
a {\it deformed} connection $\tilde{\nabla}$  on $M \times {\bf C}$:  
for any $u,v \in T_t M$, 
depending  also on  $z$, we define  
$$ 
    \tilde{\nabla}_u v := \nabla_u v + z u \cdot v,$$ 
$$
   \tilde{\nabla}_{d\over dz}v:= {\partial \over \partial z} v +E\cdot v - 
{1\over z} \hat{\mu}v,
$$ 
$$ 
  \tilde{\nabla}_{d\over dz}{d\over dz} =0,~~~~~~~
 \tilde{\nabla}_u {d\over dz} =0
$$
where 
  $E$ is the Euler vector field and 
$$
\hat{\mu}:= I-{d \over 2} -
\nabla E
$$ 
is an operator acting on $v$. In coordinates $t$: 
$$
\hat{\mu} = \hbox{diag}(\mu_1,...,\mu_n),~~~~\mu_{\alpha}=q_{\alpha} -{d 
\over 2},$$ provided that $\nabla E$ is diagonalizable. From iv) of
section \ref{definizione FM} it follows
that 
$$ 
   \eta \hat{\mu}+\hat{\mu}^T\eta=0
$$

\vskip 0.2 cm
\noindent
{\bf  Theorem} \cite{Dub1}: {\it $\tilde{\nabla}$ is flat}.
\vskip 0.2 cm 

To find a flat coordinate $\tilde{t}(t,z)$ we impose $\tilde{\nabla} d\tilde{t}
=0$, which becomes the linear system
\be
    \partial_{\alpha} \xi = z C_{\alpha}(t) \xi,
\label{systemt}
\ee
\be
     \partial_z \xi = \left[{\cal U}(t) +{\hat{\mu} \over z} \right] \xi,
\label{systemz}
\ee
where 
$ 
   \xi$ is a column vector of components $\xi^{\alpha}=\eta^{\alpha\mu}
\partial \tilde{t} /\partial t^{\mu}$, $\alpha=1,...,n$ (sum omitted), 
and $ C_{\alpha}(t)=
\bigl(
                                  c_{\alpha \gamma}^{\beta}(t)\bigr)$, $
{\cal U}:= \bigl( E^{\mu} c_{\mu \gamma}^{\beta}(t)  \bigr).
$ From the definition we have ${\cal U}^T \eta=\eta{\cal U}$. 
 The compatibility of the system is equivalent to the fact that the curvature
for $\tilde{\nabla}$ is zero. 

 We stress that the deformed connection is a natural structure on a Frobenius
 manifold. Roughly speaking, a manifold with a flat
 connection is a Frobenius manifold if the deformed connection is flat. More
 precisely, suppose that 
 $M$ is a (smooth/analytic) manifold such that $T_tM$ is a
 commutative algebra with unit element and a bilinear form $<.,.>$
 invariant w.r.t. the product. Take an arbitrary vector field $E$ and define
 the deformation $\tilde{\nabla}$ as above. Then the following is true: 

\vskip 0.2 cm
{\it $\tilde {\nabla}$ is flat if and only if $\partial_{\delta}
c_{\alpha\beta\gamma}$ is completely symmetric, $\nabla \nabla E=0$, the
product is associative and $\hbox{Lie}_E~ c = c$.}

\vskip 0.2 cm
\noindent
This statement simply translates the conditions of compatibility of
(\ref{systemt}) (\ref{systemz}). 

 \vskip 0.2 cm

Therefore if $M$ is Frobenius, then $\tilde{\nabla}$ is flat. 
Conversely, if $\tilde{\nabla}$ is flat, then 
 $M$ is a Frobenius manifold (with
$q_1=0$) provided that we  also require iii), iv)  of
section \ref{definizione FM} and $\nabla(e)=0$.


\section{ Semisimple Frobenius
 manifolds}\label{Semisimple Frobenius manifolds}  

 \noindent
{\bf Definition:} A commutative, associative algebra $A$ 
with unit element is {\it semisimple} if there is no element $a\in A$ such
that $a^k=0$ for some $k\in{\bf N}$.

\vskip 0.2 cm
{\it Any semisimple Frobenius algebra of dimension $n$ over ${\bf C}$ is
isomorphic to ${\bf C} \oplus {\bf C} \oplus ...{\bf C}$ $n$ times. The 
direct sum is orthogonal w.r.t. $<.,.>$. }

\vskip 0.2 cm 
 Therefore, $A$ has a basis $\pi_1$,...,$\pi_n$ such that
$$
   \pi_i \cdot \pi_j = \delta_{ij} \pi_i ~~~\hbox{ no sum}
$$
$$ 
  <\pi_i,\pi_j>= \eta_{ii} \delta_{ij}
$$
$\pi_1$,...,$\pi_n$ are called {\it idempotents}, determined up to permutation.  

\vskip 0.3 cm
\noindent
{\bf Definition:} A Frobenius manifold is {\it semisimple} if $T_tM$ is
semisimple at generic $t$. Such $t$ is called a semisimple point. 
\vskip 0.2 cm 

 If $A$ is semisimple, the vector 
${\cal E}:=\sum_i u_i \pi_i$, $u_i\in {\bf C}$, $u_i \neq u_j$ for $i\neq j$,  
 has $n$ distinct eigenvalues 
$u_i$ with eigenvectors $\pi_i$. Conversely, consider an algebra admitting 
 a  vector 
  ${\cal E}$ having distinct eigenvalues $u_1$,...,$u_n$ and eigenvectors 
$e_1$,...,$e_n$; from associativity and commutativity it follows that $({\cal
 E}\cdot 
 e_i)\cdot e_j= e_i\cdot({\cal E}\cdot e_j)$, namely $(u_i-u_j)e_i\cdot e_j=0
 \Rightarrow e_i\cdot e_j=0$ for $i\neq j$. Thus the algebra is
 semisimple. We have proved that: 
\vskip 0.2 cm
{\it 
    $A$ is semisimple if and only if  there exists a vector ${\cal E}$ such
    that the 
multiplication ${\cal E} \cdot$ has $n$ distinct eigenvalues.}
\vskip 0.2 cm

If follows that semisimplicity is an ``open property'' in a Frobenius
manifold $M$: 
if $t_0\in M$ is semisimple, $T_tM$ is still semisimple for any $t$
is a neighbourhood of $t_0$. 

\vskip 0.2 cm
\noindent
{\bf Theorem} \cite{Dub1}:
{\it Let $t\in M$ be a semisimple point and $\pi_1(t)$,...,$\pi_n(t)$ a basis
of idempotents 
in $T_tM$. The commutator $[\pi_i,\pi_j]=0$ and there exist local coordinates
$u_1,...,u_n$ 
around $t$ such that 
  $$ \pi_i= {\partial \over \partial u_i} 
$$
}
\vskip 0.2 cm
 The local coordinates $u_i$ are determined up to shift and permutation. 
Let ${\cal S}_n$ be the symmetric group of $n$ elements. Let  ${\bf C}^n
\backslash \hbox{ diagonals}:= \{(u_1,...,u_n)\in {\bf C}^n \hbox{ such that } 
                       u_i\neq u_j \hbox{ for } i\neq j\}
$. Finally, let $[(u_1(t),...,u_n(t))]$ be an equivalence class in $ 
    {\cal M} \to {({\bf C}^n \backslash \hbox{ diagonals} )\over {\cal S}_n}
$. 
\vskip 0.2 cm
\noindent
{\bf Theorem:} {\it 
Let ${\cal U}(t)$ be the matrix of multiplication by the Euler
vector field. Let 
$$
 {\cal M}:= \{ t\in M \hbox{ such that } \det({\cal U}(t)-\lambda)=0
 \hbox{ has $n$ distinct eigenvalues } u_1(t),...,u_n(t)\}
$$
The map 
$$ 
    {\cal M} \to {({\bf C}^n \backslash \hbox{ diagonals} )\over {\cal S}_n}
$$
defined by $    t \mapsto [(u_1(t),...,u_n(t))]$ 
is a local diffeomorphism, i.e. $u_1(t),...,u_n(t)$ are local coordinates.
 Moreover
 $$
{\partial \over \partial u_i(t)} \cdot {\partial \over \partial u_j(t)}
   = 
     \delta_{ij} {\partial \over \partial u_i(t)},~~~~
   E(t)= \sum_{i=1}^n u_i(t) {\partial \over \partial u_i(t)}
$$
}

  \vskip 0.2 cm

A point in ${\cal M}$ is semisimple, but in principle there may be semisimple
points $t$ such that  $u_i(t)=u_j(t)$. $u_1$,...,$u_n$ are called {\it
canonical coordinates}. They are defined up to permutation. 

\vskip 0.3 cm 

 We introduce the orthonormal basis $f_i:= {\pi_i \over \sqrt{<\pi_i,\pi_i>}
}$ and we define the matrix $ 
  \Psi= (\psi_{i\alpha})$ by
$$ 
   \partial_{\alpha} = \sum_{i=1}^n ~\psi_{i\alpha} f_i.
$$
 From definitions and simple computations it follows that
$$
\sqrt{<\pi_i,\pi_i>}= 
\psi_{i1},~~~~ <\partial_i,\partial_j>=\delta_{ij} \psi_{i1}^2,
$$
$$
   \partial_{\alpha} = \sum_{i=1}^n~ {\psi_{i\alpha}\over \psi_{i1}}
   \partial_i,~~~~
    \partial_i= \psi_{i1} \sum_{\alpha=1}^n \psi_{i\alpha} \eta^{\alpha\beta}
    \partial_{\beta}, 
$$
$$
   c_{\alpha\beta\gamma}= \sum_{i=1}^n
   {\psi_{i\alpha}\psi_{i\beta}\psi_{i\gamma} \over \psi_{i1}},
$$
 where    $\partial_i:={\partial\over \partial u_i}$.
Furthermore, the relation 
$$
    \Psi^T \Psi=\eta,
$$
holds true (we will prove it later). Finally we define,  
 $$
  U:= \hbox{ diag}(u_1(t),...,u_n(t))= \Psi {\cal U}(t) \Psi^{-1},~~~~
   V(u):= \Psi \hat{\mu} \Psi^{-1}
$$
$V(u)$  is skew-symmetric,  because $\eta \hat{\mu} +\hat{\mu}\eta=0$.

\vskip 0.3 cm

In the following we 
 restrict to  analytic {\it semisimple} Frobenius manifolds. The 
matrix ${\cal U}$ can be diagonalized with {\it distinct eigenvalues} on the  
open dense subset ${\cal M}$
 of $M$. 
Later it  will be convenient to introduce an alternative notation for $\Psi$,
 namely we will often denote $\Psi$ with the name 
  $\phi_0$.  
The systems (\ref{systemt}) and (\ref{systemz})  become:
\be
{\partial y \over \partial u_i}=\left[ zE_i+V_i\right]~y
\label{systemt1}
\ee
\be
{\partial y \over \partial z} = \left[U +{V\over z}\right] y,
\label{systemz1}
\ee
where the row-vector $y$ is $y:=\phi_0 ~\xi$ and $E_i$ is a 
diagonal matrix such that
$(E_i)_{ii}=1$ and all the other entries are 0; we have also  defined
  $$
V_i:= {\partial
\phi_0\over \partial u_i}~\phi_0^{-1}
$$ 
 The compatibility of the two systems is equivalent to  
$$
  \left[U,V_k\right]=\left[E_k,V\right] ~~~\Longrightarrow~~~
(V_k)_{ij}={\delta_{ki}-\delta_{kj}\over u_i-u_j}~V_{ij} 
$$
$$
  {\partial V \over \partial u_i} = \left[V_i,V \right]
$$
At $z=0$ we have a {\it fundamental matrix solution}  
\be 
  Y_0(z,u)= \left[\sum_{p=0}^{\infty}\phi_p(u)~z^p \right] 
z^{\hat{\mu}} z^{R},~~~~\phi_0(u)=\Psi(u)
\label{Y0}
\ee
  where $ {R}_{\alpha \beta}=0$  if $ \mu_{\alpha}-\mu_{\beta} 
\neq k>0$, $k\in {\bf N}$. We will return later to this point.  
The compatibility of the systems implies ${\partial R \over \partial u_i}=0$
  (we discuss this point in section \ref{Monodromy
data of a Semisimple Frobenius Manifold}) 
and then, by plugging $Y_0$ into (\ref{systemt1}) we get 
\be
   {\partial \phi_p, \over \partial u_i}= E_i \phi_{p-1} +V_i \phi_p
\label{ISOmonodromy}
\ee
Finally, let $\Phi(z,u):=\sum_{p=0}^{\infty}\phi_p(u)~z^p$. The condition
$\Phi(-z,u)^T \Phi(-z,u) =\eta$ holds (we prove it is section \ref{Monodromy
data of a Semisimple Frobenius Manifold}) and it implies
\be
\phi_0^T\phi_0=\eta,~~~~
  \sum_{p=0}^m \phi_{p}^T \phi_{m-p}=0 ~~~~\hbox{ for any } m>0
\label{SYMMEtria}
\ee

\vskip 0.3 cm
\noindent
$\bullet$ 
 Conversely, if we start from a (local) solution $V(u)$, $\Psi(u)$ of
\be 
  {\partial V \over \partial u_i}=[V_i,V],
\label{localsystem1}
\ee
\be
  {\partial \Psi \over \partial u_i} = V_i \Psi
\label{localsystem2}
\ee
 such that 
$$ 
  \hbox{det}\Psi(u) \neq 0,~~~~\Pi_{i=1}^n \psi_{i1}(u) \neq 0.
$$ 
for $u$ in a neighbourhood of some $u_0$, we can reconstruct a Frobenius
manifold locally from the formualae 
$$
  \eta=  \Psi^T \Psi,
$$
$$
   \partial_{\alpha} = \sum_{i=1}^n~ {\psi_{i\alpha}\over \psi_{i1}}
   \partial_i, ~~~~
    \partial_i= \psi_{i1} \sum_{\alpha=1}^n \psi_{i\alpha} \eta^{\alpha\beta}
    \partial_{\beta}, 
$$
$$
  <\partial_i,\partial_j>=\delta_{ij} \psi_{i1}^2,~~~~
   c_{\alpha\beta\gamma}= \sum_{i=1}^n
   {\psi_{i\alpha}\psi_{i\beta}\psi_{i\gamma} \over \psi_{i1}}.
$$
$$
  E=\sum_i u_i \partial_i,~~~~
e={\partial \over \partial t^1}=\sum_{i=1}^n {\partial \over \partial u_i}
$$
 The structure of $V(u)$ is as follows: we  fix $V(u_0)=V_0$ at
 $u_0$. Then we take an invertible solution $\Psi$ of (\ref{localsystem2}) and
 we 
 define 
$$  
 V(u)= \Psi(u)~ V_0~ \Psi(u)^{-1}. 
$$
 This solves the equation for $V$, 
as it is verified by direct substitution. It 
 is the unique solution such that the initial value is $V_0$. Let $\hat{\mu}$
 be its Jordan form and  let $\mu_1$ be the first eigenvalue (i.e. the first
 column of $\hat{\mu}$ is $(\mu_1,0,0,...,0)^T$.
$$ 
    V_0=C~\hat{\mu} ~C^{-1}
$$
where $C$ is an invertible matrix independent of $u=(u_1,..,u_n)$. Now we
observe that we can re-scale $\Psi(u) \mapsto \Psi(u) ~C$ and therefore a
solution of (\ref{localsystem1}) is 
$$
     V(u)=\Psi(u) ~\hat{\mu}~ \Psi(u)^{-1}
$$
The dimension of the manifold is 
   $$ d:=-2\mu_1$$
This corresponds to the choice of the first column
$(\psi_{11},...,\psi_{n1})^T$ 
of $\Psi(u)$.


\section{ Inverse Reconstruction of a Semisimple FM (I)}

In this section we show that it is possible to construct a local parametric
  solution of the WDVV equations in terms of the  coefficients
 of (\ref{Y0}). The result is discussed in \cite{Dub2} and it is the main
  formula which allows to reduce the problem of solving the WDVV
 equations to problems of isomonodromic deformations of linear systems of
 differential equations.  In chapter \ref{inverse reconstruction} we will 
present a more explicit
  (and computable!) 
  parametric solution of the WDVV eqs. for $n=3$  as a consequence of the 
results of this  section.  
 
As a first step we note that the condition $\tilde{\nabla} d\tilde{t}=0$ is 
satisfied both
 by a flat coordinate $\tilde{t}^{\alpha}$ and by 
 $\tilde{t}_{\alpha}:=\eta_{\alpha\beta}
 \tilde{t}^{\beta}$. Thus, we choose a fundamental matrix solution of
 (\ref{systemt}), (\ref{systemz}) of the form: 
 $$ 
    \Xi= \left( \partial^{\alpha}\tilde{t}_{\beta} \right) \equiv \left( 
                                  \eta^{\alpha \gamma} {\partial
                                  \tilde{t}_{\beta}\over \partial 
t^{\gamma}} \right)=
 \left[\sum_{p=0}^{\infty} H_p(t) z^p\right] ~z^{\hat{\mu}}~z^R,
~~~~~H_0=I,  $$
close to $z=0$. 
If we restrict to the system (\ref{systemt}) only, we can choose as a
fundamental solution 
$$ 
  H(z,t):= \sum_{p=0}^{\infty} H_p(t) z^p
$$
and so the flat coordinates of $\tilde{\nabla}$ {\it on} $M$ (not on
$M\times {\bf C}$) are 
$$
  \tilde{t}_{\alpha}= \sum_{p=0}^{\infty} h_{\alpha,p}(t) z^p$$ 
where 
\be
   h_{\alpha 0}= t_{\alpha} \equiv \eta_{\alpha \beta} t^{\beta} 
\label{ing1}
\ee
\be
  \partial_{\gamma}\partial_{\beta} h_{\alpha,p+1}= c_{\gamma \beta}^{\epsilon}
  \partial_{\epsilon} h_{\alpha,p} ,~~~~~p=0,1,2,... 
\label{ing2}
\ee
 We stress that the normalization $H_0=I$ is precisely what is necessary to
 have $\tilde{t}_{\alpha}(z=0) =t_{\alpha}$ and it corresponds exactly to
 $Y= \phi_0 \Xi$ in (\ref{Y0}). 

Observe that $ h_{\alpha 0}= t_{\alpha} \equiv \eta_{\alpha \beta} t^{\beta}$
implies 
\be \partial_{\beta} h_{\alpha ,0}=\eta_{\beta \alpha}\equiv c_{\beta \alpha 1}
\label{ing3}
\ee
 Denote by  $\nabla f :=
 \left(\eta^{\alpha \beta} \partial_{\beta} f\right) \partial_{\alpha}$  the
 gradient of the function $f$.  The following are a choice for the flat
 coordinates and for the solution of the WDVV equations
\be
t_{\alpha} = < \nabla h_{\alpha,0} ,\nabla h_{1,1}> \equiv \eta^{\mu \nu}
 \partial_{\mu}   h_{\alpha,0}\partial_{\nu}  h_{1,1}
\label{t}
\ee
$$
F(t)= {1\over 2}\left[ <\nabla h_{\alpha,1},\nabla h_{1,1}> \eta^{\alpha
\beta} <\nabla h_{\beta,0},\nabla h_{1,1}> \right.$$
\be
 \left.- < \nabla h_{1,1},\nabla h_{1,2}>-
<\nabla h_{1,3},\nabla h_{1,0}>\right]
\label{F}
\ee
To prove it, it is enough to check by direct differentiation that 
$\partial_{\alpha} t_{\beta}= \eta_{\alpha \beta} $ and $\partial_{\alpha}
\partial_{\beta} \partial_{\gamma} F(t) = c_{\alpha \beta \gamma}(t)$,  
 using (\ref{ing1}), (\ref{ing2}), (\ref{ing3}) and... some patience. 

In the following, we denote the entry $(i,j)$ of a matrix $A_k$ by
$A_{ij,k}$. We recall that $\Psi\equiv \phi_0$ and we observe that 
  $$ 
  \partial_{\mu} = \sum_{i=1}^n {\phi_{i\mu,0} \over \phi_{i1,0} }~ \partial_i 
,~~~~~ Y_{i\alpha}= {1\over \phi_{i1,0}} ~\partial_i\tilde{t}_{\alpha} 
$$
where $\partial_i= {\partial \over \partial u_i}$. From this it follows that 
$$
  {1\over \phi_{i1,0}}\partial_i h_{\alpha,p} = \phi_{i\alpha,p}$$
and thus:
\be
t_{\alpha}(u)= \sum_{i=1}^n \phi_{i\alpha,0} \phi_{i1,1}
\label{tu}
\ee
$$
F(t(u))= {1\over 2} \left[\eta^{\alpha\beta}\sum_{i=1}^n \phi_{i\alpha,1}\phi_{i1,1}
\sum_{j=1}^n \phi_{j\beta,0} \phi_{j1,1} - \sum_{i=1}^n
\phi_{i1,1}\phi_{i1,2}-\sum_{i=1}^n \phi_{i1,3}\phi_{i1,0} \right]
$$
Equivalently
\be
F(t(u))={1\over 2} \left[t^{\alpha}t^{\beta}\sum_{i=1}^n
\phi_{i\alpha,0}\phi_{i\beta,1} -\sum_{i=1}^n
\left(\phi_{i1,1}\phi_{i1,2}+\phi_{i1,3}\phi_{i1,0}\right) \right]
\label{Fu}
\ee

\vskip 0.3 cm
 It is now clear that we can locally reconstruct a Frobenius manifold from
$\phi_0(u)$, $\phi_1(u)$, $\phi_2(u)$, $\phi_3(u)$ of (\ref{Y0}). 
It is enough to know the local solutions $\phi_0(u)\equiv\Psi(u)$ and
$V(u)$ of (\ref{localsystem1}) and (\ref{localsystem2}) in order to construct
the system
(\ref{systemz1}) and is fundamental matrix (\ref{Y0}). 

\vskip 0.2 cm
 The equations (\ref{localsystem1}) and (\ref{localsystem2}) express the fact
that the dependence on $u=(u_1,...,u_n)$ of the coefficients system 
(\ref{systemz1}) is {\it isomonodromic}. We are going to describe this
property and to show how it is possible to characterize locally the
matrix coefficient $V(u)$ and the matrix $\phi_0(u)$ of such a system
from its monodromy data. We will also explain that 
the problem of solving (\ref{localsystem1}) and
(\ref{localsystem2}) is equivalent to solving a {\it boundary value problem}.

 As a consequence, any analytic semisimple Frobenius manifold can be locally
 parametrized by a set of monodromy data.


\section{ Monodromy and Isomonodromic deformation of a linear
 system}\label{secdeformata} 

 We briefly review some basic notions about the monodromy of a solution of a
 linear system of differential equations and about isomonodromic
 deformations. We make use of the results of \cite{BJL1} and \cite{JMU}. 

\subsection{ Monodromy} 

 Consider a system of differential equations whose coefficients are $N\times
 N$ matrices, meromorphic on ${\bf C}$ with a finite number of poles
 $a_1,...,a_n,\infty$. By Liouville theorem, the most general form of such a
 system is
$$ 
    {dY\over dz}=A(z)~Y,
$$
$$
  A(z)= z^{r_{\infty}-1}\left[ A_0^{(\infty)}+ {A_1^{(\infty)}\over z} +...+ 
                                 { A_{r_{\infty}}^{(\infty)}\over
                                 z^{r_{\infty}}}\right] +
$$
$$
   + \sum_{k=1}^n \left\{   {1\over (z-a_k)^{r_k+1}} \left[
A_0^{(k)}+A_1^{(k)}(z-a_k)+...+A_{r_k}^{(k)}(z-a_k)^{r_k+1}
 \right]\right\}   
$$
$$
   A_0^{(\infty)} ~\hbox{ diagonal},
~~~~r_{\nu}\geq 0 \hbox{ integers }, ~~\nu=1,...,n,\infty.
$$ 
 In the case $r_{\nu}\geq 1$ we suppose that 
$A_0^{(\nu)}$ has {\it distinct
eigenvalues}.  For $z$ sufficiently big we can find $G(z)$ holomorphic at
$\infty$ such that 
$$
   G^{-1}(z) A(z)G(z)= z^{r_{\infty}-1}\left[\Lambda_0^{(\infty)}+{\Lambda_1^{(\infty)}\over z} + 
                                 +...+{\Lambda_{r_{\infty}}^{(\infty)}\over
                                 z^{r_{\infty}}}+...\right]  
   ~~\hbox{ diagonal },~~~\Lambda_0^{(\infty)}=A_0^{(\infty)}
$$
 For $z-a_k$ sufficiently small we can find $G_k(z)$ holomorphic at
$a_k$ such that 
$$
   G_k^{-1}(z) A(z)G_k(z)= {1\over (z-a_k)^{r_k+1}} \left[
\Lambda_0^{(k)}+\Lambda_1^{(k)}(z-a_k)+...+\Lambda_{r_k}^{(k)}(z-a_k)^{r_k+1}
+... \right],~~~\hbox{ diagonal}
$$
The matrices $\Lambda_l^{(\nu)}$, $\nu=1,...,n,\infty$ are uniquely determined
by $A(z)$.

\vskip 0.3 cm
\noindent
$\bullet$ {\bf Representation of the solutions}
\vskip 0.2 cm 

It is a standard result that if $z_0$ is none of the singular points 
 $a_1,...,a_n,\infty$, there exists an invertible matrix $Y(z)$ holomorphic in
 a neighbourhood of $z_0$ which solve the system. It is called a {\it
 fundamental solution}.  Any other fundamental solution in the neighbourhood
 of $z_0$ is  $Y(z)~M$, where $M$ is an invertible matrix independent of $z$. 

   As it is well
 known, $Y(z)$ has analytic continuation along any path $\sigma$ 
 not containing the 
 singular points. The analytic continuation depends on the homotopy class of
 the path in  $CP^1\backslash\{a_1,...,a_n,\infty\}$. We fix a  base point
 $z_0$ in
 $CP^1\backslash\{a_1,...,a_n,\infty\}$ and a base of loops
 $\gamma_1,...,\gamma_n$  in the
 fundamental group $\pi(CP^1\backslash\{a_1,...,a_n,\infty\};z_0)$, starting
 at $z_0$ and encircling  $a_1,...,a_n$. Consider a loop $\gamma$:  
the solution $Y(z)$ obtained by analytic continuation along a path $\sigma$
 and the solution $Y^{\prime}(z)$  obtained by analytic continuation along
 $\gamma \cdot \sigma$ are connected by a constant matrix $M_{\gamma}$. Namely
 $Y^{\prime}(z)= Y(z)~ M_{\gamma}$. $M_{\gamma}$ is called  {\it monodromy
 matrix} for the loop $\gamma$.  Observe that $\gamma\mapsto M_{\gamma}$ is an anti-homomorphism. 
The map 
$$
    \pi(CP^1\backslash\{a_1,...,a_n,\infty\};z_0) \to GL(N,{\bf C})
$$
 is called {\it monodromy representation}.  All the monodromy matrices are
 obtained as products of the matrices $M_1,...,M_n$ corresponding to
 $\gamma_1,...,\gamma_n$ (note that
 $\gamma_1\gamma_2...\gamma_n=\gamma_{\infty}$ $\Rightarrow$ $ M_{\infty}= M_n
 ...M_2M_1$).

 \vskip 0.2 cm
\noindent 
{\it Remark:} i) If  $\pi(CP^1\backslash\{a_1,...,a_n,\infty\};z_0)$ is fixed
but we change normalization from $Y(z)$ to $Y(z)~C$, $\det C\neq 0$, 
the monodromy matrices
change by conjugation $M_i\mapsto C^{-1} M_i C$. 

ii) The monodromy matrices change by conjugation also 
 if we change the base point and the basis of loops. 

Therefore $A(z)$ determines the monodromy representation up to conjugation.

\vskip 0.3 cm     
 
In the following, a fundamental solution $Y(z)$ will be a a branch, 
namely the analytic continuation  of a fundamental solution 
  defined in a neighbourhood of $z_0$ along a path from $z_0$ to $z$.   

  We can
 choose fundamental solutions $Y^{(\infty)}_l(z)$ such they have a 
specific asymptotic
 behaviour on the sectors    
$$
  {\cal S}_l^{(\infty)}= \{ |z|>R~\hbox{ such that } {\pi\over
  r_{\infty}}(l-1) -\delta<\arg(z)< {\pi\over
  r_{\infty}}l \},~~~l=1,...,2r_{\infty}
$$
defined for $R$ sufficiently big and
 for $\delta $ sufficiently small. The behaviour is 
$$
  Y_l^{(\infty)}(z)=F_l^{(\infty)}(z) ~e^{T^{(\infty)}(z)}
$$
   $$
       F_l^{(\infty)}(z)\sim I +{F_1^{(\infty)}\over z}+{F_2^{(\infty)}\over
       z^2} + ...~~~~~~z\to\infty \hbox{ in } {\cal S}_l^{(\infty)},
$$
$$
    T^{(\infty)}={\Lambda_0^{(\infty)}\over r_{\infty}}~z^{r_{\infty}}
            + {\Lambda_1^{(\infty)}\over r_{\infty}-1}~z^{r_{\infty}-1}+...
        + 
          \Lambda_{r_{\infty}-1}^{(\infty)} z + \Lambda_{r_{\infty}} \ln(z)
$$
 In the following we denote $Y_1^{(\infty)}$ with the name $Y(z)$:
$$
      Y(z):=Y_1^{(\infty)}(z). 
$$

 For $\epsilon$ small we   consider the sectors 
$$
 {\cal S}_l^{(k)}=  \{ |z-a_k|<\epsilon~\hbox{ such that } {\pi\over
  r_k}(l-1) -\delta<\arg(z-a_k)< {\pi\over
  r_k}l \},~~~l=1,...,2r_k
$$
where we can choose fundamental solutions
$$
 Y_l^{(k)}=F_l^{(k)}(z) ~e^{T^{(k)}(z)}   
$$
$$
F_l^{(k)}(z)\sim F_0^{(k)}+F_1^{(k)}(z-a_k)+F_1^{(k)}(z-a_k)^2+...,~~~~z\to
a_k \hbox{ in }  {\cal S}_l^{(k)}
$$
$$
T^{(k)}(z)= -\left({\Lambda_0^{(k)}\over r_k (z-a_k)^{r_k}}+ 
                    {\Lambda_1^{(k)}\over (r_k-1) (z-a_k)^{r_k-1}}+...
             +   {\Lambda_{r_k-1}^{(k)}\over  (z-a_k)}
   \right)+ \Lambda_{r_k}^{(k)}\ln(z-a_k)
$$

In the case $r_k=0$, we suppose that the $A_0^{(k)}$ are diagonalizable, with
eigenvalues $\lambda_1^{(k)},..., \lambda_m^{(k)}$. 
Therefore  $\exp\{ T^{(k)}(z)\}=(z-a_k)^{\Lambda_0^{(k)}}~(z-a_k)^{R^{(k)}}$,
where the entry $R^{(k)}_{ij}\neq 0$ only if $\lambda_i^{(k)} -\lambda_j^{(k)}$
is an integer greater than zero. In the same way, if $r_{\infty}=0$, 
$\exp\{ T^{(\infty)}(z)\}= z^{\Lambda_0^{(\infty)}}~z^{R^{(\infty)}}$, where
$R^{(\infty)}_{ij}\neq 0$ only if $\lambda_j^{(\infty)} -\lambda_i^{(\infty)}$ 
is an integer greater than zero. 
 The  general case of non diagonalizability may be treated  with more
 technicalities, which are not necessary here.
If $r_{\nu}=0$, 
the asymptotic series above are convergent in a neighbourhood of the
points $a_{\nu}$, $\nu=1,...,n,\infty$ ($a_{\infty}:=\infty$). 

The matrices $F_l^{(\nu)}$, $l=0,1,2,...$ and $\nu=1,...,n,\infty$, are
uniquely constructed from the coefficients of  $A(z)$ by direct substitution of
the formal series into the system (note that $A(z)$ must be expanded close to
$a_{\nu}$ [where $a_{\infty}:=\infty$] before substituting).

\vskip 0.3 cm
\noindent
$\bullet$ {\bf Stokes' rays}

\vskip 0.2 cm
 Let $r_{\infty}\geq 1$ and 
 $\lambda_1^{(\infty)}$,..., $\lambda_n^{(\infty)}$ be the
 {\it distinct}  eigenvalues of $A_0^{(\infty)}$. We define the {\it Stokes rays} to
 be the the half lines where  the real part of 
$(\lambda_i^{(\infty)} -\lambda_j^{(\infty)}) z^{r_{\infty}}$ is
zero. Therefore, if 
$$\lambda_i^{(\infty)} -\lambda_j^{(\infty)}= |\lambda_i^{(\infty)}
-\lambda_j^{(\infty)} |
~\exp\{ i\alpha_{ij}\},
$$ 
the Stokes' rays are 
$$
 R_{ij,h}:= \left\{z=\rho e^{i\theta_{ij,h}}\hbox{ such that } \theta_{ij,h}={1\over
  r_{\infty} }\left({\pi \over 2} -\alpha_{ij}\right)+h~{\pi \over
  r_{\infty}} \right\},~~h=0,1,...,2r_{\infty}-1
      $$
The following properties are easily proved from the very definition of Stokes
rays:

\vskip 0.2 cm
i) If a fundamental solution has the asymptotic behaviour we gave above 
as $z\to \infty$ in some sector $\alpha <\arg (z) <\beta$, then the the
solution has the same behaviour on the extension of the sector up to the
nearest stokes rays (not included). 

\vskip 0.2 cm
ii) 
    If a solution has the asymptotic behaviour above, for $z\to \infty$, in a
    sector of angular amplitude greater than ${\pi \over r_{\infty}}$, then
    such solution is unique. 

 \vskip 0.2 cm

 As a consequence, we can extend the sectors ${\cal S}_l^{(\infty)}$ up to the
 nearest Stokes rays.
 Suppose we count  the distinct rays in clockwise (or counter-clockwise) order, so that
  we can label the distinct rays as $R_1$, $R_2$,... etc (if some rays
 coincide, they have the same label).
This means that the maximal extension of  ${\cal
 S}_l^{(\infty)}$ goes from one ray $R_{m}$ to the ray $R_{m+1}$ plus an angle
 ${\pi \over
 r_{\infty}}$ (rays are not included in the sector!). In the following, by
 ${\cal S}_l^{(\infty)}$ we mean the already extended sector. 

 Since two fundamental matrices are connected by the multiplication of an
 invertible matrix to the right, we have
$$
  Y_{l+1}^{(\infty)}(z)= Y_l^{(\infty)}(z)~S_l^{(\infty)},~~~z\in {\cal
  S}_l^{(\infty)} \cap {\cal
  S}_{l+1}^{(\infty)}
  ~~~~~~l=1,...,2r_{\infty} -1
$$
\vskip 0.2 cm
$$
  Y_1(z) = Y_{2r_{\infty}}(z) ~S_{2r_{\infty}}^{(\infty)},~~~z\in {\cal
  S}_{2r_{\infty}}^{(\infty)} \cap {\cal
  S}_{1}^{(\infty)}
$$
 The matrices $S_1^{(\infty)}$,...,$S_{2r_{\infty}}$ are called {\it Stokes'
 matrices}. 

\vskip 0.2 cm
The same holds true at any $a_k$, with obvious modification of the above
definitions; therefore we have a set of Stokes' matrices
$S_1^{(\nu)}$,...,$S_{2r_{\nu}}^{(\nu)}$, $\nu=1,...,n,\infty$. 
It's not our purpose to give a
full description of the structure of these matrices. We will construct
explicitly some Stokes matrices later. Here we just remark that

\vskip 0.2 cm
  a) The element on the diagonals are all equal to 1, except for 
     $$
        \hbox{diag}(S_{2r_{\nu}}^{(\nu)})= \exp\{ -2\pi i
        \Lambda_{r_{\nu}}^{(\nu)} \},~~~~\nu=1,2,...,n,\infty
$$

\vskip 0.2 cm
  b) If the entry $(S_l^{(\nu)})_{ij}\neq 0$, the $(S_l^{(\nu)})_{ji}= 0$. 

\vskip 0.2 cm
 
The connection between $Y(z):=Y_1^{(\infty)}(z)$ and the other
solutions $Y_1^{(k)}(z)$ is again given by invertible matrices $C^{(k)}$: 
$$
   Y(z)=\left\{\matrix{ Y_1^{(\infty)}(z) \cr \cr
                        Y_1^{(k)}(z) ~C^{(k)},~~~k=1,...,n \cr
                   }\right.
$$  
The above formula is to be intended as the analytic continuation
of $Y=Y_1^{(\infty)}$ along a path from a neighbourhood of $\infty$ to $z$ in a
neighbourhood of $a_k$ (imagine the path passing through the base-point
$z_0$!). 

\vskip 0.3 cm
\noindent 
$\bullet$ {\bf Monodromy}
\vskip 0.2 cm

 From the above definitions it follows that for the  counter-clockwise loop
 $\gamma_{\infty}$ defined by   $z
 \mapsto z e^{2\pi i}$ for $|z|$  as big as to 
encloses all the  singularities
 $a_1$,...,$a_n$,  the monodromy $M_{\infty}$ of $Y(z)$ is 
$$ 
  Y(z) \mapsto Y(z) M_{\infty},~~~~M_{\infty}=
  (S_1^{(\infty)}...S_{2r_{\infty}}^{(\infty)} )^{-1}
$$ 
Note: $M_{\infty}=\exp\{2\pi i \Lambda_0^{(\infty)}\}$ if $r_{\infty}=0$
    
For a counter-clockwise loop $\gamma_k$ defined by  $(z-a_k)
 \mapsto (z-a_k) e^{2\pi i}$ in a neighbourhood of $a_k$ not containing other
 singularities
$$
   Y(z) \mapsto Y(z) M_k,~~~~M_k={ C^{(k)}}^{-1} (S_1^{(k)}...S_{2r_{k}}^{(k)}
                                      )^{-1} 
                                      C^{(k)}
$$
Note:   $M_k={ C^{(k)}}^{-1}\exp\{2\pi i \Lambda_0^{(k)}\}~ C^{(k)}$
 if $r_{k}=0$.

\subsection{ Isomonodromic deformations}

 Here we quote the results of the famous paper by Jimbo, Miwa and Ueno
 \cite{JMU}. Suppose that $A(z)$ depends on some additional parameters
 $t=(t_1,...,t_q)$. This means that
 $$
    a_k=a_k(t),
~~~~  A_l^{(k)}=A_l^{(k)}(t),
~~~~ A_l^{(\infty)}= A_l^{(\infty)}(t)
$$
Again we assume $A_0^{(\infty)}$ to be diagonal. We expect that the  
 $\Lambda_l^{(k)}$'s and $F_l^{(\nu)}$'s depend on $t$. 
The fundamental matrix $Y(z)$ now depends on $t$, namely it is 
 $Y(z,t)\equiv Y_1^{(\infty)}(z,t)$. 

Let $t$ vary in a (small) open set $V$ . We say that the deformation 
is {\it isomonodromic}  if
$$
  \left.\matrix{  \Lambda_{r_{\infty}}^{(\infty)},&S_1^{(\infty)},&
  ...~,&~S_{2r_{\infty}}^{(\infty)},& C^{(\infty)}=I \cr
  \Lambda_{r_{1}}^{(1)},&S_1^{(1)},&  ...~,&~S_{2r_{1}}^{(1)},& C^{(1)} \cr
\vdots&\vdots& & \vdots & \vdots \cr
 \Lambda_{r_{n}}^{(n)},&S_1^{(n)},&  ...~,&~S_{2r_{n}}^{(n)},& C^{(n)} \cr
}\right.
$$
 are independent of $t\in V$. In the case $r_{\nu}=0$ we also require  that
 $R^{(\nu)}$ is independent of $t$.
   This implies that the monodromy matrices are
 also independent of $t$ and the differential 
$$
  \Omega(z,t):= dY(z,t)~ Y(z,t)^{-1}= \sum_{j=1}^q~ {\partial Y(z,t) \over \partial
  t_j} Y(z,t)^{-1} ~dt_j,
$$
 is single-valued and  meromorphic  in  $z$. We have not enough space here to
 describe the details. It is enough to say that 
$$
  \Omega(z,t)= \sum_{l=0}^{r_{\infty}} \Omega_l^{(\infty)}(t) ~z^l+
               \sum_{k=1}^n\left\{ 
      \sum_{l=1}^{r_k+1}{\Omega_l^{(k)}(t)\over (z-a_k)^l}  \right\}.
$$
The forms $ \Omega_l^{(\nu)}(t)$ $\nu=1,...,n,\infty$ are uniquely and
explicitly determined
by $Y(z,t)$ and therefore by $A(z,t)$. More precisely,
$$
  F^{\infty}(z,t) dT^{(\infty)}(z,t)
  F^{\infty}(z,t)^{-1}=\sum_{l=0}^{r_{\infty}} \Omega_l^{(\infty)}(t) ~z^l
  +O\left({1\over z}\right)
$$
$$
   F^{(k)}(z,t) dT^{(k)}(z,t)  F^{(k)}(z,t)^{-1} =
   \sum_{l=1}^{r_k+1}{\Omega_l^{(k)}(t)\over (z-a_k)^l} +O(1)
$$

\vskip 0.2 cm 
\noindent
{\bf Theorem} \cite{JMU}: {\it The deformation 
 $t\in V$ of ${dY\over dz}=A(z,t)Y$
is isomonodromic if and only if $Y$ satisfies
$$ 
      dY = \Omega(z,t)Y 
$$
where $\Omega(z,t)$ is uniquely determined by $A(z,t)$.
}

\vskip 0.2 cm
If we expand both sides of $\Omega(z,t)=dY(z,t)~Y(z,t)^{-1}$ at $a_k$, we take 
into account that $Y(z,t)=
\left[F_0^{(k)}(t)+O(z-a_k)\right]~e^{T^{(k)}(z,t)} ~C^{(k)}$ in the right
hand-side and we equate the  zero-order terms $(z-a_k)^0$, we find the
following equation:
$$
  dF_0^{(k)}= \theta^{(k)}(t) ~F_0^{(k)},
$$
where $\theta^{(k)}$ is a form determined by $A(z)$ and $F_0^{(k)}$.

\vskip 0.2 cm

 In the above discussion, the asymptotic expansions and the Stokes' phenomenon
are supposed to be uniform in $t\in V$. This makes it possible to exchange
asymptotic expansions and differentiation ``$d$''.  
\vskip 0.3 cm

The systems
$$ 
       {\partial Y \over \partial z} = A(z,t)~ Y,
$$
$$
      dY=\Omega(z,t)~Y,
$$
are compatible if and only if
$$
    dA={\partial \Omega \over \partial z} +[\Omega ,A]  ~~~~(\hbox{it is }
                                                               d {\partial_z}=
                                                               {\partial_z }
                                                               d),
$$
(the second compatibility 
condition $d\wedge d=0$, namely $d\Omega=\Omega \wedge
\Omega$, follows from the above (see \cite{JMU})). 

\vskip 0.3 cm 

Finally, we construct a non linear system of differential equations for  the
(entries of) $A_l^{(\nu)}(t)$, for $a_k(t)$ and for 
 $F_0^{(k)}(t)$ which ensures
that the deformation is isomonodromic. 

 If we consider the entries of  the $A_l^{(\nu)}$ and the singular points
 $a_k$ as parameters, it is proved \cite{Shi} \cite{JMU} that the fundamental
 matrices are holomorphic and the
 asymptotic expansions are uniform  in a small open set of the 
 parameters.

\vskip 0.2 cm
\noindent
{\bf Theorem} \cite{JMU}: {\it The deformation  
(for $t$ in a small open set $V$) 
is isomonodromic if and only if
$A(z,t)$ and the $F_0^{(k)}(t)$'s are  solutions of 
$$
  dA={\partial \Omega \over \partial z}+[\Omega,A],
$$
\be
   dF_0^{(k)}= \theta^{(k)}(t)~F_0^{(k)}.
\label{Miwauenoeq}
\ee
for $t\in V$. Here the forms $\Omega$ and $\theta^{(k)}$ are given  above 
as functions
of $F_0^{(k)}$ and $A_l^{(\nu)}$. 
}

\vskip 0.3 cm
Finally, we recall that in \cite{JMU} it is proved that the maximum number of
independent parameters $t$ can be chosen to be
$$
   \left.\matrix{   &  \Lambda_0^{(\infty)}, &...~,&
   \Lambda_{r_{\infty}-1}^{(\infty)} \cr
 a_1 &  \Lambda_0^{(1)}, &...~,& \Lambda_{r_{1}-1}^{(1)} \cr
  \vdots & \vdots &            & \vdots \cr
    a_n &  \Lambda_0^{(n)}, &...~,& \Lambda_{r_{n}-1}^{(n)} \cr
}\right.
$$
where  by $\Lambda_l^{(k)}$ we mean its eigenvalues.


\section{ Monodromy data of a Semisimple Frobenius Manifold}\label{Monodromy
data of a Semisimple Frobenius Manifold}

 We apply the results of the previous section to a Frobenius manifold $M$. 
In the point $t=(t^1,...,t^n)\in M$  
we consider the system (\ref{systemz}). 
For simplicity, we suppose that $\hat{\mu}$ is diagonalizable with eigenvalues
$\mu_{\alpha}$, $\alpha=1,2,...,n$. Namely 
$
\hat{\mu}= \hbox{diag}(\mu_1,...,\mu_n)
$. Let $t$ be fixed. Since $z=0$ is a fuchsian singularity
 (\ref{systemz})  has a fundamental matrix solution 
$$
  \Xi(x,t)=H(z,t) ~z^{\hat{\mu}}~z^R,
~~~ H(z,t)=\left[\sum_{p=0}^{\infty} ~H_l(t) z^p\right],~~~H_0=I,
$$
$$
R=R_1+R_2+...,~~~~ (R_k)_{\alpha \beta}\neq 0~\hbox{ only if }
\mu_{\alpha}-\mu_{\beta} =k,
$$
where the series $\left[\sum_{p=0}^{\infty} ~H_l(t) z^p\right]$ is convergent
in a neighbourhood of $z=0$. $R$ is not uniquely determined. The
ambiguity is  \cite{Dub2}, 
$$
    R\mapsto G~R~G^{-1},
$$
where
$$
     G=1+\Delta_1+\Delta_2+...,
$$
$$
    (1-\Delta_1+\Delta_2-\Delta_3+...)~\eta~(1+\Delta_1+\Delta_2+\Delta_3+...),~~~~(\Delta)_{\alpha\beta}\neq
    0 \hbox{ only if } \mu_{\alpha}-\mu_{\beta}= k
$$ 
Let's call $[R]$ such an orbit.

 In \cite{Dub2} it is proved that $R\in [R]$ is independent of $t\in M$. In the
 semisimple case this
 property of isomonodromicity follows from our general considerations about
 isomonodromic deformations. Actually, $R$ is independent 
of {\it any} $t\in M$, not only for $t$ in a small
 neighbourhhod of 
 a given $t_0\in M$. 
 Since $R$ and $\hat{\mu}$ are independent of $t\in M$, 
the following definition makes sense: 
\vskip 0.2 cm
\noindent
{\bf Definition:} $\hat{\mu}$ and $[R]$ are called {\it monodromy data of the
Frobenoius manifold} at $z=0$.
\vskip 0.3 cm

 We turn to the semisimple case through the gauge $y=\Psi\xi$. 
The system  (\ref{systemz}) becomes (\ref{systemz1}). 
A fundamental matrix solution is 
$$
   Y^{(0)}(z,u)= \Phi(z,u)~z^{\hat{\mu}}~z^R,~~~~ \Phi(z,u)=\sum_{p=0}^{\infty}
   ~\phi_p(u)z^p,~~~\phi_0(u)=\Psi(u),
$$
(the $\phi_p$'s are the matrices $F_p^{(0)}$ of section \ref{secdeformata}).
In $\pi(CP^1\backslash\{0,\infty\},z_0)$ we consider the basic loop 
$z_0\mapsto z_0 e^{2 \pi i}$. The monodromy of $Y^{(0)}$ is $e^{2\pi
i\hat{\mu}} e^{2\pi i R}$. 

\vskip 0.2 cm
\noindent
{\it Remark:}
The ``symmetries'' $\eta \hat{\mu}+ \hat{\mu}^T \eta =0$, ${\cal U}^T
\eta=\eta {\cal U}$  imply that 
$ \xi_1(-z,t)^T \eta \xi_2(z,t)$ is independent of $z$ for any 
two solutions $\xi_1(z,t)$, $\xi_1(z,t)$ of (\ref{systemz}). In particular 
$H(-z,t)^T\eta H(z,t)=\eta$. Now, $\Phi(z,u)=\Psi(u) H(z,t(u))$, therefore
$$
  \Phi(-z,u)^T~\Phi(z,u)=\eta.
$$
\vskip 0.2 cm

(\ref{systemz1})  has a formal solution 
$$
   Y_F = \left[ I + {F_1 \over z}+{F_2 \over z^2}+
... \right]~e^{z~U},~~~~z\to \infty
$$ 
where $F_j$'s are $n \times n$ matrices (the $F_j^{(\infty)}$ in our general
discussion).  We explained  
that fundamental matrix  solutions exist which have $\tilde{Y}_F$ as 
asymptotic expansion for $z \to \infty$. 
 We choose the  Stokes' rays
$$ 
   R_{rs}= \{ z=- i \rho (\bar{u_r} - \bar{u_s}), ~~~\rho>0 \},~~~ r
   \neq s. 
$$
Note that 
$\Re  ~((u_r-u_s)z)=0$ on $R_{rs}$, therefore the definition of Stokes' rays
is satisfied. Let $l$ be an oriented  
  line not containing Stokes' rays and passing through the origin, 
with  a positive 
half-line $l_{+}$ and a negative $l_{-}$. It is characterize also by the  angle
$\varphi$ between $l_{+}$ and the positive real axis.  We call $\Pi_R$
and $\Pi_L$ the half planes to the right and left of $l$ w.r.t its
orientation. $\Pi_L$ can be extended in $\Pi_R$ and $\Pi_R$  in $\Pi_L$ 
respectively, up to the first Stokes rays we encounter performing that 
extension. The resulting extended {\it overlapping} sectors will be 
called ${\cal S}_L$ and ${\cal S}_R$.
 There exist unique fundamental matrix solutions  $Y_L$ and $Y_R$ having the 
asymptotic expansion in ${\cal S}_L$ and ${\cal S}_R$ 
respectively \cite{BJL1}.  They are related by the {\it Stokes' matrix} 
matrix $S$,
such that $$
    Y_L(z,u)=Y_R(z,u)S
$$  
in the overlapping of ${\cal S}_L$ and ${\cal S}_{R}$ containing $l_{+}$. 
In the opposite overlapping region containing $l_{-} $
one can prove (as a consequence of
the skew-symmetry of $V$, see \cite{Dub1}) that the corresponding Stokes
matrix is $S^T$: namely,  
$
 Y_L(z,u)=Y_R(z,u)S^T.
$  
 The Stokes' matrix $S$ has entries  
$$ 
                 s_{ii}=1, ~~ ~~~~~ 
               s_{ij}=0~~\hbox{ if } R_{ij}\subset \Pi_R
$$
This follows from the fact that on the overlapping region $0<\arg z
<{\pi \over k}$ there are no Stokes' rays and 
$$ 
   e^{zU}~S~e^{-zU} \sim I, ~~~~z \to \infty,~~~~
\hbox{ then }~~ e^{z(u_i-u_j)}~s_{ij}\to \delta_{ij}
$$
 Moreover, $\Re \left(z(u_i-u_j)\right)>0$ to the
left of the ray $R_{ij}$, while  $\Re \left(z(u_i-u_j)\right)<
0$ to the right (the natural orientation on $R_{ij}$, from
$z=0$ to $\infty$ is understood). This implies
$$ 
   \left|e^{zu_i} \right| > \left|e^{zu_j} \right|~~~~~ \hbox{ and  }
                      e^{z(u_i-u_j)}\to \infty ~~\hbox{ as } z \to \infty 
$$
on the left, while on the right 
$$ 
   \left|e^{zu_i} \right| < \left|e^{zu_j} \right|~~~~~ \hbox{ and  }
                      e^{z(u_i-u_j)}\to 0 ~~\hbox{ as } z \to \infty 
$$

\vskip 0.2 cm

 We call  {\it central connection
 matrix} the connection matrix $C$ such  that 
$$
   Y^{(0)}(z,u)=\tilde{Y}_{R}(z,u) C   ~~~z\in \Pi_R,
$$ 
(here $C$ is the inverse of the connection matrix $C^{(0)}$ we introduced in
the general discussion).

\vskip 0.2 cm

The monodromy of 
$$
Y(z,u):=Y_R(z,u)
$$ 
at $z=0$ is therefore
$$
     M_0=  S^T S^{-1}\equiv C~e^{2\pi i \hat{\mu}}e^{2\pi i R}~C^{-1}
 $$
The ``symmetries'' $\Phi(-z,u)^T\Phi(z,u)=\eta$, $Y_R(z,u)^TY_L(z,u)=I$ (in
the overlapping region containing $l_{+}$),  $Y_R(ze^{-2\pi i},u)^TY_L(z,u)=I$ 
(in
the overlapping region containing $l_{-}$) imply 
 $$
    S= C ~e^{-i\pi R} ~e^{-i\pi \hat{\mu}}~\eta^{-1} ~C^T,
$$
$$
    S^T= C~e^{i\pi R} ~e^{i\pi \hat{\mu}}~\eta^{-1} ~C^T,
$$

We note that if $C^{\prime}$ is such that $z^{\hat{\mu}} z^R~C^{\prime}
z^{-R}z^{-\hat{\mu}}= C_0+C_1z+C_2z^2+...$ is 
 a convergent series at $z=0$, then
$Y^{(0)}C^{\prime}$ has again the form $[\sum_{p}\phi_p^{\prime}(u)]
~z^{\hat{\mu}} z^R$. In particular, $C_0\hat{\mu}=\hat{\mu}C_0$ and therefore 
$\phi_0^{\prime}=\phi_0 C_0$ is again a matrix such that
${\phi_0^{\prime}}^{-1} V 
\phi_0^{\prime}=\hat{\mu}$. Also, ${C^{\prime}}^{-1}~e^{2\pi i
\hat{\mu}}e^{2\pi 
i R} ~ C^{\prime}=e^{2\pi i \hat{\mu}}e^{2\pi
i R}$. Such matrices form a normal subgroup $C_0(\hat{\mu},R)$ 
 in the group   of
matrices commuting 
with $e^{2\pi i \hat{\mu}}e^{2 i \pi R}$. 

\vskip 0.2 cm
\noindent
{\bf Theorem} \cite{Dub2}: 
{\it If two systems 
$$ 
    {dy^{(i)}\over dz}=\left[U+{V^{(i)}\over z} \right]~y^{(i)},
    ~~~{V^{(i)}}^T=-V^{(i)},~~~i=1,2, 
$$
have the same $\hat{\mu}$ (diagonal from of $V^{(i)}$), $R$, $S$ (w.r.t. the
same line) and the same $C$ (up to $C\mapsto CC^{\prime}$, $C^{\prime}\in
C_0(\hat{\mu},R)$),  then $V^{(1)}=V^{(2)}$.}

\vskip 0.2 cm
The theorem is important to our purpose of classification of Frobenius
manifolds. It states that $\hat{\mu}$, $R$, $S$ (w.r.t. a line) and $C$
determine uniquely the system (\ref{systemz1}) 
at a fixed  $u=(u_1,...,u_n)$.

\vskip 0.3 cm
In the discussion above $u=(u_1,...,u_n)$ was fixed. Now we let if vary.

\vskip 0.2 cm
\noindent
{\bf Definition:} Let $e$ be the unit vector. We call $e$, $\hat{\mu}$, $R$, $S$, $C$  the {\it monodromy data} 
of the FM in a neighbourhood of the semisimple point $(u_1,...,u_n)$ where
they are computed.

\vskip 0.2 cm

The definition makes sense 
 because the dependence on $u$ in the coefficients of 
(\ref{systemz1}) is {\it isomonodromic}, as we are going to show. We included
 $e$ in the definition because the eigenvector of $V$ with eigenvalue
 $\mu_1=-d/2$ must be marked and it corresponds to th unity in $M$. In
 canonical coordinates it is the first column of $\Psi=\phi_0$. 

We explain why the dependence on $u$ is isomonodromic. 
 We  know 
that $\hat{\mu}$ and $R$ are independent of $t\in M$. Therefore, for the
 local change of coordinates $t=t(u)$ they are {\it locally} independent of
 $u$. Namely, they are constant if $u$ varies in a small open set.  
The system (\ref{systemz1}) is a particular example of the general systems
of section \ref{secdeformata}, Let
 $$
Y(z,u):=Y_R(z,u). 
$$ 

\vskip 0.2 cm
a) Suppose that the deformation $u$ is isomonodromic. From 
\be
  Y(z,u)=\left[I+{F_1\over z}+O\left({1\over z}\right)\right]~e^{zU},~~~~z\to
  \infty 
\label{ancora la soluzione} 
\ee
$$
  = \sum_{p=0}^{\infty} ~  \left[\phi_p(u)z^p
  \right]~z^{\hat{\mu}}z^R~C^{-1},~~~~z\to 0.
$$
we construct $\Omega(z,u)$: 
$$
  \Omega(z,u):= {\partial Y(z,u)\over u_i} ~Y^(z,u)^{-1} = 
\left\{\matrix{ 
                z E_i+ [F_1,E_i]+O\left({1\over z}\right),~~~z\to\infty\cr
\cr
                {\partial \phi_0 \over \partial u_i}~\phi_0^{-1} +O(z),~~~z\to
                0 
\cr
}\right.
$$
Therefore, $\Omega(z,u)- \left( z E_i+ [F_1,E_i]\right)\to 0$ as $z\to
\infty$ 
and it is holomorphic at $z=0$. By Liouville theorem:
$$
  \Omega(z,u)=   z E_i+ [F_1(u),E_i]
$$
We expand the two sides of 
$\Omega_0(z,u)={\partial Y(z,u)\over u_i} ~Y^(z,u)^{-1} $ near $z=0$: 
$$   z E_i+[F_1,E_i]=
  {\partial \phi_0 \over \partial u_i}~\phi_0^{-1} +O(z)
$$
At $z=0$: 
$$
     {\partial \phi_0 \over \partial u_i}=[F_1,E_i]~\phi_0
$$ 
This is precisely the equation (\ref{Miwauenoeq}). 
Finally, $F_1$ is computed from the coefficients of
(\ref{systemz1}) by substituting (\ref{ancora la soluzione}): 
$$
  \left(-{F_1^2\over z^2}+...\right)e^{zU}+\left(I+{F_1\over z} +...\right)
Ue^{zU} 
= \left[U+{V\over z}\right]\left(I+{F_1\over z} +...\right)
Ue^{zU} 
$$
From the term  ${1\over z}$ we have 
$$
    (F_1)_{ij}={V_{ij}\over u_j-u_i} 
$$
and from the term ${1\over z^2}$ we have
$$
  (F_1)_{ii}= \sum_{j\neq i} ~{V_{ij}^2\over u_j-u_i}
$$
If we put 
$$
  V_k:=[F_1,E_k]
$$
we obtain 
$$
(V_k)_{ij}={\delta_{ki}-\delta_{kj}\over u_i-u_j}~V_{ij} 
$$
which is precisely equivalent to $[U,V_k]=[E_k,V]$. 

\vskip 0.2 cm
b) From the general theory of section \ref{secdeformata} 
 we  conclude that the deformation $u$ of the system (\ref{systemz1})  is
isomonodromic if and only if 
$$ 
 {\partial Y \over \partial u_i} = [zE_i+V_i] ~Y
$$
where $V_i$ is uniquely determined by 
$$ 
[U,V_k]=[E_k,V]
$$
In particular,  the matrix $\phi_0(u)$ satisfies
$$
  {\partial \phi_0\over \partial u_i}= V_i ~\phi_0.
$$
Here we recognize precisely the properties 
 of  a semisimple Frobenius manifold!

\vskip 0.2 cm
In section \ref{secdeformata} we  learnt that the deformation is
isomonodromic if and only if $dA=\partial_z \Omega+[\Omega,z]$, $ dF_0^{(k)}=
\theta^{(k)} F_0^{(k)}$. In our case, the two conditions become 
respectively the condition of compatibility:  
$$
 {\partial V \over \partial u_i}= [V_i,V],~~~[U,V_i]= [E_i,V]
$$
  and 
$$
  {\partial \phi_0 \over \partial u_i}= V_i~\phi_0.
$$  

\vskip 0.2 cm 

 We conclude that for a semisimple Frobenius manifold, $u$ is an isomonodromic
 deformation, then $\hat{\mu}$, $R$, $S$ (for a fixed line), $C$ are
 independent of $u$, if $u$ varies in a sufficiently small open set.


\section{ Inverse Reconstruction of a Semisimple FM (II)}\label{INVERSE II}

 Let's fix $u=u^{(0)}=(u_1^{(0)},...,u_n^{(0)})$ such that $u_i^{(0)}\neq
u_j^{(0)} $ for $i\neq j$.  
Suppose we give $\mu$, $R$, an admissible line $l$, $S$ and $C$ such that 
$$
   s_{ij}\neq 0 \hbox{ only if the Stokes' ray } R_{ij}\in \Pi_L
$$
$$
     S^T S^{-1}\equiv C~e^{2\pi i \hat{\mu}}e^{2\pi i R}~C^{-1}
 $$
 $$
    S= C ~e^{-i\pi R} ~e^{-i\pi \hat{\mu}}~\eta^{-1} ~C^T,
$$
$$
    S^T= C~e^{i\pi R} ~e^{i\pi \hat{\mu}}~\eta^{-1} ~C^T.
$$
Let $D$ be a disk specified by $|z|<\rho$ for some small $\rho$. Let $P_L$ and
$P_R$ be the intersection of the external part of the disk with $\Pi_L$ and
$\Pi_R$ respectively. We denote by $\partial D_R$ and $\partial D_L$  the lines
 on the boundary of 
$D$ on the side of $P_R$ and $P_L$ respectively; we denote by $\tilde{l}_{+}$
and $\tilde{l}_{-}$
 the portion of $l_{+}$ and $l_{-}$ respectively, 
on the common boundary of $P_R$ and $P_L$.
 Let's consider the following (discontinuous) boundary value problem (b.v.p.):
``construct a piecewise holomorphic matrix function 
$$
   \Phi(z)=\left\{\matrix{\Phi_R(z),~~~z\in P_R\cr
                           \Phi_L(z),~~~z\in P_L\cr
\Phi_0(z),~~~z\in D\cr
}\right.,
$$
continuous on the boundary of $P_R$, $P_L$, $D$ respectively, such that 
$$
   \Phi_L(\zeta)=\Phi_R(\zeta)~e^{\zeta U}S e^{-\zeta U},~~~~\zeta\in
   \tilde{l}_{+}
$$
$$
  \Phi_L(\zeta)=\Phi_R(\zeta)~e^{\zeta U}S^T e^{-\zeta U},~~~~\zeta\in
   \tilde{l}_{-}
$$
$$
\Phi_0(\zeta)=\Phi_R(\zeta)~e^{\zeta U} C\zeta^{-R}\zeta^{-\hat{\mu}}
,~~~~\zeta\in \partial D_R
$$
$$
\Phi_0(\zeta)=\Phi_L(\zeta)~e^{\zeta U}S^{-1} C\zeta^{-R}\zeta^{-\hat{\mu}}
,~~~~\zeta\in \partial D_L
$$
$$
   \Phi_{L/R}(z)\to I \hbox{ if  $z\to \infty$ in $P_{L/R}$'' }. 
$$
 The reader may observe that 
 $ \tilde{Y}_{L/R}(z):=\Phi_{L/R}(z) e^{zU},
$  
$
    \tilde{Y}^{(0)}(z):=\Phi_0(z,u)z^{\hat{\mu}} z^{R}$
 have precisely the 
monodromy
 properties of the
 solutions of (\ref{systemz1}). 

\vskip 0.2 cm
\noindent
{\bf Theorem} \cite{Miwa}\cite{Malgrange}\cite{Dub2}: {\it If the above boundary value problem has solution for a given
$u^{(0)}= (u_1^{(0)},...,u_n^{(0)})$ such that $u_i^{(0)}\neq
u_j^{(0)} $ for $i\neq j$, then: 

i)  it is unique. 

ii) The solution exists and it is analytic 
for $u$ in a neighbourhood of $u^{(0)}$. 

iii) The solution has analytic continuation as a meromorphic function 
on the universal covering of ${\bf C}^n\backslash\{diagonals\}:=
\{ (u_1,...,u_n)~|~u_i\neq u_j \hbox{ for } i\neq j\}$. 
}

\vskip 0.2 cm
 Consider the solutions $\tilde{Y}_{L/R}$, $\tilde{Y}^{(0)}$ 
of the b.v.p.. We have $\Phi_R(z)=I+{F_1\over z} +O\left({1\over z^2}\right)$
as $z\to \infty$ in $P_R$. We also have $\Phi_0(z)= \sum_{p=0}^{\infty} \phi_p
z^p$ as $z\to
0$. Therefore
$$
   {\partial \tilde{Y}_R \over \partial z}~\tilde{Y}_R^{-1}=\left[U+
{[F_1,U]\over z} +O\left({1\over z^2}\right)\right],~~~~z\to \infty,
$$
$$
  {\partial \tilde{Y}^{(0)} \over \partial z}~(\tilde{Y}^{(0)})^{-1}=
{1\over z}\left[\phi_0\hat{\mu}\phi_0^{-1}+O(z)\right],~~~~z\to 0.
$$
 Since $C$ is independent of $u$ the right hand-side of the two equalities
 above are equal. Also $S$ is independent of $u$, therefore
 the matrices $\tilde{Y}$ above satisfy 
$$
  {\partial y \over \partial z} =\left[ U+{V\over
  z}\right]~y,~~~~V(u):=[F_1(u),U]\equiv \phi_0\hat{\mu}\phi_0^{-1}. 
$$
In the same way 
$${\partial  \tilde{Y}_R \over \partial
u_i}~\tilde{Y}_R^{-1}=zE_i+[F_1,E_i]+\left({1\over z}\right) ,~~~~z\to\infty 
$$
$$
 {\partial \tilde{Y}^{(0)} \over \partial u_i}~( \tilde{Y}^{(0)}
 )^{-1}={\partial 
 \phi_0 \over \partial u_i} ~\phi_0^{-1}+O(z),~~~~z\to 0.
$$
The right hand-sides are equal, therefore the $\tilde{Y}$'s satisfy
$$
 {\partial y\over \partial u_i}= [zE_i+V_i]~y,~~~~V_i(u):=[F_1(u),E_i]\equiv 
{\partial \phi_0\over \partial u_i} ~\phi_0^{-1}.
$$

 \vskip 0.2 cm
We conclude that from the solution of the b.v.p. we obtain solutions to
(\ref{systemz1}), (\ref{systemt1}). This means that we can locally reconstruct 
a Frobenius structure from the local solution of the b.v.p. and we can do
the analytic continuation of such a structure by analytically continuing the
solution of the b.v.p to the universal covering of ${\bf
C}^n\backslash\hbox{diagonals}$. In order to do this, it is enough to use the
solution of the b.v.p.  $\Phi(z)= \sum_{p=0}^{\infty} \phi_p
z^p$ as $z\to
0$ and the formulae (\ref{tu}), (\ref{Fu}), 
provided that at the given initial point $u^{(0)}$ ($u_i^{(0)}\neq u_j^{(0)}$
for $i\neq j$) the condition $ 
  \Pi_{i=1}^n \phi_{i1,0}(u^{(0)}) \neq 0.
$
is satisfied.  
 If it is not, there is a singularity in the change of
 coordinates; actually $c_{\alpha\beta\gamma} =\sum_{i=1}^n
   {\phi_{0,i\alpha}\phi_{0,i\beta}\phi_{0,i\gamma} \over \phi_{0,i1}}$ 
may diverge and 
 $dt^{\alpha}= \eta^{\alpha\beta} \sum_{i=1}^n \phi_{i\beta,0}(u)\phi_{i1,0}(u)
du_i$ are not independent. 

\subsection{ Analytic continuation}

The analytic continuation of the solution of the b.v.p. 
to the universal covering of  ${\bf
C}^n\backslash\hbox{diagonals}$ gives the analytic continuation of a 
Frobenius structure. Since $(u_1,...,u_n)$ are local coordinates, they are
defined up to permutation. 
Therefore, the analytic
continuation is described by the fundamental group  $\pi(( {\bf
C}^n\backslash\hbox{diagonals})/{\cal S}_n,u^{(0)})$, where $S_{n}$ is the
symmetric 
group of $n$ elements. This group is called {\it Braid group} ${\cal B}_n$
\cite{Braid}.   

 The local solution of the b.v.p. is obtained   
from the monodromy data. Whereas $\hat{\mu}$ and $R$ are constant on the
manifold,  $S$ and $C$ are constant only for small deformations.  
We fix the 
 line $l$ {\it once and for all} and let $u$ vary, so that also the
Stokes' rays move. From the definition of $S$, it
 follows that whenever a Stokes' ray crosses $l$ some entries of $S$ which
where zero may become non zero, and other entries must vanish. This is a
discrete jump, described by an action of the braid group. 

 In order to describe this action we first note what is the effect of
 permutations on (\ref{systemz1}), $S$ and $C$. Let 
$U = $ diag$( u_1,~u_2,~...,~u_n)$. Let $\sigma:~(1,2,..,n)$ 
$\mapsto$ $(\sigma(1),\sigma(2),...,\sigma(n))$ be a permutation. It 
is represented by an invertible matrix $P$ which acts as the gauge  
$y\mapsto P y$. The new system has matrix   
$$ P~U~P^{-1}~=~\hbox{
diag}(u_{\sigma(1)},u_{\sigma(2)},...,u_{\sigma(n)}).$$ 
$S$ and $C$ are then transformed in $PSP^{-1}$ and $PC$. 
For a suitable $P$, $PSP^{-1}$ is  upper triangular.  
As a general result \cite{BJL1}, the good permutation is the
one which puts  $u_{\sigma(1)}$, ..., $u_{\sigma(n)}$ in lexicographical
order w.r.t. the oriented line $l$. 
$P$ corresponds to a change of coordinates in the 
given chart, consisting in the permutation $\sigma$ of the
coordinates. 

 \vskip 0.2 cm 

\noindent
We  recall that the  braid group is generated by $n-1$ elementary 
braids $\beta_{12}$, $\beta_{23}$, ..., $\beta_{n-1,n}$,  with
relations: 
$$ 
    \beta_{i,i+1}\beta_{j,j+1}=\beta_{j,j+1}\beta_{i,i+1}~\hbox{ for
} i+1 \neq j,~j+1\neq i$$
$$ 
 \beta_{i,i+1}\beta_{i+1,i+2}\beta_{i,i+1}=
\beta_{i+1,i+2}\beta_{i,i+1} 
     \beta_{i+1,i+2} $$

  Let's start from $u=u^{(0)}$. If
 we move sufficiently far away from $u^{(0)}$, 
some Stokes' rays cross the
 {\it fixed} admissible line $l$. Then, we must change
 $Y_L$ and $Y_R$, $S$ and $C$.    
The motions of the points $u_1$, ..., $u_n$  
represent transformations of
the braid group. Actually, a braid $\beta_{i,i+1}$ can be represented  as an
 ``elementary'' 
 deformation 
consisting  of  a permutation of $u_i$, $u_{i+1}$ moving 
counter-clockwise (clockwise or counter-clockwise is a matter of
convention).

 Suppose $u_1$, ..., $u_n$ are already  in
lexicographical order w.r.t. $l$, so that $S$ is upper
triangular (recall that this configuration can be reached by a
suitable permutation $P$). The effect on $S$ 
of the deformation of $u_i$,
$u_{i+1}$ representing $\beta_{i,i+1}$ is the following \cite{Dub2}: 
$$ 
   S \mapsto~S^{\beta_{i,i+1}}:=A^{\beta_{i,i+1}}(S)~S~A^{\beta_{i,i+1}}(S)
$$
where 
$$ 
   \left(A^{\beta_{i,i+1}}(S) \right)_{kk}=
   1~~~~~~k=1,...,n~~~n\neq~i,~i+1
$$
$$
 \left(A^{\beta_{i,i+1}}(S) \right)_{i+1,i+1}=-s_{i,i+1}
$$
$$
 \left(A^{\beta_{i,i+1}}(S) \right)_{i,i+1}=
 \left(A^{\beta_{i,i+1}}(S) \right)_{i+1,i}=1
$$
and all the other entries are zero. 
  For the inverse braid $\beta_{i,i+1}^{-1}$ ($u_i$ and $u_{i+1}$ move
  clockwise) the representation is
 $$ 
   \left(A^{\beta_{i,i+1}^{-1}}(S) \right)_{kk}=
   1~~~~~~k=1,...,n~~~n\neq~i,~i+1
$$
$$
 \left(A^{\beta_{i,i+1}^{-1}}(S) \right)_{i,i}=-s_{i,i+1}
$$
$$
 \left(A^{\beta_{i,i+1}^{-1}}(S) \right)_{i,i+1}=
 \left(A^{\beta_{i,i+1}^{-1}}(S) \right)_{i+1,i}=1
$$
and all the other entries are zero. 
We remark that $S^{\beta}$ is still upper triangular.

The effect on $C$ is
$$ 
   C \mapsto A^{\beta_{i,i+1}} C
$$

\vskip 0.2 cm
 Not all the braids are actually to be considered.  Suppose we 
do the following gauge  $y\mapsto J y$, $J=$diag$(\pm1,...,\pm1)$,
on the system (\ref{systemz1}).  Therefore ${J}U{J}^{-1}\equiv U$
but $S$ is transformed to  ${J} S {J}^{-1}$, where some  entries
change sign. The formulae  
which define a local chart of the manifold in
terms of monodromy data (i.e. in terms of $\phi_p$, $p=0,1,2,3$) 
 are not affected by this transformation.  The analytic continuation of the
local structure on the universal covering  of $({\bf
C}^n\backslash\hbox{diagonals})/{\cal S}_n$ 
  is therefore described by the elements of  the quotient group 
\be
   {\cal B}_n/ \{\beta\in {\cal B}_n ~|~S^{\beta}=JSJ\}
\label{quozienteinutile}
\ee
Therefore (also recall that the analytic continuation is
meromorphic) we conclude \cite{Dub2}:

\vskip 0.2 cm
{\it For the given monodromy data ($e$, $\hat{\mu}$, $R$, $S$, $C$) the
local Frobenius structure
obtained from the solution of the b.v.p. extends to a
dense open subset of the manifold given by the covering of  $({\bf
C}^n\backslash\hbox{diagonals})/{\cal S}_n$ w.r.t. the covering
transformations in the  quotient (\ref{quozienteinutile}).}
\vskip 0.2 cm
{\it Let's  start from a Frobenius manifold $M$. Let ${\cal M}$ be the open
sub-manifold of points $t$ such that ${\cal U}(t)$ has distinct eigenvalues. If
we  compute its monodromy
data   ($e$, $\hat{\mu}$, $R$, $S$, $C$) at a point $u^{(0)}\in {\cal M}$ 
and we construct the Frobenius structure from the
analytic continuation of the corresponding b.v.p. on the covering of  $({\bf
C}^n\backslash\hbox{diagonals})/{\cal S}_n$ w.r.t. the quotient
(\ref{quozienteinutile}), then there is
an equivalence of Frobenius structures between this last manifold and ${\cal
M}$. 
}


\section{ Intersection Form and Monodromy Group of a Frobenius Manifold}

 The deformed flat connection was introduced as a natural structure on a
 Frobenius manifold and allows to transform the problem of solving the WDVV
 equations to a problem of isomonodromic deformations. There is a further
 natural structure on a Frobenius manifold which makes it possible to do the
 same. It is the intersection form. 

\vskip 0.2 cm
There is a natural isomorphism $\varphi:T_tM \to T_t^{*}M$ induced by
$<.,.>$. Namely,  
let $v\in T_tM$ and define $\varphi(v):=<v,.>$. This allow us to define the
product in $T_t^{*}M$ as follows: for $v,w\in T_tM$ we define $\varphi(v)\cdot
\varphi(w):=<v\cdot w,.>$. In flat coordinates $t^1,...,t^n$ the product is 
$$
    dt^{\alpha}\cdot dt^{\beta}=c^{\alpha\beta}_{\gamma}(t)~dt^{\gamma},~~~~
c^{\alpha\beta}_{\gamma}(t)=\eta^{\beta\delta}c_{\delta\gamma}^{\alpha}(t),
$$
(sums over repeated indices are omitted).

\vskip 0.2 cm
\noindent
{\bf Definition:} The {\it intersection form} at $t\in M$ is a bilinear form on
$T_t^{*}M$ defined by
$$
   (\omega_1,\omega_2):=(\omega_1\cdot\omega_2) (E(t))
$$
where $E(t)$ is the Euler vector field. In coordinates 
$$
g^{\alpha\beta}(t):=
(dt^{\alpha},dt^{\beta})= E^{\gamma}(t) c_{\gamma}^{\alpha \beta}.
$$
Recall  that $c^{\alpha\beta}_{\gamma}$ is independent of $t^1$. Let
$\tilde{t}:=(t^2,...,t^n)$. Also note that 
$c^{\alpha\beta}_1=\eta^{\alpha\beta}$; then 
$$
    g^{\alpha\beta}(t)=t^1\eta^{\alpha\beta}
    +\tilde{g}^{\alpha\beta}(\tilde{t}) 
$$
where $\tilde{g}^{\alpha\beta}(\tilde{t})$  depends only on
$\tilde{t}$. Therefore 
$$
  \det ((g^{\alpha\beta}(t)))=  \det(\eta^{-1})~ (t^1)^n+c_{n-1}~(\tilde{t})
(t^1)^{n-1}  
+c_{n-2}~(\tilde{t}) (t^1)^{n-2}+...+c_0(\tilde{t})$$
This proves that in an analytic Frobenius manifold the intersection form is
non-degenerate on an open dense subset $M\backslash
\Sigma$, where 
$$  
\Sigma:= \{ t\in M ~|~ \det((g^{\alpha\beta}(t)))= 0 \}$$ 
is called the {\it discriminant locus}. 
 On $M\backslash \Sigma$ we define
 $(~g_{\alpha\beta}(t)~):=(~g^{\alpha\beta}(t)~)^{-1} 
$. It is a result of \cite{Dub1} \cite{Dub2} that 
\vskip 0.2 cm 
{\it The metric $g_{\alpha \beta} dt^{\alpha}dt^{\beta}$ is flat on
$M\backslash \Sigma$ and its Christoffel coefficients in flat coordinates for
$<.,.>$ are 
$$
\Gamma_{\gamma}^{\alpha\beta}:=-g^{\alpha\delta}\Gamma_{\delta\gamma}^{\beta}=
\left({d+1\over 2}-q_{\beta}\right) ~c^{\alpha\beta}_{\gamma}
$$ 
}

\vskip 0.2 cm

 Let's denote by $<.,.>^{*}$ the bilinear form on $T^{*}_tM$ defined by
 $\eta^{-1}=(\eta^{\alpha\beta})$. Let $(.,.)-\lambda <.,.>^{*}$ 
be a family of bilinear forms on $T^{*}M$,
 parametrized by $\lambda\in {\bf C}$.  In flat coordinates for $\eta$: 
$$
   g^{\alpha\beta}(t)-\lambda\eta^{\alpha\beta}=
   \eta^{\alpha\beta}(t^1-\lambda)+\tilde{g}^{\alpha\beta}((\tilde{t}) 
$$
 Therefore 
$$
  \det( g^{\alpha\beta}(t)-\lambda\eta^{\alpha\beta})= \det(\eta^{-1})~
  (t^1-\lambda)^n+ c_{n-1}~(\tilde{t})(t^1-\lambda)^{n-1}+...+c_0(\tilde{t})
$$
and $(.,.)-\lambda <.,.>^{*}$ is not degenerate on the open dense subset
$M\backslash \Sigma_{\lambda}$, where
$$
\Sigma_{\lambda}:=\{ t~|~det
( g^{\alpha\beta}(t)-\lambda\eta^{\alpha\beta})=0\}. 
$$
The $n$ roots of the determinant, considered as a polynomial in
 $(t^1-\lambda)$,  are:  
\be 
       t^1-\lambda= f_1(\tilde{t}),~...,~
 t^1-\lambda= f_n(\tilde{t}),
\label{porcellini}
\ee
where $f_1$, ..., $f_n$ are functions of $\tilde{t}=(t^2,...,t^n)$. 
 
\vskip 0.2 cm
\noindent
{\bf Theorem} \cite{Dub1}: {\it 
 $(.,.)-\lambda <.,.>^{*}$ is a 
 flat pencil of metrics, namely  it is flat  and 
 its Christoffel coefficients are the sum of the
Christoffel coefficients of $g$ and $\eta$: 
$$ 
  \left(\Gamma_{(.,.)-\lambda<.,.>^{*}}\right)^{\alpha\beta}_{\gamma}=
  \left(\Gamma_{(.,.)}\right)^{\alpha\beta}_{\gamma}+0=\left({d+1\over 2}-q_{\beta}\right)
  ~c^{\alpha\beta}_{\gamma}. 
$$
}

\vskip 0.2 cm
 In order to  find a flat coordinate $x$ for the pencil we impose 
$\hat{\nabla}dx=0$, where $\hat{\nabla}$ is the connection for the metric
induced by the pencil. We skip details (see \cite{Dub1} \cite{Dub2}) and we
give the final result:
\be 
\left({\cal U}(t)-\lambda\right) {\partial \xi \over \partial \lambda} =
\left( {1\over 2} +\hat{\mu}\right) \xi
\label{gaussmanin1}
\ee
\be 
\left({\cal U}(t)-\lambda\right) \partial_{\beta} \xi +C_{\beta}\left({1\over
2}+\hat{\mu} \right)\xi=0
\label{gaussmanin2}
\ee
where the column vector  $\xi=(\xi^1,...,\xi^n)^T$ is defined by $\xi^{\alpha}=
\eta^{\alpha\beta}{\partial x \over \partial t^{\beta}}$.  ${\cal U}$ and
$C_{\beta}$ have already been defined.

\vskip 0.3 cm 
 In the semisimple case, let $u_1,...,u_n$ be local canonical coordinates,
equal to the distinct 
eigenvalues of ${\cal U}(t)$. From the definitions we have  
$$
du_i\cdot du_j= {1\over \eta_{ii}}\delta_{ij} du_i,~~~
   g^{ij}(u)= (du_1,du_j)= {u_i\over \eta_{ii}} \delta_{ij},
   ~~~~~\eta_{ii}=\psi_{i1}^2 
$$ 
Then $g^{ij}-\lambda\eta^{ij}= {u_i-\lambda\over \eta_{ii}} \delta_{ij}$ and 
$$ 
  \det((g^{ij}-\lambda\eta^{ij}))= {1\over
  \det((\eta_{ij}))}~(u_1-\lambda)(u_2-\lambda)...(u_n-\lambda).
$$ 
Namely, the roots $\lambda$ of the above polynomial are the canonical
coordinates. Now we  recall that ${\cal M}\subset M$ was defined as the 
sub-manifold of
semisimple 
points such that $u_i\neq u_j$ for $i\neq j$. 
If $\lambda$ is fixed, then the discriminant  $\Sigma_{\lambda}$
is:
$$ 
   \Sigma_{\lambda}\cap {\cal M}:=  \bigcup_{i=1}^n\{ t\in{\cal
   M}~|~u_i(t)=\lambda\} 
$$
Of course, on the component  $u_i(t)=\lambda$ we must have $u_j(t)\neq\lambda$
for any $j\neq i$. As a consequence of (\ref{porcellini}) we have a
representation for the canonical coordinates: 
$$ 
    u_i(t)= t^1-f_i(\tilde{t})
$$

 Now we perform the gauge 
 $\phi(\lambda,u):=\Psi(u(t)) \xi(\lambda,t)$, where $t=t(u)$
{\it locally}. The system (\ref{gaussmanin1}), (\ref{gaussmanin2}) becomes
$$ 
  (U-\lambda) {\partial \phi\over \partial \lambda} = \left({1\over 2}
  +V(u)\right) \phi,
$$
$$
 (U-\lambda){\partial \phi \over \partial u_i} = \left[(U-\lambda) V_i(u) -E_i
 \left({1\over 2} +V(u)\right)\right]\phi.
$$
The matrix $V_i$ is defined by 
$$ 
   V_i:= {\partial \Psi \over \partial u_i} ~\Psi^{-1} 
$$
Equivalently: 
\be 
   {\partial \phi\over \partial \lambda} = \sum_{i=1}^n{B_i\over \lambda-u_i} 
 \phi,~~~~~B_i:=-E_i\left({1\over 2}+V\right),
\label{gaussmanin1u}
\ee
\be
{\partial \phi \over \partial u_i} = \left(V_i-{B_i\over
\lambda-u_i}\right)\phi. 
\label{gaussmanin2u}
\ee

\vskip 0.3 cm 
\noindent
$\bullet$ {\bf Isomonodromy and Inverse Reconstruction} 
\vskip 0.2 cm

The compatibility of the two systems (it is enough to take ${\partial \over
\partial \lambda}{\partial \over
\partial u_i}={\partial \over
\partial u_i} {\partial \over
\partial \lambda}$ which implies  ${\partial \over
\partial u_i}{\partial \over 
\partial u_j}={\partial \over \partial u_j}{\partial \over \partial u_i}$) is
\be
   {\partial_i B_j} =[V_i,B_j]+{[B_i,B_j]\over u_i-u_j} ,~~~i\neq j
\label{CoMp1}
\ee
\be
   \sum_{j=1}^n\left(\partial_iB_j+[B_j,V_i]\right)=0
\label{CoMp2}
\ee
The first is equivalent to 
$$ 
   [U,V_i]=[E_i,V]$$
and the second to
$$ 
  \partial_i V= [V_i,V]
$$
  The compatibility conditions and   
${\partial \Psi \over \partial u_i} = V_i~\Psi$ are the  conditions of
isomonodromicity of section \ref{secdeformata}. 
Therefore, the problem of the WDVV equations is transformed into an
isomonodromy deformation problem for (\ref{gaussmanin1u}) which is equivalent
to the same problem for (\ref{systemz1}). Incidentally, we anticipate that the
systems (\ref{gaussmanin1u}) and (\ref{systemz1}) are connected by a Laplace
transform (see section \ref{Laplace transform}). The boundary value problem for
the fuchsian system  (\ref{gaussmanin1u}) is the classical {\it Riemann
Hilbert Problem}  which we are going to partly solve  in section
\ref{PaInLeVe} for $n=3$ in terms of Painlev\'e transcendents.

\vskip 0.3 cm 
\noindent
$\bullet$ {\bf Monodromy Group} 
\vskip 0.2 cm

Let $x_1(t,\lambda)=x_1(t^1-\lambda,\tilde{t}),...,
x_n(t,\lambda)=x_n(t^1-\lambda,\tilde{t})$ be flat
coordinates obtained from a fundamental matrix solution of the system
(\ref{gaussmanin1u}) (\ref{gaussmanin2u}). They are 
 multi-valued functions of $\lambda$ and of $t$; for  $\lambda=0$
 they represent the flat coordinates for the
intersection form. For a loop $\gamma$ around the 
 discriminant $\Sigma$ the monodromy  of the coordinates is linear,
namely:
$$
 (x_1(t),...,x_n(t))\to (x_1(t),...,x_n(t))M_{\gamma},~~~~M_{\gamma}\in
GL(n,{\bf C})
$$
The image of the representation 
$$ 
  \pi(M\backslash \Sigma,t_0)\to GL(n,{\bf C})
$$
is called the {\it Monodromy group of the Frobenius manifold}. To show this we
 restrict to
the semisimple case. We choose a particular 
 loop around $\Sigma$ starting and ending at
 $u_0=u(t_0)$, defined by the conditions that 
  $\tilde{t}_0$ is fixed and only  
$t^1$ varies according to the rule $t^1=t_0^1-\lambda$. The discriminant locus
$\Sigma$ is reached when one of the $u_i(t^1,\tilde{t}_0)$ becomes
zero: 
$$ 
   u_i(t)=t^1-f_i(\tilde{t})\equiv (t^1_0-\lambda)+f_i(\tilde{t}_0)=0
$$
This means that the discriminant is reached if $\lambda$ equals one of the
canonical coordinates 
$u_1(t_0)=t_0^1-f_1(\tilde{t}_0),...,u_n(t_0)=t_0^1-f_n(\tilde{t}_0)$. On the
other hand,  the coordinates $x_a(t)=x_a(t_0^1-\lambda,t_0^2,...,t_0^n)$,
$a=1,...,n$, come from the solution of (\ref{gaussmanin1u}) when the poles are
$u_1=
u_1(t_0)$, ..., $u_n=u_n(t_0)$. Hence   the monodromy group of the system
(\ref{gaussmanin1u}), which is independent of $u_0=u(t_0)$ for small
deformations of $u$ in a neighbourhood of $u_0$,   
 precisely describes the monodromy group of the Frobenius manifold. Actually,
if $\Phi=[\phi^{(1)},...,\phi^{(n)}]$ is a fundamental matrix of
(\ref{gaussmanin1u}), where $\phi^{(a)}$, $a=1,...,n$ are independent
columns, then
$$ 
   \partial_i x_a= \psi_{i1} \phi_i^{(a)},
$$
This follows from the gauge $\phi=\Psi \xi$, the definition of $\xi$ 
 and the formula $ {\partial \over
\partial u_i}= \psi_{i1}\sum_{\epsilon}
 \psi_{i\epsilon}\partial^{\epsilon}$. If $\lambda$ describes a loop around
 $u_i(t_0)$,  $\Phi$ is changed by a monodromy matrix $R_i$: 
$$ 
   \Phi\to \Phi^{\prime}:=\Phi R_i
$$
Therefore
$$
   \partial_i x^{\prime}_b(u_0)= \partial_i x_a(u_0)~ (R_i)_{ab}
$$
which can be integrated because $R_i$ is independent of $u$:
$$ (x_1^{\prime}(u_0),...,x_n^{\prime}(u_0))=
(x_1(u_0),...,x_n(u_0))~R_i+(c_1,...,c_n)
$$
where the $c_i$'s are constants. 
We can put them equal to zero. Actually, one can verify
 \cite{Dub2} that 
$$ 
  x_a(u,\lambda)={2\sqrt{2}\over 1-d} \sum_{i=1}^n (u_i-\lambda) 
\psi_{i1}(u)
  \phi_i^{(a)}(u,\lambda)
$$
If the  deformation of $u$ is such that two coordinates $u_i$, $u_j$ are
interchanged (for $i\neq j$), one can find an action of the
braid group on the $R_i$'s, and the new monodromy matrices generate the same
group. This part is analyzed in cite{Dub2} and in \cite{guz}.

\vskip 0.3 cm 

Finally,  we give  another formula for $g^{\alpha\beta}$  by
differentiating twice the 
expression 
$$
     E^{\gamma}\partial_{\gamma}F= (2-d)F +{1\over 2} A_{\alpha\beta}t^{\alpha}
     t^{\beta}+ B_{\alpha} t^{\alpha} +C
$$ 
which is the  quasi-homogeneity of $F$ up to quadratic terms. By recalling that 
$E^{\gamma}= (1-q_{\gamma})t^{\gamma}+r_{\gamma}$ and that
$\partial_{\alpha}\partial_{\beta} \partial_{\gamma} F=c_{\alpha\beta\gamma}$
we obtain 
\be 
   g^{\alpha\beta}(t)=
   (1+d-q_{\alpha}-q_{\beta})\partial^{\alpha}\partial^{\beta} F(t)
   +A^{\alpha\beta} 
\label{dall'MSRI}
\ee
where $\partial^{\alpha}=\eta^{\alpha \beta}\partial_{\beta}$, $A^{\alpha
\beta} = \eta^{\alpha \gamma}\eta^{\beta\delta} 
A_{\gamma\delta} $.


\section{ Reduction of the Equations of Isomonodromic Deformation 
to a Painlev\'e 6 Equation}\label{PAPAPA}

  The compatibility conditions
(\ref{CoMp1}),(\ref{CoMp2}) are equivalent to 
the differential equations for  $V$, $V_i$, $\Psi$ whose solution allows us to
do the inverse local reconstruction of the manifold. Solving such equations is
equivalent to solving the boundary value problem (b.v.p.), 
provided that we find a mean
to parameterize the integration constants in terms of the monodromy data
defining the b.v.p..  

For $n=3$ such 
equations were reduced \cite{Dub1} 
to a particular form of the VI Painlev\'e equation. This is the
first step towards the solution of the b.v.p. and the inverse
reconstruction. In chapter  \ref{PaInLeVe} we'll see how to parameterize the
solutions of the Painlev\'e equation in terms of monodromy data. Finally, in
chapter \ref{inverse reconstruction} we'll apply the results to the explicit
inverse reconstructions of some 3-dimensional Frobenius manifolds. 
\vskip 0.2 cm

   Let 
$$ 
   \hat{\mu}= \hbox{ diag}(\mu,0,-\mu)
$$
We consider the auxiliary system 
$$
 {\partial \Phi\over \partial \lambda}= \sum_{i=1}^n~{{\cal B}_i\over
 \lambda-u_i} 
 ~\Phi 
$$
$$
    {\partial \Phi\over \partial u_i}= \left(V_i-{{\cal B}_i\over \lambda-u_i} 
\right)~\Phi
$$
where 
   $$ {\cal B}_i:= -E_i ~V$$
The compatibility conditions expressed  in terms of $B_i$ are (\ref{CoMp1}),
(\ref{CoMp1}) where ${\cal B}_i$ is substituted by $B_i$. They are again
equivalent to  
$ [U,V_i]=[e_i,V]$ and $ 
  \partial_i V= [V_i,V]$, plus the condition $\partial_i \phi_0=V_i\phi_0$.
Therefore, we can study the isomonodromy dependence on $u$ of  
the auxiliary system instead of (\ref{gaussmanin2u}) 
(\ref{gaussmanin2u}). 

Let $n=3$.   
 We observe that 
$$ \sum_{i=1}^n~{{\cal B}_i\over
 \lambda-u_i} = -(\lambda - U)^{-1} ~V=  -(\lambda - U)^{-1} ~\phi_0 \hat{\mu}
\phi_0^{-1} 
$$
We put $X(\lambda):=
\phi_0^{-1}~\Phi(\lambda)$ and rewrite:
$$ 
   {\partial X\over \partial \lambda}=-\mu ~v(\lambda,u) ~\hbox{diag}(1,0,-1)~X
$$
where
$$
v(\lambda,u):=
 \phi_0^{-1}(\lambda-U)^{-1} \phi_0 = -\mu \sum_{i=1}^3 {{\cal A}_i \over
 \lambda-u_i},
$$
$$
{\cal A}_i:=  \phi_0^{-1}E_i \phi_0 
$$
Let $X$ be the column vector 
$$
  X=\pmatrix{ X_1 \cr X_2 \cr X_3 \cr} ~~\Rightarrow ~~~\hbox{diag}(1,0,-1)
  ~X= \pmatrix{ X_1 \cr 0 \cr -X_3 \cr}
$$
Therefore 
\be
 {\partial \over \partial \lambda} \pmatrix{X_1 \cr X_2 \cr} = 
-\mu ~ \pmatrix{ v_{11} & -v_{13} \cr
                 v_{31} & v_{33}  \cr }  \pmatrix{X_1 \cr X_2 \cr}=
-\mu ~\sum_{i=1}^3 {A_i\over \lambda -u_i} ~\pmatrix{X_1 \cr X_2 \cr}
\label{221}
\ee
where
  \be
     A_i:= \pmatrix{ \phi_{i1,0}\phi_{i3,0} & -\phi_{i3,0}^2 \cr
                      \phi_{i1,0}^2         & \phi_{i1,0}\phi_{i3,0} \cr
                    } 
,~~~~A_1+A_2+A_3 = \pmatrix{ 1 & 0 \cr 0 & -1 \cr}
\label{definite in extremis}
\ee
To write the
$A_i$'s, we have taken into account the relations 
$\phi_0^T \phi_0=\eta$ and the choice
$$
(\phi_{12,0},\phi_{22,0},\phi_{32,0})= \pm i 
(\phi_{21,0}\phi_{33,0}
-\phi_{23,0}\phi_{31,0},\phi_{13,0}\phi_{31,0}-\phi_{11,0}\phi_{33,0}, \phi_{11,0}\phi_{23,0}-
\phi_{13,0}\phi_{21,0}) 
$$
We also obtain:
\be
  {\partial \over \partial u_i}\pmatrix{X_1 \cr X_2 \cr} = \mu {A_i \over
  \lambda-u_i}~\pmatrix{X_1 \cr X_2 \cr} 
\label{222}
\ee
The second component of $X$  is obtained by quadratures: 
$$
 {\partial \over \partial \lambda}X_2= \mu~(v_{21}~X_1-v_{23}~X_3)
$$
$$
   {\partial \over \partial u_i}X_2= {\mu\over \lambda-u_i} [\phi_{i1,0}\phi_{i2,0} ~X_1
   - 
      \phi_{i2,0}\phi_{i3,0} ~X_3 ]
$$
The reduced $2\times 2$ system (\ref{221}) (\ref{222}) 
is solved (see \cite{JM1} \cite{Dub1}) by
introducing the following coordinates $q(u)$, $p(u)$ 
in the space of matrices $A_i$ modulo
diagonal conjugation: $q$ is the root of 
$$
\left(\sum_{i=1}^3 {A_i\over \lambda-u_i}\right)_{12} =0
$$
and
$$
p:= \left(\sum_{i=1}^3 {A_i\over q-u_i}\right)_{11}
$$
The entries of the $A_i$'s are re-expressed as follows:
$$
  \phi_{i1,0}~\phi_{i3,0}= -{q-u_i\over 2 \mu P^{\prime}(u_i)} \left[
                             P(q)p^2+{2\mu\over q-u_i} P(q)p +\mu^2 (q+2u_i 
                            -\sum_{j=1}^3~ u_j \right]
$$
$$
  \phi_{13,0}^2= -k{q-u_i\over  P^{\prime}(u_i)}
$$
$$
  \phi_{i1,0}^2=  -{q-u_i\over 4 \mu P^{\prime}(u_i)~k} \left[
                             P(q)p^2+{2\mu\over q-u_i} P(q)p +\mu^2 (q+2u_i 
                            -\sum_{j=1}^3~ u_j \right]^2
$$
Here $k$ is a parameter, $P(z)=(z-u_1)(z-u_2)(z-u_3)$. The compatibility of
(\ref{221}) (\ref{222}) is
\be
  {\partial q \over \partial u_i}= {P(q)\over P^{\prime}(u_i)}\left[ 
                                2p+{1\over q-u_i}\right] 
\label{pain1}
\ee
\be
    {\partial p \over \partial u_i}=-{P^{\prime}(q) p^2+(2q+u_i-\sum_{j=1}^3
    u_j)p +\mu (1-\mu)\over P^{\prime}(u_i)}
\label{pain2}
\ee
$$
   {\partial \ln k \over \partial u_i}=
 (2 \mu -1) {q-u_i\over P^{\prime}(u_i)} 
$$
In the variables 
$$ 
  x={u_3-u_1\over u_2-u_1},~~~~y={q-u_1\over u_2-u_1}
$$
the system (\ref{pain1}), (\ref{pain2}) becomes a special case of the
VI Painlev\'e  equation, with the following choice of the parameters (in the
standard notation of \cite{IN}): 
$$
  \alpha={(2\mu-1)^2\over 2},~~~\beta=\gamma=0,~~~\delta={1\over
  2}
$$
Namely:
\be
{d^2y \over dx^2}={1\over 2}\left[ 
{1\over y}+{1\over y-1}+{1\over y-x}
\right]
           \left({dy\over dx}\right)^2
-\left[
{1\over x}+{1\over x-1}+{1\over y-x}
\right]{dy \over dx}$$
$$
+{1\over 2}
{y(y-1)(y-x)\over x^2 (x-1)^2}
\left[
(2\mu-1)^2+{x(x-1)\over (y-x)^2}
\right]
,~~~~~\mu \in {\bf C}
\label{PVImu}
\ee
In the following, this equation will be referred to as $PVI_{\mu}$. From a solution
of ($PVI_{\mu}$) one can reconstruct 
$$ 
  q= (u_2-u_1)~y\left({u_3-u_1\over u_2-u_1}\right) +u_1
$$
$$
  p= {1\over 2} {P^{\prime}(u_3) \over P(q) } ~y^{\prime}\left({u_3-u_1\over
  u_2-u_1}\right) -{1\over 2} {1 \over q-u_3}
$$

From the very definition of $q$ we have:
$$
  y(x)= {x R\over x(1+R)-1},~~~~R:={(A_1)_{12}\over (A_2)_{12}}\equiv
  \left({\phi_{13,0}\over \phi_{23,0}}\right)^2 
$$


\section{Conclusion}

 A semisimple 
Frobenius manifold can be locally parametrized by a set of monodromy
 data. The local parameterization is given by the formulae of section
\ref{Semisimple Frobenius manifolds} and  by (\ref{Fu}) (\ref{tu}). 
The matrices $\phi_p$ are obtained 
  solving a
 boundary value problem or, equivalently, the 
 equations of isomonodromic deformation ${\partial V \over \partial u_i} =
 [V_i,V]$, ${\partial \phi_0 \over \partial u_i} = V_i \phi_0$. 
When $n=3$ the problem is reduced to a Painlev\'e 6 equation.



\chapter{ Quantum Cohomology of Projective spaces}\label{cap2}

 In this chapter we introduce the 
FM called quantum
 cohomology of the projective space $CP^d$ and we describe its connections to
 enumerative geometry.  In the last section we will  prove
 that the radius of convergence of the
 famous Kontsevich's solution of WDVV equations for $CP^2$ corresponds to a
 singularity in the change from flat to canonical coordinates.

 We start by introducing a
 structure of Frobenius algebra on the cohomology $H^{*}(X,{\bf C})$ of a
 closed oriented manifold $X$ of dimension $d$  such that 
$$ 
    H^{i}(X,{\bf C})=0 ~~\hbox{ for $i$ odd}
$$
then
$$ 
    H^{*}(X,{\bf C})=\otimes_{i=0}^d ~H^{2i}(X,{\bf C}).
$$
For brevity we omit ${\bf C}$ in $H$. 
We realize  $H^{*}(X)$ by classes of closed differential forms. The
unit element is a 0-form $e_1\in  H^{0}(X)$. Let denote by 
$\omega_{\alpha}$ a form in 
$H^{2q_{\alpha}}(X)$, where $q_1=0$, $q_2=1$, ...,
$q_{d+1}=d$. The product of two forms $\omega_{\alpha}$, $\omega_{\beta}$  
 is defined by the wedge product  $ \omega_{\alpha} \wedge\omega_{\beta}\in
H^{2(q_{\alpha}+q_{\beta})}(X)$ and  
the
 bilinear form is 
$$ 
  <\omega_{\alpha},\omega_{\beta}>:=\int_X ~ \omega_{\alpha}
  \wedge\omega_{\beta} \neq 0 ~\Longleftrightarrow ~q_{\alpha}+q_{\beta}=d
$$
It is not degenerate by Poincar\'e duality  and of course only 
$q_{\alpha}+q_{d-\alpha+1}=d$.

\vskip 0.3 cm

Let $X=CP^d$. 
Let $e_1=1\in H^0(CP^d)$, $ e_2\in H^2(CP^d)$,  ...,
 $e_{d+1}\in H^{2d}(CP^d)$ be a basis in $ H^{*}(CP^d)$. For a suitable
normalization we have 
$$ 
   (\eta{\alpha\beta}):=(<e_{\alpha},e_{\beta}>)= 
\pmatrix{             &   &         & & 1 \cr 
                      &   &         &1& \cr 
                      &   & \adots  & & \cr
                      & 1 &         & & \cr
                    1 &   &         & &  \cr
                   } 
$$
The multiplication is $$e_{\alpha}\wedge e_{\beta}= e_{\alpha+\beta-1}.$$ We
observe that it can also be written as 
$$ 
  e_{\alpha}\wedge e_{\beta}=c_{\alpha\beta}^{\gamma} e_{\gamma} ,~~\hbox
  { sums on $\gamma$} 
$$
where 
    $$\eta_{\alpha\delta}c_{\beta\gamma}^{\delta}:= {\partial^3 F_0(t) \over 
   \partial t^{\alpha}  \partial t^{\beta}  \partial t^{\gamma}}
$$
$$
      F_0(t):= {1\over 2} (t^1)^2 t^n +{1\over 2} t^1 \sum_{\alpha=2}^{n-1}
  t^{\alpha} t^{n-\alpha+1}
$$
$F_0$ is the trivial solution of WDVV equations. We can construct a trivial
Frobenius manifold whose points are $t:=\sum_{\alpha=1}^{d+1} t^{\alpha}
e_{\alpha}$. It has tangent space $H^{*}(CP^d)$ at any $t$. By {\it quantum
cohomology} of $CP^d$ (denoted by $QH^{*}(CP^d)$ 
we mean a Frobenius manifold whose structure is
specified by 
$$
    F(t)=F_0(t)+\hbox{ analytic perturbation}
$$
 This manifold has therefore tangent spaces $T_tQH^{*}(CP^d)= H^{*}(CP^d)$,
 with the same $<.,.>$ as above, 
 but the multiplication is a deformation, depending on $t$, 
of the wedge product (this is the
 origin of the adjective ``quantum'').

\section{ The case of $CP^2$}
 
 We restrict to  $CP^2$. In this case 
$$ 
  F_0(t)={1\over 2}\left[ (t^1)^2 t^3+t^1(t^2)^2\right]
$$ 
which generates the product for the basis $e_1=1\in H^0$, $e_2\in H^2$,
$e_3\in H^4$. The deformation was introduced by Kontsevich \cite{KM}.

\subsection{Kontsevich's solution}

 The WDVV equations for $n=3$ variables have  solutions  
 $$F(t_1, t_2,t_3)= F_0(t_1,t_2,t_3)+f(t_2,t_3)$$
where 
\be
   f_{222}f_{233}+f_{333}=(f_{223})^2 
\label{WDVV1special}
\ee
with the notation $f_{ijk}:={\partial^3 f \over \partial t_i \partial t_j 
\partial t_k}$. As for notations, the variables $t_j$ are flat coordinates in  
the Frobenius manifold associate to $F$. They should be written with upper 
indices, but we use the lower for convenience now. 

Let $N_k$ be the number of rational curves $CP^1 \to CP^2$ of degree $k$ 
through $3k-1$ generic points. Kontsevich \cite{KM} found the solution
\be
  f(t_2,t_3) = {1\over t_3} \varphi(\tau),~~~~
\varphi(\tau)=\sum_{k=1}^{\infty} A_k \tau^k,~~~
~~\tau= t_3^3~e^{t_2}
\label{Konts}
\ee
 where
 $$
 A_k= {N_k\over (3k-1)!}
$$
We note that this solution has precisely the form of the general solution of
 the WDVV eqs. for $n=3$, $d=2$ and $r_2=3$.
  If we put $\tau=e^X$ we  rewrite (\ref{WDVV1special}) as follows
$$
\Phi(X):=\varphi(e^X)=\sum_{k=1}^{\infty} ~A_k ~e^{kX},
$$
\be
    -6 \Phi+33\Phi^{\prime} -54 \Phi^{\prime \prime} -  
(\Phi^{\prime \prime})^2
      + \Phi^{\prime \prime \prime} \left(
             27+2 \Phi^{\prime} -3 \Phi^{\prime \prime}
             \right)  =0
\label{WDVVK}
\ee
 The prime stands for the derivative w.r.t $X$. 
If we fix $A_1$, the above (\ref{WDVVK}) determines the $A_k$ uniquely. Since 
$N_1=1$, we fix 
$$A_1={1\over 2}.
$$
Then  (\ref{WDVVK})  yields  the recurrence relation
\be
  A_k= \sum_{i=1}^{k-1} \left[ {A_1 A_{k-i} ~i(k-i)\bigl(
                                (3i-2)(3k-3i-2)(k+2)+8k-8
\bigr)\over 6(3k-1)(3k-2)(3k-3) }\right]
\label{recurrence}
\ee

\subsection{ Convergence of Kontsevich solution}

 The convergence of (\ref{Konts}) was studied by Di Francesco and Itzykson 
\cite{DI}. They proved that 
$$
 A_{k}= b~ a^k ~k^{-{7\over 2}}~\left(1+O\left({1\over k}\right) \right),
~~~~~k\to \infty
$$
and numerically computed 
$$
  a=0.138, ~~~~b=6.1~~~.
$$
 The result implies that $\varphi(\tau)$ converges in a neighbourhood of 
$\tau=0$ with radius of convergence ${1\over a}$. 
We will return later on the  numerical evaluation of these numbers. 

The proof of \cite{DI} is divided in two steps. The first is based on the 
  relation (\ref{recurrence}), to prove that 
$$
    A_k^{1\over k}\to a \hbox{ for } k \to \infty, ~~~~{1\over 108}
                                                              < a < {2\over 3}
$$
 $a$ is real positive because the $A_k$'s are such. It follows that  we 
can rewrite  
$$
  A_k= b a^k ~k^{\omega}~ \left(1+O\left({1\over k}\right) \right),~~~~
\omega
       \in {\bf R}
$$
The above estimate implies that $\varphi(\tau) $ has the radius of convergence 
${1\over a}$. 
The second step is the determination of $\omega$ making use of the 
differential equation (\ref{WDVVK}). Let's write 
$$
   A_k:=C_k~a^k
$$
$$
   \Phi(X)= \sum_{k=1}^{\infty} ~A_k ~e^{kX}= \sum_{k=1}^{\infty} ~C_k 
                                             ~e^{k(X-X_0)},~~~~X_0:=\ln{1\over
                                             a}
$$
The inequality ${1\over 108} < a < {2\over 3}$ implies that $X_0>0$. 
The series converges at least for  $\Re X <X_0$.  
To determine $\omega$ we divide $\Phi(X)$ into a regular part at $X_0$ and 
a singular one. Namely
$$
  \Phi(X)= \sum_{k=0}^{\infty} d_k (X-X_0)^k+ 
           (X-X_0)^{\gamma}~ \sum_{k=0}^{\infty} e_k (X-X_0)^k,~~~~
      \gamma>0 ,~~~\gamma\not\in {\bf N},
$$
$d_k$ and $e_k$ are coefficients. By substituting into (\ref{WDVVK}) we 
see that the equation is satisfied only if  $\gamma={5\over 2}$. Namely:
$$
 \Phi(X)= d_0+d_1(X-X_0)+d_2(X-X_0)^2+e_0(X-X0)^{5\over 2}+...
$$
This implies that $\Phi(X)$,  $\Phi^{\prime}(X)$ and  $\Phi^{\prime
\prime}(X)$ exist at $X_0$ but $\Phi^{\prime\prime\prime}(X)$ diverges like
\be
  \Phi^{\prime\prime\prime}(X)\asymp {1\over \sqrt{X-X_0}},~~~~X\to X_0
\label{behaviourhh}
\ee
On the other hand $
 \Phi^{\prime\prime\prime}(X)$ behaves like the series
$$ \sum_{k=1}^{\infty}~b~k^{\omega+3}~
       e^{k(X-X_0)}, ~~~~\Re (X-X_0)<0
$$
Suppose $X\in {\bf R}$, $X<X_0$. Then $\Delta:=X-X_0<0$
 and the above series is 
$$
 {b\over |\Delta|^{3+\omega}}~\sum_{k=1}^{\infty} (|\Delta|k)^{3+\omega}
 e^{-|\Delta|k}
\sim {b\over |\Delta|^{3+\omega}}~\int_{0}^{\infty} dx~ x^{3+\omega} e^{-x} 
$$
It follows from (\ref{behaviourhh}) that this
  must diverge like $\Delta^{-{1\over 2}}$, and thus $\omega=-{7\over 2}$ (the
  integral remains finite). 

\vskip 0.2 cm
As a consequence of  (\ref{WDVVK}) and of the divergence of $ \Phi^{\prime
\prime\prime}(X)$ 
$$   
             27+2 \Phi^{\prime}(X_0) -3 \Phi^{\prime \prime}(X_0)
               =0
$$

\vskip 0.2 cm

\subsection{Small Cohomology}

 We realize  $QH^{*}(CP^2)$ as the set of points $t=t^1e_1+t^2e_2+t^3e_3\in
 H^{*}(CP^2)$ (we restore the upper indices only in this formula) 
such that the tangent space at
 $t$ is again $H^{*}(CP^2)$ with the product $e_{\alpha}\cdot e_{\beta}=
 c_{\alpha\beta}^{\gamma}(t) e_{\gamma}$, where 
$$ 
  c_{\alpha\beta\gamma}=\eta_{\alpha\delta}c_{\beta\gamma}^{\delta}=
\partial_{\alpha} \partial_{\beta} \partial_{\gamma} F(t)$$
$F(t)$ being Kontsevich solution. We also observe that $ F(t_1,t_2+2\pi i,
t_3) $ differs from $F(t_1,t_2,t_3)$ by quadratic terms:
$$
   F(t_1,t_2+2\pi i, t_3)=F(t_1,t_2,t_3)+2\pi i t_1 t_2 -2\pi^2 t_1.
$$
Thus $QH^{*}(CP^2)$ is a sub-manifold of the quotient $H^{*}(CP^2)/2\pi i
H^2(CP^2,{\bf Z})$  because $e_2\in H^2(CP^2,{\bf Z})$. 

\vskip 0.2 cm
Suppose $t_3\to 0$ and $t_2\to -\infty$. From Kontsevich solution we get 
$$
  e_2\cdot e_2 = (t_3e^{t_2}+...)~e_1+\left({1\over 2}t_3^2 e^{t_2}+...\right)
  ~e_2+e_3,
$$
$$
 e_2\cdot e_3= (e^{t_2}+...)~e_1+(t_3e^{t_2}+...) ~e_2,
$$
$$
   e_3\cdot e_3=\left({t_3^2\over 2}
   e^{2t_2}+...\right)~e_1+(e^{t_2}+...)~e_2.
$$
 It follows that for $t_2\to -\infty$ we recover the ``classical'' wedge
 product. On the other hand, for $t_3=0$:
$$
  e_2\cdot e_2=e_3,~~~
   e_2\cdot e_3=q~ e_1,~~~
    e_3\cdot e_3= q~ e_2,~~~~
~~ q:=e^{t_2}
$$
The limit $t_3=0$ is called {\it small cohomology}. By the identification
$e_1\mapsto 1$, $e_2\mapsto x$, $e_3\mapsto x^2$, the algebra $T_{(t_1,t_2,0)}
QH^{*}(CP^2)$ is isomorphic to ${\bf C}[x]/(x^3=q)$.


\section{ The case of $CP^d$}\label{solutionQH}
For $d=1$ the deformation is given by 
$$ 
  F(t){1\over 2} t_1^2 t_2+e^{t_2}
$$
For any $d\geq 2$, the deformation is given by the following solution of the
WDVV equations \cite{KM} \cite{Manin}:
$$
  F(t)=F_0(t)+\sum_{k=1}^{\infty}\left[ \sum_{n=2}^{\infty}
          \tilde{\sum_{\alpha_1,...,\alpha_n}} 
          ~{N_k(\alpha_1,...,\alpha_n)\over n!} ~t_{\alpha_1}...t_{\alpha_n}
\right]e^{kt_2}
$$
where
$$ 
     \tilde{\sum_{\alpha_1,...,\alpha_n}}:=
     \sum_{\alpha_1+...+\alpha_n=2n+d(k+1) +k-3} 
$$ 
Here $N_k(\alpha_1,...,\alpha_n)$ 
is the number of rational curves $CP^1\to CP^d$ of degree $k$
through $n$ projective subspaces of codimensions
$\alpha_1-1,...,\alpha_n-1\geq2$ in general position. 
In particular,  there is one line through
two points, then 
$$
   N_1(d+1,d+1)=1
$$
Note that in Kontsevich solution $N_k=N_k(d+1,d+1)$. 

 In flat coordinates the 
 {\it Euler vector field} is 
$$ 
   E= \sum_{\alpha \neq 2} ~(1-q_{\alpha})t^{\alpha}{\partial \over 
\partial t^{\alpha} }
+~k {\partial \over \partial t^2}
$$
$$ q_1=0,~q_2=1,~q_3=2,~...,~q_{k}=k-1
$$
and 
$$\hat{\mu}= \hbox{diag}(\mu_1,...,\mu_k)=\hbox{diag}(-{d\over
2},-{d-2 \over 2},...,{d-2 \over 2},{d \over 2}),
~~~~~\mu_{\alpha}=q_{\alpha}-{d\over 2}$$


\section{ Nature of the singular point $X_0$}\label{Nature of the singular
point}

 We ask the question whether the singularity $X_0$ in Kontsevich solution for
$CP^2$   corresponds   
to the fact that two canonical coordinates $u_1$,
 $u_2$, $u_3$ merge. Actually, we know that the structure of the semisimple
manifold may become singular in such points because the solutions of the
boundary value problem (or, equivalently, of the the equations of isomonodromic
deformation) are meromorphic on the universal covering of ${\bf
C}^n\backslash\hbox{diagonals}$ and are multivalued if  $u_i-u_j$ ($i\neq j$) 
describes a loop around zero.
  In this section we restore the upper indices for the flat
coordinates $t^{\alpha}$.

The canonical coordinates can be computed from the  
intersection form. We recall that the flat metric is 
$$\eta =(\eta^{\alpha\beta}):=\pmatrix{0 & 0 & 1 \cr
                 0 & 1 & 0 \cr
                 1 & 0 & 0 \cr}
$$
The intersection form is given by the formula (\ref{dall'MSRI}):
$$
  g^{\alpha\beta}= (d+1-q_{\alpha}-q_{\beta})~\eta^{\alpha \mu} \eta^{\beta\nu}
  \partial_{\mu}
  \partial_{\nu} F+ A^{\alpha\beta}, ~~~\alpha,\beta=1,2,3,
$$
where $d=2$ and the {\it charges} are $q_1=0$, $q_2=1$, $q_3=2$. The matrix
$A^{\alpha\beta} $ appears in the action of the  Euler vector field
$$
 E:=t^1\partial_1 +3 \partial_2-t^3\partial_3
$$
on $F(t^1,t^2,t^3)$:
$$
  E(F)(t^1,t^2,t^3)=(3-d) F(t^1,t^2,t^3)+ A_{\mu \nu} t^{\mu} t^{\nu}
                   \equiv
                      F(t^1,t^2,t^3)+ 3t^1t^2
$$
Thus 
$$ 
   (A^{\alpha\beta})=(\eta^{\alpha\mu}\eta^{\beta\nu} A_{\mu \nu})
                     =\pmatrix{0 & 0 & 0 \cr
                               0 & 0 & 3 \cr
                               0 & 3 & 0 \cr}
$$
 After the above preliminaries, we are able to compute the intersection form:
$$
(g^{\alpha\beta})= \pmatrix{
                            {3\over [t^3]^3}[2\Phi - 9 \Phi^{\prime}
                            +9\Phi^{\prime\prime}  ]
                          &
                             {2\over [t^3]^2}[3\Phi^{\prime\prime} -
                            \Phi^{\prime} ] 
 &
             t^1  \cr\cr
        {2\over [t^3]^2}[3\Phi^{\prime\prime} -
                            \Phi^{\prime} ] 
&
 t^1+{1\over t^3} \Phi^{\prime\prime} 
&
   3
\cr\cr
t^1
&
3
&
-t^3
\cr
}
$$
The canonical coordinates are roots of
$$
 \det((g^{\alpha\beta}-u \eta)=0
$$
This is the polynomial
$$
u^3- \left(3t^1+{1\over t^3} \Phi^{\prime\prime}\right)~u^2
-\left(
-3[t^1]^2-2{t^1\over t^3} \Phi^{\prime\prime}+{1\over [t^3]^2} 
(9\Phi^{\prime\prime}+15\Phi^{\prime}-6\Phi)
\right)        ~u    
+P(t,\Phi)$$
where
$$
 P(t,\Phi)={1\over [t^3]^3}\left( -9t^1t^3 \Phi^{\prime\prime}+243 \Phi^{\prime\prime}- 
243\Phi^{\prime}+6\Phi
 \Phi^{\prime} \right.$$
$$\left. -9( \Phi^{\prime\prime})^2 +6t^1t^3\Phi+[t^1]^2[t^3]^2
 \Phi^{\prime\prime} -3  \Phi^{\prime}\Phi^{\prime\prime}+[t^1]^3[t^3]^3 -4
 ( \Phi^{\prime})^2 +54 \Phi-15 t^1t^3\Phi^{\prime}\right)
$$
It follows that 
$$
   u_i(t^1,t^3,X)=t^1+{1\over t^3} {\cal V}_i(X)
$$
$ {\cal V}_i(X)$ depends on $X$ through $\Phi(X)$ and derivatives. We also
observe that 
$$
  u_1+u_2+u_3=3t^1+{1\over t^3}\Phi^{\prime\prime}(X)
$$

\vskip 0.2 cm
We verify numerically that $u_i\neq u_j$ for $i\neq j$ at $X=X_0$. In order to
do this we need to compute $\Phi(X_0)$, $\Phi^{\prime}(X_0)$ ,
$\Phi^{\prime\prime}(X_0)$ in the following approximation
$$
  \Phi(X_0) \cong \sum_{k=1}^{N}~A_k~{1\over a^k},~~~
  \Phi^{\prime}(X_0) \cong \sum_{k=1}^{N}~k~A_k~{1\over a^k},~~~
~~~\Phi^{\prime\prime}(X_0) \cong \sum_{k=1}^{N}~k^2~A_k~{1\over a^k},
$$
 In our computation we fixed $N=1000$ and we computed the $A_k$,
 $k=1,2,...,1000$  exactly using
 the relation (\ref{recurrence}). Then we computed $a$ and $b$ by
 the least squares method. For large $k$, say for $k\geq N_0$, we assumed that 
\be
 A_k\cong b a^k k^{-{7\over 2}}
\label{linehh}
\ee
which implies
$$\ln(A_k~k^{7\over 2} )
                                \cong (\ln a)~k +\ln b
$$
The corrections to this law are $O\left({1\over k}\right)$. 
This is the line to fit the data $k^{7\over 2}A_k$. Let 
$$
 \bar{y}:={1\over N-N_0+1} \sum_{k=N_0}^N~\ln(A_k~k^{7\over 2}),~~~~
  \bar{k} :={1\over N-N_0+1}\sum_{N_0}^N~k.
$$
By the least squares method 
$$
  \ln a = {\sum_{k=N_0}^N ~(k-\bar{k}) (\ln(A_k~k^{7\over 2} )-\bar{y}) 
           \over \sum_{k=N_0}^N~(k-\bar{k})^2}, ~\hbox{ with error }
           \left({1\over \bar{k}^2}\right)
$$
$$
  \ln b = \bar{y} -(\ln a)~\bar{k}~\hbox{ with error }
           \left({1\over \bar{k}}\right)
$$
 For $N=1000$, $A_{1000}$ is of the order $10^{-840}$. In our computation we
 set the accuracy to $890$ digits. Here is the results, for three choices of
 $N_0$. The result should improve as $N_0$ increases, since the approximation 
 (\ref{linehh}) becomes better. 
$$
  N_0=500,~~~~a= 0.138009444...,~~~b=6.02651...
$$
$$
 N_0=700,~~~~ a= 0.138009418...,~~~b=6.03047...
$$
$$
 N_0=900,~~~~ a= 0.138009415...,~~~b=6.03062...
$$
 It follows that (for $N_0=900$) 
$$
  \Phi(X_0)=4.268908...~,~~~~\Phi^{\prime}(X_0)=5.408...~,~~~~
\Phi^{\prime\prime}(X_0)=12.25... 
$$ 
With these values we find
$$   
             27+2 \Phi^{\prime}(X_0) -3 \Phi^{\prime \prime}(X_0)=1.07...,
$$
 But the above should  vanish! The reason why this does not happen 
is that $ \Phi^{\prime \prime}(X_0)= 
\sum_{k=1}^{N}~k^2~A_k~{1\over a^k}$ converges slowly. To obtain a better
approximation we compute it numerically as 
$$
    \Phi^{\prime \prime}(X_0)={1\over 3}(27+2 \Phi^{\prime }(X_0))
                             ={1\over 3}(27+2\sum_{k=1}^{N}~k~A_k~{1\over
                             a^k})=12.60...
$$
 Substituting into $g^{\alpha\beta}$ and setting $t^1=t^3=1$ we find
$$
  u_1\approx 22.25...~,~~~u_2\approx -(3.5...)-(2.29...)i~,~~~
  u_3 =\bar{u}_2
$$
Here the bar means conjugation. Thus, with a sufficient accuracy, we have
proved that $u_i\neq u_j$ for $i\neq j$. 

\vskip 0.3 cm
 Finally, 
we prove that the singularity is a singularity for the change of coordinates 
$$
   (u_1,u_2,u_3)\mapsto (t^1,t^2,t^3)
$$
 Recall that 
$$
   {\partial u_1\over \partial t^{\alpha}}= {\psi_{i\alpha}\over \psi_{i1}} $$
This may become infinite if $\psi_{i1}=0$ for some $i$. In our case 
$$ u_1+u_2+u_3= 3 t^1+ {1\over t^3}\Phi(X)^{\prime\prime},~~~{\partial X\over
\partial 
t^1} =0,~~~ {\partial X\over \partial
t^2} =1,~~~ {\partial X\over \partial
t^3} ={3\over t^3}
$$
and we compute
$$
   {\partial\over \partial t^1}(u_1+u_2+u_3)= 3,$$
$$
   {\partial\over \partial t^2}(u_1+u_2+u_3)= {1\over
   t^3}\Phi(X)^{\prime\prime\prime}  ,$$
$$
   {\partial\over \partial t^3}(u_1+u_2+u_3)= -{1\over
   [t^3]^2}\Phi(X)^{\prime\prime} +{3\over [t^3]^2}\Phi(X)^{\prime\prime\prime}
    .
$$
Thus, we see that the change of coordinates is singular because 
 both ${\partial\over \partial t^2}(u_1+u_2+u_3)$ and 
$ {\partial\over \partial t^3}(u_1+u_2+u_3)$ behave like 
$\Phi(X)^{\prime\prime\prime}\asymp{1\over \sqrt{X-X_0}}$ for $X\to X_0$. 

I thank  A. Its and P. Bleher for suggesting me to try 
 the computations of this section.


\chapter{ Inverse Reconstruction of 2-dimensional FM}\label{rh2}

 This chapter explains in a didactic way the process of inverse reconstruction
 of a semisimple FM for $n=2$ through the formulae (\ref{tu}), (\ref{Fu}).

\section{ Exact Solution and Monodromy data}

The coefficients of the
 system $dY/dz=[U+V/z]Y$ are {\it necessarily}:  
$$ 
  V(u)=\pmatrix{ 0 & i{\sigma\over 2} \cr
                 -i{\sigma\over 2} & 0 \cr},~~~~U=\hbox{diag}(u_1,u_2)
$$
Here $u=(u_1,u_2)$, and $V$ is independent of $u$.  It has
 the diagonal form 
$$
\Psi^{-1} V \Psi = \hbox{diag}\left({\sigma\over 2}, -{\sigma \over
2}\right)=:\hat{\mu},
$$
where 
$$\Psi(u)= \pmatrix{{1\over 2 f(u)}&  f(u) \cr
                    {1\over 2 i f(u)} & i f(u) \cr
               },
~~~~\Psi^T\Psi=\eta:= \pmatrix{0 & 1 \cr  1 & 0 \cr}. 
$$

\vskip 0.2 cm 

The description of the Stokes' phenomenon for 
\be
{dY\over dz}=\left[\pmatrix{u_1 & 0 \cr 0 & u_2 \cr} +\pmatrix{ 0 & i
{\sigma\over 2} \cr -i {\sigma\over 2} & 0 \cr } \right] Y
\label{systemz1WOW2}
\ee
requires the stokes rays
$$
   R_{12}=-i\rho(u_1-u_2),~~\rho>0;~~~~R_{21}=-R_{12}
$$
Let $l$ be an oriented line through the origin, 
 not containing the Stokes' rays and having
$R_{12}$ to the left. $l=l_{+}+l_{-}$.  Then
$$
   Y_L(z,u)=Y_R(z,u)~S ~~\hbox{ from the side of } l_{+}
$$
$$
 Y_L(z,u)=Y_R(ze^{-2 \pi i },u)~S^T ~~\hbox{ from the side of } l_{-}
$$
$$
  S= \pmatrix{1 & s \cr 0 & 1 \cr},~~~~s\in {\bf C}
$$
Let $\Pi_R$ be the half plane to the right of $l$. At the origin 
\be
  Y(z,u)=\Phi(z,u) z^{\hat{\mu}} z^R,~~~~\Phi(z,u)=
                               \sum_{k=0}^{\infty} \phi_k(u)
                               z^k,~~~\phi_0=\Psi,  
\label{stai barando!}
\ee
$$
Y_0(z,u)=Y_R(z,u)~ C ~~\hbox
  { in } \Pi_R
$$
$C$ is the central connection matrix. Then
$$
   S^T S^{-1} = C e^{2 \pi i \pmatrix{{\sigma\over 2} & 0 \cr 
                                   0 & -{\sigma \over 2} \cr}} e^{2\pi i R}
 C^{-1}$$
From the trace 
$$ 
    s^2=2(1-\cos(\pi \sigma))
$$

The above monodromy data $R$, $\hat{\mu}$, $S$ define the  boundary value
 problem of section \ref{INVERSE II}.  
 The standard technique to solve a 2-dimensional boundary value problem is to
 reduce it to a system of differential equations, which is
 (\ref{systemz1WOW2}) in our case, and then to reduce the system to a second
 order differential equation. It turns out that the equation is (after a
 change of dependent and independent variables) a   Whittaker
 equation. Therefore, the solution of the b.v.p. is given in terms of
 Whittaker functions. We skip the details.   
Let $H:=u_1-u_2$; the fundamental solutions are 
$$
Y_R(z,u)=\pmatrix{  e^{i{\pi \over 2}} \left(H ~z \right)^{-{1\over 2}} 
                    e^{z{u_1+u_2\over 2}} W_{{1\over 2},{\sigma\over 2}}\left(
                   e^{-i\pi} H ~z\right)  
                 &
                    -i{\sigma \over 2}     \left(H~z \right)^{-{1\over
                    2}}  e^{z{u_1+u_2\over 2}} W_{-{1\over 2},{\sigma\over
                    2}}\left( H~z \right) \cr
                              & \cr
                      i{\sigma\over 2}  e^{i{\pi \over 2}} \left(H~z 
                         \right)^{-{1\over 2}} 
                    e^{z{u_1+u_2\over 2}} W_{-{1\over 2},{\sigma\over 2}}\left(
                   e^{-i\pi} H~z\right) 
                  &
                      \left(H~z \right)^{-{1\over
                    2}}  e^{z{u_1+u_2\over 2}} W_{{1\over 2},{\sigma\over
                    2}}\left( H~z \right)
                   \cr
}$$ 
$$
   \arg(R_{12})< \arg(z) < \arg(R_{12})
    +2\pi 
$$ 
where 
$$
 \arg(R_{12}):=-{\pi \over 2} -\arg(u_1-u_2)
$$
$$
Y_L(z,u)=
\pmatrix{\left(Y_R(z,u)\right)_{11}
& i{\sigma\over 2}  \left(H~z \right)^{-{1\over 2}} 
                    e^{z{u_1+u_2\over 2}}
     W_{-{1\over 2},{\sigma\over 2}}\left(
                   e^{-2 i\pi} H~z\right) \cr
      &   \cr
            \left(Y_R(z,u)\right)_{12}         
     & 
         - \left(H~z \right)^{-{1\over
                    2}}  e^{z{u_1+u_2\over 2}} W_{{1\over 2},{\sigma\over
                    2}}\left(e^{-2i \pi } H~z \right)
        \cr
}
$$
$$
\arg(R_{12})+\pi < \arg(z) < \arg(R_{12}) + 3 \pi.
$$
 $W_{\kappa,\mu}$ are the Whittaker functions. 

 For the choice of $Y_R$ and $Y_L$ above, also the sign of $s$
can be determined. According to our computations it is 
$$
  s= 2 \sin\left({\pi \sigma \over 2}\right)
$$
It is computed from the expansion of  $Y_R$ and $Y_L$ at $z=0$.

\vskip 0.2 cm
We stress that {\it the only monodromy data } are  $\sigma$ and the non-zero
entry of $R$.  
The purpose of this didactic  chapter is
to show that the inverse reconstruction of the manifold through (\ref{tu}) and
(\ref{Fu}) brings solutions  $F(t)$ of the WDVV eqs. {\it explicitly
parametrized 
by $\sigma$ and $R$}.


\section{ Preliminary Computations}

\vskip 0.2 cm
The functions $\phi_p(u)$ to be plugged into (\ref{tu}), (\ref{Fu}) 
may be derived from the above representations in
terms of  Whittaker functions. We prefer to proceed in a different way, namely
by imposing the conditions of isomonodromicity (\ref{ISOmonodromy}) to
the solution (\ref{stai barando!}).

\vskip 0.2 cm 

The function $f(u)$ in $\Psi(u)\equiv \phi_0(u)$ is arbitrary, but subject to the condition of 
isomonodromicity (\ref{ISOmonodromy}) for $p=0$, namely 
$$
  \partial_i \phi_0 = V_i \phi_0
$$
where 
$$ 
  V_1= {V\over u_1-u_2},~~~  V_2=-V_1,
$$
Let 
$$
{\cal U}:= \Psi^{-1} U \Psi,~~~~\mu_1:={\sigma\over 2},~~~\mu_2=-{\sigma\over
2},~~~~d=-\sigma 
$$
We will use  $h(u)$ to denote an arbitrary functions of $u$. Let's also denote
the entry $(i,j)$ of a matrix $A_k$  by  $A_{ij,k}$
or by $(A_k)_{ij}$ according to the convenience.  Let $R=R_1+R_2 +
R_3 +...$, $R_{ij,k} \neq 0$ only if $\mu_i-\mu_j=k>0$ integer.   In order to
compute $\phi_p(u)$ of (\ref{stai barando!}) 
we decompose it (and define $H_p(u)$) as follows:  
$$\phi_p(u):=\Psi~ H_p(u),~~~~p=0,1,2,...$$
Plugging the above into (\ref{systemz1WOW2}) we obtain 
\vskip 0.2 cm 
1** $$
  \phi_0=\Psi
$$

2**
$$
  \mu_1\neq \pm{1\over 2},~~~~ 
 H_{ij,1}= {{\cal U}_{ij} \over 1 +\mu_j - \mu_i}, ~~~~R=0  $$
$$
 \mu_1={1\over 2}, ~~~~ H_{12,1} = h_1(u),~~~~R_{12,1}= {\cal U}_{12}$$
$$
 \mu_1=-{1\over 2}, ~~~~ H_{21,1} = h_1(u),~~~~R_{21,1}= {\cal U}_{21}$$

3** 
   $$\mu \neq \pm 1, ~~~~H_{ij,2}={\left( {\cal U} H_1 - H_1 R_1 
    \right)_{ij}\over 2+\mu_j-\mu_i},~~~~R_2=0
$$
$$ 
  \mu=1,~~~~ H_{12,2}=h_2(u),~~~~ R_{12,2}= ({\cal U} H_1)_{12},~~~~R_1=0
$$
$$ 
  \mu=-1,~~~~ H_{21,2}=h_2(u),~~~~ R_{21,2}= ({\cal U} H_1)_{21},~~~~R_1=0
$$ 

4**
  $$ \mu\neq \pm{3\over 2} ,~~~~H_{ij,3}={\left( {\cal U} H_2 -H_1
R_2-H_2 R_1\right)_{ij} \over 3+\mu_j-\mu_i},~~~~R_3=0
$$
$$
 \mu={3\over 2} ,~~~~H_{12,3}= h_3(u), ~~~~R_{12,3}= ({\cal U}
 H_2)_{12}, 
~~~~R_1=R_2=0 
 $$
$$
\mu=-{3\over 2} ,~~~~H_{21,3}= h_3(u), ~~~~R_{21,3}= ({\cal U}
 H_2)_{21},
~~~~ R_1=R_2=0.
 $$

\vskip 0.3 cm
\noindent
We compute $t=t(u)$, $F=F(t(u))$  from the formulae
\be
t^1= \sum_{i=1}^2\phi_{i2,0}\phi_{i1,1},~~~~t^2=\sum_{i=1}^2\phi_{i1,0}
\phi_{i1,1} 
\label{t2dim}
\ee
\be
F={1\over 2} \left[t^{\alpha}t^{\beta}\sum_{i=1}^2
\phi_{i\alpha,0}\phi_{i\beta,1} -\sum_{i=1}^2
\left(\phi_{i1,1}\phi_{i1,2}+\phi_{i1,3}\phi_{i1,0}\right) \right]
\label{F2dim}
\ee
$F$ is defined up to quadratic terms. 

In the following we will compute $F=F(t)$ ($t=(t^1,t^2)$) in closed form, obtaining a solution
of the WDVV equations. The only needed ingredients are (\ref{t2dim}),
(\ref{F2dim}),
the conditions of isomonodromicity (\ref{ISOmonodromy}) and the
symmetries (\ref{SYMMEtria}) for $p=0,1,2,3$. 

\vskip 0.2 cm 
 For any value of $\sigma$  the isomonodromicity
 condition $   \partial_i \phi_0  = V_i \phi_0 $ reads 
$$
 {\partial f(u) \over \partial u_1} = -{\sigma \over 2}  {f(u)
\over u_1-u_2} 
$$ 
$$
 {\partial f(u) \over \partial u_2} = {\sigma \over 2}  {f(u)
\over u_1-u_2} 
$$
 In other words ${\partial f(u) \over \partial u_1}=- {\partial f(u) \over
 \partial u_2} $ and thus $f(u)\equiv f(u_1-u_2)$.   Let
$$
   H:=u_1-u_2
$$
Therefore 
$$ {f(H) \over d H} = -{\sigma \over 2}  {f(H)
\over H} ~~~\Longrightarrow ~~~f(H)= C~ H^{-{\sigma\over 2}},
~~~~C \hbox{ a  constant}
$$


\section{The generic case}

We start from the generic case: $\sigma$ not integer. The result of the
application of formulae (\ref{t2dim}), (\ref{F2dim}) is 
$$ 
 t^1={u_1+u_2 \over 2}
 $$
$$    
              t^2= {1\over 4 (1+\sigma)}{u_1-u_2 \over f(u)^2}
             $$
and
$$ 
 F(t(u))= {1\over 2} (t^1)^2 t^2 +{2 (1+\sigma)^3 \over (1-\sigma) (\sigma+3)
 } (t^2)^3 f(u)^4
$$
But now observe that 
$$ 
   f(u)^2= {u_1-u_2 \over 4(1+\sigma) t^2}$$
$$
   u_1-u_2=H
$$
$$
f(u)^2\equiv f(H)^2= C^2 H^{-\sigma}$$
The above three expressions imply
$$ 
  H= C_1~(t^2)^{1\over 1+\sigma}, ~~~~C_1= [4(1+\sigma)C^2]^{1\over
  1+\sigma}
$$
Therefore
$$ f(u)^4= C_2 (t^2)^{-2{\sigma\over 1+\sigma}} 
$$
and 
$$
 F(t)=  {1\over 2} (t^1)^2 t^2+ C_3 (t^2)^{\sigma+3\over
 \sigma+1} 
$$
as we wanted. Here $C_3$ ( or $C$, or $C_1$, or $C_2$) 
 is an arbitrary constant.


\section{The cases $\mu_1={3\over 2}$, $\mu_1=1$, $\mu_1=-1$}

\noindent
1) Case  $\mu_1={3\over 2}$, $\sigma=3$
\vskip 0.15 cm 
Formula (\ref{t2dim}) gives the same result  of the generic case (with
$\sigma=3$) because $h_3(u)$ appears only in $\phi_3(u)$ and does not affect $t$: 
$$ 
 t^1={u_1+u_2 \over 2}
 $$
$$    
              t^2= {1\over 4 (1+\sigma)}{u_1-u_2 \over f(u)^2}\Big|_{\sigma=3}
             $$
Although $h_3(u)$ appears in $\phi_3(u)$ it is not in $F$:
$$ 
 F(t(u))= {1\over 2} (t^1)^2 t^2 +{2 (1+\sigma)^3 \over (1-\sigma) (\sigma+3)
 } (t^2)^3 f(u)^4\Big|_{\sigma=3}
$$
We may proceed as in the generic case. Actually, now the computation of $f(u)$
is straightforward because 
$$
  R_3=\pmatrix{ 0 & -{1\over 16} (u_1-u_2)^3 f(u)^2 \cr
                0 & 0 \cr
} \equiv \pmatrix{0 & r \cr 0 & 0 \cr } ,~~~~~r=\hbox{constant} 
$$
Namely 
$$
 -{1\over 16} (u_1-u_2)^3 f(u)^2=r
$$
  On the other hand, from $t^2$ we have 
$$ 
   u_1-u_2 = 16 t^2 f(u)^2
$$
and thus
$$ 
   f(u)^4= {(-r)^{1\over 2} \over 16 (t^2)^{3\over 2}}
$$
and finally 
$$ 
 F(t)= {1\over 2} (t^1)^2 t^2 -{3\over 2} (-r)^{1\over 2} (t^2)^{3\over 2}
  \equiv  {1\over 2} (t^1)^2 t^2 +C(t^2)^{3\over 2} 
$$
 where $C$ is an arbitrary constant.

\vskip 0.3 cm
\noindent
2) Case $\mu_1=1$, $\sigma=2$.
\vskip 0.15 cm
 Again, the arbitrary function $h_2(u)$ does not appear in $t(u)$ and
 $F(t(u))$: 
$$ 
 t^1={u_1+u_2 \over 2}
 $$
$$    
              t^2= {1\over 4 (1+\sigma)}{u_1-u_2 \over f(u)^2}\Big|_{\sigma=2},
             $$
$$ 
 F(t(u))= {1\over 2} (t^1)^2 t^2 +{2 (1+\sigma)^3 \over (1-\sigma) (\sigma+3)
 } (t^2)^3 f(u)^4\Big|_{\sigma=2}.
$$
Now we proceed like in the generic case and we find the generic result with
$\sigma=2$.

\vskip 0.3 cm
\noindent
3) Case $\mu_1=-1$, $\sigma=-2$
\vskip 0.15 cm
Now the formulae (\ref{t2dim}), (\ref{F2dim}) yield 
$$ 
 t^1={u_1+u_2 \over 2}
 $$
$$    
              t^2= {1\over 4 }{u_2-u_1 \over f(u)^2}
             $$
$$
 F(t(u))= {3\over 2} (t^1)^2 t^2 -{2\over 3} (t^2)^3 f(u)^4 -t^1 h_2(u)
$$
The condition 
$$ 
0
 = \phi_0^T\phi_2-\phi_1^T\phi_1 +\phi_2^T \phi_0= \pmatrix{ 
                     2h_2(u) +{u_1^2-u_2^2\over 4 f(u)^2} & 0 \cr
                                    0                     & 0 \cr }
$$
implies 
$$
    h_2(u)= {1\over 8} {u_2^2-u_1^2 \over f(u)^2} \equiv t^1 t^2
$$ 
Therefore 
$$
  F(t(u))= {1\over 2} (t^1)^2 t^2 -{2\over 3} (t^2)^3 f(u)^4 
$$
Now we proceed as in the generic case, using $f(H)= C H^{-{\sigma\over
2}} = C H$ and we find the generic result with $\sigma=-2$.


\section{ The case $\mu_1=-{1\over 2}$}

\noindent
 We analyze the case $\mu_1=-{1\over 2}$, $\sigma=-1$. 
The formula (\ref{t2dim})  gives
$$
  t^1= { u_1+u_2 \over 2}
$$
$$
  t^2= h_1(u)
$$
By putting $h_1(u)= t^2$ we get from (\ref{F2dim}) 
$$
  F(t(u))= {1\over 2} (t^1)^2 t^2 +{1\over 16} {(u_1-t^1)^3\over f(u)^2} 
         = {1\over 2} (t^1)^2 t^2 +{1\over 16}
                          {\left({u_1-u_2\over 2}\right)^3 \over f(u)^2} 
$$
From the condition $\partial_{u_i} \Psi = V_i \Psi$ we  computed 
the differential equation for $f(u_1-u_2)=f(H)$ and we got
$f(H)= C H^{-{\sigma\over 2}} = C H^{1\over 2}$. But it is
straightforward to obtain $f(u)$ from 
$$
  R_1= \pmatrix{ 0 & 0 \cr
                 {u_1-u_2\over 4 f(u)^2} & 0 \cr
               } \equiv \pmatrix{ 0 & 0 \cr
                                  r & 0 \cr}
$$
from which 
$$ 
  f(u)^2 = {u_1-u_2 \over 4 r} = {1 \over 4 r} H
$$
The last thing we need is to determine $H$ as a function of $t^1, t^2$. We
can't use the condition $\Phi(-z)^T \Phi(z)= \eta$, because direct  
 computation shows that the following are identically  satisfied:  
$
\phi_0^T \phi_1 -\phi_1^T \phi_0=0
$, $
 \phi_0^T\phi_2-\phi_1^T\phi_1+\phi_2^T\phi_0=0
$, $
\phi_0^T\phi_3-\phi_1^T\phi_2+\phi_2^T\phi_1-\phi_3^T\phi_0=0
$. We make use of the isomonodromicity conditions 
$$
{\partial \phi_1 \over \partial u_1}= E_1 \phi_0 + V_1 \phi_1,~~~~
{\partial \phi_1 \over \partial u_2}= E_2 \phi_0 + V_2 \phi_1
$$
which become 
$$
 {\partial h_1(u) \over \partial u_1}= {1\over 4 f(u)^2},~~~~
{\partial h_1(u) \over \partial u_2}=-{\partial h_1(u) \over \partial u_1}
$$
Thus 
$$ 
h_1(u)\equiv h_1(u_1-u_2),
$$
and 
$$
  {d h_1(H) \over d H} = {r\over H} ~~\Longrightarrow~~
 t^2\equiv h_1(H) = r \ln(H) + D$$ 
$D$ a constant.  
$$ 
 f(u)^4 = {H^2 \over 16 r^2} = C e^{2 {t^2\over r}}
$$
$C$ a constant ($C= \exp(-2D)$)
We get the final result 
$$
 F(t)={1\over 2} (t^1)^2 t^2 +C e^{2 {t^2\over r}} 
$$


\section{ The case $\mu_1={1\over 2}$}

\noindent
 Let $\mu_1={1\over 2}$, $\sigma=1$. $t$ is like in the generic case
$$
  t^1={u_1+u_2 \over 2}
$$
$$
t^2={u_1-u_2 \over 8 f(u)^2}
$$
while $F$ contains $h_1(u)$ 
$$
  f(t(u))= {1\over 2} (t^1)^2 t^2 +{1\over 2} (t^2)^2 h_1(u) -3 (t^2)^3 f(u)^4
$$ 
We determine $f(u)$ as in the generic case, or better we observe that 
$$
  R_1= \pmatrix{ 0 & (u_1-u_2) f(u)^2 \cr 
                 0  & 0  \cr } 
                              \equiv \pmatrix{ 0 &  r  \cr 0 & 0 \cr}
$$ 
Thus
 $$ 
      f(u)^2 = {r \over u_1-u_2}
$$
We determine $h_1(u)$. The condition $\Phi(-z)^T \Phi(z)=\eta$ does not help,
because it  is automatically satisfied. We use the isomonodromicity conditions
$$
{\partial \phi_1 \over \partial u_1}= E_1 \phi_0 + V_1 \phi_1,~~~~~
{\partial \phi_1 \over \partial u_2}= E_2 \phi_0 + V_2 \phi_1
$$
which become 
$$
{\partial h_1(u) \over \partial u_1} = f(u)^2,~~~~~ 
  {\partial h_1(u) \over \partial u_2}=-{\partial h_1(u) \over \partial u_1}
$$
Therefore 
$$
   h_1(u)\equiv h_1(u_1-u_2)$$
We keep into account that $f(u)^2=r/H$: 
$$
   {d h_1(H) \over d H} = {r\over H}~~\Longrightarrow~~ h_1(H) 
= r \ln(H) +D
$$
$D$ a  constant. Finally, recall that 
$$
t^2= {H \over 8 f(u)^2}\equiv {H^2\over 8 r},
$$
hence 
$$ 
   f(u)^4= {r \over 8 t^2} ,
$$ 
which contributes a linear term to $F(t)$, and
$$
h_1(u)= {r\over 2} \ln (t^2) + B
$$
 $B= {r\over 2} \ln(8r) +C$ is an arbitrary constant. 
Finally
$$
F(t)= {1\over 2} (t^1)^2 t^2 +{r\over 4} (t^2)^2 \ln (t^2)
$$
as we wanted.


\section{ The case $\mu_1=-{3\over 2}$}

\noindent
 Finally, let's take $\mu_1=-{3\over 2}$, $\sigma=-3$. From (\ref{t2dim}),
 (\ref{F2dim}) we have
$$
t^1={u_1+u_2\over 2}
$$
$$
t^2= {u_2-u_1\over 8 f(u)^2}
$$
$$
F(t(u))= {3\over 4} (t^1)^2 t^2 + (t^2)^3 f(u)^4 -{1\over 2} h_3(u)
$$
$f(u)$ is obtainable as in the generic case, but it is more straightforward to
use 
$$
R_3=\pmatrix{ 0 & 0 \cr {(u_2-u_1)^3 \over 64 f(u)^2}  & 0 \cr }\equiv
\pmatrix{0 & 0 \cr r & 0 }~\Longrightarrow~ f(u)^2= {(u_2-u_1)^3\over 64 r} 
$$
To obtain $h_3(u)$ we can't rely on $\Phi(-z)^T\Phi(z)=\eta$, which turns out
to be identically satisfied. We use again the conditions 
\be
{\partial \phi_3 \over \partial u_i} = E_i \phi_2 + V_i \phi_3
\label{tmp}
\ee
It is convenient to introduce 
$$
   G(u):= {1\over 4} (t^1)^2t^2 -{1\over 2} h_3(u)
$$
The above (\ref{tmp}) becomes 
$$ 
   {\partial G \over \partial u_1} = {r \over 2(u_2-u_1)} ,~~~~ {\partial G \over \partial u_2}=-{\partial G \over \partial u_1}
$$
which implies $G(u)=G(u_1-u_2)$ and 
$$ 
{d G \over d H}= -{r\over 2 H}~~\Longrightarrow ~~ G(H)=-{r\over 2}
\ln(H) +C$$
$C$ a constant. 
Finally, recall that 
$$ 
   t^2 = {H^2\over 8 r}~~~\Longrightarrow~~ 
G(H(t^2))= {r\over 4} \ln(t^2) +C_1  
$$ 
 Thus 
$$
  F(t) = {1\over 2} (t^1)^2 t^2 +{r\over 4} \ln(t^2) 
$$
having dropped the constant terms.

\section{ Conclusions} 

The solution of the boundary value problem for the monodromy data $\sigma$
 and the non zero entry $r$ of the matrix $R$
 (Stokes matrices and other monodromy data depend on $\sigma$) was obtained
   solving the equations (\ref{ISOmonodromy}) with the constraints
 (\ref{SYMMEtria}). 

We have obtained the solutions of the WDVV eqs. which we already derived from
elementary considerations in chapter \ref{cap1}. They are parametrized
explicitly by the monodromy data:   
$$ 
\hbox{ For }\sigma \neq \pm 1, -3~~~~~
  F(t)= {1\over 2} (t^1)^2 t^2 + C (t^2)^{3+\sigma \over 1 +\sigma}
$$
where $C$ is a constant. 
$$
\sigma=-1, ~~~~~F(t)={1\over 2} (t^1)^2 t^2+ C e^{2 {t^2 \over r}}
$$
$$
\sigma=1,~~~~~F(t)={1\over 2} (t^1)^2 t^2+C (t^2)^2 \ln(t^2)
$$
$$
\sigma=-3,~~~~~F(t)={1\over 2} (t^1)^2 t^2 + C \ln(t^2)
$$


\chapter{ Stokes Matrices and Monodromy of $ QH^{*}(CP^d)$}\label{stokes}

 We briefly summarize the results of 
\cite{guz}, where  we  computed Stokes' matrices  and monodromy group for the 
Frobenius manifold given by the quantum cohomology of the projective space 
${\bf CP}^d$. The first motivation for the computation is 
 to know the monodromy data of
$QH^{*}(CP^d)$, in view of the  solution of the inverse problem. 
Actually, the global
analytic properties of the solution of the WDVV eqs. 
obtained by Kontsevich-Manin (see
section \ref{solutionQH}) are unknown and the boundary value problem may shed
light on them.

The second motivation was  to study the  links between 
quantum cohomology  and the theory of coherent sheaves, in order to prove a
long-lasting conjecture about the coincidence of the Stokes matrix with the
Gram matrix of entries $g_{ij}=\sum_k(-1)^k ~\hbox{dim} ~Ext^k({\cal
O}(i-1),{\cal O}(j-1))$, $i,j=1,2,...,d+1$.

\vskip 0.2 cm
 The main result of \cite{guz}  is the proof 
that the conjecture about coincidence of the 
Stokes matrix for quantum cohomology of 
${\bf CP}^d$  and the above  Gram matrix is true.  
We proved that the Stokes'matrix can be
reduced to the canonical form $S=(s_{ij})$ where  
$$ 
s_{ii}=1,~~~~s_{ij}=\bin{d+1}{j-i},~~~~s_{ji}=0,~~~~~~i<j,
$$
by the action of the braid group. This form is equal to the Gram matrix 
(modulo the action of the braid group ). 
In this way, we generalized to any $d$ the result obtained in \cite{Dub2} for
$d=2$.

\vskip 0.2 cm 

In \cite{guz} 
we also studied the structure of the monodromy group of the quantum cohomology
of  
${\bf CP}^d$.                                             
 We proved  that for $d=3$ the group is isomorphic to 
the subgroup of 
orientation preserving transformations in the hyperbolic triangular group
$[2,4,\infty]$.  
In \cite{Dub2} it was proved that for $d=2$ the monodromy group is isomorphic
to the  
direct product of the subgroup of orientation preserving transformations in
$[2,3,\infty]$  and the cyclic group of order 2, $C_2=\{\pm \}$. Our numerical
calculations  also suggest that for any $d$ even the monodromy group may be 
 isomorphic to the orientation preserving transformations in 
$[2,d+1,\infty]$, and for any $d$ odd to the direct product of 
 the orientation 
preserving transformations in $[2,d+1,\infty]$ by $C_2$.


%
\chapter{ Connection problem and Critical Behaviour for $PVI_{\mu}$
}\label{PaInLeVe}

 In this chapter we study the critical behaviour and the connection
 problem for  the solutions of $PVI_{\mu}$ (\ref{PVImu}). 
The purpose  is to fix some results which
 will be used in the next chapters when we  solve the inverse  
problem for the associated
 Frobenius manifold in terms of Painlev\'e transcendents.  

 We refer to the Introduction of the thesis for a general 
 discussion of   the results of this
 chapter. 

\vskip 0.4 cm 
\noindent 
I would like to thank A. Bolibruch, A. Its, M. Jimbo, M. Mazzocco, 
S. Shimomura for
discussions related to the subject of this chapter.

\section{Branches of $PVI_{\mu}$ and monodromy data - known results}

    The equation $PVI_{\mu}$ is equivalent to the equations of 
    isomonodromy deformation of the fuchsian system obtained in section 
\ref{PAPAPA} 
\be{dY\over dz}
=
\left[
{A_1(u)\over z-u_1} +{A_2(u)\over z-u_2}+
{A_3(u)\over z-u_3}\right]~Y:= A(z;u)~Y
\label{fuchsMJ}
\ee
$$u:=(u_1,u_2,u_3), 
~~~~\hbox{tr}(A_i)=\det{A_i}=0,~~~~\sum_{i=1}^3~A_i=-\hbox{diag}(\mu,-\mu)
$$
If $u$ is deformed, 
the  monodromy matrices of a solution of the system (\ref{fuchsMJ}) do not
change, provided that the deformation is ``small'' (see below). 
 The connection to   $PVI_{\mu}$ is given by
$$
    x={u_3-u_1\over u_2-u_1},~~~y(x)= {q(u)-u_1\over u_2-u_1}
$$
where $q(u_1,u_2,u_3)$ is the root of 
$$
    [A(q;u_1,u_2,u_3)]_{12}=0~~~~~~\hbox{if }\mu \neq 0
$$
The case $\mu=0$ is disregarded, because $PVI_{\mu=0} \equiv
PVI_{\mu=1}$.

The system (\ref{fuchsMJ}) has fuchsian singularities at $u_1$, $u_2$,
 $u_3$. Let us fix a branch $Y(z,u)$ of a fundamental matrix solution by
 choosing branch cuts  in the $z$  plane and a basis of loops in
  $\pi({\bf C}\backslash \{u_1,u_2,u_3\};z_0)$, where $z_0$ is a
 base-point. Let   $\gamma_i$ be a basis of loops 
encircling  counter-clockwise the point $u_i$, $i=1,2,3$. See figure
 \ref{figura8}. Then  
 $$ 
Y(z,u) \mapsto Y(z,u)M_i,~~~ i=1,2,3,
~~~~\det M_i\neq 0, 
$$
if $z$ describes a loop  $\gamma_i$. Along the 
loop $\gamma_{\infty}:=\gamma_1 \cdot \gamma_2 \cdot
\gamma_3$ we have $Y\mapsto Y M_{\infty}$. 
$M_i$ are the {\it monodromy matrices}, and they give a
representation of the fundamental group.  Of course $ 
              M_{\infty}= M_3 M_2 M_1
$. The transformations $ 
 Y^{\prime}(z,u)= Y(z,u) B ~~~~~\det(B)\neq 0
$ 
yields all the possible fundamental matrices, hence the monodromy matrices
of (\ref{fuchsMJ}) are defined up to conjugation 
$$
M_i \mapsto M_i^{\prime}= B^{-1}
M_i B.
$$  
From the standard theory of fuchsian systems it
follows that we can choose  a fundamental solution behaving as follows
$$
 Y(z;u)= \left\{ \matrix{ 
                         (I+O({1\over z})) ~z^{-\hat{\mu}}
                        z^{R}~C_{\infty},~~~~z\to \infty \cr \cr
                        G_i(I+O(z-u_i))~(z-u_i)^J~C_i,~~~~z\to
                        u_i,~~~i=1,2,3 \cr
                                     } \right.
$$
where $J=\pmatrix{ 0 & 1 \cr 0 & 0 \cr}$,  
$ \hat{\mu}= $diag($\mu,-\mu)$, $G_i J G_i^{-1}= A_i$ ($i=0,x,1$) and  
$$
     R=\left\{\matrix{ 
~~~~~~
                        0~~~~, ~~~~~~~~~~~~~~~~\hbox{ if } 2\mu \not
\in {\bf Z } \cr \cr
             \left. \matrix{  
           \pmatrix{ 0 & b \cr 0 & 0 \cr } ,~~~~~~~~\mu>0\cr
                \pmatrix{ 0 & 0 \cr b & 0 \cr } ,~~~~~~~~\mu<0 \cr 
     } \right\}
                      \hbox{ if } 2\mu\in{\bf Z} \cr   
         }\right. 
$$
$ b \in {\bf C} $ being determined by the  matrices $A_i$.  
Then $M_i=C_i^{-1} e^{2\pi i J} C_i$, $M_{\infty}= C_{\infty}^{-1}
e^{-2\pi i \hat{\mu}} e^{2\pi i R} C_{\infty}$.

The dependence of the fuchsian system on $u$ is isomonodromic. This
 means that for small deformations of $u$ the monodromy matrices do
 not change \cite{JM1} \cite{IN} and then 
$$ 
   \partial_i R = \partial_i C_1 = \partial_i C_2 =\partial_i C_3= 
\partial_i C_{\infty} =0
 $$
 ``Small'' deformation means that $x=(u_2-u_1)/(u_3-u_1)$
 can vary in the $x$-plane provided it does not describe complete
 loops around $0,1,\infty$ (in other words, we fix  
 branch cuts, like $\alpha<\arg(x)<\alpha+2\pi$ and
 $\beta<\arg(1-x)<\beta+2\pi$ for $\alpha,\beta\in{\bf R}$). 
 For ``big'' deformations, the monodromy matrices
 change according to an action of the pure braid group. This point
 will be discussed later. 

\vskip 0.2 cm 
Suppose we have a branch $y(x)$. We associate to it the fuchsian
system (we stress that we actually have a branch defined by the
cuts, thus the fuchsian system has 
monodromy matrices independent of
$x$). Therefore, to any branch of a Painlev\'e transcendent there
corresponds a monodromy representation. 

 Conversely, the problem of finding 
 the branches of the Painlev\'e
 transcendents of $PVI_{\mu}$ for given monodromy matrices (up to
 conjugation) is the problem of finding a fuchsian system (\ref{fuchsMJ})
 having the given monodromy matrices. This problem is called {\it
 Riemann-Hilbert problem}, or { \it $21^{th}$ Hilbert problem}. For a
 given $PVI_{\mu}$ (i.e. for a fixed $\mu$) there
 is a one-to-one correspondence between a monodromy representation and a
 branch of a transcendent if and only if 
the Riemann-Hilbert problem has a unique
 solution.

\vskip 0.2 cm 
\noindent
{\bf $\bullet$ Riemann-Hilbert problem (R.H.)}: find the 
 coefficients $A_i(u)$, $i=1,2,3$ from the following monodromy data:

a) the matrices  
    $$ 
        \hat{\mu} = \hbox{diag}(\mu,-\mu), ~~~~~\mu\in {\bf
        C}\backslash \{0\}
$$
\vskip 0.2 cm 
$$
     R=\left\{\matrix{ 
~~~~~~
                        0~~~~, ~~~~~~~~~~~~~~~~\hbox{ if } 2\mu \not
\in {\bf Z } \cr \cr
             \left. \matrix{  
           \pmatrix{ 0 & b \cr 0 & 0 \cr } ,~~~~~~~~\mu>0\cr
                \pmatrix{ 0 & 0 \cr b & 0 \cr } ,~~~~~~~~\mu<0 \cr 
     } \right\}
                      \hbox{ if } 2\mu\in{\bf Z} \cr   
         }\right. 
$$
$ b \in {\bf C} $.  

b)  three poles $u_1$, $u_2$, $u_3$, a base-point and a base of loops in
 $\pi({\bf C}\backslash \{u_1,u_2,u_3\};z_0)$. See figure \ref{figura8}.

\begin{figure}
\epsfxsize=12cm
\centerline{\epsffile{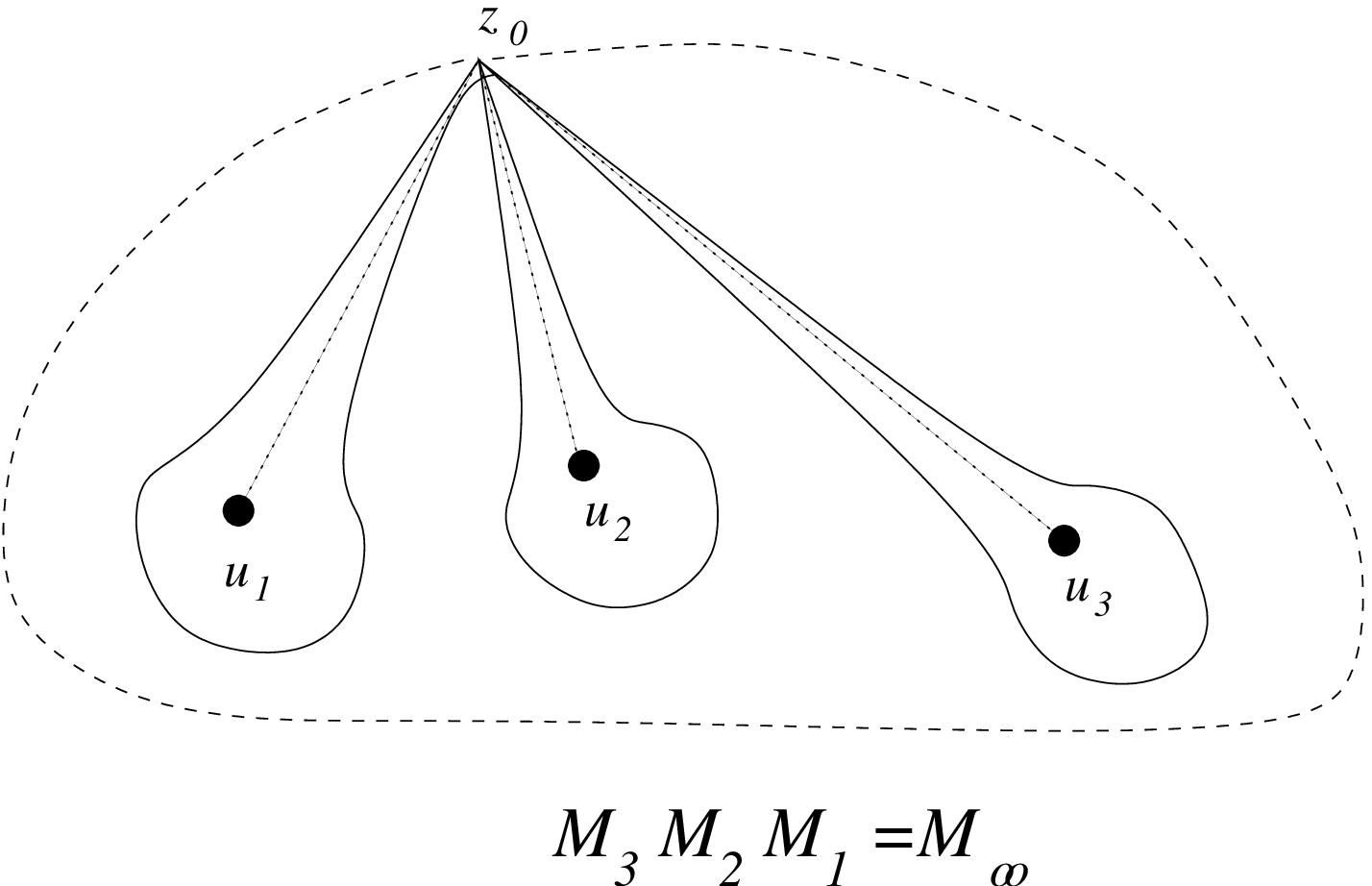}}
\caption{Choice of a basis in $\pi_0({\bf C}\backslash\{u_1,u_2,u_3\}) $}
\label{figura8}
\end{figure}

c) three  monodromy matrices $M_1$, $M_2$, $M_3$ relative to the loops
 (counter-clockwise) 
 and a matrix $M_{\infty}$ similar to $e^{-2\pi i \hat{\mu}} e^{2\pi iR}$, 
 and satisfying 
$$ 
\hbox{tr}(M_i)=2,~~~\hbox{det}(M_i)=1,~~~~i=1,2,3 $$
$$ 
  M_3~M_2~M_1 = M_{\infty}
$$
\be
M_{\infty}= C_{\infty}^{-1}~e^{-2\pi i \hat{\mu}}e^{2 \pi i  R}~C_{\infty}
\label{conn1}
\ee
where $C_{\infty}$ realizes the similitude. We also choose the indices
of the problem, namely we fix  ${1 \over 2\pi i}\log M_i$ as
follows: let
$$ 
   J:=\pmatrix{0 & 1 \cr 0 & 0 \cr}
$$ 
We require there exist 
three {\it connection matrices} $C_1$, $C_2$, $C_3$ such that 
\be
     C_i^{-1} e^{2 \pi i J} C_i = M_i,~~~~~i=1,2,3
\label{conn2}
\ee
and we look for a matrix valued meromorphic
function $Y(z;u)$  such that 
\be
   Y(z;u)= \left\{ \matrix{ 
                        G_{\infty} (I+O({1\over z})) ~z^{-\hat{\mu}}
                        z^{R}~C_{\infty},~~~~z\to \infty \cr \cr
                        G_i(I+O(z-u_i))~(z-u_i)^J~C_i,~~~~z\to
                        u_i,~~~i=1,2,3 \cr
                                     } \right.
\label{solrh}
\ee
$G_{\infty}$ and $G_i$ are  invertible matrices depending on $u$.  The
coefficient of the fuchsian system are then given by $
 A(z;u_1,u_2,u_3):= {d Y(z;u) \over dz} Y(z;u)^{-1}
$.

\vskip 0.2 cm
Recall that a $2 \times 2$ R.H. is always solvable \cite{AB}. The monodromy
matrices are considered up to the conjugation
\be
 M_i \mapsto M_i^{\prime}=B^{-1} M_i B, ~~~~~~~\det{B}\neq 0,~~~~i=1,2,3, \infty
\label{conjj}
\ee
and the coefficients of the fuchsian system itself are considered up
to conjugation $A_i \mapsto F^{-1} A_i F$ $(i=1,2,3)$, by an invertible matrix
$F$.  Actually, two conjugated
fuchsian systems admit fundamental matrix solutions with the same
monodromy, and a given fuchsian system defines the monodromy up to
conjugation (depending on the choice of  the fundamental solution). 

On the other hand, a triple of monodromy matrices $M_1$, $M_2$, $M_3$
may be realized by two fuchsian systems which are not conjugated.    
This corresponds to the fact that the solutions $C_{\infty}$, $C_{i}$
of (\ref{conn1}), (\ref{conn2})
are not unique, and the choice of different particular solutions may
give rise to   fuchsian systems which are not conjugated. If this is
the case, there is no one-to-one correspondence between monodromy matrices (up
to conjugation) and solutions of $PVI_{\mu}$. It is easy to prove that:

\vskip 0.2 cm 
{\it The
R.H. has a unique solution, up to conjugation, for $2 \mu \not \in {\bf
Z}$ or for $2 \mu \in {\bf Z}$ and $R\neq 0$}.
\footnote{ The proof  is 
 done in the following way: consider
two solutions $C$ and $\tilde{C}$ of the equations (\ref{conn1}),
(\ref{conn2}). Then 
$$ 
    (C_i \tilde{C}_i^{-1})^{-1} e^{2\pi i J} (C_i \tilde{C}_i^{-1})=
    e^{2\pi i J}
$$
$$ 
 (C_{\infty} \tilde{C}_{\infty}^{-1})^{-1} e^{-2\pi i
 \hat{\mu}}e^{2\pi i R} 
 (C_{\infty} 
\tilde{C}_{\infty}^{-1})=
    e^{-2\pi i \hat{\mu}} e^{2\pi i R}
$$
We find  
$$
C_i \tilde{C}_i^{-1}= \pmatrix{a & b \cr
                               0 & a  \cr}  ~~~~~~a,b\in{\bf C},~~~a\neq 0
$$
Note that this matrix commutes with $J$, then 
$$ 
  (z-u_i)^J C_i = (z-u_i)^J \pmatrix{a & b \cr
                               0 & a  \cr}\tilde{C}_i=\pmatrix{a & b \cr
                               0 & a  \cr}(z-u_i)^J \tilde{C_i}
$$
We also find 
$$
C_{\infty} 
\tilde{C}_{\infty}^{-1}
=
    \left\{ \matrix{ 
                     i)~~~~~~~\hbox{diag}(\alpha, \beta),~~~\alpha\beta\neq 0;
                 ~~~~~~~~~~~~~~~~~~~~~~~~~~~~~~~~~~~~ 
                    ~\hbox{if } 2\mu \not \in {\bf Z} \cr
                   ii) ~ \pmatrix{\alpha & \beta \cr
                                  0  &  \alpha \cr} ~~(\mu>0),~~~ 
 \pmatrix{\alpha & 0 \cr
                                  \beta  &  \alpha \cr} ~~(\mu<0),
~~~\alpha\neq 0, ~~~\hbox{if } 2\mu \in {\bf Z}, R\neq 0 \cr
iii)~ \hbox{ Any invertible matrix} ~~~~~~~~~~~~~~~~~~~~~~~~~~~~~~~~~~~~
~\hbox{ if }  2\mu \in {\bf Z}, R= 0 \cr
}\right.
$$
Then 
$$ i)~ z^{-\hat{\mu}}C_{\infty}= z^{-\hat{\mu}}\hbox{diag}(\alpha,\beta)
\tilde{C}_{\infty} = \hbox{diag}(\alpha,\beta)z^{\hat{\mu}}\tilde{C}_{\infty}
$$
$$
ii)~z^{-\hat{\mu}}z^{-R} C_{\infty} = ~...~= \left[\alpha I+ {1\over
z^{|2\mu|}} Q\right] 
     ~z^{-\hat{\mu}}z^{-R} \tilde{C}_{\infty}
$$
where $Q=\pmatrix{0 & \beta \cr 0 & 0\cr }$, or $Q= \pmatrix{0 &  0\cr
\beta  & 0 \cr}$. 
$$
iii)~  i)~ z^{-\hat{\mu}}C_{\infty}=~...~= \left[
              {Q_1\over z^{|2\mu|}} + Q_0 + Q_{-1}z^{|2\mu|}  
\right]
             z^{-\hat{\mu}} \tilde{C}_{\infty}
$$
where $Q_{0}=$ diag$(\alpha, \beta)$, $Q_{\pm 1}$ are respectively
upper and lower triangular (or lower and upper triangular, depending
on the sign of $\mu$), and $C_{\infty}
\tilde{C}_{\infty}^{-1}=Q_1+Q_0+Q_{-1}$  

This implies that the  two solutions $Y(z;u)$, $\tilde{Y}(z;u)$ of the
form  (\ref{solrh}) with $C$ and $\tilde{C}$ respectively, are such
that $Y(z;u)~\tilde{Y}(z;u)^{-1}$ is holomorphic at each $u_i$, while at
$z=\infty$ it is 
$$
   Y(z;u)~\tilde{Y}(z;u)^{-1} \to 
\left\{ \matrix{
                 i)~ G_{\infty} \hbox{diag}(\alpha,\beta)
G_{\infty}^{-1} \cr
ii)~ \alpha I \cr
iii)~ \hbox{ divergent }
    }
\right.
$$
Thus the two fuchsian systems are conjugated only in the   cases
$i)$ and $ii)$, because in those cases $Y\tilde{Y}^{-1}$ is
holomorphic everywhere on ${\bf P}^1$, and then it is a constant.  
In other words {\it the
R.H. has a unique solution, up to conjugation, for $2 \mu \not \in {\bf
Z}$ or for $2 \mu \in {\bf Z}$ and $R\neq 0$}.
}

\vskip 0.3 cm 

Once the R.H. is solved, the sum of the matrix coefficients $A_i(u)$ of the 
solution $A(z;u_1,u_2,u_3)=\sum_{i=1}^3
{A_i(u)\over z-u_i} $ must be diagonalized (to give
$-$ diag$(\mu,-\mu)$).
\footnote{ Note however that if $G_{\infty}=C_{\infty}=I$, then $\sum_{i=1}^3
A_i$ is  already  diagonal. 
Moreover, for $2\mu \not \in {\bf Z}$, $M_1$, $M_2$, $M_3$ and 
the choice of normalization $Y(z;u)z^{\hat{\mu}}
 \to I$ if $z\to \infty$ determine uniquely $A_1$, $A_2$,
$A_3$. Actually, for 
 any diagonal invertible matrix $D$, the matrices  $M_1^{\prime}= 
D^{-1}M_1D$, $M_2^{\prime}= D^{-1}M_2D$,
 $M_3^{\prime}= D^{-1} M_3 D $ 
 determine the coefficients $D^{-1} A_i D$, whose sum is still
 diagonal (the normalization of $Y$ is the same). }
  After that, a branch y(x)  of $PVI_{\mu}$ can be 
``computed'' from $[A(q;u_1,u_2,u_3)]_{12}=0$.   
 The fact that the R.H. has a unique solution for the given monodromy data
 (if  $2 \mu \not \in {\bf
Z}$ or  $2 \mu \in {\bf Z}$ and $R\neq 0$) means that there is a one-to-one
 correspondence between the triple  $M_1$, $M_2$, $M_3$  and the branch
 $y(x)$.

\vskip 0.2 cm 
\noindent
 We review some known results \cite{DM} \cite{M}. 
\vskip 0.15 cm 

1) One  $M_i=I$ if and only if the Schlesinger equations yield
$q(u)\equiv u_i$. This does not correspond to a solution of
$PVI_{\mu}$. 

\vskip 0.15 cm 

2) If the $M_i$'s, $i=1,2,3$, commute, then $\mu$ is integer (as it follows
from the fact that the $2\times 2$ matrices with 1's on the diagonals 
commute if and only if
they  can be simultaneously put in upper or
lower triangular form). There are solutions of
$PVI_{\mu}$ only for 
  $$ 
      M_1=\pmatrix{ 1 & i \pi a \cr
                    0 & 1       \cr},~~~M_2=\pmatrix{ 1 & i \pi  \cr
                    0 & 1       \cr},~~~M_3=\pmatrix{ 1 & i \pi(1- a) \cr
                    0 & 1       \cr},~~~~~a\neq 0,1
$$
In this case $R=0$ and $M_{\infty}=I$. For $\mu=1$ the solution is 
  $$ y(x)= {a x \over 1-(1-a)x}$$
and for other integers $\mu$ the solution is obtained from $\mu=1$ by a
birational transformation \cite{DM} \cite{M}.  In particular, these
solutions are rational. 

\vskip 0.15 cm 

3) Non commuting $M_i$'s. 

 The parameters in the space of the
 monodromy representation,  independent of conjugation of the $M_i$, are 
$$
  2-x_1^2:=\hbox{ tr}(M_1 M_2), ~~~2-x_2^2:= \hbox{ tr}(M_2 M_3),~~~
 2-x_3^2:=\hbox{ tr}(M_1M_3)
$$
The triple $(x_0, x_1, x_{\infty})$ in the introduction and in the
remaining part of this paper corresponds to
$(x_1,x_2,x_3)$.

\vskip 0.13 cm

3.1) If at least two of the $x_j$'s are zero, then one of the
$M_i$'s is $I$, or the matrices commute. We return to the case 1 or 2.
Note that in case 2 $(x_1,x_2,x_3)=(0,0,0)$.
\vskip 0.13 cm 

3.2) At most one of the $x_j$'s is zero. Namely, the triple
$(x_1,x_2,x_3)$ is {\it admissible}. In this case it is possible to
fully parameterize the monodromy using the triple $(x_1,x_2,x_3)$ :
namely, there exists a basis
such that: 
$$ 
  M_1= \pmatrix{ 1 & -x_1 \cr 
                 0 & 1     \cr} ,~~~~M_2= \pmatrix{ 1 & 0 \cr 
                 x_1 & 1     \cr},~~~~
M_3=\pmatrix{ 1 +{x_2 x_3\over x_1} & -{x_2^2\over x_1} \cr
               {x_3^2\over x_1}     &  1-{x_2 x_3\over x_1} \cr
             },
$$
if $x_1\neq 0$. If $x_1=0$ we just choose a similar parameterization
starting from $x_2$ or $x_3$. 
The relation 
$$ M_3M_2M_1 \hbox{ similar to } e^{-2\pi i \hat{\mu} }e^{2\pi i R}$$ 
implies 
$$ 
    x_1^2+x_2^2+x_3^2-x_1x_2x_3 = 4 \sin^2(\pi \mu)
$$
The signs of the $x_i$'s must be chosen in such a way that the above
relation is satisfied. 
The conjugation (\ref{conjj}) changes the triple by two signs. Thus
the true parameters for the monodromy data are classes of equivalence
of triples $(x_1,x_2,x_3)$ defined by the change of two signs. 

We distinguish three cases: 

i) $2\mu \not \in {\bf Z}$. Then there is a one to one correspondence
between monodromy data $(x_0,x_1,x_{\infty})$ 
and the branches of transcendents of $PVI_{\mu}$. 
The  solutions of \cite{DM} are included in this case: the 
 connection problem was solved for the class of
transcendents 
having the following local behaviour at the critical points $x=0,1,\infty$: 
\be
y(x)= a^{(0)} x^{1-\sigma^{(0)}}(1+O(|x|^{\delta})),~~~~x\to 0,
\label{loc1}
 \ee
\be
y(x)= 1-a^{(1)}(1-x)^{1-\sigma^{(1)}} (1+O(|1-x|^{\delta})),~~~~x\to 1,
\label{loc2}
 \ee
\be
y(x)= a^{(\infty)}
 x^{-\sigma^{(\infty)}}(1+O(|x|^{-\delta})),~~~~x\to \infty,
\label{loc3}
\ee
where $a^{(i)}$ and
 $\sigma^{(i)}$ are complex numbers such that $a^{(i)}\neq 0$ and 
 $0\leq \Re \sigma^{(i)}<1$. $\delta$ is a small positive number. 
This behaviour is true if $x$ converges
 to the critical points inside a sector in the $x$-plane with vertex
 on the corresponding critical point and finite angular width. 
The {\it connection problem} was solved finding 
  the relation among the three couples $(\sigma^{(i)},a^{(i)})$,
$i=0,1,\infty$  and the triple $(x_0,x_1,x_{\infty})$. 
In \cite{DM} all the
algebraic solutions are classified and related to the 
finite reflection groups $A_3$, $B_3$, $H_3$.

ii) For any $\mu$ half integer there is an infinite set of {\it Picard type
solutions }
(see \cite{M}), in one to one correspondence to triples of monodromy
data ($R\neq 0$) not in the equivalence class of  $(2,2,2)$. 
These solutions form a two parameter family, behave asymptotically as 
the solutions of the case
$i)$, and comprise a denumerable subclass of algebraic solutions. For any half integer $\mu\neq {1\over 2}$ there is also a one
parameter family of   
{\it Chazy solutions}.  For them the
one to one correspondence with monodromy data is lost. In fact, they
form an infinite family but  any element of the family corresponds to a triple
$(x_1,x_2,x_3)$ 
in the orbit (of the braid group) of the triple $(2,2,2)$ (this orbit
is simply obtained by changing two signs in all possible ways). 
They appear in the case $R=0$ (no other solutions of $PVI_{\mu}$ correspond
to $R=0$ and $\mu$ half integer). The result of our paper applies to the 
Picard's solutions with $x_i\neq \pm 2$. 

iii) $\mu$ integer. In this case $R\neq 0$ ( $R=0$ only in  the case 2) of
commuting monodromy matrices and $\mu$ integer). 
There is a one to one correspondence
between monodromy data $(x_0,x_1,x_{\infty})$ 
and the branches. To our knowledge,
this case 
has not yet been studied. There are relevant examples of Frobenius
manifolds where these solutions must appear, like the case of Quantum 
Cohomology of $CP^2$. In this case
$\mu=-1$, the triple $(x_1,x_2,x_3)=(3,3,3)$ (the monodromy data
coincide with the elements of the Stokes'matrix of the  corresponding
Frobenius manifold \cite{Dub2} \cite{guz}) and the real part of
$\sigma$ is equal to 1.

\vskip 0.15 cm 
 In this chapter we find the critical behaviour and we solve the
 connection problem for almost all the triples and for any $\mu \neq 0$, 
the only restriction required
 being 
$$
    x_i\neq \pm 2 ~~\Longrightarrow ~~ \sigma^{(i)} \neq 1,~~~~i=0,1,\infty
$$


\section{ Local Behaviour -- Theorem 1}\label{Local Behaviour -- Theorem 1}

A solution $y(x)$ of $PVI_{\mu}$ is a meromorphic function of
the point $x$ belonging to  the universal
covering of  ${\bf P}^1 \backslash \{0,1,\infty\}$. 
 Let $\sigma$ and $a$ denote two  complex numbers, with the
  restrictions    
$$\sigma \in \Omega:={\bf C}\backslash
  \{(-\infty,0) \cup [1,+\infty)\}~~~~a\neq 0
$$
  Our first aim is to
  show the existence of solutions of $PVI_{\mu}$ which behave like $a
  ~x^{1-\sigma}$, as $x\to 0$ along a  suitable path. 
 Thus, we  concentrate on a small punctured
neighbourhood of $x=0$, and the point $x$ can be read as a 
point in the universal covering of ${\bf C}_0:={\bf
C}\backslash \{0\}$ with $0<|x|<\epsilon$ ($\epsilon<1$). Namely, 
 $x=|x| e^{i
\arg(x)}$, where $-\infty< \arg(x) < +\infty$.

In order to specify the way  $x$ may tend to zero,  we introduce a domain 
contained in
the universal covering  of ${\bf C}_0$ (we denote the universal
covering by  $\widetilde{{\bf
C}_0}$). 
If $\sigma=0$ we define the domain 
$$D(\epsilon;\sigma=0)=\{{ x}\in \widetilde{{\bf C}_0}
   \hbox{ s.t. } |x|<\epsilon \}
$$
If $\sigma \neq 0$ we observe that 
$$
    |{x}^{\sigma}| = |x|^{\sigma^{\prime}(x)},~~\hbox{ where } 
 \sigma^{\prime}(x) := \Re \sigma- {\Im \sigma ~\arg(x) \over \log|x|}
$$
We try to define a domain where the exponent $\sigma^{\prime}(x)$ satisfies
the restriction $0\leq\sigma^{\prime}(x)<1$ for $x\to 0$. Let
$\theta_1$, $\theta_2 \in {\bf R}$, $0<\tilde{\sigma} 
<1$. The desired domain is:
$$
   D(\epsilon;\sigma;\theta_1,\theta_2,\tilde{\sigma}):= 
\{{ x}\in \widetilde{{\bf C}_0}
   \hbox{ s.t. } |x|<\epsilon ,~~ e^{-\theta_1 \Im \sigma}
   |x|^{\tilde{\sigma}} \leq |{x}^{\sigma}| 
\leq e^{-\theta_2 \Im
   \sigma} |x|^{0},~~~0<\tilde{\sigma}<1
   \},~~~0<\epsilon<1. 
$$  
 The domain can be written as 
$$|x|<\epsilon,~~~~~ 
        \Re \sigma \log|x|+\theta_2 \Im \sigma \leq \Im \sigma \arg(x)
        \leq (\Re \sigma - \tilde{\sigma}) \log|x| + \theta_1 \Im
        \sigma 
$$

 Figure \ref{figura1} shows the domains. Note that if $0\leq \Re \sigma <1$ the
domain 
contains, for $|x|$ sufficiently small, 
 any sector $\alpha<\arg (x)<\alpha+2\pi$ ($\alpha$ a real number); 
but for $\Re \sigma <0$ and
$\Re \sigma \geq 1$ it is not possible to include a sector in it when $x\to
0$.  

\begin{figure}
\epsfxsize=13cm
\centerline{\epsffile{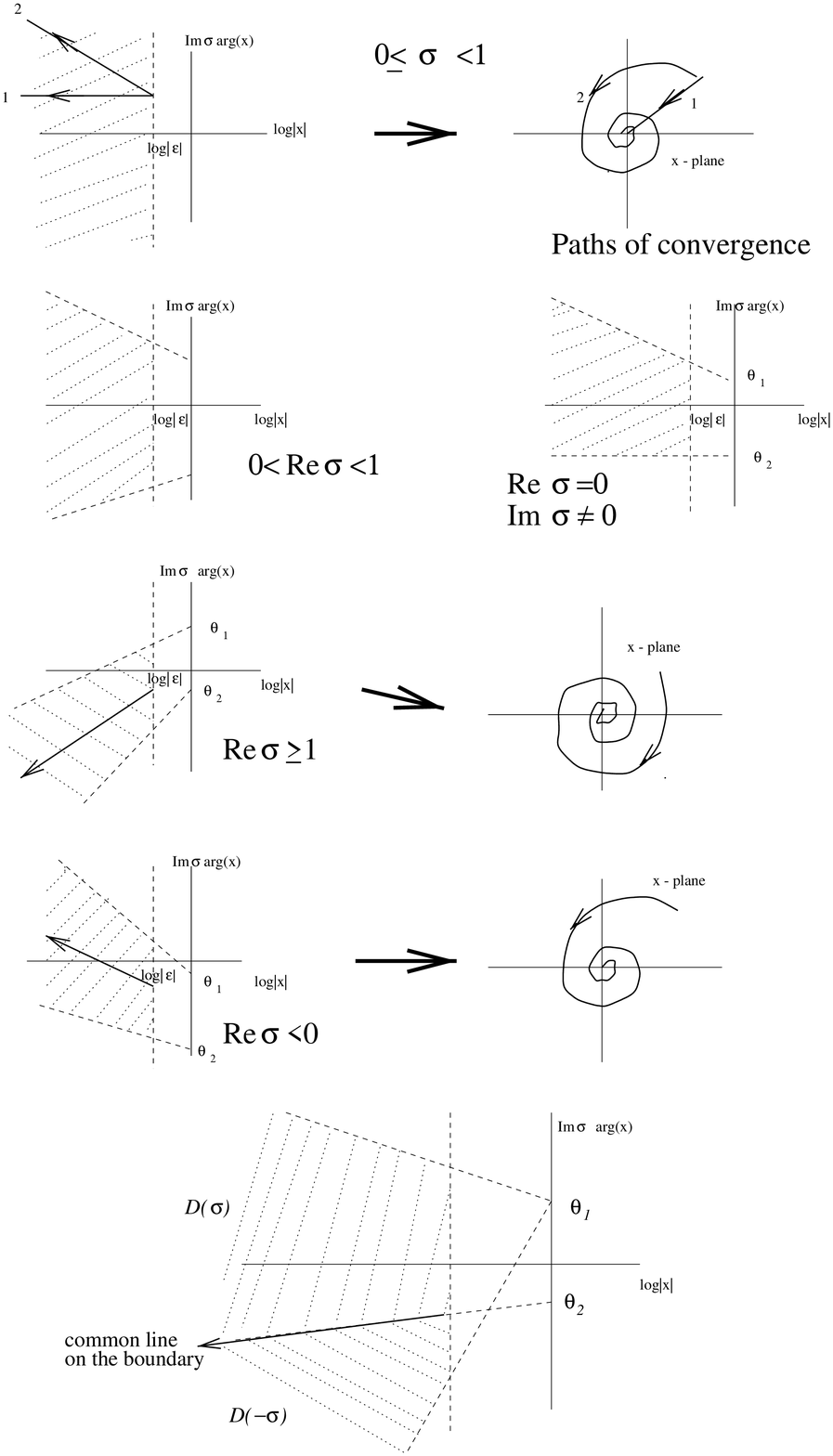}}
\caption{ We represent the domains $D(\epsilon;\sigma;\theta_1,\theta_2)$ in
the $(\ln |x|,\Im\sigma\arg(x))$-plane. Note that $
D(\epsilon;\sigma;\theta_1,\theta_2)= D(\epsilon;\sigma=0)$ for real 
$0\leq \sigma <1$. 
  We also represent some lines along
which $x$ converges to 0. These lines are also represented in the $x$-plane:
they are radial paths or spirals.}
\label{figura1}
\end{figure}

\vskip 0.3 cm 
\noindent
{\bf Theorem 1: } {\it Let $\mu \neq 0$.  For any $\sigma \in \Omega$, for
any $a \in {\bf 
C}$, $a\neq 0$, for any $\theta_1,\theta_2  \in {\bf R}$ and for any
$0<\tilde{\sigma }<1$, there exists a sufficiently small positive $\epsilon$ 
such that  the  equation $PVI_{\mu}$ has a transcendent 
$y(x;\sigma,a)$ with behaviour  
\be
    y(x;\sigma,a)=a x^{1-\sigma} \left(
1+O(|x|^{\delta})
\right)~~~~,0<\delta<1,
\label{asy0}
\ee
as $x\to 0$ in  $D(\epsilon;\sigma;\theta_1,\theta_2,\tilde{\sigma})$.  

The above local behaviour is valid with an  exception which occurs if 
 $\Im \sigma\neq 0$  and  ${x}\to 0$ along the special  paths 
 $\Im \sigma \arg(x) = \Re \sigma \log|x| + d
$, where $d$ is a constant such that the path is 
 contained in $D(\epsilon;\sigma;\theta_1,\theta_2)$: the local
behaviour  becomes: 
$$
   y({x};\sigma,a)= a(x)~{x}^{1-\sigma} (1+O(|x|^{\delta}))
,~~~~~~~~~~~~~~ \sigma\neq 0, $$
$$a({x})= a~\left(1 +{1\over 2a}~ C e^{i \alpha({x})} 
             + {1\over 16 a^2}~ C^2 e^{2 i \alpha({x})} 
 \right) = O(1), \hbox{ for }
             x\to 0
$$
where 
$$
 {x}^{\sigma}=C e^{i\alpha({x})},
~~~~~C:=e^{-d},~~~\alpha({x}):=\Re\sigma \arg(x)+
\Im \sigma \ln|x|\Bigl|_{\Im \sigma \arg(x) = \Re \sigma \log|x| + d}
$$
}

\vskip 0.3 cm
The small number $\epsilon$  depends on $\tilde{\sigma}$, $\theta_1$ and
$a$. In the following 
we may sometimes 
 omit $\epsilon$, $\tilde{\sigma}$, $\theta_i$ and write simply
$D(\sigma)$. 

 In figure \ref{figura1}  
we  draw the possible paths along which $x\to 0$. Any path is allowed if $\Im
\sigma=0$. If $\Im
\sigma\neq 0$, an example of allowed paths is the following: 
\be
  |x|<\epsilon,~~~~~~ \Im \sigma \arg(x) = 
(\Re \sigma - \Sigma) \ln|x| + b, \label{spirale}
\ee
for a suitable choice of the additive constant $b$  and for $0\leq \Sigma\leq
\tilde{\sigma}$. In general, these paths are spirals, represented in figure
\ref{figura1}. 
They include radial paths if $0\leq \Re \sigma<1$ and $  \Sigma= 
\Re \sigma$, because in this case $\arg(x)=\hbox{ constant}$.  
But there are only  spiral paths 
 whenever $\Re \sigma<0$ and $\Re \sigma \geq 1$.
The special
 paths 
 $\Im \sigma \arg(x) = \Re \sigma \log|x| + b
$,  are parallel to one of the  boundary lines of
 $D(\sigma)$ in the plane $(\ln|x|,\arg (x))$. The boundary line is $\Re
 \sigma \ln|x|+ \Im \sigma \theta_2$ and it is shared by
$ D(\sigma)$ and   $
 D(-\sigma)$ (with the same $\theta_2$).

\vskip 0.3 cm 
\noindent
{\bf $\bullet$ Restrictions on the domain 
$D(\epsilon;\sigma;\theta_1,\theta_2)$:} In  theorem 1  we can choose
$\theta_1$ arbitrarily. Apparently, if we increase $\theta_1 \Im \sigma $ the
domain $D(\epsilon;\sigma;\theta_1,\theta_2)$ becomes larger. But $\epsilon$
itself depends on $\theta_1$. In the proof of theorem 1 (section 
\ref{proof of theorem 1}) we have to impose  
$$      \epsilon^{1-\tilde{\sigma}}\leq c~e^{-\theta_1 \Im \sigma} 
$$   
where $c$ is a constant (depending on $a$). Equivalently, $\theta_1 \Im \sigma
\leq (\tilde{\sigma}-1)\ln \epsilon + \ln c$. This means that if we increase
$\Im \sigma \theta_1$ we have to decrease $\epsilon$. Therefore, for $x \in   
D(\epsilon;\sigma;\theta_1,\theta_2)$ we have: 
$$ 
   \Im \sigma \arg (x) \leq (\Re \sigma -\tilde{\sigma}) \ln|x|+ \theta_1 
\Im \sigma \leq (\Re \sigma -\tilde{\sigma}) \ln|x| + (\tilde{\sigma}-1)\ln
\epsilon  +\ln c
$$
We advise the reader to visualize a point $x$ in the plane
$(\ln|x|, \Im\sigma \arg(x))$.  With this visualization in mind, let 
$x_{\epsilon}$ be  the point   $\left\{ (\Re \sigma -\tilde{\sigma}) \ln|x|+
(\tilde{\sigma}-1)\ln 
\epsilon  +\ln c \right\} \cap \{|x|=\epsilon\}$ (see figure \ref{figure13}). 
Namely, 
$$ 
 \arg x_{\epsilon}= (\Re \sigma -1) \ln \epsilon + \ln c
$$ 
This means that the union  of the domains
$D(\epsilon(\theta_1);\sigma;\theta_1,\theta_2)$ on all values of $\theta_1$
for  given $\sigma, a, \tilde{\sigma}, \theta_2$ is
$$ 
  \bigcup_{\theta_1} D\bigl(\epsilon(\theta_1);\sigma;\theta_1,\theta_2,\tilde{\sigma}\bigr)
\subseteq B(\sigma,a;\theta_2, \tilde{\sigma})
$$ 
where 
\be
   B(\sigma,a;\theta_2, \tilde{\sigma}):= \left\{|x|< 1 ~\hbox{ such that }~ \Re \sigma \ln|x|
+\theta_2 \Im \sigma < \Im \sigma \arg (x)\leq (\Re \sigma -1) \ln|x| +\ln
c\right\}
\label{come se non bastasse!!}
\ee
The dependence on $a$ of the domain $B$ defined above is motivated by the fact
that $c $ depends on $a$ (but not on $\theta_1$, $\theta_2$).

\begin{figure}
\epsfxsize=12cm
\centerline{\epsffile{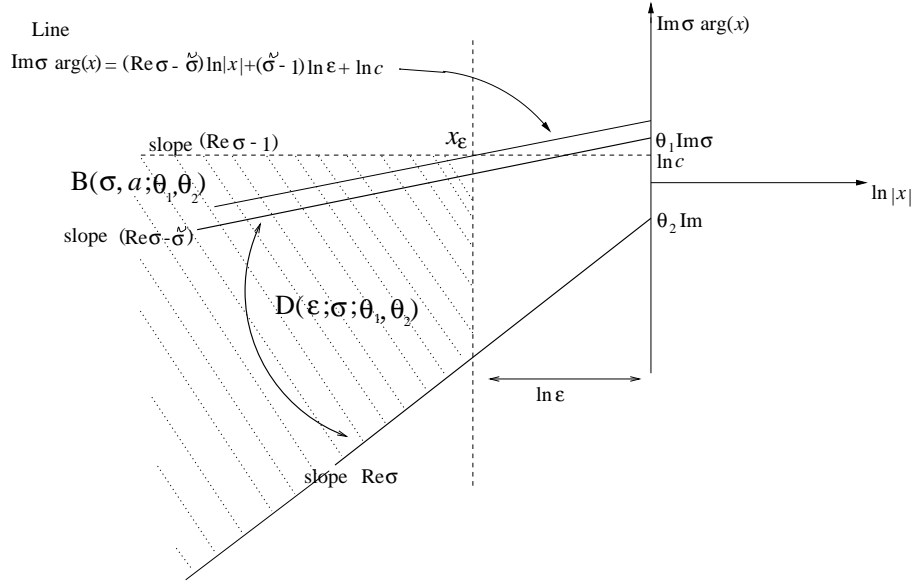}}
\caption{Construction of the domain $B(\sigma,a;\theta_2, \tilde{\sigma})$ for
$\Re \sigma =1$.}
\label{figure13}
\end{figure}

 If $0\leq \Re \sigma <1$, the above result is not a limitation on the values
 of $\arg(x)$ of the points 
 $x$ that we want to include in a given $D(\sigma;\epsilon)$ provided that
 $|x|$ is sufficiently small. Also in the case $\Re \sigma <0$ there is no
 limitation, because we can always decrease $\Im \sigma \theta_2$ without
 affecting $\epsilon$ in order to include  in
 $D(\epsilon;\sigma;\theta_1,\theta_2)$ a point $x$ such that
 $|x|<\epsilon$.   But this is not the case if $\Re
 \sigma \geq 1$.  Actually, if $x$ (in the $(\ln|x|, \Im\sigma
 \arg(x))$-plane) lies above the set  $B(\sigma,a;\theta_2, \tilde{\sigma})$ it
 never can be included in any $D(\epsilon;\sigma;\theta_1,\theta_2)$. See
 figure \ref{figure14}

\begin{figure}
\epsfxsize=12cm
\centerline{\epsffile{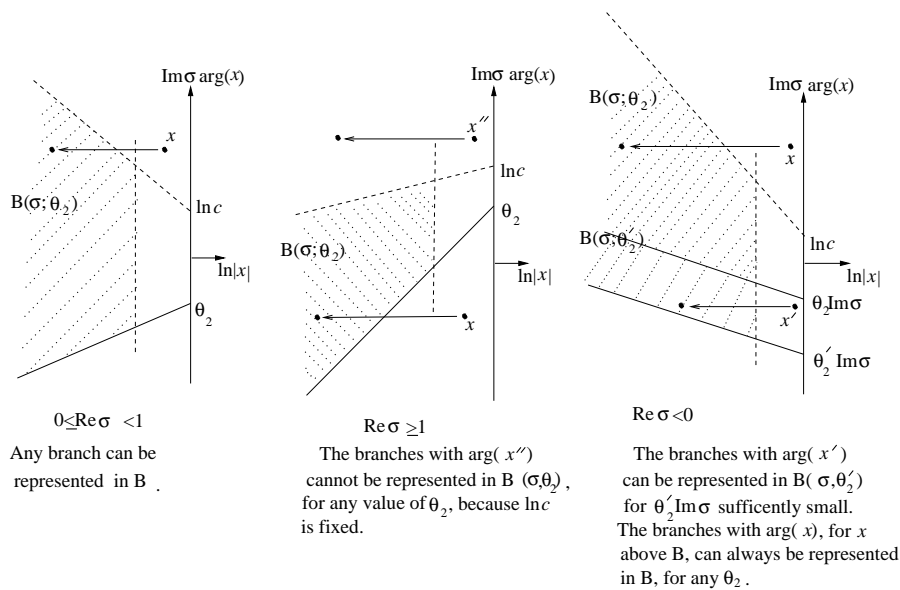}}
\caption{For $\Re \sigma\geq 1$ we can not include all values of $\arg(x)$ in
 $B$}
\label{figure14}
\end{figure}


\section{Parameterization of a branch through Monodromy Data --
Theorem 2}\label{Parametrization of a branch through Monodromy Data --
Theorem 2}

We are going to consider the fuchsian system (\ref{fuchsMJ}) for the
special choice 
$$ 
  u_1=0,~~~u_2=1,~~~u_3=x
$$ 
The labels $i=1,2,3$ will be substituted by the labels $i=0,1,x$, and
the system becomes 
$$ 
 {d Y \over dz} = \left[ {A_0(x) \over z} +{A_x(x)\over z-x}+ {A_1(x)
 \over z-1} \right]~Y
$$
Also, the triple $(x_1,x_2,x_3)$ will be denoted by
$(x_0,x_1,x_{\infty})$, as in \cite{DM}. 

We consider only  admissible triples and  $x_i\neq \pm 2$, $i=0,1,\infty$. 
We recall that two admissible triples 
are equivalent if their elements  differ just by the change of two
signs and that 
\be
x_0^2+x_1^2+x_{\infty}^2-x_0x_1x_{\infty}=4 \sin^2(\pi \mu)
\label{10ii}
\ee

 We denote by $y(x;x_0,x_1,x_{\infty})$ a branch in one to one
 correspondence with $(x_0,x_1,x_{\infty})$. Since we operate close to $x=0$,
 the branch is specified by $\alpha<\arg (x)<\alpha+2\pi$, $\alpha\in {\bf
 R}$.

\vskip 0.3 cm
\noindent
{\bf Theorem 2:} {\it Let $\mu$ be any non zero complex number. 

For any $\sigma \in
\Omega$ and for any $a\neq 0$  
 there exists a triple of monodromy data
($x_0,x_1,x_{\infty}$) uniquely determined (up to equivalence) 
by the following formulae: 
 \vskip 0.15 cm
\noindent
 i) $\sigma \neq 0, \pm 2\mu+2m$ for any $m\in {\bf Z}$. 
$$
   \left\{ \matrix{x_0=2 \sin({\pi \over 2}\sigma) \cr\cr
                  x_1=i\left({1\over
f(\sigma,\mu)G(\sigma,\mu)}~\sqrt{a}-G(\sigma,\mu) ~{1\over \sqrt{a}}\right)
                      \cr\cr
                  x_{\infty}= {1\over f(\sigma,\mu)G(\sigma,\mu)
e^{-{i\pi\sigma\over 2}} }~\sqrt{a}+ G(\sigma,\mu) e^{-i{\pi\sigma\over
2}}~{1\over \sqrt{a}}\cr
                         }\right. 
$$
where 
$$f(\sigma,\mu)={2 \cos^2({\pi \over 2} \sigma) \over \cos(\pi \sigma)- 
                  \cos(2\pi\mu)}, ~~~
~~~ G(\sigma,\mu)= {1\over 2}{4^{\sigma}\Gamma({\sigma+1\over 2})^2\over
                  \Gamma(1-\mu+ {\sigma\over
                  2})\Gamma(\mu+{\sigma\over 2})}$$
Any sign  of $\sqrt{a}$  is good (changing the sign of $\sqrt{a}$ 
                is equivalent
                  to changing the sign of both $x_1$, $x_{\infty}$).
\vskip 0.15 cm
\noindent
 ii) $\sigma =0$
$$\left\{ 
              \matrix{
                         x_0=0  \cr\cr
                        x_1^2=2\sin(\pi\mu)~\sqrt{1-a} \cr\cr
                        x_{\infty}^2=2\sin(\pi\mu)~\sqrt{a}
                        \cr
                       }
 \right.  
$$
We can take any sign of the square roots 
 
\vskip 0.15 cm
\noindent
 iii) $\sigma= \pm 2\mu +2 m$.  

\vskip 0.15 cm
  iii1) $\sigma=2\mu+2m$, $m=0,1,2,...$
 $$ 
         \left\{ 
   \matrix{ 
              x_0=2\sin(\pi\mu)\cr \cr
              x_1= -{i\over 2}{16^{\mu+m} \Gamma(\mu+m+{1\over 2})^2 \over
   \Gamma(m+1) \Gamma(2\mu+m)}~{1\over\sqrt{a}} 
                                           \cr \cr
               x_{\infty}=i~x_1~e^{-i\pi\mu}\cr
          }
                  \right.
$$

 iii2) $\sigma=2\mu+2m$, $m=-1,-2,-3,...$
$$
  \left\{ 
\matrix{
            x_0=2\sin(\pi\mu)\cr \cr
            x_1=2i{\pi^2\over \cos^2(\pi\mu)}
                {1\over 16^{\mu+m} \Gamma(\mu+m+{1\over 2})^2
\Gamma(-2\mu-m+1) \Gamma(-m)}~\sqrt{a}
\cr\cr
x_{\infty}=-ix_1 e^{i\pi\mu}
}
\right.
$$

iii3) $\sigma=-2\mu +2m$, $m=1,2,3,...$
$$
   \left\{ 
\matrix{
           x_0=-2\sin(\pi\mu) \cr\cr
         x_1= -{i\over 2}{16^{-\mu+m} \Gamma(-\mu+m+{1\over 2})^2\over
\Gamma(m-2\mu+1) \Gamma(m)}~{1\over \sqrt{a}} \cr\cr
        x_{\infty}= ix_1e^{i\pi\mu}\cr
}
\right.
$$

iii4) $\sigma= -2\mu+2m$, $m=0,-1,-2,-3,...$
$$
  \left\{
\matrix{
            x_0=-2\sin(\pi\mu) \cr\cr
            x_1= 2i {\pi^2\over \cos^2(\pi\mu)}{1\over 16^{-\mu+m}
\Gamma(-\mu+m+{1\over 2})^2\Gamma(2\mu-m) \Gamma(1-m)}  ~\sqrt{a}
\cr\cr
x_{\infty}=-ix_1 e^{-i\pi \mu}\cr
}
 \right.
$$
In all the above formulae the relation
$x_0^2+x_1^2+x_{\infty}^2-x_0x_1x_{\infty}=4 \sin^2(\pi \mu)$  is
automatically satisfied. Note that $\sigma \neq 1$ implies $x_0\neq \pm
2$. Equivalent triples (by the change of two signs) are also allowed. 

Let $x\in D(\epsilon;\sigma)$. The branch at $x$ of  $y(x;\sigma,a)$ coincides
 with 
 $y(x;x_0,x_1,x_{\infty})$. 
\footnote{
Note that we have picked up a point $x\in D(\sigma)$, therefore
 $\alpha<\arg (x) <\alpha +2\pi$ for a suitable $\alpha$. In other words, the
 branch is specified by $\alpha <\arg(x) <\alpha +2\pi$.  If we pick up
 a new point $x^{\prime}=e^{2\pi i } x \in D(\sigma)$ (provided this is
 possible),  then $y(x^{\prime};\sigma,a)$ is
 again equal to the branch  $y(x^{\prime};x_0,x_1,x_{\infty})$ at
 $x^{\prime}$, the branch being specified by
 $\alpha+2\pi <\arg (x) <\alpha +4\pi$ . 
In section \ref{Analytic Continuation of a Branch}  
we will describe in detail the problem of analytic
 continuation. We anticipate that for the loop  $x^{\prime}= x e^{2\pi
 i}$ we have the continuation 
 $y(x^{\prime};\sigma,a)=y(x^{\prime};x_0,x_1,x_{\infty})\equiv 
y(x;x_0^{\prime},x_1^{\prime},x_{\infty}^{\prime})$. 

$y(x;x_0^{\prime},x_1^{\prime},x_{\infty}^{\prime})$ is a new branch
 of  $y(x;x_0,x_1,x_{\infty})$ corresponding to the continuation above  at
 the same point $x$. In other words,  we put
 branch cuts in the 
 $x$-plane, then   $\arg (x)$ can not increase by $2\pi$ and the analytic
 continuation of a branch yields a new branch with the same $\arg(x)$ and new
 monodromy data $(x_0^{\prime},x_1^{\prime},x_{\infty}^{\prime})$. 

$y(x^{\prime};x_0,x_1,x_{\infty})$ is  the continuation of the branch 
$y(x;x_0,x_1,x_{\infty})$ in the universal
 covering of $C_0\cap \{ |x|<\epsilon\}$. It has new  
$\arg(x)$, i.e $\arg(x) \mapsto \arg(x^{\prime})
 =\arg(x) +2\pi$ and the same monodromy data.  
If $x^{\prime}$ still lies in $D(\sigma)$ we can
 represent the continuation as $y(x^{\prime};\sigma,a)=y(x^{\prime};
 x_0,x_1,x_{\infty})$, where the branch $y(x^{\prime};
 x_0,x_1,x_{\infty})$ has the branch cut specified by 
 $\alpha +2\pi <\arg(x^{\prime})<
 \alpha +4 \pi$, while $y(x;x_0,x_1,x_{\infty})$ has branch cut
$\alpha<\arg(x) <\alpha +2\pi$.   
}


\vskip 0.2 cm 

 Conversely, for any set of monodromy data ($x_0,~x_1,~x_{\infty}$)
 such that $x_0^2+x_1^2+x_{\infty}^2-x_0x_1x_{\infty}=4 \sin^2(\pi
 \mu)$, $x_i\neq \pm 2$,  there exist parameters   $\sigma$  and $a$
 obtained as follows:  
 
\vskip 0.15 cm 
\noindent
 I) Generic case 
$$
\cos(\pi \sigma)= 1-{x_0^2\over 2}
$$
$$                  %
  a={iG(\sigma,\mu)^2\over 2 \sin(\pi \sigma)} 
      \Bigl[
 2(1+e^{-i\pi\sigma})-f(x_0,x_1,x_{\infty})(x_{\infty}^2+e^{-i\pi\sigma} x_1^2)
\Bigr]  ~ f(x_0,x_1,x_{\infty})
$$
where  
$$
   f(x_0,x_1,x_{\infty}):=f(\sigma(x_0),\mu)={4-x_0^2\over 2-x_0^2-
2\cos(2\pi\mu)}={4-x_0^2\over
    x_1^2+x_{\infty}^2-x_0x_1x_{\infty}} $$
Any solution $\sigma $ of the first equation must satisfy the restriction
$\sigma \neq  \pm 2\mu+2m$ for any $m\in {\bf Z}$, otherwise we
encounter the singularities in $G(\sigma,\mu)$ and in $f(\sigma,\mu)$. 
If $x_0^2=4$ the system has
solutions $\sigma=1+2n$, $n\in {\bf Z}$, which do not belong to
$\Omega$. 

\vskip 0.15 cm 
\noindent
 II) $x_0=0$. 
           $$ \sigma=0,$$
$$
    a= {x_{\infty}^2\over x_1^2+x_{\infty}^2}.
$$
provided that $x_1\neq 0$ and $x_{\infty}\neq 0$, namely  $\mu
\not\in {\bf Z}$.  
\vskip 0.15 cm 
\noindent
III) $x_0^2=4 \sin^2(\pi\mu)$.  Then (\ref{10ii}) implies
 $x_{\infty}^2=-x_1^2 ~\exp(\pm 2\pi i \mu)$ . 
 Four cases which yield the values of $\sigma$ non included in I)
and II) must be considered 

\vskip 0.15 cm
   III1)  $x_{\infty}^2=-x_1^2 e^{- 2
\pi i \mu}$
$$\sigma=  2\mu + 2m,~~~~m=0,1,2,...$$ 
 $$ a=-{1\over 4 x_1^2} 
               { 16^{2\mu+2m}
 \Gamma(\mu+m+{1\over 2})^4\over \Gamma(m+1)^2 \Gamma(2\mu+m)^2  } 
$$

\vskip 0.15 cm
  III2) $x_{\infty}^2=-x_1^2 e^{2\pi i \mu}$
$$\sigma=2\mu+2m,~~~~m=-1,-2,-3,...$$
    $$
            a=-{\cos^4(\pi\mu)\over 4 \pi^4}
                 16^{2\mu+2m} \Gamma(\mu+m+{1\over 2})^4
\Gamma(-2\mu-m+1)^2 \Gamma(-m)^2~ x_1^2  
$$

\vskip 0.15 cm
  III3) $x_{\infty}^2=-x_1^2e^{2\pi i\mu}$
$$\sigma=-2\mu+2m,~~~~m=1,2,3,...$$
$$a=-{1\over 4 x_1^2}
         {
  16^{-2\mu+2m} \Gamma(-\mu+m+{1\over 2})^4\over \Gamma(m-2\mu+1)^2
\Gamma(m)^2  }
$$

\vskip 0.15 cm
III4)  $x_{\infty}^2=-x_1^2 e^{-2\pi i \mu}$
$$\sigma=-2\mu+2m,~~~~m=0,-1,-2,-3,...$$
$$
  a= -{\cos^4(\pi\mu)\over 4 \pi^4} 16^{-2\mu+2m}
\Gamma(-\mu+m+{1\over 2})^4\Gamma(2\mu-m)^2 \Gamma(1-m)^2 ~ x_1^2
$$

 If $x_0^2\neq 4$ ($\sigma\neq 1+2 n$, $n \in {\bf
 Z}$)
 we can always  choose (from I), II), III) ) 
 $\sigma \in \Omega$. 

 Let  $x\in D(\epsilon; \sigma)$ (then there  
exists $\alpha\in {\bf R}$ such that 
$\alpha<\arg(x) <\alpha +2\pi$). The  
branch   $y(x;x_0,x_1,x_{\infty})$ coincides at $x$ with the transcendent 
$y(x;\sigma,a)$ of theorem 1.
 }

\vskip 0.3 cm

 We stress that the proof of the theorem 
is valid also for the {\it resonant}
 case $2\mu \in {\bf Z}\backslash \{0\}$.  To read the formulae, 
it is enough to just substitute an integer for
 $2\mu$ in the above formulae i) or I). Actually, we note
 that ii), iii); II), III) cannot occur. 
Note that for $\mu$ integer the case ii), II) degenerates to
 $(x_0, x_1,x_{\infty})=(0,0,0)$
and $a$ arbitrary. This is the case in which the triple is 
 not a good parameterization for the monodromy (not admissible triple). Anyway,
we know that in this case there is a one-parameter family of rational
solutions \cite{M}, which are all obtained by a birational
transformation from the family 
$$ 
    y(x)={a x \over 1 -(1-a) x} , ~~~~~\mu =1
$$
At  $x=0$ the the behaviour is  
$y(x)= a x (1+O(x))$, and then the limit of theorem 2 for $\mu \to n \in
              {\bf Z} \backslash \{0\}$ and $\sigma=0$ yields the
above one-parameter family. Recall that $R=0$ in this case.

\vskip 0.2 cm
\noindent
{\it Remark 1: }
 The  equation
$$
 \cos(\pi \sigma) = 1-{x_0^2\over2}
$$
determines $\sigma$ up to $\sigma \mapsto \pm \sigma +2 n$, $n\in {\bf
Z}$.
\footnote{In the case $\sigma=1+i \nu$, $\nu \in {\bf R}\backslash \{0\}$, 
the freedom $\sigma \mapsto -\sigma+2$ is equivalent to $\sigma=1+i\nu \mapsto 
\sigma-2i\nu=1-i \nu$.} 
 These  different values give different coefficients $a$ and  different 
domains $D(\epsilon_n,\pm \sigma +2n)$.  Thus the same branch
$y(x;x_0,x_1,x_{\infty})$ has analytic continuations on different
domains with different asymptotic behaviours
prescribed by theorem 1. In particular, note that if $\Im \sigma \neq 0$ 
it is always possible to 
 choose
$$
0\leq \Re \sigma \leq 1.
$$
Observe however that for $0\leq \sigma<1$ the operation  $\sigma
\mapsto \pm  \sigma
+2n$ is not allowed, because $\sigma$  leaves $\Omega$. Then, our paper does 
not give any information about what happens for $\sigma=1$.

We also note that $\sigma =\sigma(x_0)$ and $a= a(\sigma;x_0,
x_1,x_{\infty})$, therefore $a(\sigma;x_0,x_1,x_{\infty})\neq a(\pm
\sigma+2n;x_0,x_1,x_{\infty})$.

\vskip 0.2 cm
\noindent
{\it Remark 2:}
The domains $D(\sigma)$ and $D(-\sigma)$, with the same $\theta_2$,
intersect along the common boundary  $  \Im \sigma \arg(x)= 
 \Re \sigma \log|x| +
 \theta_2 \Im \sigma$ (see figure \ref{figura1}). 
 The asymptotic behaviour of the analytic
 continuation of the branch  $y(x;
 x_0,x_1,x_{\infty})$ at $x$ belonging to the common boundary  is given in
 terms of $(\sigma(x_0), a(\sigma;x_0,x_1,x_{\infty}))$ and 
$(-\sigma(x_0), a(-\sigma;x_0,x_1,x_{\infty}))$ respectively.

 This implies that the two different asymptotic 
 representations 
 of theorem 1 on $D(\sigma)$ and $D(-\sigma)$ must become equal
 along the  boundary.
 Actually, from  theorem 1 it is 
 clear that along the boundary of $D(\sigma)$, the behaviour of $y(x)$ is 
 $$
 y(x) =
 A(x;\sigma,a(\sigma)) ~x~\left(1+O(|x|^{\delta})\right)$$
 where $\delta$ is a small number between 0 and 1 and 
$$
 A(x;\sigma, a(\sigma))= a (C e^{i\alpha(x;\sigma)})^{-1}+{1\over 2} +{1\over
 16 a} 
 Ce^{i\alpha(x;\sigma)} 
$$    
$$x^{\sigma}=
 Ce^{i\alpha(x;\sigma)},~~~~~C=e^{-\theta_2\Im\sigma},~~~~
~\alpha(x;\sigma)=\Re\sigma\arg(x) 
+\Im\sigma\ln|x|\bigl|_{\Im \sigma \arg(x) =
 \Re \sigma \log|x| + \theta_2 \Im \sigma }
$$   
We observe that  $\alpha(x;-\sigma)=-\alpha(x;\sigma)$. After 
the proof of theorem 2 we'll see that 
$a(\sigma) = {1 \over 16 a(-\sigma)}$: this immediately implies
that 
     $$ A(x;- \sigma, a(-\sigma))= A(x;\sigma,a(\sigma))$$
Therefore, the asymptotic behaviour, as prescribed
by theorem 1 in $D(\sigma)$ and $D(-\sigma)$, is the same along the
common boundary of the two domains. 

\vskip 0.3 cm
We end the section with the following

\vskip 0.2 cm
\noindent
{\bf Proposition:} {\it Let $y(x)\sim a x^{1-\sigma}$ as $x\to 0 $ in
a domain $D(\epsilon,\sigma)$. Then, $y(x)$ coincides with $y(x;\sigma,a)$ of
theorem 1} 
\vskip 0.2 cm
\noindent
{\it  Proof:} see section \ref{dimth2}.



 \section{ Alternative Representations of the Transcendents}\label{beyond}

  We study the critical point $x=0$ (the points  $x=1,\infty$ will be
discussed in section \ref{Singular
Points x=1, x=infty  (Connection Problem)}).

 If $ \Im \sigma\neq 0$, 
the freedom $\sigma \mapsto \pm \sigma +2n$ allows us to always reduce the 
exponent to $ 0 \leq \Re \sigma \leq 1$. So,  
the values of $\sigma$ we may restrict to  are 
$$   0 \leq \Re \sigma \leq 1 \hbox{ for } \Im \sigma\neq 0,$$
$$   0\leq \sigma <1 \hbox{ for } \sigma\in {\bf R}.  $$
In the $(\ln |x|,\Im \sigma \arg(x))$-plane we 
draw 
the domains $D(\sigma)$, $D(-\sigma)$,
$D(-\sigma+2)$, $D(2-\sigma)$, etc - see figure \ref{figure8} (left). Here 
$\theta_1$ and $\theta_2$ are the same
($\epsilon$ may not be the same, so we choose the smallest). We suppose 
that the domains and the corresponding transcendents of theorem 1 are
associated to the same triple.  So  we have different
critical behaviours predicted by theorem 1 in the different domains for 
 the same transcendent. 
Some ``small
sectors'' remain uncovered by the union of the domains (figure \ref{figure8}
(right)). If $x\to 0$ inside these sectors, we do not know the behaviour of
the transcendent. If $\Re \sigma =1$, a 
 radial path converging to $x=0$ will end up in a forbidden ``small 
sector'' (see also figure \ref{figure12} for the case $\Re \sigma=1$).  

\begin{figure}
\epsfxsize=15cm
\centerline{\epsffile{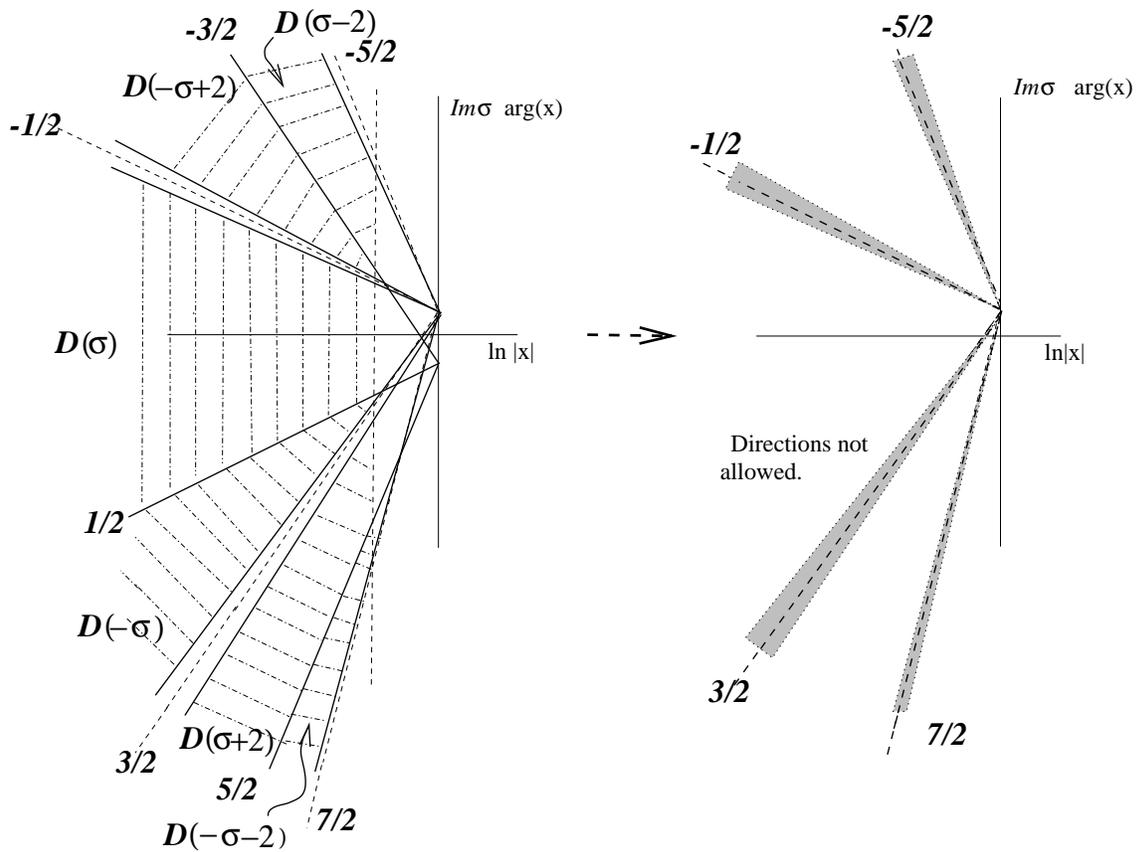}}
\caption{ Domains for $\sigma={1\over 2} +i \Im \sigma$. The numbers close to
the lines are their slopes. 
The ``small sectors''
around the dotted lines 
represented in the right figure are not contained in the union of the
domains. If $x\to 0$ along a direction which ends in one of these sectors, we
do not know the behaviour of the transcendent.}
\label{figure8}
\end{figure}

 If we draw, for the same $\theta_2$,  
 the domains $B(\sigma)$, $B(-\sigma)$, $B(-\sigma+2)$, etc, 
 defined in 
 (\ref{come se non bastasse!!}) we obtain strips in the $(\ln|x|,
 \Im\sigma\arg(x) )$-plane which are {\it certainly forbidden} to theorem 1
 (see figure \ref{figura18}).
 In the strips we know nothing about the transcendent. We guess
 that there might be poles there, as we verify in one example later.

\begin{figure}
\epsfxsize=12cm
\centerline{\epsffile{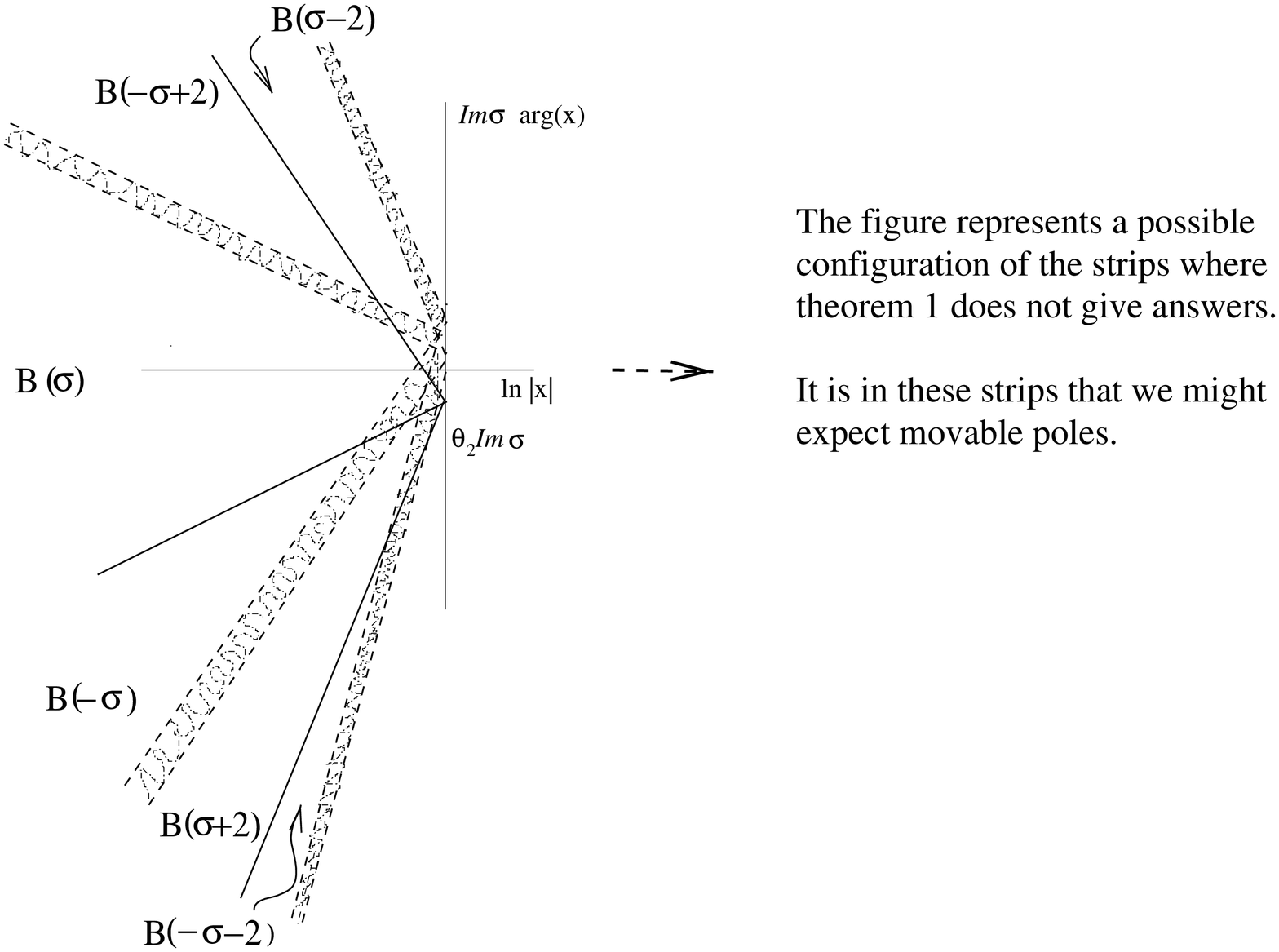}}
\caption{ }
\label{figura18}
\end{figure}


 \vskip 0.2 cm 

What is the behaviour along the directions not described by theorem 1? 
 In the very particular case $(x_0,x_1,x_{\infty}) \in \{ (2,2,2), (2,-2,-2),
(-2,-2,2), (-2,2,-2)\}$ (and so $ x_0 = \pm 2$ !), it is known that 
$PVI_{\mu=-1/2}$ has a 1-parameter family of {\it classical} 
solutions \cite{M}. The  asymptotic behaviour of a branch for {\it radial} 
convergence to the critical points 0, 1 ,$\infty$ was computed in \cite{M}:  
$$ 
   y(x)=\left\{\matrix{          
                      -\ln(x)^{-2}(1+O(\ln(x)^{-1})),~~~~~x\to 0 \cr\cr
                  1+ \ln(1-x)^{-2}(1+O( \ln(1-x)^{-1})),~~~~~x\to 1 \cr\cr
                   -x\ln(1/x)^{-2}(1+O(\ln(1/x)^{-1}),~~~~~x\to \infty \cr
     } \right.
$$ 
The branch is specified by $|\arg(x)|<\pi$, $|\arg(1-x)|<\pi$.  
This behaviour is completely different from $\sim a(x) \tilde{x}^{1-\sigma}$ 
 as $x\to 0$. Intuitively, as $x_0$ approaches the value 2,  
$1-\sigma$ approaches 0 and the decay of $y(x)\sim a x^{1-\sigma}$ becomes 
logarithmic. These solutions were called {\it Chazy solutions} in \cite{M}, 
because they can be computed as functions of solutions of the 
 Chazy equation. 

 \vskip 0.3 cm

 This section is devoted to the investigation of
the local behaviour at $x=0$ of the analytic continuation of a branch 
in the regions not described in  theorem 1.


\subsection{ Elliptic Representation}

 The transcendents of  $PVI_{\mu}$ can be represented in the elliptic form
 \cite{fuchs} 
$$ 
  y(x)= {\cal P} \left({u(x)\over 2};\omega_1(x),\omega_2(x)\right)
$$
where ${\cal P}(z;\omega_1,\omega_2)$ is the Weierstrass elliptic function of
 half-periods $\omega_1$, $\omega_2$.  u(x) solves the non-linear differential
 equation 
\be
  {\cal L}(u)= {\alpha\over x(1-x)} {\partial \over \partial u} \left(
    {\cal P}\left( {u\over
    2};\omega_1(x),\omega_2(x)\right)\right),~~~\alpha={(2\mu-1)^2\over 2}
\label{difficilissima}
\ee
where the differential linear  operator ${\cal L}$ applied to $u$ is 
$$ 
   {\cal L}(u):= 
x(1-x)~{d^2 u\over dx^2}+(1-2x)~{du\over dx}-{1\over 4}~u.  
$$
 The half-periods are two independent solutions of ${\cal L}(u)=0$ normalized
 as follows: 
$$
\omega_1(x):= {\pi \over 2} ~F\left( x\right),~~~~
\omega_2(x):= -{i\over 2} [F(x) \ln x +F_1(x)] 
$$
where $F(x)$ is the hypergeometric function 
$$
F(x):=F\left({1\over 2},{1\over 2},1;x\right)=  
\sum_{n=0}^{\infty}{ \left[\left({1 \over 2}\right)_n\right]^2
  \over (n!)^2 } x^n,
$$
and 
$$
F_1(x):=  
\sum_{n=0}^{\infty}{ \left[\left({1 \over 2}\right)_n\right]^2
  \over (n!)^2 } 2\left[ \psi(n+{1\over 2}) - \psi(n+1)\right]
x^n
$$
$$\psi(z) = 
{d \over dz}\ln \Gamma(z),~~~ \psi\left({1\over 2}\right) = -\gamma -2 \ln
2,~~~ \psi(1)=-\gamma,~~~\psi(a+n)=\psi(a)+\sum_{l=0}^{n-1} {1\over a+l}.
$$

The solutions $u$ of (\ref{difficilissima}) are in general unknown. Before
trying to study them  in general terms, we solve it in a special case:

 \vskip 0.3 cm
\noindent
{\it Example:} The equation $PVI_{\mu=1/2}$ has a two parameter family of 
solutions discovered  
 by Picard \cite{Picard} \cite{Okamoto}  \cite{M}. It is easily obtained
from (\ref{difficilissima}). Since  $\alpha=0$, $u$ solves the  hypergeometric
equation ${\cal L}(u)=0$ and has the general form 
$$
{u(x)\over 2}:= \nu_1 \omega_1(x)+\nu_2 \omega_2(x),~~~ \nu_i \in {\bf C},
~~~0\leq \Re\nu_i <2,~~~(\nu_1,\nu_2)\neq (0,0),
$$
 A branch of $y(x)$ is specified by a 
branch of $\ln x$. 
The monodromy data computed in \cite{M} are 
$$ x_0=-2\cos \pi r_1,~~~x_1= -2\cos \pi r_2,~~~x_{\infty}=-2\cos \pi r_3,$$
$$
   r_1={\nu_2\over 2},~~~r_2=1-{\nu_1\over 2},~~~r_3={\nu_1-\nu_2\over 2},~~~~
\hbox{ for } \nu_1>\nu_2$$
$$
   r_1=1-{\nu_2\over 2},~~~r_2={\nu_1\over 2},~~~
r_3={\nu_2-\nu_1\over 2},~~~~
\hbox{ for } \nu_1<\nu_2.
$$
 The modular parameter is now a function of $x$. Since we are interested in it
 when $x\to 0$ we give its expansion:  
$$ 
      \tau(x) = {\omega_2(x) \over \omega_1(x)}=  {1 \over \pi} 
                    (\arg x - i \ln|x|) +{4i\over \pi}\ln 2+O(x),~~~x \to 0.
$$
We see that $\Im \tau >0$ as $x \to 0$.  Now, if 
\be
\left|\Im {u(x)\over 4 \omega_1}\right|< \Im \tau,
\label{tau}
\ee
 we can expand the 
Weierstrass function in Fourier series.  Condition (\ref{tau}) becomes
$$
     {1\over 2} \left| \Im \nu_1 +{\Im \nu_2 \over \pi} \arg(x)-
        {\Re \nu_2 \over \pi} \ln|x|+{4\ln 2\over \pi} \Re\nu_2\right|<
        -{\ln|x| \over  
                             \pi} +{4\ln 2\over \pi} +O(x),~~~~\hbox{ as }
x\to 0
$$
namely,
\be 
   (\Re \nu_2 +2) \ln|x| - \pi \Im \nu_1-4\ln 2 ~(\Re\nu_2+2) < \Im \nu_2 \arg(x) 
< (\Re \nu_2 -2) \ln|x| - \pi \Im \nu_1-4\ln 2 ~(\Re \nu_2 -2).
\label{dominio}
\ee
or
$$\hbox{Any value of } \arg(x) \hbox{ if }  \Im \nu_2=0.
$$  
The Fourier expansion is 
$$
y(x) = {x+1\over 3} +{1 \over F(x)^2 }\left[  {1 \over \sin^2\left( 
-{1\over 2} [i \nu_2(\ln(x) +{F_1(x) \over F(x)})-\pi \nu_1] \right)}
-{1\over 3} + \right.
$$
$$\left.
 +8 \sum_{n=1}^{\infty} { x^{2n} \over e^{-2 n {F_1(x)\over 
F(x)}}
- x^{2n}} \sin^2\left( 
-{n\over 2} [i \nu_2(\ln(x) +{F_1(x) \over F(x)})-\pi \nu_1] \right) \right]
$$
$$
= {x \over 2} +(1 -{x\over 2} +O(x^2)) \left[ {1\over \sin^2\left( 
-{1\over 2} [i \nu_2(\ln(x) +{F_1(x) \over F(x)})-\pi \nu_1] \right)}+
\right.
$$
$$
\left.
 -{1 \over 4} \left[{e^{i\pi \nu_1} \over 
16^{\nu_2-1}}\right]^{-1} x^{2-\nu_2} +
O(x^2+ x^{3-\nu_2}+ x^{4 - \nu_2})  \right], ~~
~~~x\to 0 \hbox{ in the domain } 
(\ref{dominio})
$$
 As far as { \it radial} convergence is concerned, we have:

\vskip 0.2 cm
a) $0<\Re \nu_2<2 $,
$$
     {1 \over \sin^2(...)}= -{1 \over 4} \left[{e^{i\pi \nu_1} \over 
16^{\nu_2-1}} \right]~x^{\nu_2} ~\bigl( 1 +O(|x^{\nu_2}|)\bigr),$$
and so 
\be
y(x)=  \left\{-{1 \over 4} \left[{e^{i\pi \nu_1} \over 
16^{\nu_2-1}} \right]~x^{\nu_2}+ {1\over 2} x 
-{1 \over 4} \left[{e^{i\pi \nu_1} \over 
16^{\nu_2-1}}\right]^{-1} x^{2-\nu_2}\right\} ~\left(1 +O(x^{\delta}) \right),
~~~~\delta>0, \label{behaviour}
\ee
in accordance with theorem 1. We 
 can identify $1- \sigma $ with $\nu_2$ for $0<\Re \nu_2<1$, or  with 
$2-\nu_2$ for 
$1<\Re \nu_2<2$. In the case $\Re \nu_2 =1$ the three terms $x^{\nu_2}$, 
$x$, $x^{2-\nu_2}$ have the same order and we find again the behaviour of 
theorem 1 
$$
y(x) = \left\{a x^{\nu_2} +{x \over 2} +{1\over 16 a} x^{2-\nu_2}
\right\}(1+O(x^{\delta}))= a x^{\nu_2} \left\{1
+{1\over 2a} x^{-i \Im \nu_2}+ {1\over 16 a^2}
 x^{-2i \Im \nu_2}\right\}(1+O(x^{\delta})),
$$
where $a= -{1\over 4} \left[{e^{i\pi \nu_1} \over 
16^{\nu_2-1}} \right]$.

\vskip 0.2 cm
b) $\Re \nu_2 =0$. Put $ \nu_2 = i \nu$ ( namely, $\sigma =1 -i \nu$ ). The 
domain (\ref{dominio}) is now (for sufficiently small $|x|$):
\be
  2 \ln|x|-\pi \Im \nu_1-8\ln 2<\Im \nu_2 \arg(x) < -2 \ln|x| - \pi \Im
\nu_1+8\ln 2,
\label{dominio1}
\ee
or 
  $$
 2 \ln|x|+\pi \Im \nu_1-8\ln 2<\Im \sigma \arg(x) < -2 \ln|x| + \pi \Im
\nu_1+8\ln 2.
$$
For radial convergence we have 
$$
  y(x)= {1 +O(x) \over \sin^2( {\nu\over 2} \ln(x) + {\nu \over 2} {F_1(x) 
\over F(x)} + \pi \nu_1)}  + O(x). 
$$ 
This is an oscillating functions, and it may have poles. 
 Suppose for example that $\nu_1$ is real.   
Since $F_1(x)/F(x)$ is a convergent power series ($|x|<1$) with real 
coefficients and defines a  bounded function,  then $y(x)$ has a 
sequence of poles on the positive real axis, converging to $x=0$. 

If    
$\nu=0$, namely $\nu_2=0$ (and then $x_0=2$) 
$$ 
y(x) = {1\over \sin^2(\pi \nu_1)} (1+0(|x|)).
$$   

 In the domain (\ref{dominio1})  spiral convergence of $x$ to zero is also 
allowed. In that case, the non-constancy of $\arg(x)$ still gives a behaviour 
(\ref{behaviour}). 

\vskip 0.2 cm 

 The case b) in the  above example is good to 
 understand the limits of theorem 1 in giving a complete description of 
the behaviour of  Painlev\'e transcendents. Actually, theorem 1 (together with the transformation $\sigma \to - \sigma$)  yields  
the behaviour (\ref{behaviour}) 
in the domain $D(\sigma)\cup D(-\sigma)$ ($\Re \sigma=1$):
$$
 (1+\tilde{\sigma})\ln|x| + \theta_1 \Im \sigma \leq \Im \sigma \arg x 
\leq (1-\tilde{\sigma}) \ln|x| +\theta_1 \Im \sigma ,
$$ 
 where 
radial convergence to $x=0$ is not allowed.  
On the other hands, the transformations $ \sigma \to\pm( \sigma -2)$, 
gives a further domain $ D(\sigma -2 ) \cup D(-\sigma +2)$:  
$$
  (-1+ \tilde{\sigma})\ln|x| + \theta_1 \Im \sigma \leq \Im \sigma \arg x
\leq -(1+\tilde{ \sigma}) \ln|x| +\theta_1 \Im \sigma.
$$
but again it is not possible for $x$ to converge to $x=0$ along a radial
path. 
Figure \ref{figure12} 
shows $D(\sigma)\cup D(-\sigma)\cup D(2-\sigma)\cup D(\sigma -2)$. 
Note that a radial  path 
would be allowed if it were possible  to make 
$\tilde{\sigma} \to 1$ and 
the interior of the set obtained as the limit 
for 
$\tilde{\sigma}\to 1$  of  
 $D(\sigma)\cup D(-\sigma)\cup D(2-\sigma)\cup D(\sigma -2)$  has the shape of 
 (\ref{dominio1}). Actually, the intersection of the two sets is never
empty. On (\ref{dominio1}) the 
 elliptic representation predicts an oscillating
behaviour and poles along the paths not allowed by theorem 1.   
So now it must be definitely clear that the 
``limit'' of theorem 1 for $\tilde{\sigma} \to 1$ is not trivial.

\begin{figure}
\epsfxsize=15cm
\centerline{\epsffile{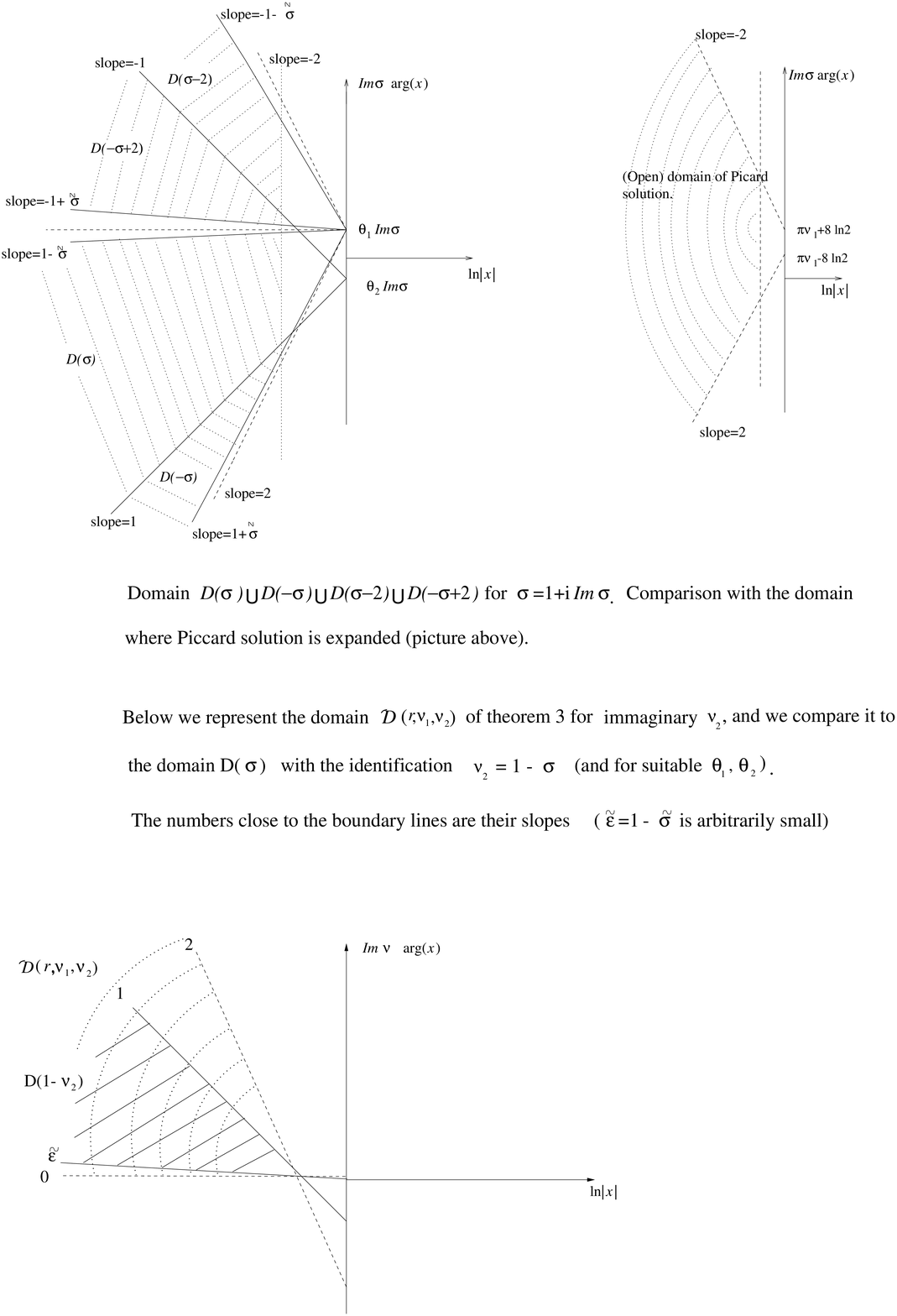}}
\caption{ }
\label{figure12}
\end{figure}

 \vskip 0.2 cm

\noindent
{\it Remark on the example:} For $\mu$ half integer all the possible 
values of $(x_0,x_1,x_{\infty})$ such that
$x_0^2+x_1^2+x_{\infty}^2-x_0x_1x_{\infty} =4$ 
are covered by Chazy and Picard's solutions, 
with the warning that for $\mu={1 \over 2}$ the image (through birational
transformations) of Chazy solutions  
is $y=\infty$. See \cite{M}.

\vskip 0.3 cm

 We turn to the general case. In section \ref{provaellittica} we will prove
 the following theorem:

\vskip 0.2 cm 
\noindent
{\bf Theorem 3:} {\it For any complex $\nu_1$, $\nu_2$ such that 
$$ 
  \nu_2\not \in (-\infty,0]\cup [2,+\infty)
$$
there exists a sufficiently small $r$ such that 
$$
  y(x)= {\cal P}(\nu_1 \omega_1(x)+\nu_2 \omega_2(x)
  +v(x);\omega_1(x),\omega_2(x))
$$
in the domain 
$$
  {\cal D}(r;\nu_1,\nu_2):= \left\{ x\in \tilde{\bf C}_0~ \hbox{ such that }
  |x|<r, \left|{e^{-i\pi \nu_1}\over 16^{2-\nu_2}} x^{2-\nu_2} \right|<r,
\left| 
{e^{i\pi \nu_1} \over 16^{\nu_2}} x^{\nu_2}\right|<r \right\}
$$
The function $v(x)$ is holomorphic  in $ {\cal D}(r;\nu_1,\nu_2)$ and has
convergent expansion 
\be
  v(x)= \sum_{n\geq 1} a_n x^n +\sum_{n\geq 0,~m\geq 1} b_{nm} x^n
  \left({e^{-i\pi \nu_1}\over 16^{2-\nu_2}} x^{2-\nu_2}\right)^m +\sum_{n\geq
  0,~m\geq 1}c_{nm} x^n \left( 
{e^{i\pi \nu_1} \over 16^{\nu_2}} x^{\nu_2}\right)^m     
\label{vdix}
\ee
where $a_n$, $b_{nm}$, $c_{nm}$ are rational functions of  $\nu_2$. 
Moreover, there exists a constant $M(\nu_2)$ depending on $\nu_2$ such that 
 $v(x)\leq M(\nu_2) \left(|x|+\left|{e^{-i\pi \nu_1}\over 16^{2-\nu_2}} x^{2-\nu_2} \right|+\left| 
{e^{i\pi \nu_1} \over 16^{\nu_2}} x^{\nu_2}\right| \right)$ in  ${\cal
D}(r;\nu_1,\nu_2)$ .
}

\vskip 0.2 cm 

 The domain ${\cal D}(r;\nu_1,\nu_2)$  is
\be
|x|<r,~~~~   \Re \nu_2 \ln|x|+ C_1-\ln r < \Im \nu_2
   \arg x < (\Re \nu_2 -2)\ln|x| +C_2 + \ln r,
\label{multiply defined}
\ee
$$
C_1:= -\bigl[4 \ln 2 \Re \nu_2 + \pi \Im \nu_1\bigr],~~~~C_2:=C_1+8\ln 2.
$$
If $\nu_2$ is real (therefore $0<\nu_2<2$), the
domain is simply $|x|<r$. 
The critical behaviour is obtained expanding  $y(x)$ in Fourier series: 
\be
   {\cal P}\left({u\over 2};\omega_1,\omega_2\right)= -{\pi^2 \over 12
   \omega_1^2} +{2\pi^2\over \omega_1^2} \sum_{n=1}^{\infty} {n e^{2\pi i n
   \tau} \over 1- e^{2\pi i n \tau}} \left(1- \cos \left(n {\pi u \over 2
   \omega_1} \right)\right) +{\pi^2\over 4 \omega_1^2} ~{1\over
   \sin^2\left( {\pi u \over 4 \omega_1}\right)} 
\label{rondine}
\ee
The expansion can be performed if $\Im \tau(x)>0$ and $\left|\Im
\left({u(x)\over 
\omega_1(x)}\right)\right|<\Im \tau$; these conditions are satisfied in ${\cal
D}(r;\nu_1,\nu_2)$. Let's put $F_1/F=-4\ln 2 +g(x)$. It follows that
 $g(x)=O(x)$.  
 Taking into account (\ref{rondine}) and  theorem 3, 
the expansion of $y(x)$ for $x\to 0$, $x\in {\cal D}(r;\nu_1,\nu_2)$,   is
$$
  y(x)= \left[{1+x\over 3} - {\pi^2\over 12 \omega_1(x)^2}\right]
 +{\pi^2\over
  \omega_1(x)^2 } \sum_{n=1}^{\infty} {n\over 1- \left({e^{g(x)}\over
  16}\right)^{2n} x^{2n}}\left\{
2\left({e^{g(x)}\over
  16}\right)^{2n} x^{2n}\right.
$$
$$\left. - e^{n(\nu_2+2)g(x)}\left[ {e^{i\pi \nu_1}\over
  16^{2+\nu_2}} x^{2+\nu_2} \right]^n e^{in\pi {v(x)\over \omega_1(x)}} - 
e^{n(2-\nu_2)g(x)} \left[{e^{-i\pi \nu_1}\over
  16^{2-\nu_2}}x^{2-\nu_2}\right]^n e^{-i\pi n {v(x)\over \omega_1(x)} } 
\right\}
$$
$$
+ {\pi^2\over 4\omega_1(x)^2}{1\over   \sin^2\left(-i{\nu_2\over 2}
  \ln x +i{\nu_2\over 2} \ln 16 +{\pi \nu_1\over 2} - i{\nu_2 \over 2} g(x)
  +{\pi v(x) \over 2 \omega_1(x)}\right) }
$$
where 
$$
{1\over \sin^2(...)}= -{4\over e^{\nu_2 g(x) } {e^{i\pi \nu_1}\over 16^{\nu_2}
} x^{\nu_2}e^{i\pi {v(x)\over \omega_1(x)}}+e^{-\nu_2 g(x) } {e^{-i\pi \nu_1}\over 16^{-\nu_2}
} x^{-\nu_2}e^{-i\pi {v(x)\over \omega_1(x)}} -2}
$$
We also observe that  $\omega_1(x) \equiv {\pi \over 2} F(x) =
{\pi \over 2} ( 1 +{1\over 4} x +O(x^2))$, 
$ 
{1+x\over 3} - {\pi^2\over 12 \omega_1(x)^2}\equiv {1+x\over 3} -{1\over 3
F(x)} =  {1\over 2} x \bigl(1+O(x)\bigr)
$, $ e^{g(x)}= 1 +O(x)$ and 
$$
e^{\pm i \pi {v(x)\over \omega_1(x)}}= 1+ O\left(|x|+\left|{e^{-i\pi
\nu_1}\over 16^{2-\nu_2}} x^{2-\nu_2} \right|+\left|  
{e^{i\pi \nu_1} \over 16^{\nu_2}} x^{\nu_2}\right| \right)
$$

In order to single out the leading terms, we observe that we are dealing with
the powers $x$, $x^{2-\nu_2}$, 	$x^{\nu_2}$ in ${\cal D}(r;\nu_1,\nu_2)$. If
$0<\nu_2<2$ (the only allowed real values of $\nu_2$) $|x^{\nu_2}|$ is leading
if $0<\nu_2<1$ and $|x^{2-\nu_2}|$ is leading if $1<\nu_2<2$. We have 

$$
 {1\over \sin^2(...)} = -4 {e^{i\pi \nu_1}\over 16^{\nu_2} } x^{\nu_2}
 \left[1+O(|x|+|x^{\nu_2}|+|x^{2-\nu_2}|)\right] 
$$ 
Thus, there exists $0<\delta<1$ (explicitly computable in terms of $\nu_2$)  
such that
$$ 
  y(x)= \left[{1\over 2} x - {1\over 4} \left[{e^{i\pi \nu_1} \over
  16^{\nu_2-1} }  \right] x^{\nu_2} - {1\over 4} \left[{e^{i\pi \nu_1} \over
  16^{\nu_2-1} }  \right]^{-1} x^{2-\nu_2}
  \right](1+O(x^{\delta}))
$$
$$
= \left\{ \matrix{ \left[{1\over 2}  - {1\over 4} \left[{e^{i\pi \nu_1} \over
  16^{\nu_2-1} }  \right]  - {1\over 4} \left[{e^{i\pi \nu_1} \over
  16^{\nu_2-1} }  \right]^{-1} 
  \right]~x~(1+O(x^{\delta})), ~~~\hbox{ if } \nu_2=1 \cr\cr
 - {1\over 4} \left[{e^{i\pi \nu_1} \over
  16^{\nu_2-1} }  \right] x^{\nu_2}(1+O(x^{\delta}))~~~\hbox{ if } 0<\nu_2<1
  \cr\cr 
- {1\over 4} \left[{e^{i\pi \nu_1} \over
  16^{\nu_2-1} }  \right]^{-1} x^{2-\nu_2}(1+O(x^{\delta}))~~~\hbox{ if }
  1<\nu_2<2 \cr
}\right.
$$ 
The behaviour is the one 
 of theorem 1 for $\sigma=0$ in the first case, $\sigma=
1-\nu_2$ in the second, $\sigma=\nu_2-1$ in the third. 

\vskip 0.2 cm

We turn to the case $\Im \nu_2\neq 0$. We consider a path contained in ${\cal
D}(r;\nu_1,\nu_2)$ of equation 
\be
  \Im \nu_2 \arg(x)= (\Re \nu_2 - {\cal V}) \ln|x|+b,~~~~0\leq {\cal V}\leq 2
\label{che il cielo ce la mandi buona}
\ee
with a suitable constant $b$. Thus $|x^{2-\nu_2}|=|x|^{2-{\cal V}} e^b$,
$|x^{\nu_2}|= |x|^{\cal V} e^{-b}$ and 
$$ 
   |x^{\nu_2}|~~\hbox{ is leading for } ~ 0\leq {\cal V}<1,
$$
$$
 |x^{\nu_2}|,~|x|,~|x^{2-\nu_2}|~~\hbox{ have the same order for } ~ 
{\cal V}=1,
$$
$$
|x^{2-\nu_2}|~~\hbox{ is leading for } ~ 1< {\cal V}\leq 2.
$$
If ${\cal V}=0$,  
$$\left| 
{e^{i\pi \nu_1} \over 16^{\nu_2}} x^{\nu_2}\right|<r,~~~~~~\hbox{ but }
~x^{\nu_2}\not\to 0~\hbox{ as } ~x\to 0.
$$
If  ${\cal V}=2$, 
  $$
 \left|{e^{-i\pi \nu_1}\over 16^{2-\nu_2}} x^{2-\nu_2} \right|<r~~~~~~\hbox
{ but } 
~x^{2-\nu_2}\not\to 0~\hbox{ as } ~x\to 0. 
$$ 
This also implies that $v(x)\not \to 0$ as $x\to 0$ 
along the paths with ${\cal V}=0$ or ${\cal V}=2$.  
Therefore we conclude that:

\vskip 0.2 cm 
a) If  $x\to 0$ in ${\cal D}(r;\nu_1,\nu_2)$ along 
(\ref{che il cielo ce la mandi buona}) for ${\cal V}\neq 0,2$, then  
$$ 
 y(x)= \left[{1\over 2} x - {1\over 4} \left[{e^{i\pi \nu_1} \over
  16^{\nu_2-1} }  \right] x^{\nu_2} - {1\over 4} \left[{e^{i\pi \nu_1} \over
  16^{\nu_2-1} }  \right]^{-1} x^{2-\nu_2}
  \right](1+O(x^{\delta})),~~~~0<\delta<1. 
$$
The term $ - {1\over 4} \left[{e^{\pi \nu_1} \over
  16^{\nu_2-1} }  \right] x^{\nu_2}(1+ \hbox{ higher orders})$ is   ${1\over
  \sin(...)^2}$. 
The three leading terms have the same order if the convergence is along a path
asymptotic to (\ref{che il cielo ce la mandi buona}) with
${\cal V}=1$. Otherwise 
$$ 
  y(x) =  - {1\over 4} \left[{e^{i\pi \nu_1} \over
  16^{\nu_2-1} }  \right] x^{\nu_2}(1+O(x^{\delta}))
$$
or 
$$ 
 y(x)= - {1\over 4} \left[{e^{i\pi \nu_1} \over
  16^{\nu_2-1} }  \right]^{-1} x^{2-\nu_2} (1+O(x^{\delta}))
$$ 
according to the path. This is the behaviour of theorem 1 with $1-\sigma=\nu_2$
or $2-\nu_2$.

Let $\nu_2=1-\sigma$ and consider the intersection ${\cal
D}(r;\nu_1,\nu_2)\cap \cal D(\sigma)$ in the $(\ln|x|,\Im \nu_2
\arg(x))$-plane.  See figure \ref{figure12}. We 
 choose $\nu_1$ such that
 $a= - {1\over 4} \left[{e^{i\pi \nu_1} \over
  16^{\nu_2-1} }  \right]$.  According to 
the proposition in  section \ref{Parametrization of a branch through
Monodromy Data -- Theorem 2},  on the intersection we identify  the 
transcendents of the elliptic representation 
 to those of theorem 1, (we do not need to
specify $\theta_1$, $\theta_2$, because the intersection is never empty).
 
Equivalently, we can choose the identification $1-\sigma=2-\nu_2$ and repeat
the argument.

The identification makes it
possible to 
investigate the behaviour of the transcendents of theorem 1 along a 
path (\ref{spirale}) with $\Sigma=1$. (\ref{spirale}) coincides with 
 (\ref{che il cielo ce la mandi buona}) for
${\cal V}=0$ if we define $1-\sigma:=\nu_2$,  or (\ref{che il cielo ce la
mandi buona}) for ${\cal V}=2$ if we define
$1-\sigma=2-\nu_2$.  We discuss the problem in the following two points:

\vskip 0.2 cm
b) If ${\cal V}=0$  the term 
$$
   {1\over \sin^2\left( -i{\nu_2\over 2} \ln x +\left[i {\nu_2\over 2} \ln 16
   +{\pi 
   \nu_1\over 2}\right] - i{\nu_2 \over 2} g(x) + {\pi v(x)\over 2\omega_1(x)} 
\right)} 
$$
$$
= -4 {e^{i\pi \nu_1}\over 16^{\nu_2}}x^{\nu_2} e^{\nu_2 g(x) + i \pi
   {v(x)\over \omega_1(x)} }~{1\over 1-2 {e^{i\pi \nu_1}\over 16^{\nu_2}}x^{\nu_2} e^{\nu_2 g(x) + i \pi
   {v(x)\over \omega_1(x)} }+\left[{e^{i\pi \nu_1}\over 16^{\nu_2}}x^{\nu_2} e^{\nu_2 g(x) + i \pi
   {v(x)\over \omega_1(x)} }\right]^2 }
$$
is {\it oscillating as $x\to 0$} 
and does not vanish. Note that there are no poles because 
the denominator does not vanish in ${\cal D}(r;\nu_1,\nu_2)$ since $\left| 
{e^{i\pi \nu_1} \over 16^{\nu_2}} x^{\nu_2}\right|<r <1 $. We prefer to keep
the trigonometric notation and write 
$$ 
y(x)= O(x) + {1\over F(x)^2} 
   {1\over \sin^2\left( -i{\nu_2\over 2} \ln x +\left[i {\nu_2\over 2} \ln 16
   +{\pi 
   \nu_1\over 2}\right] - i{\nu_2 \over 2} g(x) + { v(x)\over F(x)} 
\right)} 
$$
$$ 
  = {1+O(x)\over \sin^2\left( -i{\nu_2\over 2} \ln x +\left[i {\nu_2\over 2} \ln 16
   +{\pi 
   \nu_1\over 2}\right] + \sum_{m=1}^{\infty} c_{0m}(\nu_2) \left[{e^{i\pi \nu_1}\over 16^{\nu_2}}x^{\nu_2}\right]^m 
\right)}+O(x) 
$$
The last step is obtained taking into account the non vanishing term 
in (\ref{vdix}) and  $ {\pi v(x)\over 2\omega_1(x)}= {v(x) \over F(x)} =
v(x)(1+O(x))$.  

\vskip 0.2 cm 
c)   If ${\cal V}=2$  the  series   
$$-\sum_{n=1}^{\infty}  {n\over 1- \left({e^{g(x)}\over
  16}\right)^{2n} x^{2n}}~
e^{n(2-\nu_2)g(x)} \left[{e^{-i\pi \nu_1}\over
  16^{2-\nu_2}}x^{2-\nu_2}\right]^n e^{-i\pi n {v(x)\over \omega_1(x)} } 
$$
which appears in $y(x)$ is oscillating. Again, $y(x)$ does not vanish. 

\vskip 0.3 cm 

We apply the result to the problem of radial convergence when $\Re
\sigma=1$. We identify $\nu_2=1-\sigma$ and choose $\nu_2=i\nu$, $\nu \neq 0$
real. Let $x\to 0$ in ${\cal D}( r; \nu_1,i\nu)$ along the line $\arg(x)$=
constant (it is the line with ${\cal V}=0$). We have
$$
   y(x)= O(x)+  {1\over F(x)^2} {1\over \sin^2\left({\nu\over 2} \ln x - \nu
   \ln 16 +{\pi \over 2} \nu_1 +{\nu\over 2} g(x) +{\pi v(x) \over 2
   \omega_1(x)} \right)}
$$
We use again  (\ref{vdix}), ${1\over F(x)}=1+O(x)$, $g(x)=O(x)$ and 
 $ {\pi v(x)\over 2\omega_1(x)}= {v(x) \over F(x)} =
v(x)(1+O(x))$. We have
$$
y(x)=O(x)+{ 1+O(x)\over \sin^2\left({\nu\over 2} \ln x - \nu \ln 16 +{\pi
\nu_1\over 2} + \sum_{m=1}^{\infty} c_{om}(\nu) \left[\left({e^{i\pi \nu_1}
\over 16^{i\nu}}\right)x^{i\nu}\right]^m+O(x)\right)}
$$
$$
\equiv
O(x)+{ 1+O(x)\over \sin^2\left({\nu\over 2} \ln x - \nu \ln 16 +{\pi
\nu_1\over 2} + \sum_{m=1}^{\infty} c_{om}(\nu) \left[\left({e^{i\pi \nu_1}
\over 16^{i\nu}}\right)x^{i\nu}\right]^m\right)}
$$
The last step is possible because $\sin(f(x)+O(x))=\sin(f(x))+O(x)=
\sin(f(x))\bigl(1+O(x)\bigr)$ if $f(x)\not  \to 0$ as $x\to 0$; this is our
case for $f(x)={\nu\over 2} \ln x - \nu \ln 16 +{\pi
\nu_1\over 2} + \sum_{m=1}^{\infty} c_{om}(\nu) \left[\left({e^{i\pi \nu_1}
\over 16^{i\nu}}\right)x^{i\nu}\right]^m $ in  ${\cal D}$.

Thanks to the identification at point a), we have
extended the result of  theorem 1 when  
$\Re \sigma=1$ and the convergence to $x=0$ 
 is along a  radial path, provide that the
limitation on $\arg (x)$ imposed in   ${\cal D}( r; \nu_1,i\nu)$ is respected,
namely 
$$ 
  - \pi \Im \nu_1 - \ln r < \nu \arg (x)
$$
The above limitation is the analogous of the limitation imposed by
$B(\sigma,a;\theta_2, \tilde{\sigma})$  of (\ref{come se non bastasse!!}).

\vskip 0.2 cm
 As a last remark we observe that the coefficients in the expansion of $v(x)$
 can be computed by direct substitution of $v$ into the elliptic form of
 $PVI_{\mu}$, the right hand-side being expanded in Fourier series.


\subsection{Shimomura's Representation}

In \cite{Sh} and \cite{IKSY} S Shimomura proved the following statement for 
 the Painlev\'e VI 
equation with any 
value of 
the parameters $\alpha, \beta,\gamma,\delta$. 
\vskip 0.2 cm
{
\it For any complex number $k$ and for any $ \sigma \in {\bf C} - (-\infty,0]
-[1,+\infty)$  there is a 
sufficiently small $r$, depending  on $\sigma$, such that 
the equation $PVI_{\alpha,\beta,\gamma, \delta}$ has a 
holomorphic solution in the domain 
$$
  {\cal D}_s(r;\sigma,k)= \{{x} \in \tilde{C_0} ~|~ |x|<r, 
~|e^{-k}{x}^{1-\sigma}|<r,~|e^k {x}^{\sigma}|<r \}$$ 
with the following representation:
$$
   y({x};\sigma,k)= {1 \over \cosh^2({\sigma-1\over 2}\ln {x}
    +{k\over 2} +{v({x})\over 2})},
$$
where 
$$ 
  v({x})= \sum_{n\geq 1} a_n(\sigma) {x}^n+ \sum_{n\geq 0, 
~m\geq 1} b_{nm}(\sigma) {x}^n (e^{-k}{x}^{1-\sigma})^m 
+ \sum_{n\geq 0,~m\geq 1} c_{nm}(\sigma) {x}^n 
(e^{k}{x}^{\sigma})^m,
$$
 $a_n(\sigma), b_{nm}(\sigma), c_{nm}(\sigma)$ are rational functions of 
$\sigma$ (and they 
may actually be computed recursively by direct substitution into 
the equation), and the series defining $v(\tilde{x})$ is convergent 
(and holomorphic) in ${\cal D}(r;\sigma,k)$. Moreover, there exists a constant 
$M=M(\sigma)$ such that
\be
  |v(x)|\leq M(\sigma) ~\left(|x|+|e^{-k}{x}^{1-\sigma}|+
                     |e^{k}{x}^{\sigma}| \right).
\label{MMM}
\ee
}
\vskip 0.2 cm

The domain ${\cal D}(r;\sigma,k)$ is specified by the conditions:
\be
 |x|<r,~~~ \Re\sigma \ln|x| + [\Re k - \ln r] <
                                                    \Im \sigma \arg(x)
           < (\Re \sigma -1 )\ln |x| + [\Re k + \ln r].
\label{UFFAUFFA}
\ee
 This is an open domain in the plane $(\ln|x|,\arg(x))$. It can  be compared 
with 
the domain $D(\epsilon;\sigma,\theta_1,\theta_2)$ of theorem 1 (figure
\ref{figure10}). 
Note that (\ref{UFFAUFFA})
imposes a  limitation on $\arg(x)$. The situation is analogous to the
elliptic representation. For example, if $\Re \sigma=1$ we have 
$$ 
  \Im \sigma \arg(x)< [\Re k + \ln r],~~~~~(\ln r<0)
$$
 After the identification of the Shimomura's transcendent with those of
 theorem 1 (see point a.1) below), the above limitation turns out to be 
 the analogous of
 the limitation imposed to $D(\epsilon;\sigma;\theta_1,\theta_2)$  by
$B(\sigma,a;\theta_2, \tilde{\sigma})$  of (\ref{come se non bastasse!!}).

\vskip 0.2 cm 

 Like the elliptic representation,  Shimomura's 
 allows us to investigate what 
happens when ${x}\to 0$ along a path (\ref{spirale}) with $\Sigma=1$, 
contained in ${\cal D}_s(r;\sigma,k)$. It is a  radial path if $\Re
\sigma=1$.    Along  (\ref{spirale}) 
we have $|x^{\sigma}|=  |x|^{\Sigma} e^{-b}$.  
We suppose $\Im \sigma\neq 0$. 
\vskip 0.2 cm 

  a) $0\leq \Sigma <1$. 
 We observe that $|x^{1-\sigma} e^{-k}|\to 0$ as $x\to 0$ 
along the line.  Then:
$$
    y({x};\sigma,k)= {1 \over \cosh^2({\sigma-1\over 2}\ln x
    +{k\over 2} +{v({x})\over 2})}= {4 \over 
    x^{\sigma-1}e^ke^{v(x)}+ x^{1-\sigma}e^{-k} e^{-v(x)} +2},
$$
$$
= 4 e^{-k} e^{-v(x)} x^{1-\sigma} ~{1\over 
(1+e^{-k} e^{-v(x)} x^{1-\sigma})^2}= 4 e^{-k} e^{-v(x)} x^{1-\sigma} ~
\left( 1+e^{-v(x)} O(|e^{-k}x^{1-\sigma}|)\right) .
$$

\vskip 0.2 cm
 Two sub-cases:

\vskip 0.2 cm
a.1) $\Sigma \neq 0$. Then $|x^{\sigma} e^k| \to 0$ and $v(x) \to 0$ 
(see (\ref{MMM})). Thus
$$
     y({x};\sigma,k)= 4e^{-k} x^{1-\sigma} \left(1 + O(|x| + 
| e^k x^{\sigma}| + |e^{-k}x^{1-\sigma}|)\right)
$$
Following the proposition in section \ref{Local Behaviour -- Theorem 1}, we
identify 
$y(x;;\sigma,k)$ and $y(x;\sigma,a)$ ($a=4 e^{-k}$) 
 on  $D_s(r;\sigma,k)\cap D(\epsilon;\sigma;\theta_1,\theta_2)$, which is not
empty for any $\theta_1$, $\theta_2$.  See figure \ref{figure10}.

\begin{figure}
\epsfxsize=15cm
\centerline{\epsffile{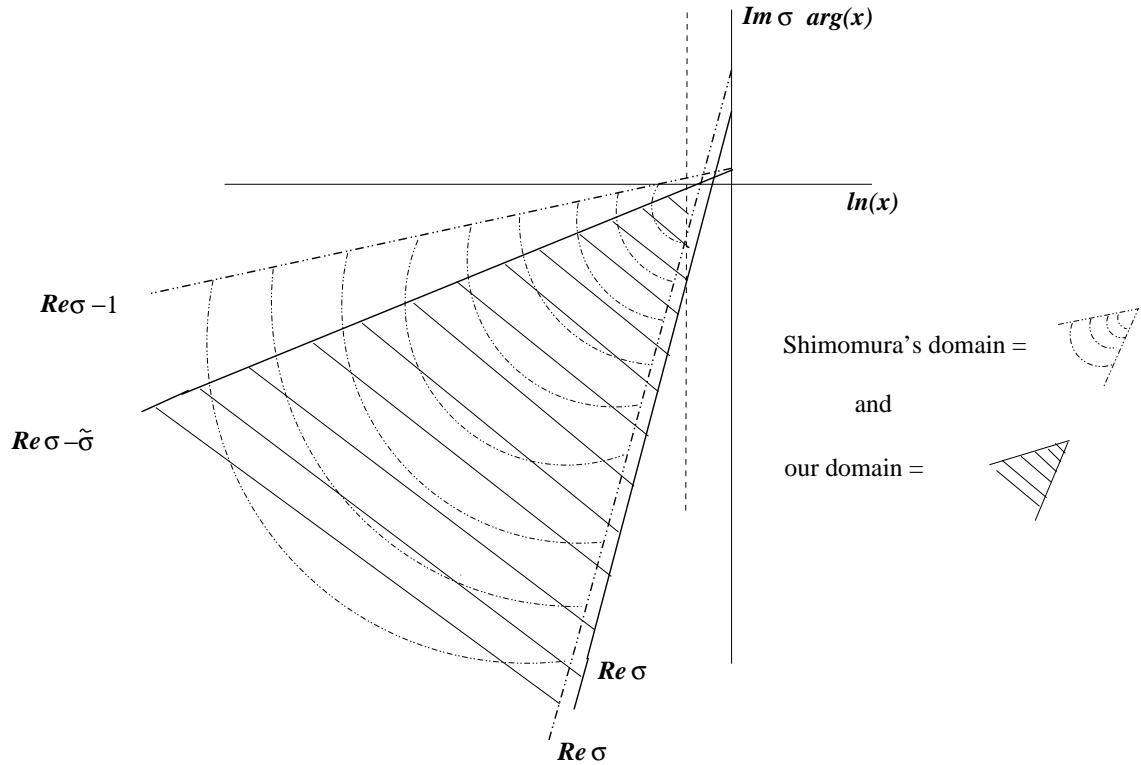}}
\caption{The domains  ${\cal
D}_s(r;\sigma,k)$ and $D(\epsilon;\sigma;\theta_1,\theta_2,\tilde{\sigma})$}
\label{figure10}
\end{figure}

\vskip 0.2 cm

a.2) $\Sigma =0$.  $|x^{\sigma} e^k| \to $ constant$<r$, so $|v(x)|$ does not
vanish.  Then 
$$
y(x)=
 a(x) x^{1-\sigma}~\left( 1+O(|e^{-k}x^{1-\sigma}|)\right),~~~~
  a(x)= 4 e^{-k} e^{-v(x)},
$$ 
which must coincide with the result of theorem 1:
 $$
 y(x) =a\left(1+{C\over 2a} e^{i\alpha(x)}+{C^2\over 16 a^2} e^{2 i \alpha(x)}
 \right)
 x^{1-\sigma}~\left(1 +
            O(|x|^{1-\sigma_1})  \right). 
$$
\vskip 0.2 cm

b) $\Sigma =1$. In this case theorem 1 fails. 
 Now   
 $|x^{1-\sigma} e^{-k}| \to $ (constant$\neq 0)<r$. Therefore $y(x)$ does not
 vanish as $x\to 0$. We keep the representation 
$$
    y({x};\sigma,k)= {1 \over \cosh^2({\sigma-1\over 2}\ln x
    +{k\over 2} +{v({x})\over 2})} \equiv 
{1 \over \sin^2(i{\sigma-1\over 2}\ln x
    +i{k\over 2} +i{v({x})\over 2}-{\pi \over 2})}
$$
 $v(x)$ does not vanish and $y({x})$ is 
oscillating 
as ${x}\to 0$, with no limit. 
We remark that like in the elliptic representation, $\cosh^2(...)$ does not
vanish in ${\cal D}_s(r;\sigma,k)$, so we do not have poles. Figure
\ref{figure11} synthesizes points a.1), a.2), b). 

\begin{figure}
\epsfxsize=15cm
\centerline{\epsffile{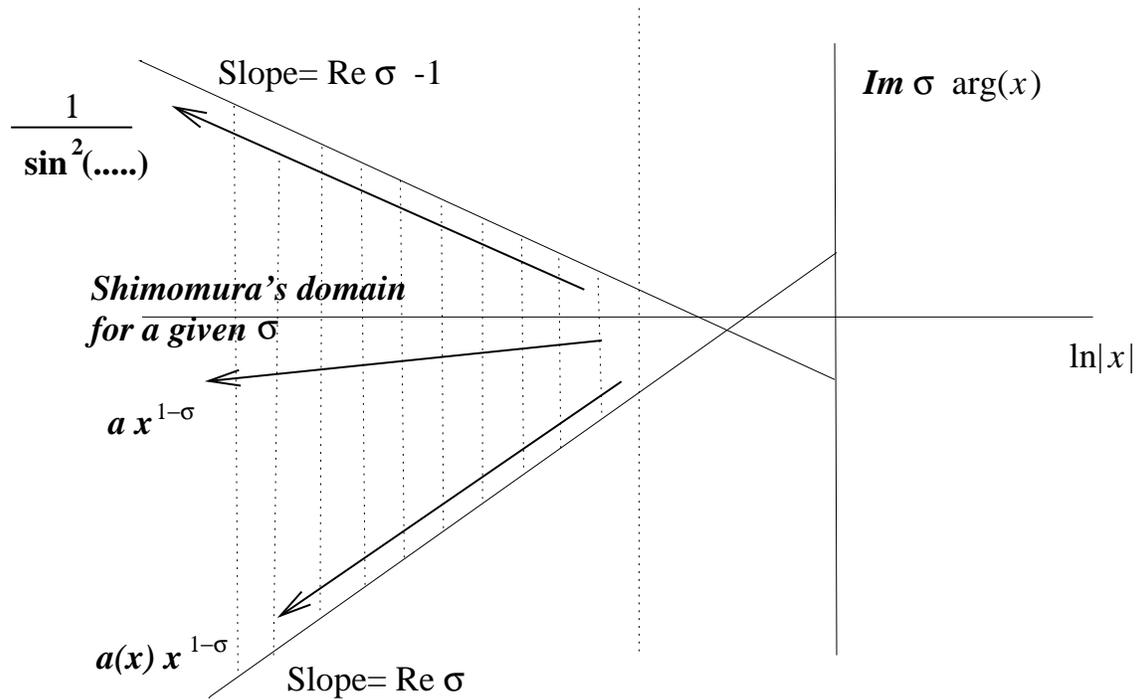}}
\caption{Critical behaviour of $y(x;\sigma,k)$ along different lines in ${\cal
D}_s(r;\sigma,k)$}
\label{figure11}
\end{figure}

\vskip 0.3 cm
\noindent
 As an application, we consider the case $\Re \sigma =1$, namely $
  \sigma = 1 -i \nu,$ $ \nu \in \bf{R}\backslash\{0\}$. 
Then, the path corresponding to $\Sigma =1$ is a {\it  radial}  path 
 in the $x$-plane  and 
$$ 
y({x}; 1 -i \nu,k)= {1 \over \sin^2 \left({\nu \over 2} \ln(x) + 
{i k \over 2} - {\pi \over 2} +i {v(x)\over 2}\right) }$$
$$
= 
 {1+O(x) \over \sin^2 \left({\nu \over 2} \ln(x) + 
{i k \over 2} - {\pi \over 2}+{i\over 2} \sum_{ 
~m\geq 1} b_{0m}(\sigma) (e^{-k}\tilde{x}^{1-\sigma})^m\right)} 
$$

\vskip 0.2 cm



\section{Analytic Continuation of a Branch}\label{Analytic Continuation of a
Branch} 

We describe the analytic continuation of the transcendent 
  $y(x;\sigma,a)$.

\begin{figure}
\epsfxsize= 14cm
\centerline{\epsffile{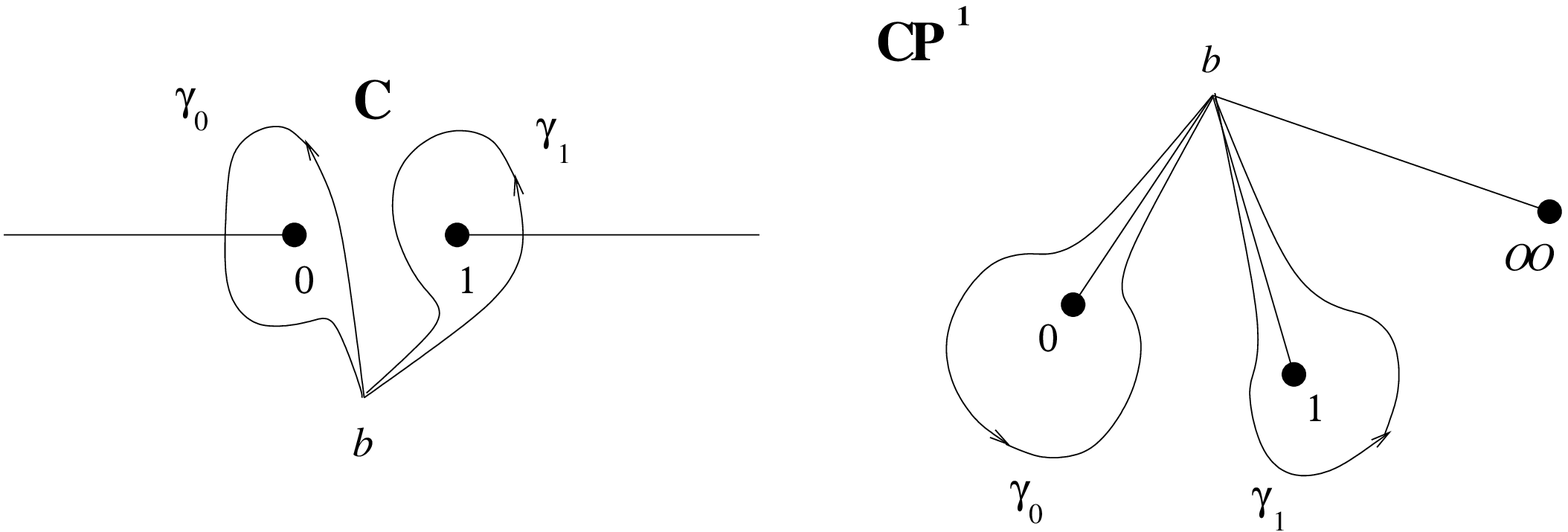}}
\caption{Base point and loops in ${\bf C}\backslash \{0,1\}$ and in ${\bf
  CP}^1\backslash \{0,1,\infty\}$. }
\label{figura2}
\end{figure}

We fix a basis in the fundamental group $\pi({\bf P}^1\backslash \{0,1,\infty\}
 , b)$, where $b$ is the base-point (figure \ref{figura2}), 
and we choose a basis
 $\gamma_0$, $\gamma_1$ of two loops around 0 and 1 respectively. 
   The analytic
 continuation of a branch 
 $y(x;x_0,x_1,x_{\infty})$  along paths encircling $x=0$ 
and $x=1$ (a loop around $x=\infty$ is homotopic  to the  product of
 $\gamma_0$, $\gamma_1$) is given by the action of the group of the 
pure braids on 
the monodromy data (see \cite{DM}).
 Namely, for a counter-clockwise loop around 0 we
 have 
to transform 
($x_0,x_1, x_{\infty}$)  by the action of the braid $\beta_1^2$, where 
$$
  \beta_1: ~~(x_0,x_1,x_{\infty}) \mapsto(-x_0,x_{\infty}-x_0 x_1,x_1)$$
$$ 
    \beta_1^2:~~ 
 ~~(x_0,x_1,x_{\infty}) \mapsto(x_0,~x_1+x_0x_{\infty}-x_1 x_0^2,~x_{\infty}-x_0 
x_1)
$$
For a counter-clockwise loop around 1 we need the braid $\beta_2^2$, where 
$$
 \beta_2: ~~(x_0,x_1,x_{\infty}) \mapsto(x_{\infty},-x_1,x_0-x_1 x_{\infty})
$$
$$
  \beta_2^2:~~(x_0,x_1,x_{\infty}) \mapsto(x_0-x_1x_{\infty},~x_1,~x_{\infty}
+x_0 x_1-x_{\infty}x_1^2)
$$
\vskip 0.2 cm
 A generic loop ${\bf P}^1\backslash \{0,1,\infty\}$ is represented by a braid
  $\beta$, which is a product of factors $\beta_1$ and $\beta_2$. 
Let $\sigma$ and $a$ be associated to $(x_0,x_1,x_{\infty})$ and let $x\in
  D(\sigma)$. At $x$, the branch $y(x;x_0,x_1,x_{\infty})$ coincides with
  $y(x;\sigma, a)$.  The braid $\beta$
 acts on $(x_0,x_1,x_{\infty})$ and
produces a new triple $(x_0^{\beta},x_1^{\beta},x_{\infty}^{\beta})$. 
We plug the new triple into the formulae of 
 theorem 2 and  we obtain the new parameters $\sigma^{\beta}$, $a^{\beta}$ for 
 the new branch $y(x;x_0^{\beta},x_1^{\beta},x_{\infty}^{\beta})$ which
  coincides, at {\it the same} point $x$, with  
  $y(x;\sigma^{\beta}, a^{\beta})$.

To further clarify the concept, consider the
transcendent $y(x;\sigma,a)$ in the point $x\in D(\sigma)$. 
Let us start at
$x$, we perform the loop $\gamma_1$ around $1$ and we go back to $x$. 
The transcendent becomes $y(x;\sigma^{\beta_2^2},a^{\beta_2^2})$. Namely
$$
   \gamma_1:~y(x;\sigma,a) \longrightarrow
   y(x;\sigma^{\beta_2^2},a^{\beta_2^2}).
$$
 
Now let $x\in D(\sigma)$. Suppose that 
$x^{\prime}:=e^{2\pi i }x\in D(\sigma)$ (this assumption is always possible if
$0\leq \Re \sigma<1$; if $\Re \sigma=1$ we may need to consider $D(\sigma)\cup
D(-\sigma)\cup D(\sigma-2)\cup D(2-\sigma)$). The loop $\gamma_0$ transforms $x\mapsto x^{\prime}$ and 
according to theorem 1 we have 
$$ 
  y(x;\sigma,a)\longrightarrow 
y(x^{\prime};\sigma,a)= a [x^{\prime}]^{1-\sigma}
  \left(1+O\left(\left|x^{\prime}\right|^{\delta}\right)\right) 
$$
$$
  = a e^{-2\pi i \sigma} x^{1-\sigma} (1+O(|x|^{\delta}))\equiv y(x;\sigma, a
  e^{-2\pi i \sigma} )
$$
This means that, {\it if we fix a branch cut for  $x$}, 
the analytic continuation of
$y(x;\sigma,a)$ starting at $x$, going around 0 with the loop $\gamma_0$ and
returning back to {\it the same }$x$ is 
\be
   \gamma_0:~y(x;\sigma,a) \longrightarrow
   y(x;\sigma,a e^{-2\pi i \sigma})
\label{analy1111}
\ee
On the other hand, if $x$ is considered as apoint on the universal covering of
${\bf C}_0 \cap \{|x|<\epsilon\}$ we simply have 
$$
   \gamma_0:~y(x;\sigma,a) \longrightarrow
   y(x^{\prime};\sigma,a)
$$
Now we note that
 the transformation of $(\sigma,a)$ according to the braid
$\beta_1$ is 
  \be
\beta_1^2:~~(\sigma,a)\mapsto(\sigma,a e^{-2\pi i \sigma}) 
\label{analy}
\ee
as it follows from the fact that 
 $x_0$ is not affected, then $\sigma$ does not change, and from the explicit 
computation of $a(x_0^{\beta_1^2},x_1^{\beta_1^2},x_{\infty}^{\beta_1^2})$
 through theorem 2  
(we will do  it at the end of  section \ref{dimth2}). Therefore
(\ref{analy1111}) is 
$$
\gamma_0:~y(x;\sigma,a) \longrightarrow
   y(x;\sigma^{\beta_1^2},a^{\beta_1^2})
$$
Thus  theorem 1 is in accordance with the
analytic continuation obtained by the action of the braid group.

\begin{figure}
\epsfxsize=15cm
\centerline{\epsffile{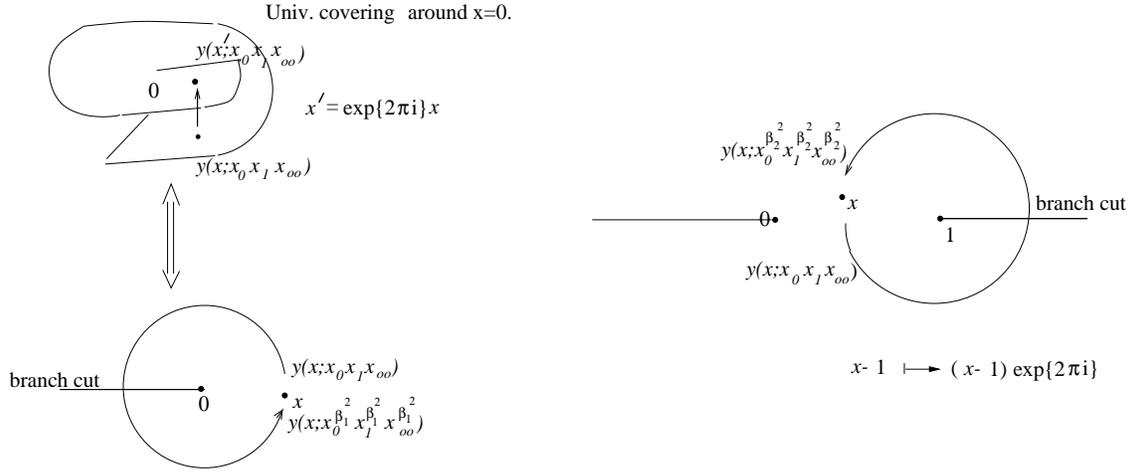}}
\caption{Analytic continuation of a branch for a loop around
$x=0$ and a loop around $x=1$. We also draw the analytic continuation on the
universal covering}
\label{figure16}
\end{figure}



\section{Singular Points $x=1$, $x=\infty$  (Connection
Problem)}\label{Singular Points x=1, x=infty  (Connection Problem)}

 In this subsection we use the notation $\sigma^{(0)}$ and $a^{(0)}$
 to denote the parameters of theorem 1 near the critical point
 $x=0$. We describe now the analogous of theorem 1 near $x=1$ and
 $x=\infty$. 
The three critical points $0$, $1$, $\infty$ are equivalent thanks
  to the symmetries of the $PVI_{\mu}$ equation. 

\vskip 0.15 cm 
 a) Let 
\be
x={1\over t} ~~~~y(x):={1\over t}~\hat{y}(t)
\label{sim1}
\ee
 Then $y(x)$
 is a solution of $PVI_{\mu}$ (variable $x$)  if and only if
 $\hat{y}(t)$ is a solution of $PVI_{\mu}$ (variable $t$). The
 singularities $0$ and $\infty$ are exchanged. Now, we can prove
 theorem 1 near $t=0$. 
 We go back to $y(x)$ and find a transcendent $y(x;\sigma^{(\infty)},
 a^{(\infty)})$ with behaviour  
\be
    y({x}; \sigma^{(\infty)},
 a^{(\infty)}) = a^{(\infty)} {{x}}^{ \sigma^{(\infty)} }\left(
 1+O({1\over |x|^{\delta}}) \right)~~~~~{x} \to \infty
\label{asy2}
\ee
in
$$ 
    D(M; \sigma^{(\infty)};\theta_1,\theta_2,\tilde{\sigma}):=\{ 
 {x}\in 
\widetilde{ {\bf C} \backslash \{\infty\} } 
\hbox{ s.t. } 
|x|>M,~ 
e^{ -\theta_1\Im\sigma^{(\infty)} }
 |x|^{-\tilde{\sigma}}
\leq |{x}^{ -\sigma^{(\infty)}}|
\leq  e^{-\theta_2\Im\sigma^{(\infty)}}
 $$
$$
0\leq
   \tilde{\sigma} < 1
\}
$$
where $M>0$ is sufficiently big and $0<\delta<1$ is small (figure
\ref{figura3}).

\begin{figure}
\epsfxsize=15cm
\centerline{\epsffile{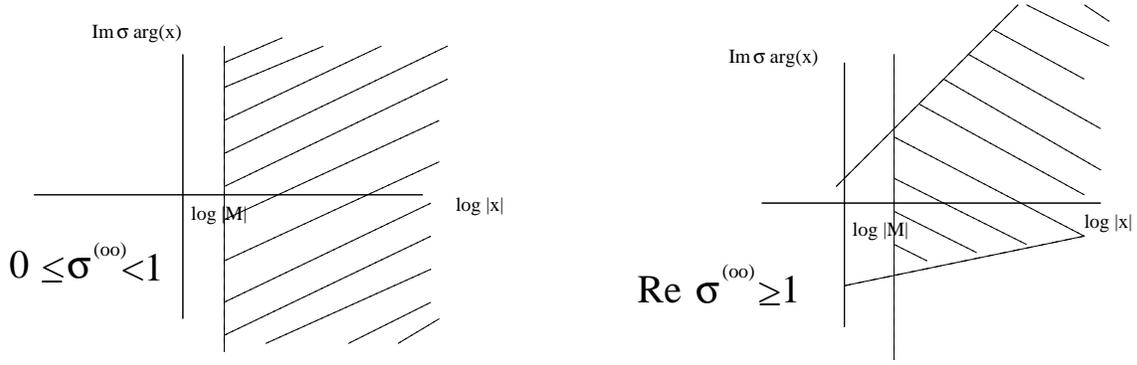}}
\caption{Some examples of the domain
$D(M;\sigma;\theta_1,\theta_2,\tilde{\sigma})$ } 
\label{figura3}
\end{figure}

\vskip 0.15 cm 
b) Let 
\be       
x=1-t,~~~~y(x)=1-\hat{y}(t)
\label{sim2}
\ee
$y(x)$ satisfies $PVI_{\mu}$ if and only if $\hat{y}(t)$ satisfies
$PVI_{\mu}$. Theorem 1 holds for $\hat{y}(t)$ near $t=0$ with coefficients
$\sigma^{(1)}$ and $a^{(1)}$.  Going back to $y(x)$ we obtain a transcendent  
 $y(x;\sigma^{(1)},a^{(1)})$ such that 
\be
y({x},\sigma^{(1)},a^{(1)})= 1-a^{(1)} (1-{x})^{1-\sigma^{(1)}}
(1+O(|1-x|^{\delta }))~~~~{x}\to 1 
\label{asy1}
\ee
in 
$$
  D(\epsilon;\sigma^{(1)};\theta_1,\theta_2,\tilde{\sigma}):= \{ {x}\in
   \widetilde{{\bf C}\backslash \{1\}} 
   \hbox{ s.t. } |1-x|<\epsilon ,~~ e^{-\theta_1 \Im \sigma}
   |1-x|^{\tilde{\sigma}} \leq |(1-{x})^{\sigma^{(1)}}| 
 \leq e^{-\theta_2 \Im
   \sigma} ,$$ 
$$0\leq \tilde{\sigma} < 1 \}$$

\vskip 0.2 cm
We choose a  triple of  monodromy
 data $(x_0,x_1,x_{\infty})$ and we compute the 
 corresponding $\sigma^{(0)} =\sigma^{(0)}
 (x_0)$, $0\leq \Re \sigma^{(0)}\leq 1$, and $a^{(0)}$. We  recall that
 $a^{(0)}$ depends on the triple but also on  $\sigma^{(0)}$, namely for the
 same triple the change $\sigma^{(0)}\mapsto \pm \sigma^{(0)} +2n$ changes
 $a^{(0)}$; so we'll write $ a^{(0)}(\pm \sigma^{(0)} +2n)$.  

 Let  $x\in
 D(\sigma^{(0)})\cup D(-\sigma^{(0)}) \cup D(2-\sigma^{(0)}) \cup
 D(\sigma^{(0)}-2)$ \footnote{This is always
 possible for any $\arg(x)$ if $|x|$ is small enough}. At $x$  
there exists  a unique branch 
 $y(x;x_0,x_1,x_{\infty})$ whose analytic continuation is 
$y\bigl(x;\sigma^{(0)},
 a(\sigma^{(0)})\bigr)$ if 
$x \in    D(\sigma^{(0)})$, or $y\bigl(x;-\sigma^{(0)},
 a(-\sigma^{(0)}))$ if $x \in    D(-\sigma^{(0)}\bigr)$, etc.

The branch $y(x;x_0,x_1,x_{\infty})$ is also defined for $x\in 
D(M;\sigma^{(\infty)}) \cup D(M;-\sigma^{(\infty)}) \cup
D(M;-\sigma^{(\infty)}) \cup 
 D(M;\sigma^{(\infty)}-2)$ 
\footnote{
 The necessity of considering the union of the four domains comes from the fact
that if $\Re \sigma^{(\infty)}=\Re \sigma^{(0)}=1$ it is not obvious 
that we can move from  $x\in D(\sigma^{(0)})$ to $x\in D(M,\sigma^{(\infty)})$ 
keeping $\arg(x)$  fixed, therefore keeping the same branch of the
transcendent.}.   
As it is proved in \cite{DM}
$$ 
   y(x;x_0,x_1,x_{\infty})= {1\over t} \hat{y}(t;x_{\infty},
-x_1,x_0-x_1x_{\infty}),~~~~x={1\over t} 
$$
From the data $(x_{\infty},
-x_1,x_0-x_1x_{\infty})$ we compute $\sigma^{(\infty)}$ and
$a^{(\infty)}(\sigma^{(\infty)})$,  $a^{(\infty)}(-\sigma^{(\infty)})$,
$a^{(\infty)}(\sigma^{(\infty)}-2)$, etc.    
 The analytic continuation of the branch  $y(x;x_0,x_1,x_{\infty})$ is then  
$y(x;\sigma^{(\infty)},a^{(\infty)}(\sigma^{(\infty)}))$ 
if $x\in D(M_1;\sigma^{(\infty)})$,
it is $y(x;-\sigma^{(\infty)},a^{(\infty)}(-\sigma^{(\infty)}))$   if $x\in
D(M_2;-\sigma^{(\infty)})$, etc.  

 The above discussion solves the connection problem between
$x=0$ and $x=\infty$.

\begin{figure}
\epsfxsize=15cm
\centerline{\epsffile{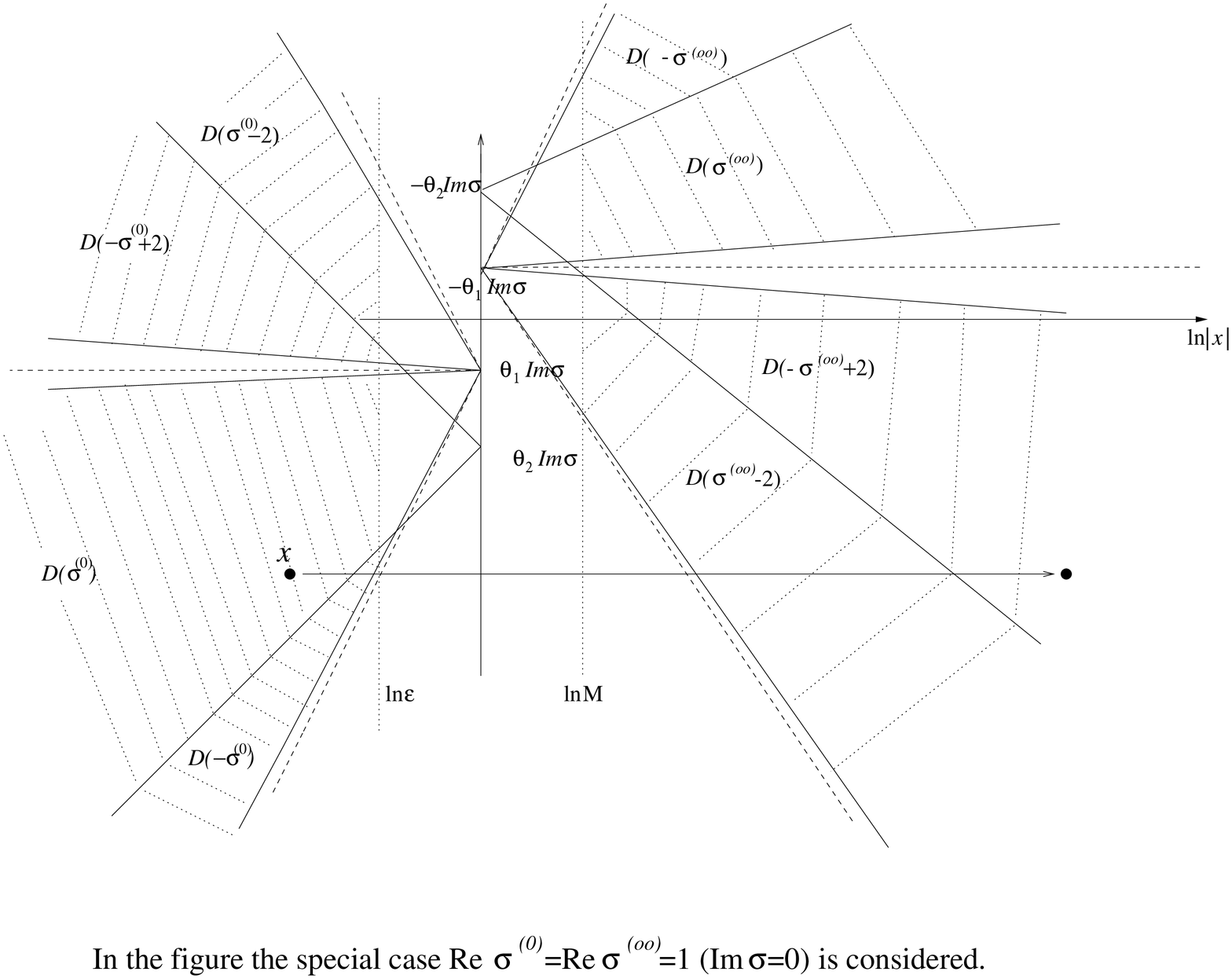}}
\caption{Connection problem for the points $x=0$, $x=\infty$}
\label{figure15}
\end{figure}

 In the same way we solve the connection problem between $x=0$ and $x=1$ by
 recalling that in \cite{DM} it is proved that a branch
 $y(x;x_0,x_1,x_{\infty})$ is 
$$ 
   y(x;x_0,x_1,x_{\infty})= 1- \hat{y}(t; x_1,x_0,x_0
   x_1-x_{\infty}),~~~~x=1-t
$$

 We repeat the same argument.  We remark again that
 for $\Re \sigma^{(0)}
 = \Re \sigma^{(1)}=1 $ it is necessary to consider the union of $
 D(\sigma^{(1)})$, $D(-\sigma^{(1)})$,
 $D(2-\sigma^{(1)})$, $D(\sigma^{(1)}-2)$ to include all
 possible values of $\arg(1-x)$.


\section{ Proof of Theorem 1}\label{proof of theorem 1}

In order to prove theorem 1 we have to recall the connection between
$PVI_{\mu}$ and Schlesinger equations for 2$\times$2 matrices $A_0(x)$,
$A_{x}(x)$, $A_{1}(x)$ 
\be 
\left.\matrix{
               {d A_0\over dx}= {[A_x,A_0]\over x} \cr
\cr               
               {dA_1\over dx}={[A_1,A_x]\over 1-x} \cr
\cr  
               {d A_{x} \over d x} = {[A_x,A_0]\over x}+{[A_1,A_x]\over 1-x}\cr
}\right.
\label{sch}
\ee   
They are the analogous of (\ref{CoMp1}) (\ref{CoMp2}). 
We look for solutions   satisfying  
$$A_0(x)+A_x(x)+A_{1}(x)=\pmatrix{-\mu&0\cr
                                 0 & \mu \cr}:=-A_{\infty}
~~~\mu\in{\bf C},~~~2\mu\not\in{\bf Z}
$$
$$ \hbox{ tr}(A_i)=\det(A_i)=0$$
Now let  
$$ 
A(z,x):={A_0\over z}+{A_x\over z-x}+{A_1\over z-1}
  $$
We explained that  
$y(x)$ is a solution of $PVI_{\mu}$ if and only if $A(y(x),x)_{12}=0$.

\vskip 0.15 cm
 
 The system (\ref{sch}) is a particular case of the system 
\be 
\left. \matrix{
{d A_{\mu}\over dx}=\sum_{\nu=1}^{n_2} [A_{\mu},B_{\nu}]~f_{\mu \nu}(x)
\cr\cr
{dB_{\nu}\over dx}= -{1\over x} \sum_{ \nu^{\prime}=1 }^{ n_2 }[B_{\nu},
B_{\nu^{\prime}}]+
\sum_{\mu=1}^{n_1}[B_{\nu},A_{\mu}]~g_{\mu \nu}(x) +     
\sum_{\nu^{\prime}=1}^{n_2} [B_{\nu},B_{\nu^{\prime}}] ~h_{\nu
\nu^{\prime} }(x) 
}
\right.
\label{sch1}
\ee
where the functions $f_{\mu\nu}$, $g_{\mu \nu}$, $h_{\mu \nu}$ are
meromorphic with poles at $x=1,\infty$ and $\sum_{\nu} B_{\nu} +
\sum_{\mu} A_{\mu}=-A_{\infty}$ (here the subscript $\mu$ is a label,
not the eigenvalue of $A_{\infty}$ !). System (\ref{sch}) is obtained
for $f_{\mu \nu}=g_{\mu \nu}=b_{\nu}/(a_{\mu}-xb_{\nu})$, $h_{\mu
\nu}=0$, $n_1=1$, $n_2=2$, $a_1=b_2=1$, $b_1=0$ 
and $B_1=A_0$, $B_2=A_x$, $A_1=A_1$. 
\vskip 0.15 cm 
 We prove the analogous result of \cite{SMJ}, page 262, for 
 the domain $D(\epsilon;\sigma)$:

\vskip 0.3 cm
\noindent
 {\bf Lemma 1: } {\it  
                  Consider  matrices $B_{\nu}^0$ ($\nu=1,..,n_2$), 
$A_{\mu}^0$ ($\mu=1,..,n_1$) and   $\Lambda$, independent of $x$ and such
                  that 
$$
   \sum_{\nu} B_{\nu}^0 + \sum_{\mu} A_{\mu}^0 = -A_{\infty}
$$
$$
   \sum_{\nu} B_{\nu}^0=\Lambda, ~~~~\hbox{ eigenvalues}(\Lambda)=
   {\sigma \over 2},~-{\sigma\over 2},~~~\sigma \in \Omega
$$
  Suppose that  $f_{\mu \nu}$, $g_{\mu \nu}$, $h_{\mu \nu}$ are holomorphic
if  $|x|<\epsilon^{\prime}$, for some small $\epsilon^{\prime}<1 $.

For any $0<\tilde{\sigma}<1$ and  $\theta_1,\theta_2$ real there exists a
sufficiently small $0<\epsilon<\epsilon^{\prime}$ such that 
the system  (\ref{sch1}) has holomorphic 
solutions 
 $A_{\mu}({x})$, $B_{\nu}({x})$   in
$D(\epsilon;\sigma;\theta_1,\theta_2,\tilde{\sigma})$   satisfying:
$$
   || A_{\mu}({x})-A_{\mu}^0||\leq C~ |x|^{1-\sigma_1}$$
$$
           || {x}^{-\Lambda} B_{\nu}({x})~ {x}^{\Lambda} 
-B_{\nu}^0||\leq C~ |x|^{1-\sigma_1}
$$          
            Here $C$ is a positive  constant and $\tilde{\sigma}<\sigma_1<1$
}

\vskip 0.3 cm
\noindent
{\it Important remark:} There is no need 
to assume here that $2 \mu \not \in {\bf Z}$. 
The theorem holds true for any value of $\mu$.  If in the system (\ref{sch1}) 
the functions $f_{\mu \nu}$, $g_{\mu \nu}$, $h_{\mu \nu}$ are chosen in such a 
way to yield Schlesinger equations for the fuchsian system of $PVI_{\mu}$, 
the assumption $2 \mu \not \in {\bf Z}$ is still not necessary, provided the 
matrix  $R$ is considered as a monodromy datum independent 
of the deformation parameter $x$.

\vskip 0.3 cm
\noindent
{\it Proof:} Let $A(x)$ and $B(x)$ be $2\times 2$ matrices holomorphic
on $D(\epsilon;\sigma)$ (we understand $\theta_1,\theta_2,\tilde{\sigma}$ which
have been chosen) and such that 
$$ 
  ||A(x) ||\leq C_1,~~~~||B(x)||\leq C_2~~~\hbox{ on }
  D(\epsilon;\sigma)
$$
Let $f(x)$ be a holomorphic function for $|x|<\epsilon^{\prime}$. 
Let $\sigma_2$ be a real number such that $\tilde{\sigma}<\sigma_2<1$. 
Then, there exists a 
 sufficiently small $\epsilon<\epsilon^{\prime}$ such that 
for $x\in D(\epsilon;\sigma)$ we have:
$$ 
    ||x^{\pm\Lambda}~A(x)~x^{\mp\Lambda}||\leq C_1 |x|^{-\sigma_2}
$$
$$ 
    ||x^{\pm\Lambda}~B(x)~x^{\mp\Lambda}||\leq C_2 |x|^{-\sigma_2}
$$
$$
\left|\left| x^{-\Lambda}~\int_{L(x)}ds
~A(s)~s^{\Lambda}~B(s)~s^{-\Lambda}~f(s)~x^{\Lambda}  \right|    \right|
\leq C_1C_2 ~|x|^{1-\sigma_2}
$$
$$
\left|\left| x^{-\Lambda}~\int_{L(x)}ds
~s^{\Lambda}~B(s)~s^{-\Lambda}~A(s)~f(s)~x^{\Lambda}   \right|   \right|
\leq C_1C_2 ~|x|^{1-\sigma_2}
$$
where $L(x)$ is a path in $D(\epsilon;\sigma)$ joining $0$ to $x$. To
prove the estimates, we observe that  
$$
  ||x^{\Lambda}||=||x^{\hbox{diag}({\sigma\over 2},-{\sigma\over
  2})}||= \max \{ |x^{\sigma}|^{1\over 2},
  |x^{\sigma}|^{-{1\over 2}} \}\leq e^{{\theta_1\over 2}
\Im \sigma}~|x|^{-{\tilde{\sigma}\over 2}},
~~~~~~
\hbox{ in }~D(\epsilon;\sigma)
$$
Then 
$$
    || x^{\Lambda} ~A(x) ~x^{-\Lambda}||\leq ||
    x^{\Lambda}||~||A(x)||~ ||~x^{-\Lambda}||\leq 
             e^{\theta_1 \Im \sigma}~C_1 ~|x|^{-\tilde{\sigma}}
   $$
$$
  = \left(    e^{\theta_1 \Im \sigma}
  ~|x|^{\sigma_2-\tilde{\sigma}} \right)~C_1~|x|^{-\sigma_2}
$$
Thus, if $\epsilon$ is small enough (we require
$\epsilon^{\sigma_2-\tilde{\sigma}} \leq e^{-\theta_1 \Im \sigma}$) we
obtain $  ||x^{\Lambda}~A(x)~x^{-\Lambda}||\leq C_1 |x|^{-\sigma_2}$.

\begin{figure}
\epsfxsize=9cm
\centerline{\epsffile{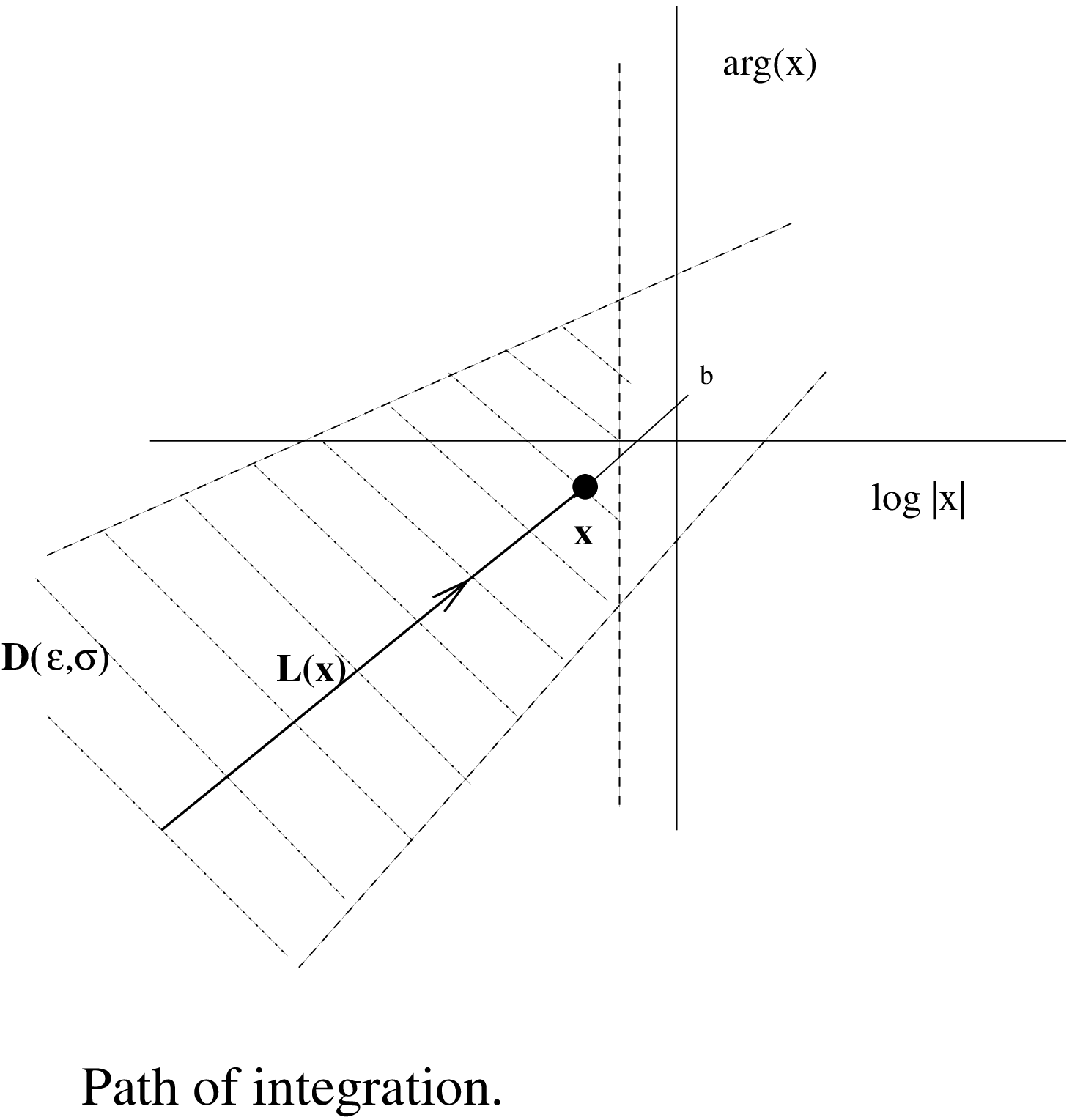}}
\caption{}
\label{figura4}
\end{figure}

We turn to the integrals. The integrands are holomorphic on
$D(\epsilon;\sigma)$, then they do not depend on the choice of the path, but
only on the point $x$. We choose a real number $\sigma^{*}$ such that $
 0\leq \sigma^{*}\leq \tilde{\sigma}$ and we choose the path 
$$
      \arg(s)=a \log|s| +b ,~~~~a={\Re \sigma - \sigma^{*}\over \Im
      \sigma} , ~~~\hbox{ or
      } \arg(s)=b\equiv \arg(x) ~\hbox{ if }~ \Im \sigma =0
$$
where $b$ is chosen appropriately such that $L(x)$ stays in
$D(\epsilon;\sigma)$. See figure \ref{figura4}. Then
we compute 
$$
\left|\left| x^{-\Lambda}~\int_{L(x)}ds
~A(s)~s^{\Lambda}~B(s)~s^{-\Lambda}~f(s)~x^{\Lambda}  \right|    \right|
= \left|\left| \int_{L(x)} ds ~ 
   x^{-\Lambda}~A(s)~x^{\Lambda} ~\left({s\over x}\right)^{\Lambda}  
~B(s)~\left( {s\over x}   \right)^{-\Lambda} ~f(s)~                
           \right|\right|
$$
$$
\leq 
     e^{2 \theta_1\Im\sigma} x^{-\tilde{\sigma}}
     ~C_1C_2~\max_{|x|<\epsilon} |f(x)| ~ \int_{L(x)}|ds| ~{|s|^{-\sigma^*} 
\over |x|^{-\sigma^*}}  
$$
where by $|ds|$ I mean $d\vartheta ~|ds(\vartheta)/d\vartheta|$, and 
$\vartheta$ is a parameter on the curve $L(x)$. 
The last step in the inequality  follows from 
$$
  \left| \left| \left({s\over x}  \right)^{\Lambda} \right|\right|
= 
 \left| \left|
\hbox{diag}({s^{\sigma\over 2}\over x^{\sigma\over
2}},{s^{-{\sigma\over 2}} \over
x^{-{\sigma\over 2}}}) \right|\right| =  
\max_L \left\{{|s^{\sigma\over 2}|\over |x^{\sigma\over 2}|},
{|s^{-{\sigma\over 2}}|\over
|x^{-{\sigma\over 2}}|}  
\right\}  
$$
and the observation  that, on $L$, $|s^{\sigma}|= 
|s|^{\bar{\sigma}}e^{-b\Im\sigma}$.  Thus 
$$   
\max \left\{{|s^{\sigma\over 2}|\over |x^{\sigma\over 2}|},
{|s^{-{\sigma\over 2}}|\over
|x^{-{\sigma\over 2}}|}  
\right\}=
\max \left\{  {|s|^{\sigma^*\over 2}\over |x|^{\sigma^*\over 2}},
{|s|^{-{\sigma^*\over 2}}\over |x|^{-{\sigma^*\over 2}}}  \right\}=
 {|s|^{-{ \sigma^*\over 2}}\over |x|^{-{\sigma^*\over 2}}}, 
~~~~~~|s|\leq |x|
$$
The parameter $s$ on $L(x)$ is 
$$ 
     s(\vartheta):=e^{\vartheta -b \over a} ~e^{i \vartheta } 
$$
where $\vartheta\in(-\infty,\arg(x)]$ or $\in [\arg(x),+\infty)$.
Then 
$$ 
   \int_{L(x)}|ds| ~ |s|^{-\sigma^*}= \int_{L(x)} d \vartheta
   ~\left| {ds \over d\vartheta}
\right|~|s(\vartheta)|^{-\sigma^*}
= {|a|\over 1-\sigma^*}
\left|i+{1\over a}   \right| 
   e^{{\arg(x)-b \over a} (1-\sigma^*)}= {|a|\over 1-\sigma^*}
\left|i+{1\over a}   \right| ~|x|^{1-\sigma^*}
$$
Then, the initial integral is less or equal to 
$$
     e^{ \theta_1\Im \sigma}~\max_{|x|<\epsilon}|f(x)| 
~C_1C_2~ \hbox{constant}~|x|^{1-\tilde{\sigma}}
$$
Now, we write   
$|x|^{1-\tilde{\sigma}}=|x|^{\sigma_2-\tilde{\sigma}}~
 |x|^{1-\sigma_2}$ and we obtain, for sufficiently small $\epsilon$: 
$$
     e^{ \theta_1 \Im \sigma}~\max_{|x|<\epsilon}|f(x)| 
~C_1C_2~ \hbox{constant}~|x|^{1-\tilde{\sigma}}
\leq C_1 C_2 ~|x|^{1-\sigma_2}
$$
We remark that for $\sigma=0$ the above estimates are still valid. Actually
$
||x^{\Lambda}||\equiv||x^{\pmatrix{0&1 \cr 0 & 0\cr}}||$ diverges like
$|\log x|$, $|| x^{\Lambda} A(x) x^{-\Lambda}||$ are less or equal to 
$ C_1~|\log(x)|^2$, and finally  $|| x^{-\Lambda}~\int_{L(x)}ds
~A(s) $ $ s^{\Lambda}~B(s)~s^{-\Lambda}~f(s)~x^{\Lambda}  ||
$ is less or  equal to $C_1 C_2 \max|f| ~|\log(x)|^2 \int_{L(x)} |ds|
~|\log s|^2$.  We can choose $L(x)$ to be a radial path $s=\rho \exp(i
\alpha)$, $0<\rho <|x|$, $\alpha$ fixed. Then the integral is $|x|(
\log|x|^2- 2 \log|x| +2 + \alpha^2)$. The factor $|x|$ does the job,
because we rewrite it as $|x|^{\sigma_2}~|x|^{1-\sigma_2}$ (here
$\sigma_2$ is any number between 0 and 1) and we
proceed as above to choose $\epsilon$ small enough in such a way
that $\bigl(\max|f|~ |x|^{\sigma_2} \times$ function diverging like
$\log^2|x|\bigr) \leq 1$. 

The estimates above are in a sense enough to prove the
 lemma. 

We solve the Schlesinger equations by successive
 approximations, as in \cite{SMJ}: let
 $\tilde{B}_{\nu}(x):=x^{-\Lambda}B_{\nu}(x)x^{\Lambda} $. The
 Schlesinger equations are re-written as 
$$ 
 {d A_{\mu}\over dx}=\sum_{\nu=1}^{n_2}
 [A_{\mu},x^{\Lambda}\tilde{B}_{\nu} x^{-\Lambda}]~f_{\mu \nu}(x)
$$
$$
{d\tilde{B}_{\nu}\over dx}= {1\over x} [\tilde{B}_{\nu},
\sum_{\mu}x^{-\Lambda}( A_{\mu}(x)-A_{\mu}^0)x^{\Lambda} ]+
\sum_{\mu=1}^{n_1}[\tilde{B}_{\nu},x^{-\Lambda}A_{\mu}x^{\Lambda} 
]~g_{\mu \nu}(x) +     
\sum_{\nu^{\prime}=1}^{n_2} [\tilde{B}_{\nu},\tilde{B}_{\nu^{\prime}}] ~h_{\nu
\nu^{\prime} }(x) 
$$
Then, by successive approximations:
$$
 A_{\mu}^{(k)}(x)= A_{\mu}^0+ \int_{L(x)}ds~\sum_{\nu}
 [A_{\mu}^{(k-1)}(s),s^{\Lambda} \tilde{B}_{\nu}^{(k-1)}(s)
 s^{-\Lambda}]~f_{\mu \nu}(s) 
$$
$$
\tilde{B}_{\nu}^{(k)}(x)= B_{\nu}^0+ \int_{L(x)}ds~ \left\{ 
         {1\over s} [\tilde{B}_{\nu}^{(k-1)}(s),\sum_{\mu}
         s^{-\Lambda} ( A_{\mu}^{(k-1)}(s)-A_{\mu}^0)s^{\Lambda}] +
         \right. $$
$$
\left. +
\sum_{\mu}[\tilde{B}_{\nu}^{(k-1)}(s),s^{-\Lambda}A_{\mu}^{(k-1)}(s)
x^{\Lambda} ]~g_{\mu
\nu}(s)+\sum_{\nu^{\prime}}[\tilde{B}_{\nu}^{(k-1)}(s),
 \tilde{B}_{\nu^{\prime}}^{(k-1)}(s)]~h_{\nu\nu^{\prime}}
\right\}
$$
The functions $A_{\mu}^{(k)}(x)$, $B_{\nu}^{(k)}(x)$ are holomorphic in
$D(\epsilon;\sigma)$, by construction. 
Observe that $||A_{\mu}^0||\leq C$, $||B_{\nu}^0||\leq C$ for some
constant $C$.   We claim that for $|x|$ sufficiently small 
\be
\left.\matrix{
||A_{\mu}^{(k)}(x)-A_{\mu}^0||\leq C|x|^{1-\sigma_1}
\cr\cr
\left|\left| x^{-\Lambda}\left(A_{\mu}^{(k)}(x)-A_{\mu}^0 \right)x^{\Lambda} 
 \right|\right|\leq C^2 |x|^{1-\sigma_2}
\cr\cr
|| \tilde{B}_{\nu}^{(k)}(x)-B_{\nu}^0||\leq C|x|^{1-\sigma_1}
\cr
}\right.
\label{succ1}
\ee
where  $ \tilde{\sigma}<\sigma_2<\sigma_1<1$. 
Note that the above inequalities imply $||A_{\mu}^{(k)}||\leq 2C$,
$||\tilde{B}_{\nu}^{(k)}||\leq 2C$. Moreover we claim that 
\be
\left. \matrix{
||A_{\mu}^{(k)}(x)-A_{\mu}^{(k-1)}(x)||\leq C~\delta^{k-1}~|x|^{1-\sigma_1}
\cr\cr
\left|\left| x^{-\Lambda}\left(A_{\mu}^{(k)}(x)-A_{\mu}^{(k-1)}(x)
 \right )x^{\Lambda} 
 \right|\right|\leq C^2~\delta^{k-1}~ |x|^{1-\sigma_2}
\cr\cr
|| \tilde{B}_{\nu}^{(k)}(x)-B_{\nu}^{(k-1)}(x)||\leq
C~\delta^{k-1}|x|^{1-\sigma_1}
\cr
}
\right.
\label{succ2}
\ee
where $0<\delta<1$. 

For $k=1$ the above inequalities are proved using the simple methods
used in the estimates at the beginning of the proof. Then we proceed
by induction, still using the same estimates. We leave this technical
point to the reader, but we give at least one example of how to
proceed. As an example, we prove the $(k+1)^{th}$ step of the first of
(\ref{succ2}):
$$
||A_{\mu}^{(k+1)}(x)-A_{\mu}^{(k)}(x)||\leq C~\delta^k~|t|^{1-\sigma_1}
$$
supposing the  $k^{th}$ step of (\ref{succ2}) is true.  
Let us proceed using the integral equations:
$$
||A_{\mu}^{(k+1)}(x)-A_{\mu}^{(k)}(x)||=
\left|\left| \int_{L(x)}ds~\sum_{\nu=1}^{n_2} \left( 
A_{\mu}^{(k)}s^{\Lambda}\tilde{B}_{\nu}^{(k)}s^{-\Lambda} -
A_{\mu}^{(k-1)} s^{\Lambda} \tilde{B}_{\nu}^{(k-1)} s^{-\Lambda}
+\right.\right. \right. $$
$$
\left.\left.\left. +
s^{\Lambda} \tilde{B}_{\nu}^{(k-1)} s^{-\Lambda}
A_{\mu}^{(k-1)}-s^{\Lambda} \tilde{B}_{\nu}^{(k)} s^{-\Lambda} A_{\mu}^{(k)}
\right)
         ~f_{\mu\nu}(s)
  \right|\right| \leq
$$
 $$
\leq 
 \int_{L(x)}d|s|~\sum_{\nu=1}^{n_2} \left|\left|
A_{\mu}^{(k)}s^{\Lambda}\tilde{B}_{\nu}^{(k)}s^{-\Lambda} -
A_{\mu}^{(k-1)} s^{\Lambda} \tilde{B}_{\nu}^{(k-1)} s^{-\Lambda}
\right|\right|~|f_{\mu\nu}(s)|+  $$
$$ 
+\int_{L(x)}d|s|~\sum_{\nu=1}^{n_2}
\left| \left|
s^{\Lambda} \tilde{B}_{\nu}^{(k-1)} s^{-\Lambda}
A_{\mu}^{(k-1)}-s^{\Lambda} \tilde{B}_{\nu}^{(k)} s^{-\Lambda} A_{\mu}^{(k)}
  \right|\right|~|f_{\mu\nu}(s)|
$$
Now we estimate 
$$
 \left|\left|
A_{\mu}^{(k)}s^{\Lambda}\tilde{B}_{\nu}^{(k)}s^{-\Lambda} -
A_{\mu}^{(k-1)} s^{\Lambda} \tilde{B}_{\nu}^{(k-1)} s^{-\Lambda}
\right|\right|\leq
$$
$$
\leq \left|\left|  
A_{\mu}^{(k)}s^{\Lambda}\tilde{B}_{\nu}^{(k)}s^{-\Lambda}-A_{\mu}^{(k-1)}s^{
\Lambda} \tilde{B}_{\nu}^{(k)} s^{-\Lambda} \right|\right| +
$$
$$ +\left| \left|
A_{\mu}^{(k-1)}s^{
\Lambda} \tilde{B}_{\nu}^{(k)} s^{-\Lambda} -
A_{\mu}^{(k-1)} s^{\Lambda} \tilde{B}_{\nu}^{(k-1)} s^{-\Lambda}
\right|\right|
$$
$$
 \leq || A_{\mu}^{(k)}-A_{\mu}^{(k-1)}||
 ~||s^{\Lambda}\tilde{B}_{\nu}^{(k)} s^{-\Lambda}||+
 ||A_{\mu}^{(k-1)}|| ~ || s^{\Lambda}||~||\tilde{B}_{\nu}^{(k)} 
-\tilde{B}_{\nu}^{(k-1)} || ~ ||s^{-\Lambda}||
$$
     By induction then:
$$
  \leq 
         \left( C~\delta^{k-1}~|s|^{1-\sigma_1} \right) 
~2C~e^{\theta_1\Im \sigma}|s|^{-\tilde{\sigma}} + 
      2 C~\left(C~\delta^{k-1}~|s|^{1-\sigma_1}\right) ~e^{\theta_1\Im \sigma}
|s|^{-\tilde{\sigma}}
$$ 
The other term is estimated in an analogous way. Then 
$$ 
  ||A_{\mu}^{(k+1)}-A_{\mu}^{(k)}||\leq \hbox{constant}~
8n_2C^2 \max |f_{\mu \nu}| ~\delta^{k-1}~e^{\theta_1\Im \sigma}
  |x|^{1-\tilde{ \sigma}}~|x|^{1-\sigma_1}
$$ 
We choose $\epsilon$  small enough to have $ \hbox{
constant}~ 8n_2C\max|f|
~e^{\theta_1\Im \sigma} |x|^{1-\tilde{\sigma}}\leq \delta$. Note that the choice of
$\epsilon$  is independent of
$k$.  In the case $\sigma=0$,  $|x|^{1-\tilde{\sigma}}$ is substituted by
$|x| (\log^2|x| +O(\log|x|))$. 

\vskip 0.2 cm

 The inequalities  (\ref{succ1}) (\ref{succ2}) 
 imply the convergence of the successive
 approximations to a solution of the Schlesinger equations,
 satisfying the  assertion of the lemma, plus the
 additional inequality
$$ ||   x^{-\Lambda}(A_{\mu}(x)-A_{\mu}^0)x^{\Lambda}||\leq C^2
|x|^{1-\sigma_2} $$

 \rightline{$\Box$}
\vskip 0.3 cm

 We observe that we imposed
\be
     \epsilon^{\sigma_2-\tilde{\sigma}} \leq c e^{-\theta_1\Im \sigma },~~~
   \epsilon^{1-\tilde{\sigma}} \leq c e^{-\theta_1\Im \sigma },
\label{condizioni su epsilon}
\ee
where $c$ is a constant constructed in the theorem ( $c$ is proportional to 
${1-\tilde{\sigma}\over 8 n_2 C}$ and 
$C=\max\{||A_{\mu}^0||,||B_{\nu}^0||\}$ ). 

\vskip 0.3 cm  

We turn to the  case in which we are concerned: we consider three
matrices $A_0^0,A_x^0,A_1^0$ such that 
$$ 
   A_0^0+A_x^0=\Lambda,~~~A_0^0+A_x^0+A_1^0=\hbox{diag}(-\mu,\mu)$$
$$ \hbox{tr}(A_i^0)=\det(A_i^0)=0, ~~~i=0,x,1
$$
\vskip 0.3 cm
\noindent
{\bf Lemma 2: } {\it Let $r$ and $s$ be two complex numbers not equal
  to 0 and
  $\infty$. 
  Let $T$ be the matrix which brings  $\Lambda$ to the  Jordan form:
$$
   T^{-1} \Lambda T = \left\{ \matrix{
                 \hbox{ diag}({\sigma \over 2}, -{\sigma \over2})
                                   ,~~~~\sigma \neq 0 \cr\cr
                         \pmatrix{0&1\cr
                                  0&0\cr},~~~~\sigma =0 \cr
}\right.
                     $$
The general solution of} 
$$A_0^0+A_x^1+A_1^0=\pmatrix{-\mu&0\cr 
                             0 & \mu \cr},~~~~\hbox{
                             tr}(A_i)=\det(A_i)=0,~~~~
A_0^0+A_x^0=\Lambda$$
{\it 
is the following:
\vskip 0.15 cm
For $\sigma\neq 0, \pm 2\mu$:
$$
\Lambda = {1\over 8 \mu} \pmatrix{
                                   -\sigma^2-(2\mu)^2 &
                                            (\sigma^2-(2\mu)^2)r \cr 
                            {(2\mu)^2-\sigma^2 \over r} &
                                   \sigma^2+(2\mu)^2 \cr }
~~~~A_1^0={\sigma^2 -(2\mu)^2\over 8\mu}\pmatrix{ 1  &  -r \cr
                                                  {1\over r} & -1 \cr}
$$
$$
  A_0^0=T \pmatrix{ {\sigma\over 4} & {\sigma \over 4} ~s \cr
                        -{\sigma\over 4}{1\over s} &  -{\sigma\over
                        4} \cr }~T^{-1},
~~~~A_x^0=T \pmatrix{ {\sigma\over 4} & -{\sigma \over 4} ~s \cr
                        {\sigma\over 4}{1\over s} &  -{\sigma\over
                        4} \cr }~T^{-1}
$$
where 
$$ T=\pmatrix{ 1 & 1 \cr
                {(\sigma + 2\mu)^2\over \sigma^2-(2\mu)^2} {1\over r}
                &
                  {(\sigma - 2\mu)^2\over \sigma^2-(2\mu)^2} {1\over
                r}\cr 
}$$
\vskip 0.15 cm 
For $\sigma=-2\mu$: $A_0^0$ and $A_x^0$ as above, but 
\be 
     \Lambda=\pmatrix{-\mu & r \cr 
                       0 & \mu \cr}~~~~~
    A_1^0=\pmatrix{0 & -r \cr
                   0 & 0 \cr }
       ~~~~~
   T = \pmatrix{1 & 1\cr 
                0 & {2 \mu \over r} \cr
           }
\label{cucu}
\ee
or 
\be
       \Lambda=\pmatrix{-\mu & 0 \cr 
                       r & \mu \cr}~~~~~
    A_1^0=\pmatrix{0 & 0 \cr
                   -r & 0 \cr }
       ~~~~~
   T = \pmatrix{1 & 0\cr 
                 -{r \over 2\mu }&1 \cr
           }
\label{cucu2}
\ee           
\vskip 0.15 cm
For $\sigma=2\mu$: $A_0^0$ and $A_x^0$ as above, but 
\be
     \Lambda=\pmatrix{-\mu & r \cr 
                       0 & \mu \cr}~~~~~
    A_1^0=\pmatrix{0 & -r \cr
                   0 & 0 \cr }
       ~~~~~
T= \pmatrix{1 & 1\cr 
                 { 2\mu\over r }&0 \cr
           }
\label{cucu1}
\ee
or 
\be
 \Lambda=\pmatrix{-\mu & 0 \cr 
                       r & \mu \cr}~~~~~
    A_1^0=\pmatrix{0 & 0 \cr
                   -r & 0 \cr }
       ~~~~~
 T = \pmatrix{0 & 1\cr 
                1 & -{r\over 2\mu} \cr
           }
\label{cucu3}
\ee
\vskip 0.15 cm 
For $\sigma=0$:
$$
   A_0^0=T \pmatrix{ 0 & s \cr 
                     0 & 0 \cr }T^{-1}
~~~~~A_x^0=T \pmatrix{0 & 1-s \cr 0 & 0 \cr }T^{-1}
$$
$$
\Lambda= \pmatrix{ -{\mu \over 2} & -{\mu^2 \over 4} r \cr
                   {1\over r}  & {\mu \over 2} \cr
}
~~~~~~A_1^0=\pmatrix{ -{\mu \over 2} & {\mu^2 \over 4} r \cr
                       -{1\over r} & {\mu \over 2} \cr
}
$$
$$
T= \pmatrix{1 & 1 \cr
             -{2\over \mu r} & -2{\mu+2\over \mu^2} {1\over r}\cr}
$$
}

\vskip 0.2 cm
\noindent
We leave the proof as an exercise for the reader. $\Box$

\vskip 0.3 cm 

We are ready to prove theorem 1:

\vskip 0.3 cm 
\noindent 
{\bf Theorem 1:} { \it The solutions of $PVI_{\mu}$, corresponding to
the solutions of Schlesinger  equations (\ref{sch}) 
obtained in lemma 1, have the
following behaviour for ${x}\to 0$ along a path $\Im \sigma
\arg(x) = (\Re \sigma - \Sigma) \log|x| +b \Im \sigma$, $0\leq \Sigma
\leq \tilde{\sigma}$, contained in  
 $D(\epsilon;\sigma,\theta_1,\theta_2)$: 
$$
   y(x)= a(x)~{x}^{1-\sigma} (1+O(|x|^{\delta}))
,~~~~~~~~~~~~~~ \sigma\neq 0$$
$$
y({x})=s {x}~(1+O(|x|^{\delta})),~~~~~~~~~~~~~~\sigma=0
$$
where $0<\delta<1$ is a small number, and $a(x)$ can be computed as a
function of $s$. Namely 
$$ a(x) = -{ 1 \over 4 s}$$ 
along any path, except
for the paths $\Im \sigma \arg(x) = \Re \sigma \log|x| + b \Im \sigma $, 
along which  ${{x}}^{\sigma}= C e^{i\alpha(x)}$ ( $C$ is a constant
=$ |x^{\sigma}| $ and   $\alpha(x)$ is the 
(real)  phase). In this case 
$$a(x)= -{1\over 4s}~\left(1 -2~s~ C e^{i \alpha(x)} 
             + s^2~ C^2 e^{2 i \alpha(x)}  \right) = O(1), \hbox{ for }
             x\to 0
$$
}

\vskip 0.2 cm 
\noindent
{\it Proof:}  
$y(x)$ can be computed in terms of the $A_i(x)$ from $A(y(x),x)_{12}=0$:
$$ 
  y(x)={x (A_0)_{12}\over (1+x) (A_0)_{12}+(A_x)_{12}
 +x(A_1)_{12}}\equiv {x (A_0)_{12}\over x(A_0)_{12} -(A_1)_{12} + x (A_1)_{12}}
$$ 
$$
   = -x {(A_0)_{12} \over (A_1)_{12}} ~{1\over 1-
   x(1+{(A_0)_{12}\over (A_1)_{12}}    )}
$$
As a consequence of lemmas 1 and 2 it follows that $|x~(A_1)_{12}|\leq
c ~|x|~(1+O(|x|^{1-\sigma_1}))$ and $|x~(A_0)_{12}|\leq
c ~|x|^{1-\tilde{\sigma}}~(1+O(|x|^{1-\sigma_1}))$, where $c$ is a
constant.  Then 
$$
y(x) =  -x {(A_0)_{12} \over (A_1)_{12}}~(1
+O(|x|^{1-\tilde{\sigma}}))
$$
From lemma 2 we find, for $\sigma \neq 0, \pm 2\mu$: 
$$
(A_0)_{12}= -r {\sigma^2 - 4 \mu^2 \over 32 \mu} \left[ 
   {x^{-\sigma} \over s}
   (1+O(|x|^{1-\sigma_1}))+s~x^{\sigma}~(1+O(|x|^{1-\sigma_1})) -2
   (1+O(|x|^{1-\sigma_1}))  \right]
$$
$$ 
  (A_1)_{12}=-r {\sigma^2 - 4 \mu^2 \over 8 \mu}
  (1 +O(|x|^{1-\sigma_1}))
$$
Then  (recall that $\tilde{\sigma}<\sigma_1$)
$$
   y(x) = -{x\over 4} \left[    {x^{-\sigma} \over s}
   (1+O(|x|^{1-\sigma_1}))+s~x^{\sigma}~(1+O(|x|^{1-\sigma_1})) -2
   (1+O(|x|^{1-\sigma_1}))\right]~\bigl( 1+O(|x|^{1-\sigma_1}) \bigr)
                             $$
Now $x\to 0$ along a path 
$$
\Im \sigma \arg(x)= (\Re \sigma -\Sigma) \log|x|+b \Im \sigma$$ 
for a suitable $b$ and $0\leq \Sigma\leq\tilde{\sigma}$. Along this path
we rewrite  $x^{\sigma}$ in terms of its absolute value $|x^{\sigma}|=C
|x|^{\Sigma} $ ($C=e^{-b\Im\sigma}$) and its real phase $\alpha(x)$ 
$$
  x^{\sigma}= C~|x|^{\Sigma} ~e^{i \alpha(x)},~~~~\alpha(x)=
\Re\sigma\arg(x)+\Im\sigma \ln |x|\Bigl|_{
\Im \sigma \arg(x)= (\Re \sigma -\Sigma) \log|x|+b \Im \sigma}.$$
Then
$$
 y(x)  =  -{x^{1-\sigma} \over 4} \left[ 
                                     {1\over s} -2 C e^{i \alpha(x)}
                                     |x|^{\Sigma}(1+O(|x|^{1-\sigma_1})) + 
s C^2 e^{2 i
                                     \alpha(x)} ~|x|^{2\Sigma}
                                     ~(1+O(|x|^{1-\sigma_1}) )
                                                  \right]~\left( 
                     1 + O(|x|^{1-\sigma_1} )\right)
$$
For $\Sigma \neq 0$ the above expression becomes
$$ 
  y(x)= -{1 \over 4s} x^{1-\sigma} ~\left( 1 + O(|x|^{1-\sigma_1}) + 
                       O(|x|^{\Sigma})
                       \right)
$$
We collect the two  $O(..)$ contribution in $O(|x|^{\delta})$ where
$\delta =   \min
\{ 1-\sigma_1, \Sigma \}$ is a small number between 0 and 1. We take
the occasion here to remark that in the case of real $0<\sigma<1$,
 if we consider $x\to 0$ along a
radial path (i.e. $\arg(x) = b $), then $\Sigma = \tilde{\sigma}
= \sigma$  and  thus:
$$   
  y(x) = \left\{ \matrix{ { 1 \over 4s} x^{1-\sigma} (1 +O(|x|^{\sigma})) \hbox{ for }
                                               0<\sigma <{1 \over 2}
                                               \cr
 \cr
 { 1 \over 4s} x^{1-\sigma} (1 +O(|x|^{1-\sigma_1})) \hbox{ for }
                                               {1\over 2}<\sigma <1 
                                               \cr
} \right.
 $$
Finally, along the path with $\Sigma=0$ we have:
$$
 y(x) = -{x^{1-\sigma} \over 4}~\left({1\over s} -2 C e^{i \alpha(x)} 
             + s~ C^2 e^{2 i \alpha(x)}  \right)~\left(1 +
             O(|x|^{1-\sigma_1})  \right) 
$$

We let the reader verify that also in the cases $\sigma =\pm 2 \mu$
the behaviour of $y(x)$ is as above (use the matrices (\ref{cucu})and
(\ref{cucu1}) -- the reason why we disregard the matrices
(\ref{cucu2}), (\ref{cucu3}) will be
clarified at the end of the proof of theorem 2) and that  for $\sigma=0$ 
$$ 
  y(x) = s ~x~ \bigl(1 +O(|x|^{1-\sigma_1}) \bigr)
$$ 
For $\sigma=0$, we recall that  $0<\sigma_1<1$ 
is arbitrarily small. 
 $\Box$

\vskip 0.2 cm 
In  the proof of lemma 1 we imposed (\ref{condizioni su epsilon}). Hence, the  
reader may observe that $\epsilon$ depends on  $\tilde{\sigma}$, $\theta_1$
and on  $||A_0^0||$, $||A_x^0||$, $||A_1^0||$; thus it depends also 
 on  $s$ ($\Rightarrow$ on $a$).    The second inequality  in
 (\ref{condizioni su epsilon})
has been used in section \ref{Local Behaviour -- Theorem 1}
 to construct the domain (\ref{come se non
bastasse!!}).  


\section{ Proof of theorem 2}\label{dimth2}

We are interested in lemma 1 when  
$$
  f_{\mu \nu}=g_{\mu \nu}={b_{\nu}\over a_{\mu}
  -xb_{\nu}},~~~~h_{\mu\nu}=0
$$
$$
  a_{\mu},b_{\nu}\in{\bf C},~~~~~a_{\mu}\neq 0~~~ \forall \mu=1,...,n_1
$$
Equations (\ref{sch1}) are the isomonodromy deformation equations for the
fuchsian system 
$$
   {dY\over dz} = \left[\sum_{\mu=1}^{n_1}{A_{\mu}(\tilde{x})\over
   z-a_{\mu}}+\sum_{\nu=1}^{n_2} {B_{\nu}(\tilde{x})\over z- xb_{\nu}}
   \right] ~Y
$$
As a corollary of lemma 1, for a fundamental matrix solution $Y(z,x)$
of the fuchsian system  the limits 
$$
   \hat{Y}(z):=\lim_{{x}\to 0} Y(z,{x})
$$
$$
 \tilde{Y}(z):=\lim_{{x}\to 0} ~{x}^{-\Lambda} 
Y({x}z,{x})
$$
exist when  ${x}\to 0$ in $D(\epsilon;\sigma)$. They satisfy 
$$
  {d \hat{Y}\over dz} = \left[\sum_{\mu=1}^{n_1}{A_{\mu}^0\over
   z-a_{\mu}}+{\Lambda \over z}
   \right] ~\hat{Y}
$$
$$
  {d\tilde{Y}\over dz} = \sum_{\nu=1}^{n_2} {B_{\nu}(x)\over z- b_{\nu}}
    ~\tilde{Y}
$$
In  our case, the last three  systems reduce to 
\be
  {dY\over dz} = \left[{A_0({x})\over z} +{A_x({x})\over z-x}+
{A_1({x})\over z-1}
   \right] ~Y
\label{stofucs}
\ee
\be
{d \hat{Y}\over dz} = \left[{A_1^0 \over z-1}+{\Lambda \over z}
   \right] ~\hat{Y}
\label{systemhat}
\ee
\be
  {d\tilde{Y}\over dz} = \left[{A_0^0\over z}+{A_x^0\over z-1}   \right]
    ~\tilde{Y}
\label{systemtilde}
\ee
Before taking the limit ${x}\to 0$, let us choose
\be
Y(z,{x})= \left(I+O\left({1\over z}\right)\right) 
~z^{-A_{\infty}}z^R ,~~~~~~z\to \infty
\label{solutionzx}
\ee 
and define as above
$$ 
   \hat{Y}(z):=\lim_{{x}\to 0} Y(z,{x}), 
~~~ \tilde{Y}(z):=\lim_{{x}\to 0} ~\tilde{x}^{-\Lambda} 
Y({x}z,{x})
$$
For the system (\ref{systemhat}) we choose a fundamental
matrix solution   normalized as follows
\be
   \hat{Y}_N(z)= \left(I+O\left({1\over z}\right)\right) ~z^{-A_{\infty}}z^R,
~~~~z\to \infty
\label{solutionhat}
\ee
$$
  =   (I+O(z))~z^{\Lambda}~\hat{C}_0,~~~~z\to 0 
$$
$$
   = \hat{G}_1 (I+O(z-1))~(z-1)^J ~\hat{C}_1,~~~~z \to 1 
$$
Where $\hat{G}_1^{-1} A_1^0 \hat{G}_1=J$, $J=\pmatrix{0&1\cr 0&0\cr}$. 
$\hat{C}_0$,  $\hat{C}_1$ are {\it connection matrices}. Note that $R$
is the same of (\ref{solutionzx}), since it is independent of  $x$. For  
(\ref{systemtilde}) we choose a fundamental matrix solution normalized
as follows
 \be
   \tilde{Y}_N(z)= \left(I+O\left({1\over z}\right)\right) ~z^{\Lambda},
~~~~z\to \infty
\label{solutiontilde}
\ee
$$
  =  \tilde{G}_0 (I+O(z))~z^J~\tilde{C}_0,~~~~z\to 0 
$$
$$
   = \tilde{G}_1 (I+O(z-1))~(z-1)^J ~\tilde{C}_1,~~~~z \to 1. 
$$
Here $\tilde{G}_0^{-1}A_0^0 \tilde{G}_0 =J$, $\tilde{G}_1^{-1}A_x^0
\tilde{G}_1 =J$.

Now we prove that 
$$
    \hat{Y}(z) = \hat{Y}_N(z)
$$ 
\be     
 \tilde{Y}(z)=
 \tilde{Y}_N(z)~\hat{C}_0
 \label{lim}
\ee
The proof we give here uses the  technique of the proof of
 Proposition 2.1. in \cite{Jimbo}. The (isomonodromic) dependence of
 $Y(z,x)$ on $x$ is given by 
$$
  {dY(z,x)\over dx}= - {A_x(x) \over z-x} Y(z,x):=F(z,x) Y(z,x)
$$
Then 
$$ 
  Y(z,x)= \hat{Y}(z) + \int_0^x dx_1~ F(z,x_1) Y(z,x_1)
$$
The integration is on a path $\arg(x)=a \log|x| +b$,
$a={\Re\sigma-\sigma^* \over \Im \sigma}$ ($0\leq \sigma^*\leq
\tilde{\sigma}$),  or $\arg(x)=0$ if $\Im
\sigma=0$. The path is contained 
in $D(\sigma)$ and  joins $0$ and $x$,  like $L(x)$ 
in the proof of theorem 1 (figure 10). By successive approximations we have:
$$ 
   Y^{(1)}(z,x)= \hat{Y}(z) + \int_0^x dx_1 F(z,x_1) \hat{Y}(z)
$$
$$
Y^{(2)}(z,x)= \hat{Y}(z) + \int_0^x dx_1 F(z,x_1) Y^{(1)}(z,x_1)
$$
$$
 \vdots
$$
$$
Y^{(n)}(z,x)= \hat{Y}(z) + \int_0^x dx_1 F(z,x_1) Y^{(n-1)}(z,x_1)
$$
$$
   = 
    \left[ I + \int_0^x dx_1 \int_0^{x_1} dx_2...\int_0^{x_{n-1}} dx_n
        F(z,x_1)F(z,x_2)...F(z,x_n)\right]~\hat{Y}(z)
$$
 Performing integration like in the proof of theorem 1 we evaluate 
$
  || Y^{(n)}(z,x) - Y^{(n-1)}(z,x)||$. Recall that $\hat{Y}(z)$ has 
singularities at $z=0$, $z=x$. Thus, if $|z|>|x|$ we obtain 
 $$
|| Y^{(n)}(z,x) - Y^{(n-1)}(z,x)||\leq {M C^n \over \Pi_{m=1}^n
(m-sigma^*)} |x|^{n-\sigma^*},$$
where $M$ and $C$ are constants. 
Then $Y^{(n)}= \hat{Y}+ (Y^{(1)}-\hat{Y})+...+
(Y^{(n)}-Y^{(n-1)})$ converges for $n\to \infty$ 
 uniformly in $z$ in every compact set contained in $\{z~|~|z|>|x|\}$
and uniformly in $x\in D(\sigma)$. We can exchange limit and
integration, thus obtaining $
Y(z,x)= \lim_{n\to \infty} Y^{(n)}(z,x)$. Namely
$$
   Y(z,x)= U(z,x) \hat{Y}(z),$$
$$
   U(z,x)= I + \sum_{n=1}^{\infty} \int_0^x dx_1 \int_0^{x_1} dx_2...
         \int_0^{x_{n-1}} dx_n F(z,x_1)F(z,x_2)...F(z,x_n)
$$
being the convergence of the series uniformly in  $x\in D(\sigma)$
and in $z$ in every compact set contained in $\{z~|~|z|>|x|\}$. Of
course 
$$ 
   U(z,x)= I + O\left({1\over z}\right) 
 \hbox{ for } x\to 0 \hbox{ and } Y(z,x) \to \hat{Y}(z)
$$
But now observe that 
$$ 
  \hat{Y}(z) = U(z,x)^{-1} Y(z,x) = \left(I+O\left({1\over z}\right)
  \right) ~\left(I+O\left({1\over z}\right) \right) z^{-A_{\infty}}
  z^R, ~~~~ z\to \infty
$$
Then 
    $$ \hat{Y}(z) \equiv \hat{Y}_N(z)
$$
Finally, for $z\to 1$, 
$$
  Y(z,x) = U(x,z) \hat{Y}_N(z)= U(x,z) ~\hat{G}_1 (I +O(z-1)) (z-1)^J
  \hat{C}_1
$$
$$
  = G_1(x)(I+O(z-1))(z-1)^J \hat{C}_1$$
This implies 
             $$   C_1\equiv \hat{C}_1$$
and then
\be 
 M_1 =
\hat{C}_1^{-1} e^{2\pi i J} \hat{C}_1
\label{MONODROMY1}
\ee   
 Here we have  chosen a monodromy representation for
(\ref{stofucs}) by fixing  a base-point and a basis in the fundamental
group of 
${\bf P}^1$ as in figure \ref{figura567}. 
 $M_0$, $M_1$, $M_x$, $M_{\infty}$ are 
the monodromy matrices for the solution (\ref{solutionzx}) 
corresponding to the loops $\gamma_i$
$i=0,x,1,\infty$.  $M_{\infty}M_1M_xM_0=I$. The result
 (\ref{MONODROMY1}) may also be proved simply observing that $M_1$ becomes
$\hat{M}_1$ as $x\to 0$ in
$D(\sigma)$ because  the system (\ref{systemhat}) 
is obtained from (\ref{stofucs}) when
$z=x$ and $z=0$ merge and the singular point $z=1$ does not move. 
$x$ may converge to $0$ along spiral paths (figure \ref{figura567}). 
We recall that the braid
$\beta_{i,i+1}$ changes the monodromy matrices of ${dY\over dz} = \sum_{i=1}^n
{A_i(u)\over z-u_i} Y$ according to $M_i\mapsto M_{i+1}$, $M_{i+1}\mapsto
M_{i+1}M_i M_{i+1}^{-1}$, $M_k \mapsto M_k$ for any $k\neq i, i+1$; therefore, 
if $\arg (x)$ increases of $2\pi$ as $x\to 0$ in (\ref{stofucs}) we have 
$$M_0 \mapsto M_x,~~~~M_x\mapsto M_x M_0 M_x^{-1},~~~M_1\mapsto M_1$$ 
If follows that $M_1$ does not change and then   
 \be 
          M_1\equiv\hat{M}_1=\hat{C}_1^{-1} e^{2\pi i J} \hat{C}_1
\label{mono1}
\ee
where $\hat{M}_1$ is the monodromy matrix of   (\ref{solutionhat}) for
the loop $\hat{\gamma}_1$ in
the basis of figure \ref{figura567}.

\vskip 0.2 cm 

Now we turn  to $\tilde{Y}(z)$. Let $\tilde{Y}(z,x):= x^{-\Lambda}
Y(xz,x)$, and by definition $\tilde{Y}(z,x) \to \tilde{Y}(z)$ as $x\to
0$. In this case 
$$ 
  {d \tilde{Y}(z,x) \over dx} = \left[{ x^{-\Lambda} (A_0+A_x)
  x^{\Lambda} - \Lambda \over x} + { x^{-\Lambda } A_1 x^{\Lambda}
  \over x - {1\over z}}\right] \tilde{Y}(z,x) := \tilde{F}(z,x)
  \tilde{Y}(z,x)
$$ 
Proceeding by successive approximations as above we get 
$$ 
   \tilde{Y}(z,x) = V(z,x) \tilde{Y}(z),$$
$$ V(z,x)=
I+\sum_{n=1}^{\infty} \int_0^x dx_1...\int_0^{x_{n-1}} dx_n
\tilde{F}(z,x_1)... \tilde{F}(z,x_n) \to I \hbox{ for } x\to 0 
$$ 
uniformly in $x\in D(\sigma)$ and in $z$ in every compact subset of
$\{ z~|~ |z|<{1\over |x|} \}$. 

Let's investigate the behaviour of $\tilde{Y}(z)$ as $z\to \infty$ and
compare it to the behaviour of $\tilde{Y}_N(z)$. First we note that
$$ x^{-\Lambda}\hat{Y}_N(xz) = x^{-\Lambda} (I+O(xz)) (xz)^{\Lambda}
\hat{C}_0 \to z^{\Lambda} \hat{C}_0 \hbox{ for } x\to 0. 
$$
Then 
$$
    \left[ x^{-\Lambda} Y(xz,x) \right] \left[x^{-\Lambda}
    \hat{Y}_N(xz)  \right]^{-1} = x^{-\Lambda} U(xz,x) x^{\Lambda}
   \to \tilde{Y}(z) \hat{C}_0^{-1} z^{-\Lambda}.
$$ 
On the other hand, from the properties of $U(z,x)$ we know that
$x^{-\Lambda} U(xz,x) x^{\Lambda}$ is holomorphic in every compact
subset 
of $\{z~|~ |z|>1\}$ and $x^{-\Lambda} U(xz,x) x^{\Lambda}=
I+O\left({1\over z}\right)$ as $z\to \infty$. Thus 
$$
\tilde{U}(z) := \lim_{x\to 0} x^{-\Lambda} U(xz,x) x^{\Lambda}$$ 
exists uniformly in every compact subset of  $\{z~|~ |z|>1\}$ and 
$$ 
  \tilde{U}(z) = I +O\left({1\over z}  \right) ,~~~~z\to \infty
$$
Then 
$$ 
\tilde{Y}(z) = \tilde{U}(z) z^{\Lambda} \hat{C}_0 \equiv
\tilde{Y}_N(z) \hat{C}_0,
$$
as we wanted to prove. Finally, the above result implies 
$$
    Y(z,x) = x^{\Lambda} V({z\over x},x) \tilde{Y}_N({z\over x} )
    \hat{C}_0 
$$
$$
= \left\{ \matrix{
                         x^{\Lambda} V({z\over x} ,x) \tilde{G}_0 \left( I
                         + O(z/x)\right) x^{-J} z^J \tilde{C}_0
                         \hat{C}_0 ~=~ G_0(x) (I+O(z)) z^{J}  \tilde{C}_0
                         \hat{C}_0,~~~z\to 0
                                               \cr
                      \cr
                          x^{\Lambda} V({z\over x} ,x)
                         \tilde{G}_1\left( 
                                             O\left({z\over x}-1)\right) 
                                             \right) \left({z\over x}
                         -1\right)^J \tilde{C}_1 \hat{C}_0 
                     ~   =  ~G_x(x)(I+O(z-x)) (z-x)^J
                         \tilde{C}_1\hat{C}_0,~~~z\to x \cr 
} 
\right. 
$$ 
Let $\tilde{M}_0$, $\tilde{M}_1$ denote the monodromy matrices of
$\tilde{Y}_N(z)$ in the basis of figure 13. The above result implies: 
\be
M_0=
 \hat{C}_0^{-1}\tilde{C}_0^{-1} e^{2\pi i J} \tilde{C}_0 \hat{C}_0=\hat{C}_0^{-1}\tilde{M}_0 \hat{C}_0
\label{mono2}
\ee
\be
M_x=
 \hat{C}_0^{-1}\tilde{C}_1^{-1} e^{2\pi i J} \tilde{C}_1 \hat{C}_0=\hat{C}_0^{-1}\tilde{M}_1 \hat{C}_0
\label{mono3}
\ee
The same result may be obtained observing that from  
\be
  {d (x^{-\Lambda}Y(xz,x))\over d z}= \left[
{x^{-\Lambda}A_0 x^{\Lambda}\over z} +{x^{-\Lambda}A_x
x^{\Lambda}\over z-1}+
{x^{-\Lambda}A_1 x^{\Lambda}\over z-{1\over x}}
\right]~x^{-\Lambda}Y(xz,x)
\label{prov}
\ee
we obtain the system (\ref{systemtilde})  as $z={1\over x}$ and $z=\infty$
merge (figure \ref{figura567}). The singularities  $z=0$, $z=1$, $z=1/x$
of (\ref{prov}) correspond to $z=0$, $z=x$, $z=1$ of
(\ref{stofucs}). The poles $z=0$ and $z=1$ of (\ref{prov}) do not move
as $x\to 0$ and  ${1\over x}$  converges to $\infty$, in general  along
spirals. At any turn of the spiral the system (\ref{prov}) has new
monodromy matrices according to the action of the braid group
$$ M_1 \mapsto M_{\infty},~~~M_{\infty}\mapsto M_{\infty} M_1
M_{\infty}^{-1}$$ 
but  
$$ M_0 \mapsto M_0,~~~~M_{x}\mapsto M_x$$
Hence, the limit 
$\tilde{Y}(z)$ still has monodromy $M_0$ and $M_x$ at $z=0,x$.   Since 
$\tilde{Y}=\tilde{Y}_N
\hat{C}_0$  we conclude that $M_0$ and $M_x$ are
(\ref{mono2}) and  (\ref{mono3}).

\begin{figure}
\epsfxsize=15cm
\centerline{\epsffile{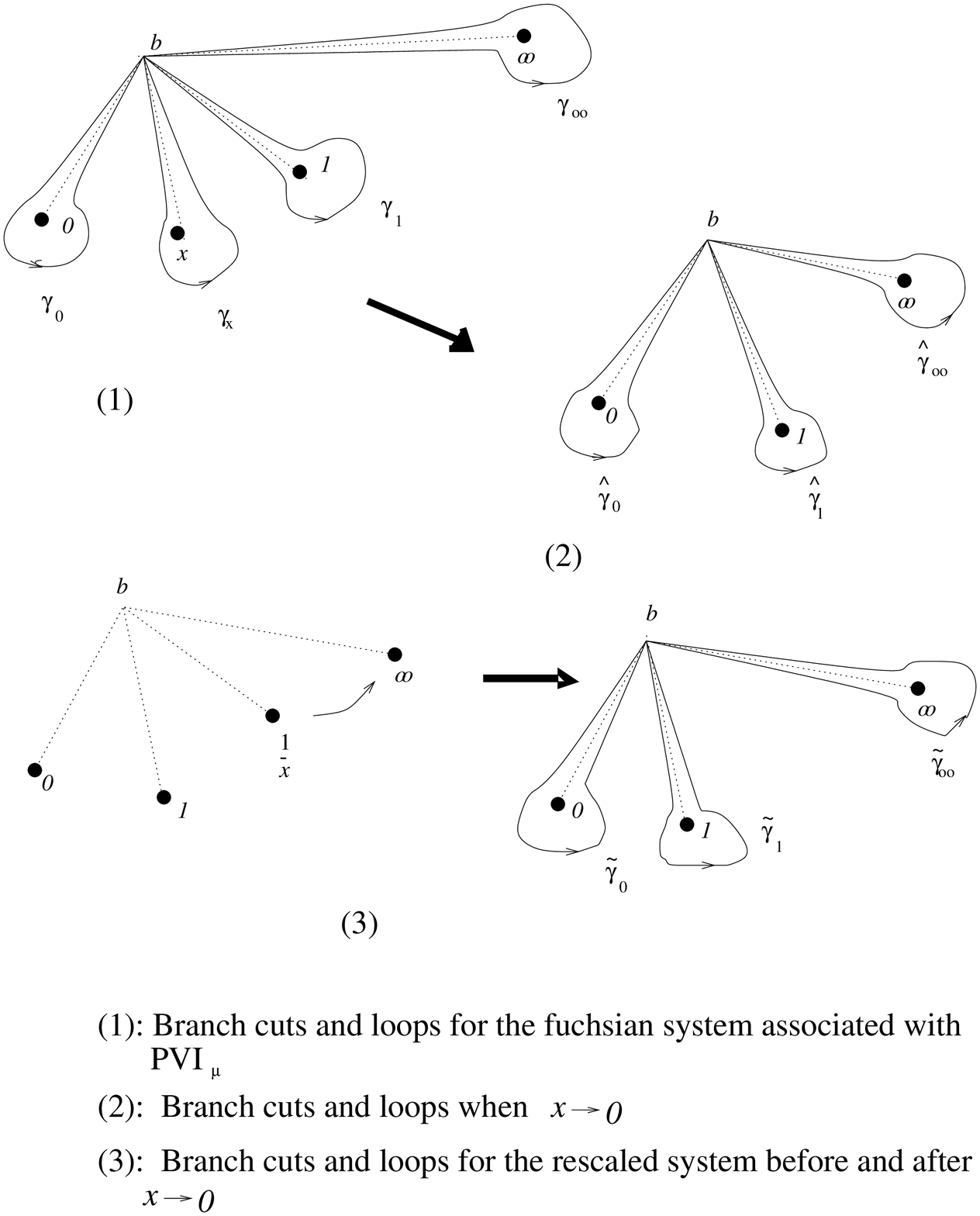}}
\caption{}
\label{figura567}
\end{figure}

\vskip 0.3 cm 

In order to find the parameterization $y({x};\sigma,a)$ in terms
of $(x_0,x_1,x_{\infty})$ we have to compute the monodromy matrices
$M_0$, $M_1$, $M_{\infty}$ in terms of $\sigma$ and $a$, and then take
the traces of their products. In order to do this we use the formulae 
(\ref{mono1}), (\ref{mono2}),(\ref{mono3}). In fact, the matrices
$\tilde{M}_i$ ($i=0,1$) and $\hat{M}_1$ can be computed explicitly
because a $2\times 2$ fuchsian system with three singular points 
can be reduced to the
hypergeometric equation, whose monodromy is completely known. 

\vskip 0.3 cm
\noindent 
{\bf Lemma 3: } {\it 
The Gauss hypergeometric equation 
\be
   z(1-z)~ {d^2 y \over dz^2} +[\gamma_0-z(\alpha_0+\beta_0+1)]~{dy\over dz}
   -\alpha_0 \beta_0 ~y=0
\label{hyper1}
\ee
 is equivalent to the system 
\be
{d\Psi\over dz}= \left[{1\over z}\pmatrix{0 & 0 \cr
                                       -\alpha_0 \beta_0 & -\gamma_0 \cr}
                                +{1\over z-1}\pmatrix{0&1\cr
                                                     0 & \gamma_0
                                       -\alpha_0-\beta_0
                                                     \cr} \right]~\Psi
\label{hyper2}
\ee
where $\Psi=\pmatrix{y \cr (z-1) {dy\over dz}\cr}$. 
}

\vskip 0.3 cm 
\noindent
{\bf  Lemma 4: } { \it Let $B_0$ and $B_1$ be matrices  of eigenvalues $0,
1-\gamma$,  and  $0, \gamma-\alpha-\beta-1$ respectively, such
that }
$$
  B_0+B_1= \hbox{ diag}(-\alpha, -\beta),~~~\alpha\neq \beta$$
{\it Then
$$B_0=\pmatrix{
                         {\alpha(1+\beta-\gamma)\over \alpha-\beta} &
                         {\alpha (\gamma-\alpha-1)\over \alpha-\beta}~r
                         \cr
               {\beta(\beta+1-\gamma)\over \alpha-\beta}~{1\over r} & 
               {\beta(\gamma-\alpha-1)\over \alpha-\beta} \cr
}
 $$
$$
   B_1=\pmatrix{  
                   {\alpha(\gamma-\alpha-1)\over \alpha-\beta} &
                   -(B_0)_{12} \cr
                   -(B_0)_{21} & {\beta(\beta+1-\gamma)\over
                   \alpha-\beta} \cr
}
$$
for any  $r\neq 0$. 
}
\vskip 0.3 cm 
\noindent
We leave the proof as an exercise.       
The following lemma connects lemmas 3 and 4:
\vskip 0.3 cm
\noindent
{\bf Lemma 5: } {\it The system (\ref{hyper2}) with 
$$\alpha_0=\alpha,~~~\beta_0=\beta+1,~~~\gamma_0=\gamma,~~~\alpha\neq \beta$$ 
is gauge-equivalent to
the system 
\be
 {d X\over dz}= \left[{B_0\over z}+{B_1\over z-1}\right] ~X
\label{hyper3}
\ee
where $B_0$, $B_1$ are given in lemma 4. This means that there exists
a matrix 
$$
   G(z):=\pmatrix{1 & 0 \cr 
 {\alpha- \beta\over (\beta+1-\gamma)\beta}~\left[\alpha ~z
 +{\alpha(\beta+1-\gamma)\over \alpha-\beta}\right]~{1\over r} & 
  z~{ \alpha-\beta \over (\beta+1-\gamma)\beta}~{1\over r}\cr
}
$$
such that $X(z)=G(z) ~\Psi(z)$. 
 It follows that (\ref{hyper3}) and the corresponding hypergeometric
equation (\ref{hyper1}) have the same fuchsian singularities
0,1,$\infty$ and the same monodromy group.
}
\vskip 0.3 cm 
\noindent
{\it Proof:}  By direct computation. $\Box$

\vskip 0.3 cm
\noindent
 Note that the form of $G(z)$ ensures that if
$y_1$, $y_2$ are independent solutions of the hypergeopmetric
equation, then a fundamental matrix of (\ref{hyper3}) may be chosen to
be 
$X(z)=
\pmatrix{ y_1(z) & y_2(z) \cr
           *  & *\cr}
$
\vskip 0.3 cm
Now we compute the monodromy matrices for the systems
(\ref{systemhat}), (\ref{systemtilde}) by reduction to an
hypergeometric equation. 
We first study the case  $\sigma \not \in {\bf Z}$. 
 Let us start with (\ref{systemhat}). With the
gauge 
$$ Y^{(1)}(z):=z^{-{\sigma\over 2}}~\hat{Y}(z)$$ 
we transform (\ref{systemhat}) in
\be 
      {d Y^{(1)}\over dz}=\left[ {A_1^0\over
      z-1}+{\Lambda-{\sigma\over2 }I\over z} \right]~Y^{(1)}
\label{novosis}
\ee
We identify the matrices $B_0$, $B_1$ with    
$\Lambda-{\sigma\over2 }I$ and  $A_1^0$ with eigenvalues 0, $-\sigma$
and 0, 0 
respectively. Moreover $A_1^0 +\Lambda-{\sigma\over2 }I=$
diag($-\mu-{\sigma\over 2},\mu-{\sigma\over 2})$. Thus:
$$
  \alpha=\mu+{\sigma\over 2},~~~\beta=-\mu+{\sigma\over
  2},~~~\gamma=\sigma+1$$
$$ \alpha-\beta= 2\mu \neq 0 ~~~\hbox{ by hypothesis}$$
The parameters of the correspondent hypergeometric equation are 
$$ 
 \left\{ \matrix{                                         
\alpha_0=\mu+{\sigma\over 2} \cr
\beta_0=1-\mu+{\sigma\over
  2} \cr
      \gamma_0=\sigma+1 \cr
}\right.
$$
 From them  we deduce the nature of two linearly independent solutions
  at $z=0$. Since $\gamma_0\not \in {\bf Z}$   ($\sigma \not \in
  {\bf Z}$   ) the solutions are expressed in terms of
  hypergeometric functions. On the other hand, the effective 
  parameters at $z=1$ and $z=\infty$ are respectively:
$$ 
 \left\{ \matrix{
\alpha_1:=\alpha_0=\mu+{\sigma\over 2} \cr
\beta_1 := \beta_0=1-\mu+{\sigma\over
  2} \cr
     \gamma_1:=\alpha_0+\beta_0-\gamma_0+1 =1 \cr
}\right.
$$
$$ 
 \left\{ \matrix{
\alpha_{\infty}:=\alpha_0=\mu+{\sigma\over 2} \cr
\beta_{\infty}:= \alpha_0-\gamma_0+1=\mu-{\sigma\over
  2} \cr
     \gamma_{\infty}=\alpha_0-\beta_0+1 =2\mu \cr
}\right.
$$
Since $\gamma_1=1$, at least one  solution has a
logarithmic singularity at $z=1$. Also note that 
$\gamma_{\infty}=2\mu$, therefore  logarithmic singularities appear
at $z=\infty$ if $2\mu \in {\bf Z}\backslash \{0\}$. 

\vskip 0.15 cm 
For the derivations which follows,  we use the notations of the
fundamental paper by Norlund \cite{Norlund}. 
 To derive the connection formulae we use
the paper of Norlund when logarithms are involved. Otherwise, in the
generic case, any textbook of special functions (like \cite{Luke}) may
be used.

\vskip 0.2 cm
 \noindent
{\bf First case: $\alpha_0, \beta_0 \not \in {\bf Z}$}. This means 
$$ 
  \sigma \neq \pm 2\mu +2m ,~~~~m\in{\bf Z}
$$

We can choose the following independent solutions of the
hypergeometric equation: 

\vskip 0.15 cm
\noindent
At $z=0$ 
$$ 
  y_1^{(0)}(z)=F(\alpha_0,\beta_0,\gamma_0;z)$$
\be
  y_2^{(0)}(z)=
  z^{1-\gamma_0}~F(\alpha_0-\gamma_0+1,\beta_0-\gamma_0+1,2-\gamma_0;z) 
\label{ffirr1}
\ee
where $F(\alpha,\beta,\gamma;z)$ is the well known hypergeometric
  function (see \cite{Norlund}).

\vskip 0.15 cm
\noindent
At $z=1$
$$ 
 y_1^{(1)}(z)==F(\alpha_1,\beta_1,\gamma_1;1-z)$$
$$
y_2^{(1)}(z)=g(\alpha_1,\beta_1,\gamma_1;1-z)
$$
Here $g(\alpha,\beta,\gamma;z)$ is a logarithmic solution introduced in
\cite{Norlund}, and $\gamma\equiv \gamma_1=1$. 

\vskip 0.15 cm
\noindent
At $z=\infty$, we consider first the case $2\mu \not \in {\bf Z}$, while the
resonant case will be considered later. Two independent solutions are:  
$$
 y_1^{(\infty)}=z^{-\alpha_0}
 ~F(\alpha_{\infty},\beta_{\infty},\gamma_{\infty} ;{1\over z})$$
$$
y_2^{(\infty)}=z^{-\beta_0}~F(\beta_0,\beta_0-\gamma_0+1,\beta_0-\alpha_0
+1; {1\over z})
$$
Then, from the connection formulas between $F(...;z)$ and $g(...;z)$ 
of \cite{Luke} and \cite{Norlund} we
derive 
$$
  [y_1^{(\infty)},y_2^{(\infty)}]=[y_1^{(0)},y_2^{(0)}] ~C_{0\infty}
$$
$$ 
  C_{0\infty}= \pmatrix{ 
                 e^{-i\pi \alpha_0}
                 {\Gamma(1+\alpha_0-\beta_0)\Gamma(1-\gamma_0) \over
                 \Gamma(1-\beta_0) \Gamma(1+\alpha_0-\gamma_0)}& 
                 e^{-i\pi \beta_0}{
                 \Gamma(1+\beta_0-\alpha_0)\Gamma(1-\gamma_0)\over
                 \Gamma(1-\alpha_0) \Gamma(1+\beta_0-\gamma_0)} \cr
                 e^{i\pi(\gamma_0-\alpha_0-1)}{\Gamma(1+\alpha_0-\beta_0)
\Gamma(\gamma_0 -1)\over \Gamma(\alpha_0)\Gamma(\gamma_0-\beta_0)}&
                      e^{i\pi(\gamma_0-\beta_0-1)}{\Gamma(1+\beta_0-\alpha_0)
\Gamma(\gamma_0-1)\over \Gamma(\beta_0)\Gamma(\gamma_0-\alpha_0)} \cr
}
$$

\vskip 0.2 cm
$$
[y_1^{(0)},y_2^{(0)}]=[y_1^{(1)},y_2^{(1)}] ~C_{01}
$$
$$
   C_{01}=\pmatrix{
                   0 & -{\pi \sin(\pi(\alpha_0+\beta_0))\over
                   \sin(\pi\alpha_0) \sin(\pi
                   \beta_0)}{\Gamma(2-\gamma_0)\over
                   \Gamma(1-\alpha_0) \Gamma(1-\beta_0) } \cr
                   -{ \Gamma(\gamma_0) \over \Gamma(\gamma_0
                   -\alpha_0) \Gamma(\gamma_0-\beta_0)} & -{
                   \Gamma(2-\gamma_0) \over \Gamma(1-\alpha_0)
                   \Gamma(1- \beta_0)} \cr
}
$$
We observe that 
$$Y^{(1)}(z)= \left(I+{F\over
z}+O\left({1\over z^2}\right)\right)~z^{\hbox{diag}(-\mu-{\sigma\over2},\mu-{\sigma\over 2})},~~~~z\to
\infty
$$
$$ 
 =
\hat{G}_0(I+O(z))~z^{\hbox{diag}(0,-\sigma)}~\hat{G}_0^{-1}\hat{C_0},~~~~z\to
0
$$
$$= \hat{G}_1(I+O(z-1))~(z-1)^J~\hat{C}_1,~~~~z\to 1
$$
where $\hat{G}_0\equiv T$ of lemma 2; namely $\hat{G}_0^{-1}\Lambda
\hat{G}_0=$ diag$({\sigma\over 2},-{\sigma\over 2})$. 
By direct substitution in the
differential equation we compute the coefficient $F$ 
$$
   F= -\pmatrix{ (A_1^0)_{11} & {(A_1^0)_{12}\over 1-2\mu} \cr
                  {(A_1^0)_{21}\over 1+2\mu} &
         (A_1^0)_{22}
}
,~~~\hbox{where}
~~A_1^0={\sigma^2 -(2\mu)^2\over 8\mu}\pmatrix{ 1  &  -r \cr
                                                  {1\over r} & -1 \cr}
$$
Thus, from the asymptotic behaviour of the hypergometric function
($F(\alpha,\beta, \gamma;{1\over z})\sim 1$, $z\to \infty$ ) we
derive 
$$
  Y^{(1)}(z)= \pmatrix{ y_1^{(\infty)}(z) & r{\sigma^2-(2\mu)^1\over
  8\mu(1-2\mu) } ~y_2^{(\infty)} \cr
  * & * \cr
}
$$
From  
\be
 Y^{(1)}(z) \sim \pmatrix{ 1 & z^{-\sigma} \cr
 * & * \cr}
~\hat{G}_0^{-1} \hat{C}_0,~~~z\to 0
\label{ffirr2}
\ee
we derive 
 $$
 Y^{(1)}(z)= 
          \pmatrix{ y_1^{(0)}(z) & y_2^{(0)}(z) \cr
* & * \cr
} 
\hat{G}_0^{-1} \hat{C}_0
$$
Finally, observe that $\hat{G}_1=\pmatrix{ a & {a\over \omega} +b r \cr
                                 {a\over r} & b \cr
}
$ for arbitrary $a,b\in{\bf C}$, $a\neq 0$, and $\omega:= 
{\sigma^2-(2\mu)^2\over 8\mu}$. We recall 
 that $y_2^{(1)}=g(\alpha_1,\beta_1,1;1-z)\sim
                                 \psi(\alpha_1)+\psi(\beta_1)-2\psi(1)
                                 -i\pi +\log(z-1)$, $|\arg(1-z)|<
\pi$, as $z\to 1$. We can choose  $a=1$ and a suitable $b$, in such a way
                                 that the asymptotic behaviour of
                                 $Y^{(1)}$ for $z\to 1$ is precisely
                                 realized by 
$$
  Y^{(1)}(z)= \pmatrix{ y_1^{(1)}(z) & y_2^{(1)}(z) \cr
                              * & * \cr
}~\hat{C}_1
$$
Therefore  we conclude that the connection matrices are: 
$$
  \hat{C}_0= \hat{G}_0 \pmatrix{
(C_{0\infty})_{11} & r~{\sigma^2-(2\mu)^2\over 8\mu
(1-2\mu)}(C_{0\infty})_{12}\cr
(C_{0\infty})_{21}& r~ {\sigma^2-(2\mu)^2\over 8\mu
(1-2\mu)}(C_{0\infty})_{22}\cr
}
$$
$$
\hat{C}_1= C_{01}~(\hat{G}_0^{-1}\hat{C}_0)= C_{01}   
\pmatrix{
            (C_{0\infty})_{11} & r~{\sigma^2-(2\mu)^2\over 8\mu
(1-2\mu)}(C_{0\infty})_{12} \cr
(C_{0\infty})_{21}& r~ {\sigma^2-(2\mu)^2\over 8\mu
(1-2\mu)}(C_{0\infty})_{22} \cr
}
$$

\vskip 0.2 cm
 It's now time to consider the more complicated resonant case $2\mu \in {\bf
Z}\backslash \{0\}$. The behaviour of $Y^{(1)}$ at $z=\infty$ is 
$$
   Y^{(1)}(z) = \left(I +{F\over z} +O\left({1\over z^2} \right) \right) 
                   ~z^{\hbox{diag}(-\mu -{\sigma\over 2}, \mu -{\sigma\over
                   2}) } z^R$$
$$ 
   R= \pmatrix{ 0 & R_{12} \cr 
                0 & 0      \cr}, ~~~\hbox{ for } \mu ={1\over 2},1,{3\over
                2}, 2, {5\over 2},...
$$
$$ 
   R= \pmatrix{ 0 & 0 \cr 
                R_{21} & 0      \cr}, ~~~\hbox{ for } \mu =-{1\over 2},-1,
                - {3\over2}, -2, -{5\over 2},...
$$
and the entry $R_{12}$ is determined by the entries of $A^{0}_1$. For example,
if $\mu={1\over 2}$ we can compute $R_{12}= (A^0_1)_{12}= -r {\sigma^2-1 \over
4} $ (and $F_{12}$ arbitrary); if $\mu=-{1\over 2}$ we have $ R_{21}=
(A^0_1)_{21}=-{1\over r}
{\sigma^2 -1 \over 4}$ (and $F_{21}$ arbitrary); if $\mu =1$  we have 
$R_{12}= -r {\sigma^2(\sigma^2-4)\over 32}$.

Since $\sigma\not \in {\bf Z}$, $R\neq 0$. 
This is true for any $2\mu  \in {\bf Z}\backslash
\{0\}$. Note that the $R$ computed here coincides (by
isomonodromicity) to the $R$ of the system(\ref{stofucs}).

There is a logarithmic solution at $\infty$. 
Only
$C_{0\infty}$ 
and thus
$\hat{C}_0$  and $\hat{C}_1$ change with respect to the non-resonant
case. 
We will see in a while that such matrices
disappear in the computation of tr($M_i M_j)$, $i,j=0,1,x$. Therefore, it is
not necessary to know them explicitly. Actually, it was not necessary to
compute them also in the non resonant case, the only important matrix
to know being $C_{01}$, which is not affected by resonance of
$\mu$. This is the reason why the formulae of theorem 2 hold true also
in the resonant case.


\vskip 0.3 cm
\noindent
{\bf Second case:} $\alpha_0, \beta_0 \in {\bf Z}$, namely
$$ 
\sigma=\pm 2\mu+2m ,~~~~~~m\in {\bf Z}
$$

The formulae are almost identical to the first case, but $C_{01}$
changes. To see this, we need to distinguish four  cases.

i) 
  $\sigma=2\mu+2m$, $m=-1,-2,-3,...$. We choose
$$
   y_2^{(1)}(z)= g_0(\alpha_1,\beta_1,\gamma_1;1-z)
$$
Here $g_0(z)$ is another logarithmic solution of \cite{Norlund}. 
Thus 
$$C_{01}= \pmatrix{ {\Gamma(-m)\Gamma(-2\mu-m+1)\over \Gamma(-2\mu-2m)}
& 0 \cr
 0 & -{\Gamma(1-2\mu -2m)\over \Gamma(1-m-2\mu)\Gamma(-m)} \cr
}
$$
As usual, the matrix is computed from the connection formulas between
the hypergeometric functions and $g_0$ that the reader can find in
\cite{Norlund}. 

ii) $\sigma=2\mu+2m$, $m=0,1,2,...$. We choose
$$
y_1^{(2)}=g(\alpha_1,\beta_1,\gamma_1;1-z)
$$
Thus
$$C_{01}= \pmatrix{0 & {\Gamma(m+1)\Gamma(2\mu+m)\over
\Gamma(2\mu+2m)}\cr
-{\Gamma(2\mu+2m+1)\over \Gamma(2\mu+m)\Gamma(m+1)} & 0 
\cr}
$$

iii)  $\sigma=-2\mu+2m$, $m=0,-1,-2,...$. We choose
$$
   y_2^{(1)}(z)=  g_0(\alpha_1,\beta_1,\gamma_1;1-z)
$$
Thus
$$C_{01}= \pmatrix{ {\Gamma(1-m)\Gamma(2\mu-m)\over \Gamma(2\mu-2m)}
& 0 \cr
 0 & -{\Gamma(1+2\mu -2m)\over \Gamma(2\mu-m)\Gamma(1-m)} \cr
}
$$

iv)  $\sigma=-2\mu+2m$, $m=1,2,3,...$. We choose
 $$
  y_2^{(1)}(z)= g(\alpha_1,\beta_1,\gamma_1;1-z)
$$
Thus
$$C_{01}= \pmatrix{0 & {\Gamma(m)\Gamma(m+1-2\mu)\over
\Gamma(2m-2\mu)}\cr
-{\Gamma(2m+1-2\mu)\over \Gamma(m+1-2\mu)\Gamma(m)} & 0 
\cr}
$$

Note that this time $F=\pmatrix{0 & {r\over 1-2\mu} \cr
 0 & 0 \cr}$ in the case $\sigma=\pm 2\mu$ (i.e. $m=0$) because
 $A_1^0$ has a special form in this case.  Then in
 $\hat{C}_0$ the elements $ {\sigma^2-(2\mu)^2\over 8\mu
(1-2\mu)}(C_{0\infty})_{12}$, $ {\sigma^2-(2\mu)^2\over 8\mu
(1-2\mu)}(C_{0\infty})_{22}$ must be substituted, for $m=0$,  with 
$ {1\over 
1-2\mu}(C_{0\infty})_{12}$, 
$ {1\over 
1-2\mu}(C_{0\infty})_{22}$.

\vskip 0.3 cm 
We turn to the system (\ref{systemtilde}). Let $\tilde{Y}$ be the
fundamental matrix (\ref{solutiontilde}). With the gauge 
$$ Y^{(2)}(z):=\hat{G}_0^{-1}~\left(\tilde{Y}(z)\hat{G}_0 \right)$$ 
we have 
$$
  {d Y^{(2)} \over d z} = \left[ {\tilde{B}_0 \over
  z}+{\tilde{B}_1\over z-1} \right] Y^{(2)}
$$
  $$ \tilde{B}_0=\hat{G}^{-1} A_0^0 \hat{G}_0 = \pmatrix{{\sigma\over
4} & {\sigma\over 4} s \cr
-{\sigma \over 4 s} & -{\sigma \over 4} \cr
  }$$
   $$ \tilde{B}_1=\hat{G}^{-1} A_x^0 \hat{G}_0 = \pmatrix{{\sigma\over
4} & -{\sigma\over 4} s \cr
{\sigma \over 4 s} & -{\sigma \over 4} \cr
  }$$
This time then 
$$
\left\{\matrix{\alpha_0=-{\sigma\over 2} \cr
                \beta_0={\sigma\over 2} +1 \cr
\gamma_0=1
}\right.  
$$
$$
\left\{\matrix{\alpha_1=-{\sigma\over 2} \cr
                \beta_1={\sigma\over 2} +1\cr
\gamma_1=1
}\right.  
$$
$$
\left\{\matrix{\alpha_{\infty}=-{\sigma\over 2} \cr
                \beta_{\infty}={\sigma\over 2} \cr
\gamma_{\infty}=\sigma
}\right.  
$$
If follows that both at $z=0$ and $z=1$ there are logarithmic
solutions. We skip all the derivation of the connection formulae,
which is done as in the previous cases, with some more technical
complications. Before giving the results we observe that 
$$
   Y^{(2)}(z)=\left(I+O\left({1\over z}\right)\right)
 z^{\hbox{diag}({\sigma\over
   2},-{\sigma\over 2})} ,~~~~z\to \infty
$$
$$
  = \hat{G}_0^{-1}\tilde{G}_0~(1+O(z))z^J~C_0^{\prime},~~~~z\to 0
$$
$$
 = \hat{G}_0^{-1}\tilde{G}_1~(1+O(z-1))(z-1)^J~C_1^{\prime},~~~~z\to 1
$$
where 
$$  
   C_i^{\prime}:= \tilde{C}_i\hat{G}_0,~~~~~i=0,1
$$
Then
 $$
        \tilde{M}_0=\hat{G}_0 ~(C_0^{\prime})^{-1} ~
\pmatrix{1&2\pi i \cr
0& 1}
~  C_0^{\prime}~\hat{G}_0^{-1}
$$
$$
        \tilde{M}_1=\hat{G}_0 ~(C_1^{\prime})^{-1} ~
\pmatrix{1&2\pi i \cr
0& 1}
~  C_1^{\prime}~\hat{G}_0^{-1}
$$
Then, the connection problem may be solved computing
$C_i^{\prime}$. The result is

$$
C_0^{\prime} = 
\pmatrix{
           (C_{0\infty}^{\prime})_{11} & {\sigma\over \sigma+1}{s\over
                                                                      4}
           (C_{0\infty}^{\prime})_{12}             \cr
      (C_{0\infty}^{\prime})_{21}&{\sigma\over \sigma+1}{s\over
                                       4} (C_{0\infty}^{\prime})_{22}\cr
            }
$$
\vskip 0.2 cm
$$
   C_1^{\prime}=C_{01}^{\prime} C_0^{\prime}$$

\noindent
where 
     $$ 
(C_{0\infty}^{\prime})^{-1}= 
\pmatrix{ 
{\Gamma(\beta_0-\alpha_0)\over
\Gamma(\beta_0)\Gamma(1-\alpha_0)} e^{i\pi\alpha_0} & 0 \cr
 { \Gamma(\alpha_0-\beta_0)\over
\Gamma(\alpha_0)\Gamma(1-\beta_0)} e^{i\pi\beta_0}&
-{\Gamma(1-\alpha_0)\Gamma(\beta_0)\over
\Gamma(\beta_0-\alpha_0+1)}~e^{i\pi \beta_0} \cr
}
$$
$$
C_{01}^{\prime}= \pmatrix{
0 & -{\pi \over \sin(\pi\alpha_0)} \cr
                                 -{\sin(\pi\alpha_0)\over \pi} &
                                 -e^{-i\pi \alpha_0} \cr
                          }$$

\vskip 0.3 cm 

The case $\sigma\in{\bf Z}$ interests us only if $\sigma=0$ (otherwise
$\sigma\not \in \Omega$). We observe that the system (\ref{systemhat})
is precisely the system for $Y^{(2)}(z)$ with the substitution
$\sigma\mapsto -2\mu$. In the formulae for $x_i^2$, $i=0,1,\infty$ 
 we only need
$C_{01}$, which is obtained from $C_{01}^{\prime}$ with
$\alpha_0=\mu$. 

As for the system (\ref{systemtilde}), the gauge
$Y^{(2)}=\hat{G}_0^{-1}\tilde{Y} \hat{G}_0$   
yields  $\tilde{B}_0=\pmatrix{0&s\cr 0 & 0\cr}$, $\tilde{B}_1= \pmatrix{
0 & 1-s \cr
0 & 0 \cr}$. Here $\hat{G}_0$ is the matrix such that $\hat{G}_0^{-1}
\Lambda \hat{G}_0= \pmatrix{0 & 1 \cr 0 & 0 \cr}$. The behaviour of
$Y^{(2)}(z)$ is now:  
$$
 Y^{(2)}(z)=(I+O({1\over z}))~z^J ~~~~~~~~~~~z\to \infty
$$  
$$
  = \tilde{\tilde{G}}_0^{-1}~(1+O(z))z^J~C_0^{\prime},~~~~z\to 0
$$
$$
 = \tilde{\tilde{G}}_1~(1+O(z-1))(z-1)^J~C_1^{\prime},~~~~z\to 1
$$ 
Here  $\tilde{\tilde{G}}_i$ is  the matrix that puts 
$\tilde{B}_i$ in Jordan form, for $i=0,1$.  
$Y^{(2)}$ can be computed explicitly: 
$$ Y^{(2)}(z)= \pmatrix{1 & s\log(z)+(1-s)\log(z-1) \cr
                       0 & 1\cr
}
$$
If we choose   $\tilde{\tilde{G}}_0=$diag$(1,1/s)$, then  
$$
  C_0^{\prime}=\pmatrix{1&0\cr
                        0 & s\cr}
$$
In the same way we find 
$$
  C_1^{\prime}=\pmatrix{1&0\cr
                        0 & 1-s\cr}
$$

\vskip 0.3 cm 
To prove theorem 2 it is now enough to compute 
$$
   2- x_0^2= \hbox{ tr}(M_0M_x)\equiv  \hbox{ tr}(e^{2\pi i J} 
(C_{01}^{\prime})^{-1} e^{2\pi i J} 
C_{01}^{\prime} )
$$
$$
2-x_1^2=  \hbox{ tr}(M_xM_1)\equiv  
 \hbox{ tr}((C_1^{\prime})^{-1}e^{2\pi i J} C_1^{\prime} C_{01}^{-1}
 e^{2\pi i J} C_{01})
$$

$$
2-x_{\infty}^2=  \hbox{ tr}(M_0M_1)\equiv  
 \hbox{ tr}((C_0^{\prime})^{-1}e^{2\pi i J} C_0^{\prime} C_{01}^{-1}
 e^{2\pi i J} C_{01})
$$
Note the remarkable simplifications obtained from the cyclic property
of the trace (for example, $\hat{C}_0$, $\hat{C}_1$   and $\hat{G}_0$ 
disappear). 
The fact that
$\hat{C}_0$ and 
$\hat{C}_1$  disappear implies that the formulae of theorem 2 
are derived for any
$\mu\neq 0$, including the resonant cases. Thus, the connection
formulae in the resonant case $2\mu \in {\bf Z}\backslash \{0\}$ are
the same of the non-resonant case. 
The final result of the computation of the traces is: 

\vskip 0.15 cm
\noindent
 I) Generic case: 
\be
   \left\{ \matrix{2(1-cos(\pi \sigma))=x_0^2 \cr\cr
                  {1\over
                  f(\sigma,\mu)}\left(2+F(\sigma,\mu)~s+{1\over
                  F(\sigma,\mu)~s}\right) =x_1^2
                         \cr\cr
                  {1\over f(\sigma,\mu)}\left(
                         2-F(\sigma,\mu)e^{-i\pi\sigma}~s-{1 \over
                  F(\sigma,\mu)
e^{-i\pi\sigma}~s }
                         \right) =x_{\infty}^2\cr
                         }\right. 
\label{20ii}
\ee
where   
$$
f(\sigma,\mu)={2 \cos^2({\pi \over 2} \sigma) \over \cos(\pi \sigma)- 
                  \cos(2\pi\mu)}\equiv{4-x_0^2\over
    x_1^2+x_{\infty}^2-x_0x_1x_{\infty}}, ~~~
~~~ F(\sigma,\mu)=f(\sigma,\mu) {16^{\sigma}\Gamma({\sigma+1\over 2})^4\over
                  \Gamma(1-\mu+ {\sigma\over
                  2})^2\Gamma(\mu+{\sigma\over 2})^2}$$

\vskip 0.15 cm 
\noindent
 II) $\sigma\in 2{\bf Z}$, $x_0=0$. 
           $$
             \left\{
                     \matrix{ 2(1-cos(\pi \sigma))=0\cr\cr
                              4\sin^2(\pi\mu) ~(1-a)=x_1^2 \cr\cr 
                             4 \sin^2(\pi\mu) ~a=x_{\infty}^2               
                 }
 \right.
$$

\vskip 0.15 cm 
\noindent
III) $x_0^2=4 \sin^2(\pi\mu)$.  Then (\ref{10ii}) implies
 $x_{\infty}^2=-x_1^2 ~\exp(\pm 2\pi i \mu)$ . 
 Four cases which yield the values of $\sigma$ non included in I)
and II) must be considered 

\vskip 0.15 cm
   III1)  $x_{\infty}^2=-x_1^2 e^{- 2
\pi i \mu}$
$$\sigma=  2\mu + 2m,~~~~m=0,1,2,...$$ 
 $$ s=
               { \Gamma(m+1)^2 \Gamma(2\mu+m)^2 \over 16^{2\mu+2m}
 \Gamma(\mu+m+{1\over 2})^4}~x_1^2 
$$

\vskip 0.15 cm
  III2) $x_{\infty}^2=-x_1^2 e^{2\pi i \mu}$
$$\sigma=2\mu+2m,~~~~m=-1,-2,-3,...$$
    $$
            s={\pi^4\over \cos^4(\pi\mu)}\left[
                 16^{2\mu+2m} \Gamma(\mu+m+{1\over 2})^4
\Gamma(-2\mu-m+1)^2 \Gamma(-m)^2~ x_1^2 \right]^{-1} 
$$

\vskip 0.15 cm
  III3) $x_{\infty}^2=-x_1^2e^{2\pi i\mu}$
$$\sigma=-2\mu+2m,~~~~m=1,2,3,...$$
$$s=
         {
\Gamma(m-2\mu+1)^2 \Gamma(m)^2 \over   16^{-2\mu+2m} \Gamma(-\mu+m+{1\over 2})^4}~ x_1^2 
$$

\vskip 0.15 cm
III4)  $x_{\infty}^2=-x_1^2 e^{-2\pi i \mu}$
$$\sigma=-2\mu+2m,~~~~m=0,-1,-2,-3,...$$
$$
  s=             {\pi^4\over \cos^4(\pi\mu)} \left[16^{-2\mu+2m}
\Gamma(-\mu+m+{1\over 2})^4\Gamma(2\mu-m)^2 \Gamma(1-m)^2 ~ x_1^2
\right]^{-1} 
$$

\vskip 0.2 cm 
 To compute $\sigma$ and $s$ in  the generic case $I$), with $x_0^2\neq 4$, 
 we solve the system (\ref{20ii}). It has two unknowns  and three equations and
we need to prove that it is 
compatible. Actually, the first
equation $ 2(1-cos(\pi \sigma))=x_0^2$ has always solutions. Let us choose
a solution 
$\sigma_0$  ($\pm \sigma_0 + 2 n$, $\forall n\in {\bf Z}$
are also solutions). Substitute it in the last two equations. We need
to verify they are compatible.  
Instead of $s$ and ${1\over s}$ write $X$ and $Y$. We have the linear
system in two variable $X$, $Y$
$$ 
                 \pmatrix{ F(\sigma_0) & {1 \over F(\sigma_0)} \cr\cr
                            F(\sigma_0)~e^{-i \pi \sigma_0}& {1\over
                           F(\sigma_0)} e^{-i\pi\sigma_0} \cr
}
 \pmatrix{X\cr \cr Y\cr}=
                          \pmatrix{
                                    f(\sigma_0) ~x_1^2-2 \cr
\cr
                                    2 - f(\sigma_0) ~x_{\infty}^2 \cr
}$$
 The system has a unique solution if and only if $2 i \sin(\pi \sigma_0)=\det 
     \pmatrix{ F(\sigma_0) & {1 \over F(\sigma_0)} \cr\
                           F(\sigma_0)~e^{-i \pi \sigma_0}& {1\over
                           F(\sigma_0)} e^{-i\pi\sigma_0} \cr
}\neq 0$. This happens for $\sigma_0 \not\in {\bf Z}$. The condition
                           is
                            not restrictive, because for $\sigma$ even
                           we
                           turn to the case $II$) of the theorem 2,
 and $\sigma$ odd is not in $\Omega$.  The solution
                           is then 
$$ X= {2 (1+e^{-i \pi\sigma_0})-f(\sigma_0)(x_1^2+x_{\infty}^2 e^{-i \pi
\sigma_0})\over F(\sigma_0)(e^{-2 \pi i \sigma_0} -1)}$$
$$  
   Y=
   F(\sigma_0)~{f(\sigma_0)e^{-i\pi\sigma_0}(e^{-i\pi\sigma_0}x_1^2+x_{\infty}^2)-2 e^{-i\pi\sigma_0}(1+e^{-i\pi\sigma_0})\over e^{-2\pi i \sigma_0}-1}
$$
Compatibility of the system means that $X~Y\equiv 1$. This is
verified  by direct computation:
$$ 
   XY= {e^{-i\pi\sigma}~\left[
   2(1+e^{-i\pi\sigma})-(x_1^2+x_{\infty}^2e^{-i\pi\sigma})f(\sigma) 
\right]~\left[
           (x_1^2e^{-i\pi\sigma}+x_{\infty}^2)f(\sigma)- 
 2(1+e^{-i\pi\sigma}) 
\right] \over (e^{-2i\pi\sigma}-1)^2}$$
$$
    = { 8 \cos^2({\pi\sigma\over
    2})(x_1^2+x_{\infty}^2)f(\sigma)
    -4(4-\sin^2({\pi\sigma\over
    2}))-((x_1^2+x_{\infty}^2)^2-x_0^2x_1^2x_{\infty}^2) f(\sigma)^2
    \over -4 \sin^2(\pi\sigma)}$$
Using the relations $\cos^2({\pi \sigma \over 2})=1-x_0^2/4$,
    $\cos(\pi\sigma) =1-x_0^2/2$ and $f(\sigma)= {4-x_0^2\over
    x_1^2+x_{\infty}^2-x_0x_1x_{\infty}}$ we obtain 
$$
 XY= {1\over x_0^2}~\left(
    -2(x_1^2+x_{\infty}^2)f(\sigma)+4+{(x_1^2+x_{\infty}^2)^2- 
        (x_0x_1x_{\infty})^2 \over
    x_1^2+x_{\infty}^2-x_0x_1x_{\infty}} f(\sigma)  
                         \right)
$$
$$ 
  ={1 \over x_0^2}~\left(4 -(x_1^2+x_{\infty}^2-x_0x_1x_{\infty})f(\sigma)
  \right)= {1\over x_0^2} \left(4-(4-x_0^2)\right)=1$$

It follows from this construction that for any $\sigma$ solution of
the first equation of (\ref{20ii}), there always exists a unique $s$
which solves the last two equations. Recall that 
 $$ a =-{1 \over 4 s}
$$

\vskip 0.2 cm 

To complete the proof of theorem 2 (points $i)$, $ii)$, $iii)$), 
 we just have to compute the square
roots of the $x_i^2$ ($i=0,1,\infty$) in such a way that (\ref{10ii}) is
satisfied.  For example, the square root of I) satisfying (\ref{10ii}) is 
$$
   \left\{ \matrix{x_0=2 \sin({\pi \over 2}\sigma) \cr\cr
                  x_1={1\over \sqrt{f(\sigma,\mu)}}
                    \left(\sqrt{F(\sigma,\mu)~s}+{1\over
 \sqrt{F(\sigma,\mu)~s}
                 }\right)
                      \cr\cr
                  x_{\infty}={i\over
 \sqrt{f(\sigma,\mu)}}\left(\sqrt{F(\sigma,\mu)~s}~e^{-i{\pi \sigma\over
 2}}-{1\over \sqrt{F(\sigma,\mu)~s}~e^{-i{\pi \sigma\over 2}}}\right)\cr
                         }\right. 
$$
which yields i), with $F(\sigma,\mu)=f(\sigma,\mu) (2 G(\sigma,\mu))^2$. . 

\vskip 0.2 cm

We remark that in case II) only $\sigma =0$ is in $\Omega$. If $\mu$ integer
in II), 
the formulae  give $(x_0,x_1,x_{\infty})=(0,0,0)$. The
triple is not admissible, and direct computation gives $R=0$ for the
system (\ref{novosis}). This is the case of commuting monodromy
matrices with a 1-parameter family of rational solutions of
$PVI_{\mu}$. $\Box$
\footnote{
{\it Remark:} In order to solve the R.H. for the monodromy data
$(x_0,x_1,x_{\infty})$ we  choose branch cuts in the $x$-plane. When $x$ is
small we just need to fix $\alpha<\arg(x)<\alpha+2\pi$, $\alpha\in{\bf R}$.
 Let $A(z,x;x_0,x_1,x_{\infty})$ be the matrix solution of the R.H. 
Consider the loop $x\mapsto x^{\prime}=x e^{2\pi i}$; the 
analytic continuation of $A$ along the loop 
is $A(z,x^{\prime};x_0,x_1,x_{\infty})\equiv
A(z,x; x_0^{\beta_1^2},x_1^{\beta_1^2},x_{\infty}^{\beta_1^2})$. In other
words,  $A(z,x^{\prime};x_0,x_1,x_{\infty})$ coincides with 
 the solution of the
R.H. obtained with the same branch cut $\alpha<\arg(x)<\alpha+2\pi$ and 
new monodromy data
$(x_0^{\beta_1^2},x_1^{\beta_1^2},x_{\infty}^{\beta_1^2})$ transformed by the
action of the braid group.  When we write
$A(z,x^{\prime};x_0,x_1,x_{\infty})$ we are considering $A$ as a function on
the universal covering of  ${\bf C}_0\cap \{|x|<\epsilon\}$; when we write 
 $ A(z,x; x_0^{\beta_1^2},x_1^{\beta_1^2},x_{\infty}^{\beta_1^2})$ we are
considering the solution of the R.H. as a ``branch''. 

In the proof of theorem 2 we start from a point $x\in D(\sigma)$ 
and we take the limits
$x\to 0$ of the system and of the rescaled system. At $x$ we assign the
monodromy $M_0, M_1, M_x$ characterized by $(x_0,x_1,x_{\infty})$ and then we
take the limit proving the theorem. 
If we had start from another point $x^{\prime}=x e^{2\pi i}\in D(\sigma)$
(provided that this is possible for the given $D(\sigma)$ and $x$) we  take
 the same monodromy $M_0, M_1, M_x$, because what we were
doing is  the limit, for $x\to 0$ in $D(\sigma)$, of
$A(z,x;x_0,x_1,x_{\infty})$ considered as a function defined on the universal
covering of ${\bf C}_0\cap \{|x|<\epsilon\}$.     
}

 \vskip 0.3 cm 
\noindent
{\it  Proof of remark 2}

We prove that $a(\sigma)={1\over 16 a(-\sigma)}$, namely 
$s(\sigma)={1\over s(-\sigma)}$ for $a=-{1\over 4s}$. 
 Given monodromy data $(x_0, x_1,
 x_{\infty})$ the parameter  $s$
 corresponding to $\sigma$ is uniquely determined by  
 $$
  {1\over
                  f(\sigma)}\left(2+F(\sigma)~s+{1\over
                  F(\sigma)~s}\right) =x_1^2
  $$
 $$                      
                  {1\over f(\sigma)}\left(
                         2-F(\sigma)e^{-i\pi\sigma}~s-{1 \over
                  F(\sigma)
 e^{-i\pi\sigma}~s }
                          \right) =x_{\infty}^2$$
   We observe that $f(\sigma)=f(-\sigma)$ and that  
 the properties of the Gamma function 
 $$
   \Gamma(1-z)\Gamma(z)= {\pi \over \sin(\pi z)},~~~~~
  \Gamma(z+1)=z\Gamma(z)$$
 imply
 $$F(-\sigma)= {1 \over F(\sigma)}$$
 Then the value of $s$ corresponding to $-\sigma$ is (uniquely) determined by 
$${1\over
                   f(\sigma)}\left(2+{s\over F(\sigma)}+{F(\sigma)\over
                  s}\right) =x_1^2
$$                      
 $$
     {1\over f(\sigma)}\left(
                         2-{s \over F(\sigma)e^{-i\pi\sigma}}-{
                  F(\sigma)
 e^{-i\pi\sigma}\over s }
                         \right) =x_{\infty}^2
$$
 We conclude that $s(-\sigma)= -{1\over s(\sigma)}$.

\vskip 0.3 cm 

The last remark concerns the choice of (\ref{cucu}), (\ref{cucu1})
instead of (\ref{cucu2}), (\ref{cucu3}). The reason is that at $z=0$
the system  (\ref{novosis})
has solution  
corresponding to (\ref{ffirr1}). This is true for any
$\sigma \neq 0$ in $\Omega$, also for $\sigma \to \pm 2\mu$. This is
equivalent to the behaviour (\ref{ffirr2}), which is obtainable from
the $\hat{G}_0 = T$ of (\ref{cucu}), (\ref{cucu1}) but not of 
(\ref{cucu2}), (\ref{cucu3}). 

\rightline{$\Box$}

 \vskip 0.3 cm 
\noindent
{\it  Proof of formula (\ref{analy}) }

 We are ready to prove formula (\ref{analy}), namely: 
 $$\beta_1^2:~~(\sigma,a)\mapsto(\sigma,a e^{-2\pi i \sigma}) $$
 For $\sigma=0$ we have $x_0=0$ and $\beta_1^2:(0,x_1,x_{\infty})
 \mapsto ( 0,x_1,x_{\infty})$. Thus 
$$ 
a={x_{\infty}^2\over x_1^2+x_{\infty}^2}\mapsto{x_{\infty}^2\over
 x_1^2+x_{\infty}^2}\equiv a$$

\noindent
 For $\sigma=\pm 2 \mu +2 m$, we consider the example
 $\sigma=2 \mu+2m$, $m=0,1,2,...$. The other cases are analogous. We have $s=
 x_1^2~H(\sigma)=-x_{\infty}^2H(\sigma)e^{2\pi i \mu}$, where the
 function $H(\sigma)$ is explicitly given in theorem 2, $III$). Then
$$ 
  \beta_1:~~s= -x_{\infty}^2 H(\sigma)e^{2\pi i \mu} \mapsto -x_1^2
   H(\sigma) e^{2\pi i \mu}= -s e^{2\pi i \mu}
$$
Then 
$$ \beta_1^2:~~s\mapsto se^{4\pi i \mu}~~~ \Longrightarrow~~~ a\mapsto
ae^{-4\pi i \mu} \equiv a e^{-2\pi i \sigma}$$

\vskip 0.15 cm
\noindent 
 For  the generic case $I$) ($\sigma\not\in{\bf Z}$, $\sigma \neq \pm
 2\mu +2m$) recall that 
$$
\left\{ \matrix{
                 F(\sigma)~s+{1\over F(\sigma)~s} = x_1^2 f(\sigma)-2
                 \cr\cr
               F(\sigma)e^{-i\pi\sigma}~s+{1\over F(\sigma)e^{-i\pi
                 \sigma}~s } = 2-x_{\infty}^2 f(\sigma)
}\right.
$$
has a unique solution $s$. 
Also observe that $\beta_1: x_{\infty}\mapsto x_1$. Then 
  the transformed parameter $\beta_1:~s\mapsto s^{\beta_1}$
satisfies the equation 
$$ 
     F(\sigma)e^{-i\pi\sigma}~s^{\beta_1}+{1\over F(\sigma)e^{-i\pi
                 \sigma}~s^{\beta_1} }=  2-x_{1}^2 f(\sigma)
                $$
$$ \equiv  
                 -\left(    F(\sigma)~s+{1\over F(\sigma)~s} \right)$$
Thus $s^{\beta_1}=-e^{i\pi\sigma}~ s$. This implies 
$$
  \beta_1^2:~s\mapsto s e^{2 \pi i \sigma}~~~\Longrightarrow~~~ a \mapsto
  a~e^{-2\pi i \sigma}
$$

\rightline{$\Box$}

\vskip 0.3 cm

 We still have to prove the following 

\vskip 0.2 cm
\noindent
{\bf Proposition:} {\it Let $y(x)\sim a x^{1-\sigma}$ as $x\to 0 $ in
a domain $D(\epsilon,\sigma)$. Then, $y(x)$ coincides with $y(x;\sigma,a)$ of
theorem 1} 
\vskip 0.2 cm
\noindent
{\it  Proof:} 
 Observe that both $y(x)$ and $y(x;\sigma,a)$ have the same
asymptotic behaviour for $x\to 0$ in $D(\sigma)$. Let  
 $A_0(x)$, $A_1(x)$, $A_x(x)$ be the matrices constructed from
 $y(x)$ and $A^*_0(x)$, $A^*_1(x)$, $A^*_x(x)$ constructed from
 $yx;\sigma,a)$ by means of   
 the formulae (\ref{definite in extremis}) of
section \ref{PAPAPA} and the formulae which give $\phi_0$ in terms of $y(x)$
 in  section \ref{PAPAPAnew}.  
 It follows that $A_i(x)$ and $A^*_i(x)$, $i=0,1,x$, 
 have the same asymptotic behaviour as
 $x\to 0$. This is the
behaviour of lemma 1 of section \ref{proof of theorem 1} 
(adapted to our case). From the proof of theorem 2 if follows that $A_0(x)$,
$A_1(x)$, $A_x(x)$  and $A^*_0(x)$, $A^*_1(x)$, $A^*_x(x)$ produce   the same
 triple 
 $(x_0,x_1,x_{\infty})$. 
The solution of the Riemann-Hilbert problem for such a triple 
is unique and therefore $A_i(x)\equiv  A^*_i(x)$, $i=0,1,x$. We conclude that 
$y(x)\equiv y(x;\sigma,a)$. $\Box$


\section{ Proof of Theorem 3}\label{provaellittica}

 To start with, we
 derive the elliptic form for the general Painlev\'e 6 equation. We
follow \cite{fuchs}. 
We put
\be
u=\int_{\infty}^y ~{d\lambda \over \sqrt{\lambda(\lambda-1)(\lambda-x)}}
\label{elliptint1}
\ee
We recall that
$$
  {du \over dx}= {\partial u \over \partial y }~{dy \over dx} + {\partial u
  \over \partial x } ={1\over \sqrt{y(y-1)(y-x)} } ~{dy \over dx} + {\partial u
  \over \partial x } 
$$
from which we compute
$$
{d^2 u\over dx^2}+{2x-1\over x(x-1)} {du\over dx} +{u\over 4 x(x-1)}=$$
$$
  = {1\over  \sqrt{y(y-1)(y-x)}}\left[{d^2 y\over dx^2}+\left({1\over x}
  +{1\over x-1}+{1\over y-x}\right) ~{dy \over dx}-{1\over 2}\left({1\over y}
  +{1\over y-1}+{1\over y-x}\right)~\left({dy \over dx}\right)^2 \right]
$$
$$
+{\partial^2 u \over \partial x^2}+{2x-1\over x(x-1)}~{\partial u \over
  \partial x}+ {u\over 4x(x-1)} 
$$
By direct calculation we have:
$$
{\partial^2 u \over \partial x^2}+{2x-1\over x(x-1)}~{\partial u \over
  \partial x}+ {u\over 4x(x-1)} = -{1\over 2} {\sqrt{y(y-1)(y-x)} \over x(x-1)}
  ~{1\over (y-x)^2}
$$
 Therefore, $y(x)$ satisfies the Painlev\'e 6 equation if and only if 
\be
 {d^2 u\over dx^2}+{2x-1\over x(x-1)}~ {du\over dx} +{u\over 4
x(x-1)}={\sqrt{y(y-1)(y-x)} \over 2 x^2 (1-x)^2}~\left[2 \alpha+2\beta{x\over
y^2} +\gamma {x-1\over (y-1)^2} +\left(\delta-{1\over2}\right){x(x-1)\over
(y-x)^2} \right]
\label{ellipainleve1}
\ee

\vskip 0.2 cm
We invert the function $u=u(y)$ by observing that we are dealing with an
elliptic integral. Therefore, we write
$$ 
  y=f(u,x)$$
where $f(u,x)$ is an elliptic function of $u$. This implies that 
$$
  {\partial y\over \partial u}=\sqrt{y(y-1)(y-x)}
$$
 The above equality allows us to  rewrite (\ref{ellipainleve1}) in the
 following way: 
\be
x(1-x)~{d^2 u\over dx^2}+(1-2x)~{du\over dx}-{1\over 4}~u = {1\over 2x(1-x)}
 ~{\partial \over \partial u}\psi(u,x),
\label{ellipainleve2}
\ee
where 
$$
\psi(u,x):= 2 \alpha f(u,x)-2\beta{x\over f(u,x)} +2\gamma{1-x\over f(u,x)-1}
+(1-2\delta){x(x-1)\over f(u,x)-x}
$$

\vskip 0.2 cm
The last step concerns the form of $f(u,x)$. We observe that 
$4\lambda(\lambda-1)(\lambda-x)$ is not in Weierstrass canonical form. We change variable:
$$\lambda= t +{1+x\over 3},$$
and we get the Weierstrass form:
$$
    4\lambda(\lambda-1)(\lambda-x)= 4 t^3-g_2 t-g_3,~~~~ 
 g_2={4\over 3}(1-x+x^2),~~~~g_3:={4\over
    27}(x-2)(2x-1)(1+x)
$$
Thus
$$
  {u\over 2}=\int_{\infty}^{y-{1+x\over 3}}~{dt\over  \sqrt{4 t^3-g_2 t-g_3}} 
$$
which implies 
$$
   y(x)= {\cal P}\left({u\over 2}; \omega_1,\omega_2\right)+{1+x\over 3}
$$
 We still need to explain what are the {\it half periods} $\omega_1$,
 $\omega_2$. In order to do that, we first observe that the Weierstrass form is 
$$
   4 t^3-g_2 t-g_3= 4(t-e_1)(t-e_2)(t-e_3)
$$
where
$$ 
    e_1= {2-x\over 3},~~~e_2={2x-1\over 3},~~~e_3=-{1+x\over 3}.
$$
Therefore
$$
   g:=\sqrt{e_1-e_2}=1,~~~{\kappa}^2:={e_2-e_3\over e_1-e_3} =
   x,~~~{{\kappa}^{\prime}}^2 := 1- {\kappa}^2=1-x
$$
and the half-periods are 
$$
\omega_1= {1\over g} \int_0^1 ~{d\xi\over \sqrt{(1-\xi^2)(1-\kappa^2 \xi^2)}} 
= \int_0^1 ~{d\xi\over \sqrt{(1-\xi^2)(1-x \xi^2)}} ={\bf K}(x)
$$
$$
 \omega_2= {i\over g} \int_0^1 ~{d\xi\over
\sqrt{(1-\xi^2)(1-{\kappa^{\prime}}^2  
\xi^2)}} 
= i~\int_0^1 ~{d\xi\over \sqrt{(1-\xi^2)(1-(1-x) \xi^2)}} =i{\bf
K}^{\prime}(1-x) 
$$
The elliptic integrals ${\bf K}(x)$ and ${\bf K}^{\prime}(1-x)$ are known:
$$
{\bf K}(x)= {\pi \over 2} ~F\left({1\over 2},{1\over 2}, 1; x\right)
$$
$$
{\bf K}^{\prime}(1-x)= {\pi \over 2}F\left({1\over 2},{1\over 2}, 1;1- x\right) 
$$
where $F(x)$ is the hypergeometric function 
$$
F\left({1\over 2},{1\over 2},1;x\right)=  
\sum_{n=0}^{\infty}{ \left[\left({1 \over 2}\right)_n\right]^2
  \over (n!)^2 } x^n,
$$
${\bf K}(x)$ and ${\bf K}^{\prime}(1-x)$ 
 are two linearly independent solutions of the hypergeometric equation 
$$
 x(1-x) \omega^{\prime\prime}+(1-2x) \omega^{\prime} -{1\over 4} \omega=0.
$$
Observe that for $|\hbox{arg}(x)|<\pi$: 
$$-\pi F\left({1\over 2},{1\over 2}, 1;1- x\right) =  F\left({1\over
2},{1\over 2},1;x\right)  \ln(x) +
F_1(x)  
$$ 
where 
$$
F_1(x):=  
\sum_{n=0}^{\infty}{ \left[\left({1 \over 2}\right)_n\right]^2
  \over (n!)^2 } 2\left[ \psi(n+{1\over 2}) - \psi(n+1)\right]
x^n,~~~~\psi(z) = 
{d \over dz}\ln \Gamma(z).
$$
 Therefore $\omega_2(x)= -{i\over 2}[F(x)\ln(x)+F_1(x)]
$ where $F(x)$ is a abbreviation for $ F\left({1\over 2},{1\over 2},
1;x\right)$. The series of $F(x)$ and $F_1(x)$ 
 converge for $|x|<1$. Incidentally, we observe that 
$$
  y(x)= {\cal P}\left(u(x)/ 2; \omega_1(x),\omega_2(x)\right)-e_3= 
{1\over \hbox{sn}^2\left(u(x)/2,\kappa^2=x\right)}
$$

\vskip 0.3 cm 
\noindent
{\it Proof of theorem 3:}  
We let $x\to 0$; if $ \Im \tau>0$  and
\be 
\left|\Im \left({u\over 4\omega_1}\right)\right|<\Im \tau
\label{speriamobene}
\ee
we expand the elliptic function in Fourier series 
(\ref{rondine}).  The first condition $\Im \tau>0$ is always satisfied for
$x\to 0$ because 
$$ 
    \Im \tau(x)= -{1\over \pi} \ln|x| +{4\over \pi} \ln 2 +O(x), ~~~~x\to 0.
$$
We look for a solution $u(x)$ of (\ref{difficilissima}) of the form 
$$
  u(x)= 2 \nu_1 \omega_1(x)+2 \nu_2 \omega_2(x) + 2 v(x)
$$
where $v(x)$ has to be determined from (\ref{difficilissima}). We look for
a holomorphic solutions $v(x)$, bounded if $x\to 0$. 
 We observe that  
$$
   {u(x)\over 4 \omega_1(x)} = {\nu_1\over 2} +{\nu_2\over 2} \tau(x) +
   {v(x)\over 2\omega_1(x)} 
$$
$$
   = {\nu_1\over 2} +{\nu_2\over 2} \left[-{i\over \pi} \ln x - {i\over \pi}
 {F_1(x)\over F(x)}\right]+{v(x)\over 2\omega_1(x)},~~~~\Bigl(\hbox{note that }
 {F_1(x)\over F(x)} =-4\ln 2 +O(x) \hbox{ as } x\to 0\Bigr).
$$ 
 Thus, if   $x\to 0$ the 
 condition (\ref{speriamobene}) becomes 
\be
(2+\Re \nu_2) \ln|x|-{\cal C}(x,\nu_1,\nu_2) -8 \ln 2  < \Im \nu_2
\arg(x)< (\Re\nu_2-2) \ln|x|- {\cal C}(x,\nu_1,\nu_2)  +8 \ln 2.
\label{speriamobene1}
\ee
where $ {\cal C}(x,\nu_1,\nu_2)=
[\Im {\pi v\over 2\omega_1}+4\ln 2 \Re \nu_2 + \pi \Im
\nu_1 +O(x)]$. 
We expand the derivative of ${\cal P}$ appearing in (\ref{difficilissima})
$$
   {\partial \over \partial u} {\cal P}\left({u\over 2};
   \omega_1,\omega_2\right) = \left({\pi\over \omega_1}\right)^3
   \sum_{n=1}^{\infty} {n^2 e^{2\pi i n \tau}\over 1-e^{2\pi i n\tau}} \sin
   \left( {n
  \pi u \over 2 \omega_1}\right) - \left({\pi \over 2
   \omega_1}\right)^3{\cos\left({\pi u 
   \over 4\omega_1}\right) \over \sin^3\left({\pi u \over 4 \omega_1}\right)}
$$
$$
= {1\over 2i}  \left({\pi\over \omega_1}\right)^3\sum_{n=1}^{\infty} {
n^2 e^{2\pi i n \tau}\over 1-e^{2\pi
i n\tau}} \left(e^{in{\pi u \over 2 \omega_1}}-e^{-in{\pi u \over 2 \omega_1}}
\right)+4 i  \left({\pi \over 2
   \omega_1}\right)^3 {e^{i {\pi u \over 4 \omega_1}} + e^{-i {\pi u \over 4
   \omega_1}}\over \left( e^{i {\pi u \over 4 \omega_1}} - e^{-i {\pi u \over 4
   \omega_1}}  \right)^3 }
$$
Now we come to a crucial step in the construction: we collect $ e^{-i {\pi u \over 4
   \omega_1}} $ in the last term, which becomes
$$ 
4 i  \left({\pi \over 2
   \omega_1}\right)^3 {e^{4 \pi i { u \over 4 \omega_1}} + e^{2\pi i
 { u \over 4
   \omega_1}}\over \left( e^{2 \pi i { u \over 4 \omega_1}} - 1  \right)^3 }.
$$
The denominator {\it does not vanish} if $ \left|e^{2 \pi i { u \over 4
\omega_1}}\right|<1$. 
From now on, this condition is added to (\ref{speriamobene})
and reduces the domain (\ref{speriamobene1}). The expansion of $ {\partial
\over \partial u} {\cal P}$ becomes
$$
 {\partial \over \partial u} {\cal P}\left({u\over 2};
   \omega_1,\omega_2\right) ={1\over 2i}  \left({\pi\over
   \omega_1}\right)^3\sum_{n=1}^{\infty} {n^2 e^{i \pi   n
   \left[-\nu_1+(2-\nu_2) \tau-{v\over 2\omega_1}\right]}\over 1-e^{2\pi i n
   \tau}} \left(e^{i \pi n \left[ \nu_1+\nu_2 \tau +{v\over 2 \omega_1}
   \right]} -1\right)
$$
$$
+ 4i \left({\pi \over 2
   \omega_1}\right)^3{ e^{2\pi i \left[  \nu_1+\nu_2 \tau +{v\over 2 \omega_1}
   \right]}+ e^{\pi i  \left[  \nu_1+\nu_2 \tau +{v\over 2 \omega_1}
   \right]} \over \left(
 e^{\pi i  \left[  \nu_1+\nu_2 \tau +{v\over 2 \omega_1}
   \right]} -1\right)^3}
$$
Let's write ${F_1(x)\over F(x)}=-4\ln 2 +g(x)$, where $g(x)=O(x)$ is a power
series starting with $x$. We have 
$$e^{i \pi C\tau}= {x^C\over 16^C}
~e^{C~g(x)}= {x^C\over 16^C} (1+O(x)), ~~~x\to 0,~~
 ~~~~\hbox{ for any } C\in {\bf C}.
$$ 
Hence
$$
{\partial \over \partial u} {\cal P}\left({u\over 2};
   \omega_1,\omega_2\right)= {\cal F}\left(x, {e^{-i\pi \nu_1} \over
   16^{2-\nu_2} } 
   x^{2-\nu_2} e^{-i\pi {v\over 2\omega_1}},{e^{i\pi \nu_1}\over 16^{\nu_2}
   }x^{\nu_2} e^{i \pi {v\over 2\omega_1}} \right)
$$
where 
$$ 
  {\cal F}(x,y,z)=  {1\over 2i} \left({\pi \over 
   \omega_1(x)}\right)^3 \sum_{n=1}^{\infty} {n^2 e^{n(2-\nu_2)g(x)}\over 
                  1- \left[{1\over 16} e^{g(x)}\right]^{2n} x^{2n}} y^n(z^n-1)
   +4i\left({\pi \over 2\omega_1(x)} \right)^3 {z^2+z\over (z-1)^3}
$$
The series converges for sufficiently small $|x|$ and for $|y|<1$, $|yz|<1$;
this is precisely (\ref{speriamobene}). However, we require that the last term
is holomorphic, so we have to further impose $|z|<1$. On the resulting domain
$|x|<r<1$, $|y|<1$, $|z|<1$ ${\cal F}(x,y,z)$ is holomorphic and satisfies 
$$
  {\cal F}(0,0,0)=0. 
$$
The condition  $|y|<1$, $|z|<1$ is $\left|{e^{-i\pi \nu_1} \over
   16^{2-\nu_2} } 
   x^{2-\nu_2} e^{-i\pi {v\over 2\omega_1}}\right|<1$, $\left|,{e^{i\pi \nu_1}
\over 16^{\nu_2}
   }x^{\nu_2} e^{i \pi {v\over 2\omega_1}}\right|<1$  , namely
\be
\Re \nu_2 \ln |x| -{\cal C}(x) < \Im \nu_2 \arg(x)< (2-\Re\nu_2) \ln |x|
   -{\cal C}(x) +8 \ln 2,
\label{speriamobene2}
\ee
which is more restrictive that (\ref{speriamobene1}).

The function ${\cal F}$ can be decomposed as follows: 
$$
 {\cal F}={\cal F}\left(x, {e^{-i\pi \nu_1} \over
   16^{2-\nu_2} } 
   x^{2-\nu_2},{e^{i\pi \nu_1}\over 16^{\nu_2}
   }x^{\nu_2} \right)+
$$
$$
+\left[ {\cal F}\left(x, {e^{-i\pi \nu_1} \over
   16^{2-\nu_2} } 
   x^{2-\nu_2} e^{-i\pi {v\over 2\omega_1}},{e^{i\pi \nu_1}\over 16^{\nu_2}
   }x^{\nu_2} e^{i \pi {v\over 2\omega_1}} \right)- {\cal F}\left(x, {e^{-i\pi \nu_1} \over
   16^{2-\nu_2} } 
   x^{2-\nu_2},{e^{i\pi \nu_1}\over 16^{\nu_2}
   }x^{\nu_2} \right)\right]
$$
$$
  =: {\cal F}\left(x, {e^{-i\pi \nu_1} \over
   16^{2-\nu_2} } 
   x^{2-\nu_2},{e^{i\pi \nu_1}\over 16^{\nu_2}
   }x^{\nu_2} \right)+{\cal G} \left(x, {e^{-i\pi \nu_1} \over
   16^{2-\nu_2} } 
   x^{2-\nu_2},{e^{i\pi \nu_1}\over 16^{\nu_2}
   }x^{\nu_2},v(x) \right)
$$
The above defines ${\cal G}(x,y,z,v)$. As a function of its arguments it 
 is holomorphic for $|x|$,
$|y|$, $|z|$,  $|v|$ less then some $r^{\prime}<1$. Moreover
$$ 
  {\cal G}(0,0,0,v)={\cal G}(x,y,z,0)=0.
$$

 Let us put $u=u_0+2v$, where $u_0= 2 \nu_1
 \omega_1+2\nu_2\omega_2$. Therefore ${\cal L}(u_0)=0$ and ${\cal L}(u_0+2v)=
 {\cal L}(u_0)+{\cal L}(2v)\equiv 2{\cal L}(v)$. Hence (\ref{difficilissima})
 becomes 
\be
   {\cal L}(v)= {\alpha \over 2x(1-x)} ({\cal F}+{\cal G}).
\label{zazen}
\ee
 We put 
$$
   w:= x v^{\prime} ~~~(\hbox{where } v^{\prime}= {d v\over dx}),
$$ 
and the equation (\ref{zazen}) becomes 
$$
   w^{\prime} = {1\over x} \left[ {\alpha \over 2(1-x)^2} {\cal F} + {x
   (w+{1\over 4} v) \over 1-x} +{\alpha \over 2(1-x)^2 } {\cal G} \right]
$$
Now, let us define
$$
    \Phi(x,y,z):=  {\alpha \over 2(1-x)^2} {\cal F}(x,y,z),$$
$$
   \Psi(x,y,z,v,w):= {x
   (w+{1\over 4} v) \over 1-x} +{\alpha \over 2(1-x)^2 } {\cal G}(x,y,z,v).
$$
They are holomorphic for $|x|,|y|,|z|,|v|,|w|$ small (say less then
$r^{\prime}<1$) and they are such that
$$
   \Phi(0,0,0)=0,~~~~\Psi(0,0,0,v,w)=\Psi(x,y,z,0,0)=0.
$$

Our initial equation (\ref{difficilissima}) becomes the system 
$$ 
    x {d v \over dx} = w,
$$
$$
    x {dw\over dx}= \Phi (x, {e^{-i\pi \nu_1}\over 16^{2-\nu_2}} x^{2-\nu_2}, 
{e^{i\pi \nu_1} \over 16^{\nu_2}} x^{\nu_2} ) + \Psi(x, {e^{-i\pi \nu_1}\over 16^{2-\nu_2}} x^{2-\nu_2}, 
{e^{i\pi \nu_1} \over 16^{\nu_2}} x^{\nu_2},v(x),w(x)).
$$
A system having the structure of the system above, with some slight changes in
the arguments of $\Phi$ and $\Psi$, 
has  been studied by S. Shimomura in 
\cite{IKSY}. He reduced it to a system of integral equations 
$$
w(x)=\int_{L(x)} {1\over t} \left\{
 \Phi (t, {e^{-i\pi \nu_1}\over 16^{2-\nu_2}} t^{2-\nu_2}, 
{e^{i\pi \nu_1} \over 16^{\nu_2}} t^{\nu_2} ) + \Psi(t, {e^{-i\pi \nu_1}\over 16^{2-\nu_2}} t^{2-\nu_2}, 
{e^{i\pi \nu_1} \over 16^{\nu_2}} t^{\nu_2},v(t),w(t))\right\}~dt
$$
$$
v(x)=\int_{L(x)}{1\over s}~\int_{L(s)} {1\over t}\left\{
  \Phi (t, {e^{-i\pi \nu_1}\over 16^{2-\nu_2}} t^{2-\nu_2}, 
{e^{i\pi \nu_1} \over 16^{\nu_2}} t^{\nu_2} ) + \Psi(t, {e^{-i\pi \nu_1}\over 16^{2-\nu_2}} t^{2-\nu_2}, 
{e^{i\pi \nu_1} \over 16^{\nu_2}} t^{\nu_2},v(t),w(t))\right\}~dtds
$$
and he solved it by successive approximations, with the initial condition
$v_0=w_0=0$. The path $L(x)$ is a path connecting $x$ to 0 in
(\ref{speriamobene2}), like the path considered in the proof of theorem 1.  

We refer the reader to \cite{IKSY} and to the last of \cite{Sh}; 
for reasons of space we just take the result:

\vskip 0.2 cm 
{\it 
For any complex $\nu_1$,  $\nu_2$  such that 
$$\nu_2 \not\in (-\infty,0]\cup[2,+\infty)$$ there exists a sufficiently small
$r<1$ such that 
the system has a solution $v(x)$ holomorphic in 
$$
  {\cal D}(r;\nu_1,\nu_2):= \left\{ x\in \tilde{\bf C}_0~ \hbox{ such that }
  |x|<r, \left|{e^{-i\pi \nu_1}\over 16^{2-\nu_2}} x^{2-\nu_2} \right|<r,
\left| 
{e^{i\pi \nu_1} \over 16^{\nu_2}} x^{\nu_2}\right|<r \right\}
$$
 with an expansion convergent in ${\cal D}(r;\nu_1,\nu_2)$ 
$$
  v(x)= \sum_{n\geq 1} a_n x^n +\sum_{n\geq 0,~m\geq 1} b_{nm} x^n
  \left({e^{-i\pi \nu_1}\over 16^{2-\nu_2}} x^{2-\nu_2}\right)^m +\sum_{n\geq
  0,~m\geq 1}c_{nm} x^n \left( 
{e^{i\pi \nu_1} \over 16^{\nu_2}} x^{\nu_2}\right)^m     
$$
where $a_n$, $b_{nm}$, $c_{nm}$ are rational functions of  $\nu_2$. 
Moreover, there exists a constant $M(\nu_2)$ depending on $\nu_2$ such that 
 $v(x)\leq M(\nu_2) \left(|x|+\left|{e^{-i\pi \nu_1}\over 16^{2-\nu_2}} x^{2-\nu_2} \right|+\left| 
{e^{i\pi \nu_1} \over 16^{\nu_2}} x^{\nu_2}\right| \right)$ in  ${\cal
D}(r;\nu_1,\nu_2)$ . 
}

\vskip 0.2 cm 
We conclude that: 
\vskip 0.2 cm 
\noindent
{\bf Theorem 3:} 
{\it  for any $\nu_1$ $\nu_2$ such that 
$$ 
\nu_2\in \left({\bf C}-\left\{\right
(-\infty,0]\cup[2,+\infty)\}\right),
$$
there exists a sufficiently small $r$ such that 
$$
  y(x)= {\cal P}(\nu_1 \omega_1(x)+\nu_2 \omega_2(x)
  +v(x);\omega_1(x),\omega_2(x))
$$
in the domain ${\cal D}(r;\nu_1,\nu_2)$,  
where $v(x)$ is given above. }

\vskip 0.2 cm 
 
  ${\cal D}(r;\nu_1,\nu_2)$ is more explicitly written in the form 
 \ref{multiply defined} which makes it evident that it is contained in
 (\ref{speriamobene2}).


\chapter{ Reconstruction of 3-dimensional FM }\label{inverse reconstruction}

 Chapter \ref{rh2} was a didactic exposition of the procedure of inverse
 reconstruction of $F(t)$ through (\ref{tu}), (\ref{Fu}). Monodromy data were
 assigned, the functions $\phi_p$ were computed directly from
 (\ref{systemz1WOW2}) and the conditions of isomonodromicity
 (\ref{ISOmonodromy}). This procedure is an alternative to the direct solution
 of the boundary value problem. 

 Here we follow the same procedure, but the situation is highly
 non-trivial. The solutions of 
        $$ 
           \partial_i \phi_0=V_i \phi_0
$$
   $$ 
        \partial V_i = [V_i,V]
$$
where $
V= \phi_0~
        \hbox{diag}(\mu,0,-\mu)~\phi_0^{-1} 
$ and $  (V_k)_{ij}=([\delta_{ki}-\delta_{kj}]~V_{ij} )/( u_i-u_j)$, was partly
given in section \ref{PAPAPA} in terms of Painlev\'e trascendents; in 
section \ref{PAPAPAnew} we'll
 give the  explicit, computable 
 solution of the above system in terms of the transcendent. Its 
 dependence on the Stokes' matrix is as follows: the
Painlev\'e  transcendents are parametrized by a triple $(x_0,x_1,x_{\infty})$
of monodromy data  of a $2\times 2$ fuchsian system (see  chapter
 \ref{PaInLeVe}) and therefore $\phi_0(u)=\phi_0(u;x_0,x_1,x_{\infty})$
locally.  Since the fuchsian  system  is the $2
 \times 2$ reduction of a $3\times 3$ fuchsian system (see section
\ref{PAPAPA}) 
which is connected to
 (\ref{systemz1}) by Laplace transform (see \cite{Dub2}),  
the Stokes' matrix of the corresponding (\ref{systemz1}) is expressed in terms
 of the triple $(x_0,x_1,x_{\infty})$ itself (see \cite{DM}):
$$
S=\pmatrix{ 1 & x_{\infty}  & x_0\cr
              0 & 1           & x_1 \cr
              0 & 0           & 1  \cr }, ~~~~x_0^2+x_1^2+x_{\infty}^2
              -x_0x_1x_{\infty} = 4 \sin^2(\pi \mu).
$$

 After the determination of $\phi_0$, we compute 
 $\phi_1(u)$, $\phi_2(u)$ and $\phi_3(u)$ by direct substitution
of the solution (\ref{Y0}) 
 into the differential equation (\ref{systemz1}). In generic cases this is
enough, but for resonant values of $\mu$ the matrix $R$ in (\ref{Y0}) does
not vanish. Some entries of $\phi_p$ are indeterminate and we can fix them
thanks to the  higher order conditions of isomonodromicity
 (\ref{ISOmonodromy}) and the condition (\ref{SYMMEtria}), for $p=1,2,3$. The
non 
zero entries of $R$ appear as parameters in the $\phi_p$'s. 

The final hard problem is to obtain the closed form $F=F(t)$ from the
parametric equations (\ref{tu}), (\ref{Fu}). 

At the end of the procedure we get $F=F(t)$ in terms of the monodromy data
$S$, $\mu$, $R$.


\section{ Computation of $\phi_0$ and $V$ in terms of  Painlev\'e
transcendents}\label{PAPAPAnew}

Let $n=3$. We can bring $\eta$ to the form:
$$
  \eta=\pmatrix{ 0 & 0 & 1 \cr
                 0 & 1 & 0 \cr
                 1 & 0 & 0 \cr}
$$
 Let 
$$
  V(u)= \pmatrix{ 0 & -\Omega_3 & \Omega_2 \cr
                  \Omega_3 & 0 & -\Omega_1 \cr
                  -\Omega_2 & \Omega_1 & 0 \cr
                }
$$
which is similar to 
$$ \hat{\mu}=
\hbox{diag}(\mu,0,-\mu),~~~~\mu=\sqrt{-(\Omega_1^2+\Omega_2^2+\Omega_3^2)} 
\hbox{ constant},~~~~\mu=-{d\over 2}.  $$
By simple linear algebra we find the eigenvalues and eigenvectors of
$V$. $\phi_0$ is precisely the matrix whose columns are the eigenvectors;
imposing also the condition 
$$ \phi_0^T \phi_0 =\eta 
$$
we find
$$
   \phi_0= \pmatrix{ {i\over \sqrt{2} \mu} {\Omega_1\Omega_2 - \mu
   \Omega_3\over (\Omega_1^2+\Omega_3^2)^{1\over 2} } ~G(u) 
& 
  {\Omega_1\over i\mu} 
&
   {i\over \sqrt{2} \mu} {\Omega_1\Omega_2 + \mu
   \Omega_3\over (\Omega_1^2+\Omega_3^2)^{1\over 2} }~{1\over G(u)}
\cr
-{i\over \sqrt{2} \mu}(\Omega_1^2+\Omega_3^2)^{1\over 2} ~G(u)
& {\Omega_2\over i \mu} 
&
-{i\over \sqrt{2} \mu}(\Omega_1^2+\Omega_3^2)^{1\over 2}~{1\over G(u)}
\cr
 {i\over \sqrt{2}\mu} {\Omega_2\Omega_3 +\mu
   \Omega_1\over (\Omega_1^2+\Omega_3^2)^{1\over 2} } ~G(u)
&
{\Omega_3\over i \mu}
&
 {i\over \sqrt{2}\mu} {\Omega_2\Omega_3 -\mu
   \Omega_1\over (\Omega_1^2+\Omega_3^2)^{1\over 2} } ~{1\over G(u)}
\cr
}
$$
where $G(u)$ is so far an arbitrary function of $u=(u_1,u_2,u_3)$. To
determine it we impose the condition 
\be
 {\partial \phi_0 \over \partial u_i} = V_i(u) \phi_0
\label{equationphi0}
\ee
We observe that $\phi_{i2,0}=\Omega_i/(i\mu)$, $i=1,2,3$, 
 and that (\ref{equationphi0})
for the $\phi_{i2,0}$'s is equivalent to the equation $\partial_i V
 =[V_i,V]$. In 
 particular, the last equation implies $\sum_i \partial_i V= \sum_i u_i
 \partial_i V =0$. Thus $V(u_1,u_2,u_3)\equiv V(x)$, where 
$$ x ={u_3-u_1\over u_2-u_1} 
$$
Finally,   $\partial_i V =[V_i,V]$ becomes:
$$
{d \Omega_1\over dx} = {1 \over x} ~ \Omega_2 \Omega_3  $$
$$
{d \Omega_2\over dx} = {1 \over 1-x} ~ \Omega_1 \Omega_3 $$
\be
{d \Omega_3\over dx} = {1 \over x( x-1)} ~ \Omega_1 \Omega_2  
\label{equationV}
\ee
 
\vskip 0.2 cm
The equations (\ref{equationphi0}), (\ref{equationV}) are reduced to
$PVI_{\mu}$. The product and squares of the entries  of $\phi_0$ are 
expressed in terms of a Painlev\'e transcendent $y(x)$ in section
\ref{PAPAPA}, where we followed \cite{Dub1}. Now we work out the explicit
expressions for the entries.   Let
 $$ 
  H:=u_2-u_1.
$$
The reader may verify that the
following entries of $\phi_0$ satisfy (\ref{equationphi0}), provided that
$y=y(x)$ is a Painlev\'e transcendent and 
$$
   k=k(x,H):= {k_0~\hbox{exp}\left\{ (2\mu-1)~ \int^x d\zeta~
                                                     {y(\zeta)-\zeta\over
   \zeta(\zeta-1) } \right\} \over H^{2\mu-1}},~~~~k_0\in {\bf
                                                     C}\backslash\{0\}. 
$$
$\phi_0$ is a function of $(x,H)$: 
$$
\phi_{13,0}= i{\sqrt{k} \sqrt{y} \over \sqrt{H}  \sqrt{x}}
$$
$$
\phi_{23,0}=i {\sqrt{k} \sqrt{y-1} \over \sqrt{H}\sqrt{1-x}}
$$
$$
\phi_{33,0}=- {\sqrt{k} \sqrt{y-x} \over \sqrt{H} \sqrt{x} \sqrt{1-x}}
$$
$$
\phi_{12,0}={1\over \mu} { \sqrt{y-1}\sqrt{y-x} \over \sqrt{x}} \left[
{A\over (y-1)(y-x)} + \mu\right]
$$
$$
\phi_{22,0}={1\over \mu} { \sqrt{y}\sqrt{y-x} \over \sqrt{1-x}} \left[
{A\over y(y-x)} + \mu\right]
$$
$$
\phi_{32,0}={i\over \mu} { \sqrt{y}\sqrt{y-1} \over \sqrt{x}\sqrt{1-x}} \left[
{A\over y(y-1)} + \mu\right]
$$
$$
 \phi_{11,0}= {i\over 2 \mu^2} {\sqrt{H} \sqrt{y} \over \sqrt{k(x)} \sqrt{x} }
 \left[ A \left(B +{2\mu\over y}\right)+\mu^2(y-1-x) \right]
$$
$$
 \phi_{21,0}= {i\over 2 \mu^2} {\sqrt{H} \sqrt{y-1} \over \sqrt{k(x)} \sqrt{1-x} }
 \left[ A \left(B +{2\mu\over y-1}\right)+\mu^2(y+1-x) \right]
$$
$$
 \phi_{31,0}= 
-{1\over 2 \mu^2} {\sqrt{H} \sqrt{y-x} \over \sqrt{k(x)} \sqrt{x} 
\sqrt{1-x} }
 \left[ A \left(B +{2\mu\over y-x}\right)+\mu^2(y-1+x) \right]
$$
where  
$$
 A=A(x):= {1\over 2} \left[ {dy\over dx} x(x-1)-y(y-1)\right]
,~~~~B=B(x):= {A\over y(y-1)(y-x)}
$$
 An equivalent way to write is
$$
  \phi_0 = \pmatrix{ {E_{11}\over f} & E_{12} &  E_{13} f \cr\cr
                      { E_{21}\over f} & E_{22} &  E_{23} f \cr\cr
 {E_{31}\over f} & E_{32} &  E_{33} f \cr
}
$$
where
$$
   f=f(x,H):=
 i~ {\sqrt{k}\sqrt{y-1}\over \sqrt{H} \sqrt{1-x}}\equiv
 {\sqrt{k}\sqrt{y-1}\over \sqrt{H} \sqrt{x-1}}
$$
$$
   E_{i2}={\Omega_i\over i\mu} ,~~~~~i=1,2,3
$$
$$E_{11}= { \Omega_1 \Omega_2-\mu\Omega_3\over 2 \mu^2},~~~E_{13}=
-{\Omega_1\Omega_2 +\mu \Omega_3\over \Omega_1^2+\Omega_3^2} 
$$
$$
E_{21}=-{\Omega_1^2+\Omega_3^2 \over 2\mu^2} ,~~~~~~~~~~E_{23}=1
$$
$$
E_{31}={\Omega_2\Omega_3 +\mu\Omega_1\over 2\mu^2}
,~~~E_{33}=-{\Omega_2\Omega_3-\mu\Omega_1\over \Omega_1^2+\Omega_3^2}
$$
and 
$$
\Omega_1=i { \sqrt{y-1}\sqrt{y-x} \over \sqrt{x}} \left[
{A\over (y-1)(y-x)} + \mu\right]
$$
$$
\Omega_2=i { \sqrt{y}\sqrt{y-x} \over \sqrt{1-x}} \left[
{A\over y(y-x)} + \mu\right]
$$
$$
\Omega_3=- { \sqrt{y}\sqrt{y-1} \over \sqrt{x}\sqrt{1-x}} \left[
{A\over y(y-1)} + \mu\right]
$$

 The branches (signs) in the square roots above are arbitrary. A change of the
 sign of one root (for example of $\sqrt{H}$) implies a change of two signs in
 $(\Omega_1,\Omega_2,\Omega_3)$, or the change
 $(\phi_{i,0},\phi_{i3,0})\mapsto -(\phi_{i,0},\phi_{i3,0})$. The reader may
 verify that all these changes do not affect the equations for $\phi_0$ and
 $V$. We only remark that in the definition of $f(x,H)$ we {\it chose}
 $\sqrt{1-x}= i\sqrt{x-1}$.

\vskip 0.3 cm

Conversely, given a solution $(\Omega_1,\Omega_2,\Omega_3)$ of
(\ref{equationV}), a corresponding solution of $PVI_{\mu}$ is
\be
  y(x)= {x R(x) \over
  x~[1+R(x)]~-1},~~~~R(x):=\left({\phi_{13,0}\over\phi_{23,0}}  \right)^2= \left(\Omega_1 \Omega_2+\mu
  \Omega_3 \over \mu^2 +\Omega_2^2\right)^2
\label{yR}
\ee

\vskip 0.2 cm
 The reader may 
verify directly that the above formulae solve (\ref{equationphi0}),
  (\ref{equationV}) and satisfy the equations of section \ref{PAPAPA}.


\section{ Explicit Computation of the Flat Coordinates and of $F$ for $n=3$}

Let $t=(t^1,t^2,t^3)$ with higher indices. 
We compute the parametric form $t=t(x,H)$ and $F=F(x,H))$ using 
\be
   t^1=\sum_{i=1}^3 \phi_{i3,0}\phi_{i1,1},~~~~
 t^2=\sum_{i=1}^3 \phi_{i2,0}\phi_{i1,1},~~~~ t^3=\sum_{i=1}^3
 \phi_{i1,0}\phi_{i1,1}, 
\label{ttu}
\ee
\be
F={1\over 2} \left[t^{\alpha}t^{\beta}\sum_{i=1}^3
\phi_{i\alpha,0}\phi_{i\beta,1} -\sum_{i=1}^3
\left(\phi_{i1,1}\phi_{i1,2}+\phi_{i1,3}\phi_{i1,0}\right) \right]
\label{ffu}
\ee
The final pourpose is to obtain a closed form $F=F(t)$.  
We recall that $\mu_1=\mu$, $\mu_2=0$, $\mu_3=-\mu$. 
 Let us compute $\phi_1$, $\phi_2$, $\phi_3$. We decompose
$$
   \phi_p:=\phi_0 ~H_p,~~~p=0,1,2,...$$
namely, $H_i$ appears in the fundamental matrix 
$$
  \Xi(z,u)= (I+H_1~z+H_2~z^2+H_3~z^3+...)z^{\hat{\mu}}z^R, ~~~~z\to 0
$$
which is solution to the equation
$$
  {d \xi \over dz}= \left[{\cal U} +{\hat{\mu} \over z} \right]
  ~\xi,~~~~~{\cal U}= \phi_0^{-1}~U~\phi_0
$$
By plugging $\Xi(z,u)$ into the equation we find the $H_i$'s. We give now 
the explicit expression for the entries of the $H_i$'s. The
``generic'' expression is valid whenever $\mu$ does not have one of the
special values listed below; in the following $h_{ij}^{(k)}(u)$ are arbitrary 
functions of $u=(u_1,u_2,u_3)$ to be determined later:

\vskip 0.2 cm
\noindent
$H_1$; generic case:
$$
~H_{ij,1}={{\cal U}_{ij} \over 1+\mu_j-\mu_i},~~~R_1=0
$$
$
  \mu={1\over 2}:$
$$H_{13,1}=h_{13}^{(1)}(u),~~~~~R_{13,1}={\cal U}_{13},
$$
$$
                     H_{ij,1}={{\cal U}_{ij} \over 1+\mu_j-\mu_i}\hbox{ if }
                     (i,j)\neq(1,3)
$$
$
\mu=-{1\over 2}:$
$$H_{31,1}=h_{31}^{(1)}(u),~~~~~R_{31,1}={\cal U}_{31},
$$
$$
                     H_{ij,1}={{\cal U}_{ij} \over 1+\mu_j-\mu_i}\hbox{ if }
                     (i,j)\neq(3,1)
$$
$
\mu=1:$,
$$H_{12,1}=h_{12}^{(1)}(u),~~~~~R_{12,1}={\cal U}_{12},~~
                           $$
$$
H_{23,1}=h_{23}^{(1)}(u),~~~~~R_{23,1}={\cal U}_{23},
$$
$$                     H_{ij,1}={{\cal U}_{ij} \over 1+\mu_j-\mu_i}\hbox{ if }
                     (i,j)\not\in \{(1,2),~(2,3)\}
$$
$
\mu=-1:$,
$$H_{21,1}=h_{21}^{(1)}(u),~~~~~R_{21,1}={\cal U}_{21},~~
                           $$
$$
H_{32,1}=h_{32}^{(1)}(u),~~~~~R_{32,1}={\cal U}_{32},
$$
$$                     H_{ij,1}={{\cal U}_{ij} \over 1+\mu_j-\mu_i}\hbox{ if }
                     (i,j)\not\in \{(2,1),~(3,2)\}
$$

\vskip 0.2 cm
\noindent
$H_2$; let ${\cal U}_2:= {\cal U}~H_1 -H_1~R_1$.

\noindent
Generic case: 
$$ 
   H_{ij,2}={ {\cal U}_{ij,2}\over 2+\mu_j-\mu_i},~~~~~R_2=0$$
$
\mu=1$:
$$
     H_{13,2}=h_{13}^{(2)}(u),~~~~~R_{13,2}={\cal U}_{13,2}$$
$$
    H_{ij,2}={ {\cal U}_{ij,2}\over 2+\mu_j-\mu_i} \hbox{ if } (i,j)\neq (1,3)
$$
$\mu=-1$: 
$$
     H_{31,2}=h_{31}^{(2)}(u),~~~~~R_{31,2}={\cal U}_{31,2}$$
$$
    H_{ij,2}={ {\cal U}_{ij,2}\over 2+\mu_j-\mu_i} \hbox{ if } (i,j)\neq (3,1)
$$
$\mu=2$: 
$$
     H_{12,2}=h_{12}^{(2)}(u),~~~~~R_{12,2}={\cal U}_{12,2}$$
$$
     H_{23,2}=h_{23}^{(2)}(u),~~~~~R_{23,2}={\cal U}_{23,2}$$
$$
    H_{ij,2}={ {\cal U}_{ij,2}\over 2+\mu_j-\mu_i} \hbox{ if } (i,j)\not\in\{ 
 (1,2),~(2,3)\}
$$
   $\mu=-2$: 
$$
     H_{21,2}=h_{21}^{(2)}(u),~~~~~R_{21,2}={\cal U}_{21,2}$$
$$
     H_{23,2}=h_{32}^{(2)}(u),~~~~~R_{32,2}={\cal U}_{32,2}$$
$$
    H_{ij,2}={ {\cal U}_{ij,2}\over 2+\mu_j-\mu_i} \hbox{ if } (i,j)\not\in\{ 
 (2,1),~(3,2)\}
$$

\noindent
$H_3$; let ${\cal U}_3:= {\cal U}~H_2 -H_2~R_1-H_1~R_2$.

\noindent
Generic case: 
$$ 
   H_{ij,3}={ {\cal U}_{ij,3}\over 3+\mu_j-\mu_i},~~~~~R_3=0$$
$
\mu={3\over 2}$; 
$$
     H_{13,3}=h_{13}^{(3)}(u),~~~~~R_{13,3}={\cal U}_{13,3}$$
$$
    H_{ij,3}={ {\cal U}_{ij,3}\over 3+\mu_j-\mu_i} \hbox{ if } (i,j)\neq (1,3)
$$
$\mu=-{3\over 2}$: 
$$
     H_{31,3}=h_{31}^{(3)}(u),~~~~~R_{31,3}={\cal U}_{31,3}$$
$$
    H_{ij,3}={ {\cal U}_{ij,3}\over 3+\mu_j-\mu_i} \hbox{ if } (i,j)\neq (3,1)
$$
$\mu=3$: 
$$
     H_{12,3}=h_{12}^{(3)}(u),~~~~~R_{12,3}={\cal U}_{12,3}$$
$$
     H_{23,3}=h_{23}^{(3)}(u),~~~~~R_{23,3}={\cal U}_{23,3}$$
$$
    H_{ij,3}={ {\cal U}_{ij,3}\over 3+\mu_j-\mu_i} \hbox{ if } (i,j)\not\in\{ 
 (1,2),~(2,3)\}
$$
   $\mu=-3$: 
$$
     H_{21,3}=h_{21}^{(3)}(u),~~~~~R_{21,3}={\cal U}_{21,3}$$
$$
     H_{32,3}=h_{32}^{(3)}(u),~~~~~R_{32,3}={\cal U}_{32,3}$$
$$
    H_{ij,3}={ {\cal U}_{ij,3}\over 3+\mu_j-\mu_i} \hbox{ if } (i,j)\not\in\{ 
 (2,1),~(3,2)\}
$$


\subsection{ The generic case $\mu \neq \pm{1\over 2},~\pm 1,~\pm{3\over 2},~\pm 2,~\pm 3$}

Let $\mu \neq \pm{1\over 2},~\pm 1,~\pm{3\over 2},~\pm 2,~\pm 3$ and  let 
$$
  \phi_0 = \pmatrix{ {E_{11}\over f} & E_{12} &  E_{13} f \cr\cr
                      { E_{21}\over f} & E_{22} &  E_{23} f \cr\cr
 {E_{31}\over f} & E_{32} &  E_{33} f \cr
}
$$ 
 $E_{ij}=E_{i,j}(x)$ and $f(x,H)$ have been previously defined. A hard 
computation gives
 us the entries of $H_1$, $H_2$, $H_3$, then $\phi_1$, $\phi_2$, $\phi_3$ and
 finally $t$ and $F$ from (\ref{ttu}), (\ref{ffu}). Being the computation very
 hard and long, we  omit it and just give the result:
$$
t^1=u_1+a(x)~H$$
$$
t^2={1\over 1+\mu}~ b(x)~{H\over f(x,H)}
$$
$$
t^3= {1\over 1+2\mu}~c(x)~{H\over f(x,H)^2}  
$$
$$
F=F_0(t)+  \left[{a_1(x)~c(x)^2\over 2(1-2\mu)(3+2\mu)} +{(\mu+4)~b(x)~b_1(x)~c(x)\over
2(1-\mu)(2+\mu)(3+2\mu)} +{b(x)^2~(b_2(x)-a(x))\over (2+\mu)(3+2\mu)}\right]
~{H^3\over f(x,H)^2}
$$
$$
F_0(t):={1\over 2} t^1(t^2)^2+{1\over 2} (t^1)^2 t^3
$$
where 
$$a(x):=E_{21}E_{23}+x~E_{31}E_{33}$$
$$b(x):=E_{22}E_{21}+x~E_{32}E_{31}$$
$$b_1(x):=E_{23}E_{22}+x~E_{33}E_{32}$$
$$a_1(x):=E_{23}^2+x~E_{33}^2$$
$$b_2(x):=E_{22}^2+x~E_{32}^2$$
$$c(x):=E_{21}^2+x~E_{31}^2$$
They depend rationally on $x$, $y(x)$, ${d y(x) \over dx}$. 
Note that $F-F_0$ is independent of $u_1$, namely it is independent of
$t^1$.


\subsection{ The case of the Quantum Cohomology of projective spaces:
$\mu=-1$}

 Let $\mu=-1$. This is a non-generic case, corresponding to the Frobenius
 Manifold called the {\it Quantum Cohomology} of $CP^2$. In this case 
 the unknown functions $h_{21}^{(1)}$, $h_{32}^{(1)}$, $h_{31}^{(2)}$ have to
 be determined. 
It is known that 
$$R_1=\pmatrix{0&0&0 \cr 
               3 & 0 & 0 \cr
                0 & 3 & 0 \cr
         } ,~~~~~R_2=0$$
The direct computation gives
$$
 R_1=\pmatrix{0&0&0 \cr 
               b(x) ~H~f^{-1} & 0 & 0 \cr
                0 & b(x) ~H~f^{-1}  & 0 \cr
         },~~~~~~~ R_2=\pmatrix{0&0&0 \cr 
               0 & 0 & 0 \cr
                 b(x) ~H~f^{-1}(h_{21}^{(1)}-h_{32}^{(1)})&0  & 0 \cr
         } $$
which implies
\be
  f(x,H)= {H\over 3} ~b(x),
\label{f}
\ee
$$
h_{21}^{(1)}=h_{32}^{(1)}
$$
$h_{32}^{(1)}$ is determined  using the differential
equation :
$$
{\partial \phi_1 \over \partial u_i} = E_i \phi_0 + V_i \phi_1
$$
which implies:
$$
{\partial  h_{32}^{(1)}\over \partial u_1}= {E_{12} E_{11}\over f},
~~~~ 
{\partial  h_{32}^{(1)}\over \partial u_2}= {E_{22} E_{21}\over f},
~~~~
{\partial  h_{32}^{(1)}\over \partial u_3}= {E_{32} E_{31}\over f},
$$ 
and thus 
$$
{\partial  h_{32}^{(1)}\over \partial u_1}+{\partial  h_{32}^{(1)}\over
\partial u_2}+{\partial  h_{32}^{(1)}\over \partial u_3}=0
$$
because:
$$
E_{12} E_{11}+E_{22} E_{21}+E_{32} E_{31}=0
$$
as it follows from $\phi_0^T\phi_0=\eta$. Therefore  
  $h_{32}^{(1)}$ is a function of $x=(u_3-u_1)/(u_2-u_1)$ and
  $H=u_2-u_1$. Keeping into account (\ref{f}) and the relations:
$$
   {\partial x\over \partial u_1}={x-1\over H},~~~{\partial x\over \partial
   u_2} =-{x\over H},~~~~{\partial x\over \partial
   u_3} ={1\over H}$$
$${\partial H\over \partial u_1}=0,~~~{\partial H\over \partial u_2}=1,~~~
{\partial H\over \partial u_3}=-1,
$$
   we obtain
$$
     {\partial h_{32}^{(1)}\over \partial x}= {3\over x+
  {E_{21}E_{22}\over E_{31}E_{32}}
                                                 },
~~~~
{\partial h_{32}^{(1)}\over \partial H}=3
$$
which are integrated as follows:
\be
  h_{32}^{(1)}= 3 \ln(H) +3\int^x
 d\zeta{1\over \zeta +{E_{21}E_{22}\over E_{31}E_{32}} }.
\label{FlaTt1}
\ee

\vskip 0.2 cm
Before determining $h_{31}^{(2)}$ it is worth computing 
$t$ through \ref{ttu}.  Explaining the details is too long and tedious, so we
give just the final result:
$$
t^1=u_1+a(x)~H,$$
\be
t^2= h_{32}^{(1)},
\label{FlaTt2}
\ee
\be
t^3=-c(x)~{H\over f^2}= -9~{c(x)\over b(x)^2}~ {1\over H}
\label{FlaTt3}
\ee
We observe that  $h_{31}^{(2)}$ does not appear in $t$. We also observe that
both  $t^1$ and $t^3$ coincide with the limits for $\mu\to -1$ 
 of the same coordinates computed in the generic case. Instead, such a limit
does not exist for $t^2$.

\vskip 0.2 cm

Now we turn to the differential equation 
$$
 {\partial \phi_2 \over \partial u_i}= E_i \phi_1+ V_i \phi_2
$$
which gives three differential equations 
 for  $h_{31}^{(2)}$. Since we already know $t$, I
write the coefficients of the equation in terms of $t$:
$$
  {\partial h_{31}^{(2)}\over \partial 
u_i} = t^1~{\partial t^3 \over \partial u_i}+
t^2~{\partial t^2 \over \partial u_i}+
t^3~{\partial t^1 \over\partial  u_i},~~~~i=1,2,3.
$$
which are immediately integrated:
$$
  h_{31}^{(2)} = {1\over 2} (t^2)^2+t^1t^3
$$.

\vskip 0.2 cm
Finally, I give the result of the (hard) computation of $F$ through \ref{ffu},
 without explaining further details:
\be
   F= F_0(t)+\left[
                    {1\over 6} a_{1}(x)c(x)^2+ {3\over 4} b(x)b_1(x)c(x) + 
(b_2(x)-a(x))b(x)^2
\right]~{H^3\over f^2}
\label{formulaQH}
\ee
Remarkably, this coincides with the limit, for $\mu\to -1$, of the generic
case.


\section{ $F(t)$ in closed form }

\noindent
1)Generic case  $\mu \neq \pm{1\over 2},~\pm 1,~\pm{3\over 2},~\pm 2,~\pm 3$
\vskip 0.2 cm

 If we keep into account the dependence of $f(x,H)$ and $k(x,H)$ 
on $H$, we see that both $t$ and $F-F_0$  
 can be factorized in a part
 depending only on $x$ and another one depending only on $H$
$$
 t^2(x,H)= \tau_2(x) ~H^{1+\mu},
$$
$$
  t^3(x,H)=\tau_3(x)~H^{1+2\mu}
$$
$$
F(x,H)=F_0(t)+ {\cal F}(x)~H^{3+2\mu}
$$
where $\tau_2(x)$, $\tau_3(x)$ and ${\cal F}(x)$ are explicitly given as
rational functions of $x$, $y(x)$, ${d y(x) \over dx}$  by the
formulae of the previous sections.   
Hence the ratio
$$
{t^2\over
(t^3)^{{1+\mu\over 1+2\mu}}}
$$
is independent of $H$. This is actually the crucial point, because now 
the closed form $F=F(t)$ must be:
$$
F(t)=F_0(t)+ (t^3)^{{3+2\mu\over 1+2\mu}}\varphi\left({t^2\over
(t^3)^{{1+\mu\over 1+2\mu}}}\right)
$$
where the function $\varphi(\zeta)$ has to be determined.

\vskip 0.3 cm
\noindent
{\it Remark:} Of course, also 
   $$
F(t)=F_0(t)+ (t^2)^{{3+2\mu\over 1+\mu}}~\varphi_1\left({t^2\over
(t^3)^{{1+\mu\over 1+2\mu}}}\right)
$$
or 
$$
F(t)=F_0(t)+ {(t^2)^4\over t^3}~\varphi_2\left({t^2\over
(t^3)^{{1+\mu\over 1+2\mu}}}\right)
$$
are  okay.

\vskip 0.3 cm
\noindent
{\it Remark:} The above forms of $F$ can also be obtained 
by imposing quasi-homogeneity (namely, (\ref{WDVV3b})).

\vskip 0.3 cm

 We'll obtain closed forms $F=F(t)$ in the following way: Suppose that the
 entries of $S$ are a triple $(x_0,x_1,x_{\infty})$.  
\vskip 0.2 cm
\noindent
i) First we choose a critical point $x=0,1,\infty$ of $PVI_{\mu}$ and we
 expand  $y(x; x_0,x_1,x_{\infty})$ close to it up to any desired order: this
 is the transcendent $y(x;\sigma,a)$ of theorem 1 of chapter  \ref{PaInLeVe}. 
 The
 coefficients of the expansion, which are rationals in $a$ and $\sigma$, 
 are therefore classical functions of the triple
 (actually, $a$ and $\sigma$ 
are rational, trigonometric or $\Gamma$ functions of the monodromy data; 
we refer  to chapter \ref{PaInLeVe}  for details). 
The most efficient way to do the expansion 
 is to start by
computing the expansions of $\Omega_1(x)$, $\Omega_2(x)$, $\Omega_3(x)$ (we
have already given the explicit connection between $y(x)$ and the
$\Omega_i$'s). The algorithm used  is an expansion 
of  the $\Omega_i$'s in a small parameter [see the appendix \ref{appendiceA}]. 
It turns out that the effective variable in the expansion is a variable 
$s\to 0$ 
$$
  s:= \left\{\matrix{ 
                x \hbox{ if } x\to 0 \cr\cr
                      1-x \hbox{ if } x \to 1 \cr\cr
                       {1\over x} \hbox{ if } x \to \infty
\cr 
} \right.$$
\vskip 0.2 cm
\noindent
ii) We plug the above expansions into $\tau_i(x)$ and ${\cal F}(x)$, 
obtaining an
expansion in $s$. In particular 
$$
{t^2\over
(t^3)^{{1+\mu\over 1+2\mu}}}\equiv {\tau_2(x(s))\over
\tau_3(x(s))^{{1+\mu\over 1+2\mu}}}
$$
is expanded.

\vskip 0.2 cm
\noindent
iii) One of the following cases may occur 
  $${\tau_2\over
\tau_3^{{1+\mu\over 1+2\mu}}}\to \left\{ \matrix{
  0 \cr\cr
          \infty \cr \cr
           \zeta_0 \cr \cr
        \hbox{ no limit} \cr} \right.
\hbox{ for } s \to 0
$$
where $\zeta_0$ is a non-zero complex number. If the limit does not exist, the
problem becomes complicated. This may actually occur for particular values of
the monodromy data (we'll see later that this is the case of the Quantum
Cohomology of $CP^2$, provided that we take $e^{t^2} ~(t^3)^3$ instead of
$t^2~(t^3)^{-{1+\mu\over 1+2\mu}}$). If the limit exists, 
we have a small quantity $X=X(s)\to 0$ as $s\to 0$; 
in the three cases above $X$ is 
$$
  X:= \left\{ \matrix{
  {\tau_2\over
\tau_3^{{1+\mu\over 1+2\mu}}}\cr\cr
  \left({\tau_2\over
\tau_3^{{1+\mu\over 1+2\mu}}}\right)^{-1}\cr\cr
{\tau_2\over
\tau_3^{{1+\mu\over 1+2\mu}}}-\zeta_0  } \right.
$$

\vskip 0.2 cm
\noindent
iv) We invert the series $X=X(s)$ and find a series $s=s(X)$ for $X\to
0$. Thus we can rewrite $\tau_2=\tau_2(X)$,  $\tau_3=\tau_3(X)$, ${\cal
F}={\cal F}(X)$.

\vskip 0.2 cm
\noindent
v) We compute $H$ as a series in $X$ and as a function of $t^3$. Now $t^3$ 
becomes the variable; namely:
    $$ 
       H = H(X,t^3) = \left[{t^3 \over \tau_3(X)}\right]^{1\over 1+2\mu}
$$

\vskip 0.2 cm
\noindent
vi) By substituting $H(X,t^3)$ into $F-F_0={\cal F}(X) H^{3+2\mu}$ we obtain 
a power series for $F-F_0$ in the small variable $X$. 
 In other words, we obtain 
$\varphi(\zeta) $ as a power series in $\zeta$ or ${1\over \zeta}$ or $\zeta
-\zeta_0$.  

\vskip 0.2 cm
\noindent
vii) Finally, we simply re-express $X$ in term of the {\it variables} $t^2$
and $t^3$ and that's all. We get the closed form $F(t)$ as a power series
whose coefficients are classical functions of the monodromy data.


\section{ $F(t)$ from Algebraic Solutions of $PVI_{\mu}$}

 We refer to \cite{DM} for the algebraic solutions of $PVI_{\mu}$. 
The Stokes' matrix of the manifold is 
$$
  S=\pmatrix{ 1 & x_{\infty} & x_0 \cr
              0 & 1          & x_1 \cr
              0 & 0          &  1  \cr
}
$$
and in \cite{DM}  branches of the algebraic solutions of $PVI_{\mu}$ are
reconstructed from the above monodromy data. The construction fits into the
general framework of chapter \ref{PaInLeVe}, 
but in \cite{DM} the monodromy data
corresponding to algebraic solutions are  carefully analyzed and parametric
simple forms for the transcendents are given. In particular, the Stokes'
matrices coincide with the Stokes' matrices of the Coxeter groups $A_3$,
$B_3$, $H_3$ (or to their images with respect to the action of the braid
group).  

 Algebraic solutions  are
 defined up to the equivalence relations given by symmetries of
 $PVI_{\mu}$. There are five equivalence classes. We
 compute $F(t)$ in closed form for the representatives of classes given
 below.

\vskip 0.3 cm
\noindent
I) Tetrahedron ($A_3$), $\mu =-{1\over 4}$

\vskip 0.2 cm
i) $x\to 0$ 
$$ (x_0,x_1,x_{\infty})=(0,-1,-1)$$
$$
  y(x)= {1\over 2} x+O(x^2),~~~~x\equiv s\to 0$$
We apply the procedure above. I computed $y$ up to order $x^{15}$. The small
variable is 
$$
  X= {t^2\over (t^3)^{3\over 2}}\to 0 \hbox{ for } x\to 0
$$
and the final result is 
 $$
   F-F_0= {4\over 15} k_0^4 ~(t^3)^5-k_0^2~ (t^3)^5~X^2 +O(X^{m}),~~~~X\to 0
$$
$$
   = {4\over 15} k_0^4~ (t^3)^5-k_0^2~(t^2)^2(t^3)^2 + O\left(
\left[ {t^2\over (t^3)^{3\over 2}}\right]^{m}\right)
$$
$k_0$ is the arbitrary integration constant in $k(x,H)$. I checked the result
up to $m=16$.  Note that different
solutions $F(t)$ corresponding to different 
values of $k_0$ are connected by symmetries of the WDVV equations \cite{Dub1}.

\vskip 0.2 cm
ii) $x\to 1$ 
$$
       (x_0,x_1,x_{\infty})=(-1,0,-1)
$$
$$
  y(x(s))= 1-{1\over 2} s +O(s^2),~~~~~s=1-x\to 0
$$
$X$ is like in i) and 
$$
   F-F_0= {4\over 15} k_0^4 ~(t^3)^5+k_0^2 ~(t^3)^5~X^2 +O(X^{m}),~~~~X\to 0
$$
$$
   = {4\over 15} k_0^4~ (t^3)^5+k_0^2~(t^2)^2(t^3)^2 + O\left(
\left[ {t^2\over (t^3)^{3\over 2}}\right]^{m}\right)
$$
Here $k_0^2$ has the opposite sign w.r.t. the previous case. The result is
checked up to $m=16$.  

\vskip 0.2 cm
ii) $x\to \infty$
$$
       (x_0,x_1,x_{\infty})=(-1,-1,0)
$$
$$
  y(x(s))= {1\over 2}\left[1+O(s)\right],~~~~~s={1\over x}
  \to 0
$$
Again, $X$ is as in $i)$ and $ii)$, and the result is precisely as in $i)$.

\vskip 0.2 cm
iii) 
Now one example for different monodromy data $x_0,x_1,x_{\infty}$ and $x\to
0$  
$$(x_0,x_1,x_{\infty})=(1,1,1)$$
$$
  y(x)= {4^{2\over 3}\over 50}  ~x^{2\over 3} (1+O(x^{\delta})), ~~~~0<\delta<1,~~~~x=s\to 0
$$
determined by the formulae in chapter \ref{PaInLeVe} 
or by the explicit parametric
form for $y(x)$ in \cite{DM}. This time the computation of the expansion of
$y(x)$ is harder than before, because of the fractional exponent. I omit any
detail. The final result is:
$$
   {t^2\over(t^3)^{3\over 2}}\to \zeta_0=-{72\over 25}\sqrt{2}~k_0,~~~~~x\to 0 
$$
$$
 X= \left[{t^2\over(t^3)^{3\over 2}}-\zeta_0 \right],~~~~~x\to 0 
$$
$$
F-F_0= {119751372\over 1953125} k_0^4~ (t^3)^5-{314928 \sqrt{2}
 \over 15625} k_0^3~
(t^3)^5~X^2+ {2187\over 625} k_0^2~ (t^3)^5 ~X^4+ O(X^m),~~~~X\to 0
$$
I've checked it up to $m=14$. Substituting $X$ as a function of $t^2$ and
$t^3$ we obtain 
$$
  F-F_0= {4\over 15}\alpha^2 ~(t^3)^5 - \alpha~ (t^2)^2(t^3)^2 +O(X^m),~~~~~
\alpha=- {2187\over 625 }k_0^2.
$$

\vskip 0.3 cm
\noindent
II) Cube ($B_3$), $\mu =-{1\over 3}$.  

The computations are similar to those for the case $A_3$. I just give a few
details, namely only the case $x\to 0$ and 
$$
  (x_0,x_1,x_{\infty})=(0,-1,-\sqrt{2})
$$
$$
  y(x)= {2\over 3} x+ O(x^2),~~~~~x\to 0
$$
Now 
$$
  X={t^2\over( t^3)^2}
$$
and the final result is 
$$
 F-F_0= {512\over 8505} k_0^6~(t^3)^7 -{16\over 27} k_0^3 ~(t^2)^2 (t^3)^3-{2i 
\sqrt{2}\over 9} k_0^{3\over 2}~ (t^2)^3 t^3 + O(X^m)
$$
for any $m$ I have checked.

\vskip 0.3 cm
\noindent
III) Icosahedron ($H_3$), $\mu=-{2\over 5}$. We give one example:
$$
  (x_0,x_1,x_{\infty})=(0,1,{1+\sqrt{5}\over 2})$$
$$
  y(x)= {3+\sqrt{5}\over 5+\sqrt{5}}x+ O(x^2), ~~~~~x\to 0
$$
Now
$$ 
  X={t^2\over (t^3)^3} 
$$
The final result is
$$
 F-F_0= {18\over 55} \alpha^4~(t^3)^{11} +{ 9\over 5}\alpha^2~(t^2)^2(t^3)^5+ 
      \alpha ~(t^2)^3(t^3)^2+ O(X^m)
$$
where 
$
\alpha= -512\sqrt{5}/[3 
(\sqrt{5}-5)^{5\over 2}(\sqrt{5}+5)^{5\over 2}]~k_0^{5\over 2}$.

\vskip 0.2 cm
\noindent
{\it Remark:} In I), II), III) we have recovered the polynomial solutions
(\ref{A3}), (\ref{B3}), (\ref{H3}).

\vskip 0.3 cm
\noindent
IV) Great dodecahedron ($H_3$), $\mu= -{1\over 3}$. Just the final result,
which is a power series:
$$F-F_0=(t^3)^7 \left[A_0+\sum_{k=2}^{\infty}~A_k \left({t^2\over
(t^3)^2}\right)^k\right] $$
$$
=
{512\over 8505}k_0^6~(t^3)^7-{16\over 27} k_0^3~(t^2)^2(t^3)^3+{4 i \sqrt{5}
\over 9} k_0^{3\over 2} ~t^3(t^2)^3 + {1\over 8} ~{(t^2)^4\over t^3} +
{i \sqrt{5} \over 80 k_0^{3\over 2}} ~ {(t^2)^5\over (t^3)^3} -{1\over 64
k_0^3} ~{(t^2)^6\over (t^3)^5} +....
$$
 the series can be computed up to any order in $X=t^2/(t^3)^2\to 0$. We
 obtained it for $x\to 0$ and $(x_0,x_1,x_{\infty})=(0,{1+\sqrt{5}\over
 2},{\sqrt{5}-1\over 2})$, for $x\to 1$ and  $(x_0,x_1,x_{\infty})=({1+\sqrt{5}\over
 2}, 0, {\sqrt{5}-1\over 2})$, for $x\to \infty $ and
   $(x_0,x_1,x_{\infty})=({1+\sqrt{5}\over 2},  {\sqrt{5}-1\over 2},0)$.

\vskip 0.3 cm
\noindent
V) Great icosahedron ($H_3$), $\mu=-{1\over 5} $. The final result is a power
series:
$$
 F-F_0 =(t^3)^{13\over 3} \left[A_0+\sum_{k=2}^{\infty}~A_{k}\left({t^2\over
(t^3)^{4\over 3}}\right)^k\right]
$$
$$
=
{54\over 284375} \alpha^4 ~(t^3)^{13\over 3} -{3\over 125} \alpha^2 ~ (t^2)^2
(t^3)^{5\over 3} +{i\over 15} \alpha~(t^2)^3 (t^3)^{1\over 3} +
{1\over 72} ~{(t^2)^4\over t^3} +{i\over 108\alpha} ~{(t^2)^5 \over
(t^3)^{7\over 3}}+....
$$
for $X= t^2/(t^3)^{4\over 3}\to 0$. Here $k_0= 270^{1\over 5} / 30~\alpha^{6\over
5}$. We checked this expansion for $x\to 0$ and
$(x_0,x_1,x_{\infty})=(0,-1,{1-\sqrt{5}\over2})$, $x\to \infty$ and  
 $(x_0,x_1,x_{\infty})=(-1,{1-\sqrt{5}\over2},0)$. 
\vskip 0.3 cm


\section{ Closed form for $QH^{*}(CP^2)$} 

 In this case the factorization is $t^3=\tau_3(x) H^{-1} $ and $F-F_0={\cal
F}(x) H$ , but 
$$t^2= h_{32}^{(1)}= \ln(H^3) +\int^x(...) $$
This implies that 
$$
   e^{t^2} (t^3)^3 
$$
is independent of $H$ and that 
$$
 F(t)=F_0(t)+ {1\over t^3} ~\varphi\left( e^{t^2} (t^3)^3 \right)
$$
or
$$
 F(t)=F_0(t)+e^{t^2\over 3}~ \varphi_1\left( e^{t^2} (t^3)^3 \right)
$$
The situation is more complicated now, 
because the  behaviour of $y(x)$ for
$x\to 0$  is not like $a~x^{1-\sigma}(1~+$ higher order terms)
 as in the previous section. The same holds for $x\to 1$ and $x\to
\infty$.  The reason for this is that 
the monodromy data are in the
orbit w.r.t  the action of the braid group of the triple
$(x_0,x_1,x_{\infty})= (3,3,3)$. Hence, the monodromy data are real and their
absolute value is greater than 2, so $\Re \sigma^{(i)}=1$. For the data
$(3,3,3)$ we have $\sigma^{(i)}=1-i\nu$, $\nu = -{2\over \pi}
\ln\left({3+\sqrt{5} \over 2}\right)$. This case corresponds to an oscillatory
transcendent   if $x$ converges to the 
critical points along radial directions.  Moreover, we do not even know the
behaviour of the transcendent and  if it has poles for some values of
$\arg(x)$. 
 We recall that the effective
parameter is $s=x$ or $1-x$ or ${1\over x}$. It turns out that
$e^{t^2}~(t^3)^3 $ has no limit as $s\to 0$.  

\vskip 0.3 cm
 In the following, we compute $F(t)$ in closed form starting from the
expansion of $y(x)$ and the $\Omega_i(x)$'s close to the non-singular point 
$x_c= \exp\{-i{\pi\over 3}\}$ and we obtain the expansion of $F(t)$ due to
Kontsevich
\be 
  F(t)=F_0(t)+ {1\over t^3}~\sum_{k=1}^{\infty} ~{N_k\over (3k-1)!}
  ~\left[(t^3)^{3} e^{t^2}\right]^k.
\label{DiFranc}
\ee

\vskip 0.2 cm
 Our result is interesting because it allows us to obtain such a relevant
 expression starting from the isomonodromy deformation theory applied to
Frobenius manifolds.

On the other hand, it  is not completely satisfactory. Since the Frobenius
manifold can in principle be reconstructed from its monodromy data, 
we should be able 
to express the coefficients $N_k$  as functions of the monodromy data. 
The critical behaviour of $y(x)$ close to a critical point
depends upon two parameters which are classical functions of the monodromy
data (actually, they depend on $(x_0,x_1,x_{\infty})$ through algebraic
operations and trigonometric and $\Gamma$ functions). Thus, we
should  compute  $F(t)$ from the local behaviour of
$y(x)$ close to a critical point. The choice of a non-singular point $x_c$ is
not satisfactory because the expansion of
$y(x)$ close to $x_c$ depends  
on two parameters (initial data $y(x_c)$, $y^{\prime}(x_c)$),
which in general are not classical (known) functions of the two parameters
upon which the critical behaviour close to a critical point depends. Thus,
they are not classical functions of the monodromy data. This is
due to the  fact that in
 general the Painlev\'e transcendents are not classical functions.

\vskip 0.3 cm

 We now compute $F(t)$ in closed form. We expand $y(x)$ close to 
$$ 
  x_c = e^{-i{\pi\over 3}}
$$

 This choice comes from the knowledge of the structure of $QH^{*}(CP^2)$ at
 the ``classical'' point $t^1=t^3=0$ \cite{Dub2} (se also \cite{guz}). 
Namely, we know that 
$$
  {\cal U} = \pmatrix{  0 & 0 & 3 q \cr
                        3 & 0 & 0    \cr
                        0 & 3 & 0    \cr},~~~~~q:=e^{t^2}
$$
with eigenvalues
 $$ 
     u_1= 3q^{1\over 3} ,~~~u_2= 3q^{1\over 3}e^{-i{2\pi\over 3}},~~~
                     u_3= 3q^{1\over 3}e^{i{2\pi\over 3}}.
 $$
The matrix $\phi_0$ is 
$$
  \phi_0= \pmatrix{ q^{-{1\over 3}}  &  1   &  q^{1\over 3} \cr
       q^{-{1\over 3}} e^{-i{\pi\over 3}} & -1 & q^{1\over 3}e^{i{\pi \over
                   3}}\cr 
                     q^{-{1\over 3}} e^{i{\pi\over 3}} & -1 & q^{1\over
                   3}e^{-i{\pi \over 
                   3}}\cr}
$$
Thus
$$
  x_c= {u_3-u_1\over u_2-u_1} =e^{-i{\pi\over 3}} 
$$

\vskip 0.2 cm
\noindent
{\it Remark:} $x_c$ 
lies at the intersection of the ``spherical'' neighourhhods of the
critical points  $x=0$, $x=1$,
$x=\infty$, defined by the condition that each neighbourhood does not contain
another critical point in its interior. Also the point $\bar{x}_c=e^{i{\pi
\over 3}} $ is at the other intersection. Any
permutation of $u_1$, $u_2$, $u_3$ yields  $x_c$ or $\bar{x}_c$. The analysis
which follows can be repeated at the point $\bar{x}_c$.
 
\vskip 0.2 cm
\noindent
{\it Remark:} 
We understand that the choice of $x_c$ is not satisfactory, because we have 
 to rely on the knowledge of $QH^{*}(CP^2)$ at the point $t^1=t^3=0$, and not
 only on the monodromy data (which actually were computed in chapter
\ref{stokes} 
 from this knowledge itself!).

\vskip 0.2 cm

 In order to compute $y(x)$ we start from $\Omega_1(x)$, $\Omega_2(x)$,  
$\Omega_3(x)$. We look for a regular expansion
$$
   \Omega_i(x)= \sum_{k=0}^{\infty} ~\Omega_i^{(k)} ~(x-x_c)^k,~~~~i=1,2,3.
$$
We need the initial conditions $\Omega_i^{(0)}$. We can compute them using 
$$ 
   \Omega_i = i\mu~\phi_{i2,0}
$$
 which implies 
$$
\Omega_1^{(0)}= -{i\over \sqrt{3}},~~~
 \Omega_2^{(0)}= {i\over \sqrt{3}}, ~~~\Omega_3^{(0)}= {i\over \sqrt{3}},
$$ 
 Then we plug the expansion into (\ref{equationV}) and we compute the
 coefficients at any desired order. Finally, we obtain $y(x)$ from
 (\ref{yR}). We skip  the details and we just give the first terms of the
 expansions, the ``effective'' small variable  being $s:=x-x_c \to 0$:
$$
 \Omega_1= -{i\sqrt{3}\over 3}-\left({1\over 6} +{1\over 6}i\sqrt{3}\right) s+
 {1 \over 9}i\sqrt{3}s^2+\left({1\over 18} -{1 \over
 18}i\sqrt{3}\right)s^3 
 -\left({5\over 36} -{5\over 108} i\sqrt{3}\right) s^4+...
$$ 
$$
\Omega_2={i\sqrt{3}\over 3}+\left({1\over 6} -{1\over 6}i\sqrt{3}\right)
s- {1 \over 9}i\sqrt{3}s^2-\left({1\over 18} +{1 \over
18}i\sqrt{3}\right)s^3
-\left({5\over 36} 
+{5\over 108} i\sqrt{3}\right)
s^4+...$$
$$
\Omega_3= {i\sqrt{3}\over 3}-{1\over 3}
s+{2 \over 9}i\sqrt{3}s^2+{4\over 9}s^3-{13\over 54 }i\sqrt{3} 
s^4+...
$$
$$
y(x)= {1\over 2} - {1\over 6} i\sqrt{3}+{1\over 3}s -{1\over 3} i\sqrt{3} s^2
-{1\over 3}s^3+{1\over
9}i\sqrt{3}s^4
+{13\over 45}s^5 -{37\over 135} i \sqrt{3}s^6
-{17\over 27} s^7 +...
$$
 Once we have $\Omega_1$, $\Omega_2$, $\Omega_3$, we can compute the
 $E_{ij}$'s and finally the flat coordinates $t(x,H)$ and $F(x,H)$. At low 
 orders: 
$$
 t^1=u_1+\left[{1\over 2} -{1\over 6} i\sqrt{3} +{1\over 3} s
 +O(s^2) \right]~H
$$
$$
  t^3= \left[-9s+O(s^2)\right]~H^{-1}
$$
$$
 q= \exp(t^2)= {1\over 143} i \sqrt{3}~ q_0 
~\left[ 1+i\sqrt{3}s+O(s^2)\right]~H^3
$$
where $q_0$ is an arbitrary integration constant (recall that $t^2$ is obtained
by integration). 
$$
 F= {1\over 6} i \sqrt{3} s^2 + {1\over 6} s^3 -{1\over 18}  i
 \sqrt{3} s^4 +O(s^5)
$$
The following quantity is independent of $H$ 
$$
  X:=t^3 q^{1\over 3} $$
For example, if we take the cubic root 
$\left(-{1\over 6} +{1\over 18} i \sqrt{3}\right)~q_0^{1\over 3}$ 
 of $ {1\over 143} i \sqrt{3}~ q_0$ we compute 
$$X=q_0^{1\over 3}~\left[ \left({3\over 2} - 
{1\over 2} i \sqrt{3}\right)s 
                      -\left( {1\over 2} +{1\over 2} i \sqrt{3} \right)
s^2 +O(s^3)\right]
$$
Another choice of the cubic root does not affect the final result (actually,
we will see that  $F-F_0$ is a series in $X^3$).  
$X$ is the small parameter that tends to 0 as $s\to 0$. 
We invert the series and find $s=s(X)$, and then we find
$H=\tau_3(X)/t^3$ as a series in $X\to 0$. 
Finally, the non cubic term in $F$ is computed:
$$
  F-F_0= {1\over t^3} \left[{1\over 2~q_0} X^3 + {1\over 120~q_0^2} X^6 +
  {1\over 3360~q_0^3 }
  X^9 +{31 \over 1995840~q_0^4} X^{12} + {1559\over 1556755200~q_0^5 } X^{15}
  +O(X^{17})\right] 
$$
We obtained this expansion through the expansions of the  $\Omega_i$'s and of
$y(x)$ at order 16. 
If we put $q_0=1$, this is exactly Kontsevich's solution, with 
$$
 N_1=1,~~~N_2=1,~~~N_3= 12,~~~N_4= 620,~~~N_5=87304.
$$

 \vskip 0.2 cm
 Though not completely satisfactory to our theoretical purposes, the above is
 a  procedure to compute
 Gromov-Witten invariants, which is alternative to the usual procedure
 consisting in  the direct substitution of
 the expansion (\ref{DiFranc}) in the WDVV equations. 

\vskip 0.3 cm
\noindent
{\bf Remark:} We observe that for the very special case of $QH^{*}(CP^2)$
$$
      {\partial^2 \over \partial (t^2)^2} ~(F-F_0) = u_1+u_2+u_3 - 3 t^1
                   $$
This follows from the computation of the intersection form of the Frobenius
manifold $QH^{*}(CP^2)$ in terms of $F$ and by recalling that its eigenvalues
are $u_1$, $u_2$, $u_3$ (see section \ref{Nature of the singular
point}). 
 Therefore, 
$$
      {\partial^2 \over \partial (t^2)^2} ~(F-F_0)= u_1+u_2+u_3-3u_1-3~a(x)~H
$$
$$
   =  H(1+x-3a(x)),~~~~x={u_3-u_1\over u_2-u_1} ,~~~H=u_2-u_1.
$$
 The above formula allows to compute $F-F_0$ faster than (\ref{formulaQH}).

\vskip 0.3 cm 

 The formulae (\ref{FlaTt1}), (\ref{FlaTt2}), (\ref{FlaTt3}) and
 (\ref{formulaQH}) are  completely explicit as rational functions of
 $y(x)$ and ${dy\over dx}$. Of $y(x)$ we partly know the behaviour close to a singular point
 (chapter \ref{PaInLeVe}), but the information we have on it when $x$ tends to
 the point along a radial path is not complete. We know that 
$$
   y(x)= {\cal P}(\nu_1 \omega_1(x)+\nu_2 \omega_2(x)
   +v(x);\omega_1(x),\omega_2(x))  
$$
and we know the form of $v(x)$ for a limited domain of $x$ which include
radial paths with some limitations on $\arg(x)$ (namely $\arg(x)$ must be
greater or less of some angle if $x\to 0$).
 Therefore, for a limited range of $\arg(x)$ we may write the asymptotic
expansion for the parametric solution (\ref{Fu}) (\ref{tu}) of the WDVV eqs. 
 Still, the problem of inversion to obtain a closed form $F(t)$ is very hard.


\appendix
\chapter{}\label{appendiceA}
  
We present a procedure to compute the expansion of the painlev\'e
transcendents of $PVI_{\mu}$ and of the solutions of 
$$
{d \Omega_1\over ds} = {1 \over s} ~ \Omega_2 \Omega_3  $$
$$
{d \Omega_2\over ds} = {1 \over 1-s} ~ \Omega_1 \Omega_3 $$
\be
{d \Omega_3\over ds} = {1 \over s( s-1)} ~ \Omega_1 \Omega_2  
\label{system}
\ee
close to the critical point $s=0$. Here we use the notation $s$ instead of
$x$. 
 The system (\ref{system}) determines the matrix $V(u_1,u_2,u_3) \equiv
V(s)$ 
\be
   V(s)= \pmatrix{  0  &   -\Omega_3     &    \Omega_2  \cr
                    \Omega_3  &  0   &  -\Omega_1  \cr
                   -\Omega_2  & \Omega_1  &  0   \cr
                  }
, ~~~~~~~s={u_3-u_1\over u_2-u_1}.
\label{V}
\ee
The corresponding solution of the $PVI_{\mu}$ equation is 
\be 
  y(s)= { - s R(s) \over 1 - s (1 +R(s)},~~~~ R(s) := \left[ {\Omega_1
  \Omega_2 + \mu \Omega_3 \over \mu^2 + 
  \Omega_2^2} \right]^2 
\label{PVI}
\ee

\section{ Expansion with respect to a Small Parameter}

We want to study the behaviour of the solution of (\ref{system}) for
$s \to 0$. Let 
               $$s:=\epsilon ~z$$ 
where $\epsilon$ is the small
parameter. The system (\ref{system}) becomes: 
$$
{d \Omega_1\over dz} = {1 \over z} ~ \Omega_2 \Omega_3  $$
$$
{d \Omega_2\over dz} = {\epsilon \over 1-\epsilon z} ~ \Omega_1 \Omega_3 $$
\be
{d \Omega_3\over dz} = {1 \over z( \epsilon z -1)} ~ \Omega_1 \Omega_2  
\ee
The coefficient of the new system are holomorphic for $\epsilon \in E:=
\{ \epsilon\in {\bf C}~|~|\epsilon| \leq \epsilon_0 \}$ and for $0<|z|<{1\over
|\epsilon_0|}$, in  
particular for $z\in D:=\{z\in {\bf C} |~R_1\leq |z| \leq R_2 \}$,
 where $R_1$ and $R_2$ are independent of $\epsilon$  and satisfy 
$0<R_1<R_2<{1\over \epsilon_0}$. 

 We will use the small parameter expansion as a formal way to compute the
 expansions of the $\Omega_j's$  for $s\to 0$, the only justification being
 that in the cases we apply it we find expansions in $s$ which we already know
 they are convergent.   To our knowledge, there is no rigorous 
  justification of the (uniform)  
convergence of the expansions for the $\Omega_j$'s in terms of the variable 
$s$  restored
 after the small parameter expansion in powers of $\epsilon$. 

 For $\epsilon \in E$ and $z\in D$ we can expand the fractions as
 follows: 
$$
{d \Omega_1\over dz} = {1 \over z} ~ \Omega_2 \Omega_3  $$
$$
{d \Omega_2\over dz} = \epsilon ~\sum_{n=0}^{\infty} z^n \epsilon^n 
 ~ \Omega_1 \Omega_3 $$
\be
{d \Omega_3\over dz} = -{1 \over z}~\sum_{n=0}^{\infty} z^n \epsilon^n
 ~ \Omega_1 \Omega_2  
\label{systpower}
\ee
and we look for a solution  expanded in powers of $\epsilon$:
\be
   \Omega_j(z,\epsilon)=\sum_{n=0}^{\infty}
   \Omega_j^{(n)}(z)~\epsilon^n~,~~~~
~~j=1,2,3.
\label{solution}
\ee

\vskip 0.2 cm
  We find the $\Omega_j^{(n)}$'s substituting (\ref{solution}) into
  (\ref{systpower}). Et order $\epsilon^0$ we find
$$ 
{\Omega_2^{(0)}}^{\prime}=0 \Longrightarrow  \Omega_2^{(0)}={i \sigma
\over 2}
$$
$$ 
{\Omega_1^{(0)}}^{\prime}={1\over z}~ \Omega_2^{(0)}\Omega_3^{(0)} 
$$
$$ 
 {\Omega_3^{(0)}}^{\prime}=-{1\over z}~ \Omega_2^{(0)}\Omega_1^{(0)} 
$$
where $\sigma$ is so far an arbitrary constant, and the prime denotes
the derivative w.r.t. $z$.  Then we solve the linear system for
$\Omega_1^{(0)}$ and $\Omega_3^{(0)}$ and find
$$ 
   \Omega_1^{(0)}= \tilde{b}~ z^{-{\sigma \over 2}}+\tilde{a} 
~ z^{\sigma\over 2}$$
$$
   \Omega_3^{(0)}= i\tilde{b}~z^{-{\sigma \over 2}}-i\tilde{a}
~ z^{\sigma\over 2}$$
where $\tilde{a}$ and $\tilde{b}$ are integration constants. 
The higher orders are
$$
\Omega_2^{(n)}(z)=\int^z d \zeta ~
\sum_{k=0}^{n-1}~\zeta^k\sum_{l=0}^{n-1-k}~\Omega_1^{(l)}(\zeta)
~\Omega_3^{(n-1-k-l)} (\zeta)
$$
$$ {\Omega_1^{(n)}}^{\prime}={1\over z}~ \Omega_2^{(0)}\Omega_3^{(n)}+
A_1^{(n)}(z)$$ 
$$
 {\Omega_3^{(n)}}^{\prime}=-{1\over z}~
\Omega_2^{(0)}\Omega_1^{(n)}+A_3^{(n)}(z)
$$
where 
$$
A_1^{(n)}(z)= {1\over z} \sum_{k=1}^n ~\Omega_2^{(k)}(z)~\Omega_3^{(n-k)}(z)
$$
$$
A_3(z)= -{1\over z} \left[
                          \sum_{k=1}^n~\Omega_2^{(k)}(z)~ \Omega_1^{(n-k)}(z) +
\sum_{k=1}^n ~z^k\sum_{l=0}^{n-k} ~\Omega_1^{(l)}(z)~ \Omega_2^{(n-k-l)}(z)
                          \right] 
$$
The system for $\Omega_1^{(n)}$, $\Omega_3^{(n)}$ is closed and
non-homogeneous. By variation of parameters we find the particular
solution 
$$
\Omega_1^{(n)}(z)= {z^{\sigma/2} \over \sigma} ~\int^{z} d \zeta~
\zeta^{1-{\sigma\over 2}} R_1^{(n)}(\zeta)~-{z^{-{\sigma/2}} \over
\sigma} ~\int^z \zeta^{1+{\sigma\over 2}} R_1^{(n)}(\zeta)
$$
$$
\Omega_3^{(n)}(z)= {z \over i \sigma /2} \left(
{\Omega_1^{(n)}}(z)^{\prime}-A_1^{(n)}(z)   \right)
$$
where
$$ R_1^{(1)}(z)= {1\over z} A_1^{(n)}(z) + {i \sigma \over 2 z} 
A_3^{(n)}(z) + {A_1^{(n)}}(z)^{\prime}
$$

\vskip 0.3 cm 

Result: 
$$
\Omega_j(s)= s^{-{\sigma \over 2}} \sum_{k,~q=0}^{\infty} b_{kq}^{(j)} 
~s^{k +(1-\sigma)q} + s^{{\sigma \over 2}} \sum_{k,~q=0}^{\infty} a_{kq}^{(j)} 
~s^{k +(1+\sigma)q}~,~~~~~~~j=1,3
$$
\be
\Omega_2=  \sum_{k,~q=0}^{\infty} b_{kq}^{(2)} 
~s^{k +(1-\sigma)q} + \sum_{k,~q=0}^{\infty} a_{kq}^{(2)} 
~s^{k +(1+\sigma)q}
\label{smallserie}
\ee
The coefficients $a_{kq}^{(j)}$ and $b_{kq}^{(j)}$ contain
$\epsilon$. 
In fact, they are functions of $a:=\tilde{a} 
\epsilon^{-{\sigma\over 2}}$, $b:=\tilde{b} \epsilon^{{\sigma\over
2}}$.

\vskip 0.3 cm
\section{ Solution by Formal Computation}

 Consider the system (\ref{system}) and expand the fractions as $s\to
 0$. We find
$$
{d \Omega_1\over ds} = {1\over s} ~ \Omega_2 \Omega_3  $$
$$
{d \Omega_2\over ds} = \sum_{n=0}^{\infty} s^n ~ \Omega_1 \Omega_3 $$
\be
{d \Omega_3\over ds} = - {1 \over s}\sum_{n=0}^{\infty} s^n
                       ~ \Omega_1 \Omega_2  
\ee
  We can look for a solution written in formal series 
$$
\Omega_j(s)= s^{-{\sigma \over 2}} \sum_{k,~q=0}^{\infty} b_{kq}^{(j)} 
~s^{k +(1-\sigma)q} + s^{{\sigma \over 2}} \sum_{k,~q=0}^{\infty} a_{kq}^{(j)} 
~s^{k +(1+\sigma)q}~,~~~~~~~j=1,3
$$
$$
\Omega_2=  \sum_{k,~q=0}^{\infty} b_{kq}^{(2)} 
~s^{k +(1-\sigma)q} + \sum_{k,~q=0}^{\infty} a_{kq}^{(2)} 
~s^{k +(1+\sigma)q}
$$
 Plugging the series into the equation we find solvable relations
 between the coefficients and we can determine them. For example, the
 first relations  give 
$$
  \Omega_2 = {i \sigma \over 2} + \left({i [b_{00}^{(1)}]^2 \over 1 -\sigma}
  s^{1-\sigma}+...\right) - \left({i [a_{00}^{(1)}]^2 \over 1 +\sigma}
  s^{1+\sigma}+...\right)$$
$$
\Omega_1= (b_{00}^{(1)}s^{-{\sigma\over 2}}+...)+ (a_{00}^{(1)}
s^{\sigma \over 2} +...)$$
$$
\Omega_3= (ib_{00}^{(1)}s^{-{\sigma\over 2}}+...)+ (-ia_{00}^{(1)}
s^{\sigma \over 2} +...)$$
All the  coefficients determined by successive relations  
are functions of $\sigma$, $b_{00}^{(1)}$, $
a_{00}^{(1)}$. These are the three parameters on which the solution of
 (\ref{system}) must depend.  We can identify $b_{00}^{(1)}$ with $b$ and
$a_{00}^{(1)}$ with $a$.

 
\vskip 0.3 cm 
\section{ The Range of $\sigma$}

 The above computations make sense if $\sigma$ is not an {\it odd} 
integer, otherwise some 
 coefficients of the expansions for the $\Omega_j$'s diverge (see for
 example the first terms of $\Omega_2$ in the preceding section). 

Moreover, the expansion in the small parameter yields the following
approximation at order 0 for $\Omega_2$:
$$
  \Omega_2 \approx {i \sigma \over 2} \equiv \hbox{ constant}
$$
The approximation at order 1 contains powers $z^{1-\sigma}$,
$z^{1+\sigma}$. If we assume that the  approximation at order 0 in
$\epsilon$ is actually the limit of $\Omega_2$ as $s=\epsilon~z \to 0$,
than we need 
$$ 
  -1 < \Re \sigma <1
$$
Of course, this makes sense if $s\to 0$ along a radial path (i.e. within
a sector of amplitude less than $2 \pi$). 

The ordering of the expansion (\ref{smallserie}) is
 somehow conventional: namely, we could transfer some terms multiplied
 by $s^{{\sigma\over 2}}$ in the series multiplied by
 $s^{-{\sigma\over 2}}$, and conversely.  I report
 the first terms: 
$$
\Omega_1(s) = b s^{-{\sigma\over 2}} \left(1 -{b^2\over
(1-\sigma)^2}s^{1-\sigma}+{\sigma^2\over 4(1-\sigma)}s+ {a^2\over
(1+\sigma)^2}s^{1+\sigma} 
+...\right)$$
$$
+as^{\sigma\over 2}\left(1+{b^2\over
(1-\sigma)^2}s^{1-\sigma}+{\sigma^2\over 4(1+\sigma)}s-  {a^2\over
(1+\sigma)^2}s^{1+\sigma} 
+... \right)
$$
\vskip 0.2 cm 
$$
\Omega_3(s)=ibs^{-{\sigma\over 2}}\left(
           1-{b^2\over (1-\sigma)^2}s^{1-\sigma}+{\sigma(\sigma-2) \over
           4(1-\sigma)} s +{a^2\over (1+\sigma)^2}s^{1+\sigma} +...
\right)$$
$$
-ias^{\sigma\over 2} \left( 
           1+{b^2\over (1-\sigma)^2} s^{1-\sigma} +{\sigma (\sigma+2) \over
           4(1+\sigma)} s - {a^2 \over (1+\sigma)^2} s^{1+\sigma} +...
\right)
$$
\vskip 0.2 cm 
$$
\Omega_2(s)= i~{\sigma\over 2}+i {b^2\over (1-\sigma)^2} s^{1-\sigma}
  -i {a^2\over 1+\sigma} s^{1+\sigma} +...
$$
Note that the dots do not mean higher order terms. There may be 
 terms bigger than
 those written above (which are computed through the expansion in the small
 parameter up to order $\epsilon$) depending on the value of $\Re \sigma $
 in $(-1,1)$. 

 \vskip 0.2 cm
Finally, we note that we  can always assume:
$$
   0\leq \Re \sigma <1, 
$$ 
because that would not affect the expansion of the solutions but for
the change of two signs. With this in mind, 
the expansions above are:
$$ 
\Omega_1= b s^{-{\sigma \over 2}}( 1 +O(s^{1-\sigma})) +
                    a s^{\sigma \over 2} (1+ O(s))
$$
$$
\Omega_3= ib s^{-{\sigma \over 2}}( 1 +O(s^{1-\sigma})) 
           - i a s^{\sigma \over 2} (1+ O(s))$$
$$ 
\Omega_2= { i \sigma\over 2}(1 +O(s^{1-\sigma}))
$$

 They  give a Painlev\'e transcendent 
with the behaviour of theorem 1 of chapter \ref{PaInLeVe}.


\vskip 0.3 cm
\section{ Small Parameter Expansion in one Non-generic Case: ``Chazy
Solution''} 

  We need to investigate what happens if $\Re \sigma=1$. 

We do it only for $\sigma=1$. If we perform the small parameter
expansions as before, we find the same $\Omega_j^{(0)}$ than
before. But due to the exponent $z^{-1/2}$ the integration for
$\Omega_2^{(1)}$ gives 
$$
 \Omega_2^{(1)} = -{i\over 2} \tilde{a}^2 z^2 + i \tilde{b}^2\ln(z)
$$
 In this way, we find for $\Omega_1$ and $\Omega_3$ an
 expansion in power of $\epsilon$ with  coefficients which are
 polynomials in $\ln(z)$; also the powers $z^{-1/2}$, $z^{1/2}$,... ,
 $z^{n/2}$, $n>0$ appear in the coefficients. $\Omega_2$ is an 
 expansion in power of $\epsilon$ with coefficients which are
 polynomials in $\ln(z)$ and $z$. 

It is not obvious how to recombine $z$ and
$\epsilon$ when logarithms appear. We can put $\epsilon=1$. Anyway, we
see that the first correction to the constant ${i \over 2}$ in
$\Omega_2$ is  $\ln(s)$, which is {\it not a correction} to the
constant when $s\to 0$, because it diverges. 

 This makes us not trust the validity of the expansion
 (\ref{solution}) for $\sigma=1$.

\vskip 0.2 cm

 We can try an expansion which already contains logarithms of
 $\epsilon$. The experience with Chazy solutions to the Painlev\'e VI
 eq. suggests to choose: 
\be
   \Omega_j(z,\epsilon)= \sum_{k=-1}^{+\infty}  \sum_{n=0}^{+\infty}
    \Omega_{j,n}^{(k)}(z) ~{\epsilon^{2 k +1 \over 2} \over (\ln
 \epsilon)^n}
,~~~~j=1,3
\label{chazy1}
\ee
\be
 \Omega_2=  \sum_{k=0}^{+\infty}  \sum_{n=0}^{+\infty} 
\Omega_{2,n}^{(k)}(z) ~{\epsilon^{k} \over (\ln
 \epsilon)^n}
\label{chazy2}
\ee
Then we substitute in (\ref{systpower}) and we equate powers of
$\epsilon$ and $\ln \epsilon$. I omit long computations. The
requirement that we could re-compose the powers of $\ln \epsilon$ and
$\ln z$ appearing in the expansion in the form $\ln (z \epsilon)$
imposes very strong relations on the integration constants. The result
to which we are led  when we solve the equations for the   
coefficients $\Omega_{j,n}^{(k)}$ equating powers up to ${1\over( \ln
\epsilon)^n}$ and $\epsilon^{-1/2}$ is : 
$$ s = z \epsilon, ~~~~ \epsilon \to 0$$
$$
 \Omega_1 = {i \over  s^{1\over2}(\ln(s) + C)} +O\left({1\over (\ln
 \epsilon)^n}\right)
               + O(\epsilon^{1\over 2}) ,
$$
$$
 \Omega_3 = {-1 \over  s^{1\over2}(\ln(s) + C)} +O\left({1\over (\ln
 \epsilon)^n}\right)
               + O(\epsilon^{1\over 2}) ,
$$
$$ 
\Omega_2 = {i \over 2} + {i \over \ln s + C } + O\left({1\over (\ln
 \epsilon)^n}\right)
               + O(\epsilon^{1\over 2}).  
$$ 
  The substitution of these formulas in (\ref{PVI}) gives the
  asymptotic behaviour for $s\to 0$ of the Chazy solutions. We remark that the
  above $\Omega_j$'s imply
$$
 \Omega_1^2+\Omega_2^2+\Omega_3^2= -{1\over 4} +o(1)
$$
therefore, $\mu^2={1\over 4}$, an then we have only Chazy solutions!


\backmatter



\end{document}